\documentclass[reqno,a4paper]{amsart} 
\usepackage[english]{babel}

\usepackage{graphicx,amsmath,amsfonts,amssymb, amsthm, mathrsfs, tikz,tikz-cd, amscd,amsbsy,amstext,mathtools, enumerate,enumitem, stmaryrd, latexsym, empheq, verbatim,a4wide, ifthen,fancyhdr,tensor, xfrac, epigraph, dynkin-diagrams,tabularray, etoolbox, patchcmd} 
\usetikzlibrary{shapes.geometric}

\mathtoolsset{showonlyrefs,showmanualtags}

\usepackage{mathpazo}
\usepackage[all]{xy} 

\usepackage{hyperref}
\hypersetup{anchorcolor=black, linktocpage=true,
	colorlinks=true,
	citecolor=MyDarkBlue,
	linkcolor=MyDarkBlue,
	urlcolor=black,
	breaklinks=true,
	plainpages=true
}

\usepackage[hyperpageref]{backref}

\usetikzlibrary{positioning,shapes,shadows,arrows}

\usepackage{color}
\definecolor{MyDarkBlue}{rgb}{0.15,0.25,0.45}

\usepackage[utf8x]{inputenc}

\usepackage[toc,page]{appendix}
\usepackage[mathscr]{euscript} 
\allowdisplaybreaks

\newif\ifpersonal
\personaltrue 
\ifpersonal

\newcommand{\todo}[1]{\textcolor{red}{(Todo: #1)}}
\else
\newcommand\ast{\personal}[1]{\ignorespaces}
\newcommand\ast{\todo}[1]{\ignorespaces}
\fi

\newcommand{\calA}{\mathcal{A}}

\newcommand{\calC}{\mathcal{C}}

\newcommand{\calE}{\mathcal{E}}
\newcommand{\calF}{\mathcal{F}}
\newcommand{\calG}{\mathcal{G}}
\newcommand{\calH}{\mathcal{H}}

\newcommand{\calK}{\mathcal{K}}
\newcommand{\calL}{\mathcal{L}}
\newcommand{\calM}{\mathcal{M}}

\newcommand{\calP}{\mathcal{P}}
\newcommand{\calQ}{\mathcal{Q}}

\newcommand{\calS}{\mathcal{S}}
\newcommand{\calT}{\mathcal{T}}

\newcommand{\calV}{\mathcal{V}}

\newcommand{\calW}{\mathcal{W}}

\newcommand{\D}{\mathbb{D}}
\newcommand{\E}{\mathbb{E}}
\newcommand{\I}{\mathbb{I}}

\newcommand{\R}{\mathbb{R}}
\newcommand{\N}{\mathbb{N}}
\newcommand{\C}{\mathbb{C}}
\newcommand{\Z}{\mathbb{Z}}
\newcommand{\Q}{\mathbb{Q}}
\newcommand{\T}{\mathbb{T}}
\newcommand{\PP}{\mathbb{P}}
\newcommand{\A}{\mathbb{A}}
\newcommand{\G}{\mathbb{G}}
\newcommand{\F}{\mathbb{F}}

\newcommand{\LL}{\mathbb{L}}
\newcommand{\V}{\mathbb{V}}
\newcommand{\Y}{\mathbb{Y}}

\newcommand{\scrA}{\mathscr{A}}
\newcommand{\scrB}{\mathscr{B}}
\newcommand{\scrC}{\mathscr{C}}
\newcommand{\scrD}{\mathscr{D}}
\newcommand{\scrE}{\mathscr{E}}
\newcommand{\scrF}{\mathscr{F}}
\newcommand{\scrG}{\mathscr{G}}
\newcommand{\scrH}{\mathscr{H}}
\newcommand{\scrK}{\mathscr{K}}
\newcommand{\scrO}{\mathscr{O}}
\newcommand{\scrT}{\mathscr{T}}
\newcommand{\scrS}{\mathscr{S}}

\newcommand{\scrP}{\mathscr{P}}
\newcommand{\scrU}{\mathscr{U}}
\newcommand{\scrV}{\mathscr{V}}
\newcommand{\scrW}{\mathscr{W}}
\newcommand{\scrX}{\mathscr{X}}
\newcommand{\scrY}{\mathscr{Y}}
\newcommand{\scrZ}{\mathscr{Z}}

\newcommand{\sfE}{\mathsf{E}}
\newcommand{\sfH}{\mathsf{H}}
\newcommand{\sfG}{\mathsf{G}}

\newcommand{\sfL}{\mathsf{L}}
\newcommand{\sfP}{\mathsf{P}}

\newcommand{\sfS}{\mathsf{S}}
\newcommand{\sfT}{\mathsf{T}}

\newcommand{\sfU}{\mathsf{U}}
\newcommand{\sfV}{\mathsf{V}}

\newcommand{\sfp}{\mathsf{p}}
\newcommand{\sfs}{\mathsf{s}}
\newcommand{\sft}{\mathsf{t}}

\newcommand{\fraka}{\mathfrak{a}}
\newcommand{\frakg}{\mathfrak{g}}
\newcommand{\frakm}{\mathfrak{m}}

\newcommand{\bfC}{\mathbf{C}}
\newcommand{\bfD}{\mathbf{D}}

\newcommand{\bfT}{\mathbf{T}}
\newcommand{\bfF}{\mathbf{F}}
\newcommand{\bfH}{\mathbf{H}}
\newcommand{\bfM}{\mathbf{M}}
\newcommand{\bfP}{\mathbf{P}}

\newcommand{\bfU}{\mathbf{U}}

\newcommand{\bff}{\mathbf{f}}

\newcommand{\bfv}{\mathbf{v}}

\newcommand{\ch}{\mathsf{ch}}

\newcommand{\bfLambda}{\mathbf{\Lambda}}
\newcommand{\bfDelta}{\mathbf{\Delta}}

\newcommand{\id}{{\mathsf{id}}}
\newcommand{\Spec}{\mathsf{Spec}}

\DeclareMathOperator*{\colim}{colim}

\newcommand{\sfrel}{\mathsf{rel}}
\newcommand{\Hom}{\mathsf{Hom}}
\newcommand{\calHom}{\mathcal{H}\mathsf{om}}
\newcommand{\Ext}{\mathsf{Ext}}
\newcommand{\calExt}{\mathcal{E}\mathsf{xt}}

\newcommand{\End}{\mathsf{End}}

\newcommand{\Pic}{\mathsf{Pic}}
\newcommand{\Ind}{\mathsf{Ind}}
\newcommand{\Pro}{\mathsf{Pro}}
\newcommand{\Fun}{\mathsf{Fun}}
\newcommand{\FunL}{\mathsf{Fun}^{\mathsf L}}
\newcommand{\FunR}{\mathsf{Fun}^{\mathsf{R}}}

\newcommand{\PSh}{\mathsf{PSh}}

\DeclareMathOperator*{\fcolim}{``colim''}
\DeclareMathOperator*{\flim}{``lim''}
\DeclareMathOperator*{\fbigoplus}{``\bigoplus''}

\newcommand{\Top}{\mathsf{Top}}
\newcommand{\Cat}{\mathsf{Cat}}
\newcommand{\Alg}{\mathsf{Alg}}
\newcommand{\CAlg}{\mathsf{CAlg}}
\newcommand{\Aff}{\mathsf{Aff}}
\newcommand{\dAff}{\mathsf{dAff}}
\newcommand{\dSch}{\mathsf{dSch}}
\newcommand{\dSchqcqs}{\mathsf{dSch}^{\mathsf{qcqs}}}
\newcommand{\St}{\mathsf{St}}
\newcommand{\PrL}{\mathsf{Pr}^{\mathsf{L}}}

\newcommand{\PrLR}{\mathsf{Pr}^{\mathsf{L},\mathsf{R}}}

\newcommand{\Map}{\mathsf{Map}}

\newcommand{\sfHilb}{\mathsf{Hilb}}

\newcommand{\bfQuot}{\mathbf{Quot}}

\newcommand{\Mod}{\textrm{-} \mathsf{Mod}}
\newcommand{\Modd}{\mathsf{Mod}}

\newcommand{\qcqs}{\mathsf{qcqs}}

\newcommand{\naive}{\mathsf{nv}}

\newcommand{\catVect}{\mathsf{Vect}}
\newcommand{\catPerf}{\mathsf{Perf}}
\newcommand{\catAPerf}{\mathsf{APerf}}

\newcommand{\catCoh}{\mathsf{Coh}}
\newcommand{\IndCoh}{\mathsf{IndCoh}}
\newcommand{\catCohps}{\mathsf{Coh}_{\mathsf{ps}}}

\newcommand{\bfSet}{\mathbf{Set}}
\newcommand{\dASp}{\mathsf{dASp}}

\newcommand{\catCohb}{\mathsf{Coh}^{\mathsf{b}}}

\newcommand{\catQCoh}{\mathsf{QCoh}}

\newcommand{\catDb}{\mathsf{D}^\mathsf{b}}

\newcommand{\ev}{\mathsf{ev}}
\newcommand{\coev}{\mathsf{coev}}

\newcommand{\bfRep}{\mathbf{Rep}}
\newcommand{\bfPerf}{\mathbf{Perf}}
\newcommand{\bfPerfps}{\mathbf{Perf}_{\mathsf{ps}}}
\newcommand{\bfPerfpsext}{\mathbf{Perf}_{\mathsf{ps}}^{\,\mathsf{ext}}}

\newcommand{\bfPic}{\mathbf{Pic}}

\newcommand{\bfFlagPerf}{\mathbf{FlagPerf}}
\newcommand{\bfFlagPerfps}{\mathbf{FlagPerf}_{\mathsf{ps}}}
\newcommand{\bfFlagCoh}{\mathbf{FlagCoh}}
\newcommand{\bfFlagCohps}{\mathbf{FlagCoh}_{\mathsf{ps}}}

\newcommand{\bfCoh}{\mathbf{Coh}}
\newcommand{\bfCohps}{\mathbf{Coh}_{\mathsf{ps}}}

\newcommand{\bfCohext}{\mathbf{Coh}^{\mathsf{ext}}}
\newcommand{\bfCohpsext}{\mathbf{Coh}_{\mathsf{ps}}^{\mathsf{ext}}}

\newcommand{\bfMap}{\mathbf{Map}}

\newcommand{\trunc}[1]{\tensor*[^{\mathsf{cl}}]{#1}{}}
\newcommand{\op}{^{\mathsf{op}}}
\newcommand{\opp}{\mathsf{op}}
\newcommand{\an}{^{\mathsf{an}}}
\newcommand{\Sym}{\mathsf{Sym}}
\newcommand{\dSt}{\mathsf{dSt}}

\newcommand{\dGeom}{\mathsf{dGeom}}

\newcommand{\dGeomqs}{\dGeom^{\mathsf{qs}}}
\newcommand{\dGeomqsconn}{\dGeom^{\mathsf{qs},\mathsf{conn}}}
\newcommand{\dGeomqcqs}{\dGeom^{\qcqs}}

\newcommand{\LambdadSt}{\Lambda\textrm{-}\dSt}
\newcommand{\LambdadGeomqs}{\Lambda\textrm{-}\dGeom^{\mathsf{qs}}}
\newcommand{\Corr}{\mathsf{Corr}}

\newcommand{\TwoSeggraded}[1]{2_{#1}\textrm{-}\mathsf{Seg}}

\newcommand{\sfSegal}{\mathsf{Segal}}
\newcommand{\Spc}{\mathsf{Spc}}
\newcommand{\Sp}{\mathsf{Sp}}

\newcommand{\tiltedheart}{\tensor*[^{\upsilon}]{\scrC}{^\heartsuit}}

\newcommand{\BettiD}{\tensor*[^{\mathsf{top}}]{\bfD}{}}
\newcommand{\motD}{\tensor*[^{\mathsf{mot}}]{\bfD}{}}
\newcommand{\sfDM}{\mathsf{DM}}
\newcommand{\Shhyp}{\mathsf{Sh}^{\mathsf{hyp}}}

\newcommand{\ind}{\mathsf{ind}}

\newcommand{\proD}{\tensor*[^{\mathsf{pro}}]{\bfD}{}}

\newcommand{\HBM}{\mathsf{H}^{\mathsf{BM}}}

\newcommand{\HBMD}{\mathsf{H}^{\bfD}}
\newcommand{\HDnaive}{\tensor*[^{\mathsf{nv}}]{\sfH}{_{\bfD}^\ast}}

\newcommand{\HBMDGamma}{\mathsf{H}^{\bfD,\Gamma}}
\newcommand{\HBMDGammaext}{\tensor*[^{\mathsf{ext}}]{\mathsf{H}}{^{\bfD,\Gamma}_\ast}}
\newcommand{\HBMDGammaextzero}{\tensor*[^{\mathsf{ext}}]{\mathsf{H}}{^{\bfD,\Gamma}_0}}
\newcommand{\HBMDnaive}{\tensor*[^{\mathsf{nv}}]{\mathsf{H}}{^{\bfD}_\ast}}
\newcommand{\HBMDext}{\tensor*[^{\mathsf{ext}}]{\mathsf{H}}{^{\bfD}_\ast}}
\newcommand{\HBMDnaiveshort}[1]{\tensor*[^{\mathsf{nv}}]{\mathsf{H}}{^{\bfD}_{#1}}}

\newcommand{\HBMGamma}{\mathsf{H}^{\mathsf{BM},\Gamma}}

\newcommand{\Hmotnaive}{\tensor*[^{\mathsf{nv}}]{\mathsf{H}}{^{\mathsf{mot}}_\ast}}

\newcommand{\CBMD}{\mathsf{C}^{\bfD}_\bullet}
\newcommand{\CDnaive}{\tensor*[^{\mathsf{nv}}]{\mathsf{C}}{_{\bfD}^{\bullet}}}

\newcommand{\CBMDGamma}{\mathsf{C}^{\bfD,\Gamma}_\bullet}
\newcommand{\CBMDGammaext}{\tensor*[^{\mathsf{ext}}]{\mathsf{C}}{^{\bfD,\Gamma}_\bullet}}
\newcommand{\CBMDnaive}{\tensor*[^{\mathsf{nv}}]{\mathsf{C}}{^{\bfD}_\bullet}}

\newcommand{\CBMDnaiveshort}{\tensor*[^{\mathsf{nv}}]{\mathsf{C}}{^{\bfD}_{\bullet}}}

\newcommand{\CBM}{\mathsf{C}^{\mathsf{BM}}_\bullet}

\newcommand{\Cmot}{\mathsf{C}^{\mathsf{mot}}_\bullet}

\newcommand{\sfSH}{\mathsf{SH}}

\newcommand{\PsCpt}{\mathsf{PsCpt}}

\newcommand{\PrLomega}{\mathsf{Pr}^{\mathsf{L},\omega}}
\newcommand{\PrLotimes}{\mathsf{Pr}^{\mathsf{L},\otimes}}
\newcommand{\dCAlg}{\mathsf{dCAlg}}

\newcommand{\pr}{\mathsf{pr}}

\makeatletter
\newcommand{\ostar}{\mathbin{\mathpalette\make@circled\star}}
\newcommand{\make@circled}[2]{%
\ooalign{$\m@th#1\smallbigcirc{#1}$\cr\hidewidth$\m@th#1#2$\hidewidth\cr}%
}
\newcommand{\smallbigcirc}[1]{%
\vcenter{\hbox{\scalebox{0.77778}{$\m@th#1\bigcirc$}}}%
}
\makeatother

\newcommand{\triend}{\parbox{2mm}{\hfill} \hfill\mbox{\hspace{0.2mm}}\hfill$\triangle$}
\newcommand{\ocend}{\parbox{2mm}{\hfill} \hfill\mbox{\hspace{0.2mm}}\hfill$\oslash$}

\patchcmd{\part}{\normalfont}{\large\scshape}{}{}

\numberwithin{section}{part}
\numberwithin{equation}{section}

\newtheorem{theorem}{Theorem}[section]

\newtheorem{proposition}[theorem]{Proposition}
\newtheorem{lemma}[theorem]{Lemma}
\newtheorem{corollary}[theorem]{Corollary}

\newtheorem{definition*}{Definition}
\newtheorem{corollary*}{Corollary}
\newtheorem*{theorem*}{Theorem}
\newtheorem*{remark*}{Remark}
\newtheorem*{proposition*}{Proposition}
\newtheorem*{conjecture*}{Conjecture}
\newtheorem*{warning*}{Warning}

\newtheorem{theoremintroduction}{Theorem}

\newtheorem{corollaryintroduction}{Corollary}[section]

\theoremstyle{remark}
\newtheorem{ex}[theorem]{Example}
\newenvironment{example}{\begin{ex}}{\triend\end{ex}}

\newtheorem{exampleintroduction}[corollaryintroduction]{Example}

\newtheorem{war}[theorem]{Warning}
\newenvironment{warning}{\begin{war}}{\triend\end{war}}
\newtheorem{variant}[theorem]{Variant}

\theoremstyle{remark}
\newtheorem{rem}[theorem]{Remark}
\newenvironment{remark}{\begin{rem}}{\triend\end{rem}}
\newtheorem{remarkintroduction}[corollaryintroduction]{Remark}

\theoremstyle{definition}
\newtheorem{defin}[theorem]{Definition}
\newenvironment{definition}{\begin{defin}}{\ocend\end{defin}}
\newtheorem{construction}[theorem]{Construction}
\newtheorem{recollection}[theorem]{Recollection}

\newtheorem{definitionintroduction}[corollaryintroduction]{Definition}

\newtheorem{assum}{Assumption}

\newenvironment{assumption}{\begin{assum}}{\ocend\end{assum}}

\newtheorem{notat}[theorem]{Notation}
\newenvironment{notation}{\begin{notat}}{\ocend\end{notat}}

\title[COHAs, their categorification, and their representations via torsion pairs]{Cohomological Hall algebras, their categorification, and their representations via torsion pairs}

\author[D.-E.~Diaconescu]{Duiliu-Emanuel Diaconescu}
\address[Duiliu-Emanuel Diaconescu]{New High Energy Theory Center - Serrin Building, Rutgers, The State University Of New Jersey, 126 Frelinghuysen Rd., Piscataway, NJ 08854-8019, USA}
\curraddr{}
\email{\href{mailto:duiliu@physics.rutgers.edu}{duiliu@physics.rutgers.edu}}

\author[M.~Porta]{Mauro Porta}
\address[Mauro Porta]{Institut de recherche mathématique avancée (IRMA), Université de Strasbourg, France}
\curraddr{}
\email{\href{mailto:porta@math.unistra.fr}{porta@math.unistra.fr}}

\author[F.~Sala]{Francesco Sala}
\address[Francesco Sala]{Università di Pisa, Dipartimento di Matematica, Largo Bruno Pontecorvo 5, 56127 Pisa (PI), Italy}
\address{Kavli IPMU (WPI), UTIAS, The University of Tokyo, Kashiwa, Chiba 277-8583, Japan}
\curraddr{}
\email{\href{mailto:francesco.sala@unipi.it}{francesco.sala@unipi.it}}

\thanks{The work of D.-E.~D was partially supported by NSF grant DMS-180241, while the work of F.~S. was partially supported by JSPS KAKENHI Grant Number JP21K03197.}



\begin{document}

\begin{abstract}

	This paper addresses one long-standing limitation of the theory of cohomological and K-theoretical Hall algebras (COHAs/KHAs): they only produce ``positive parts'' of \textit{whole} algebras. In the special case of preprojective COHAs/KHAs associated to quivers, these whole algebras are known as Yangians, and their existence in the general case is due to Maulik-Okounkov. For geometric COHAs/KHAs arising from smooth surfaces, producing Yangians is an important open problem of the theory.
	
	We propose an approach to this problem based on torsion pairs and tiltings. More specifically, given a stable $\infty$-category $\scrC$ equipped with a $t$-structure $\tau$ and a torsion pair $\upsilon = (\scrT,\scrF)$ on $\scrC^\heartsuit$, we construct, under some favorable hypothesis: 
	\begin{enumerate}
		\item a COHA/KHA attached to $\scrT$; 
		\item a left and right representation of this COHA/KHA attached to $\scrF$.
	\end{enumerate}
	
	Together, they generate a \textit{whole} algebra; the operators arising from the left (resp.\ right) action correspond to positive (resp.\ negative) operators. In fact, we work equally at the categorical level, thereby obtaining categorifications of the corresponding KHAs and of its representations.
	
	In the quiver case, we recover the action of the two-dimensional COHA of a quiver on the cohomology of Nakajima quiver varieties within our framework and we provide its categorification. Besides the quiver case, we also apply our framework to two explicit torsion pairs on a smooth projective complex surface, and we investigate the corresponding Hall algebras and their representations associated to them. In these cases, our formalism recovers the action of affine Yangian of $S$ on the cohomology of moduli spaces of Gieseker-stable sheaves on $S$, and provides a categorification of it. Finally, we slightly modify our method to construct representations of the COHA and the CatHA of zero-dimensional sheaves on $S$ on the (motivic) Borel-Moore homology and bounded derived category of the (derived) moduli stack of Pandharipande-Thomas stable pairs on surfaces and ``dually'' on (derived) relative Hilbert schemes of points. Contrary to the previous examples, the action comes from extensions of flags rather than extensions of sheaves (as it is more commonly encountered in the theory of cohomological Hall algebras).

	These results are achieved via several preliminary foundational results. Most importantly, we need to adapt Khan's motivic Borel-Moore homology framework to take into account stacks that are not necessarily quasi-compact. We also show that to any pair $(\scrC,\tau)$ satisfying certain natural conditions one can functorially attach a COHA. This construction does not only generalize and categorify the known examples of two-dimensional COHAs of quivers and surfaces, but also provide new examples of COHAs and CatHAs associated to $t$-structures underlying Bridgeland stability conditions (in the sense of \cite{BLMNPS_Stability}).

\end{abstract}

\maketitle\thispagestyle{empty}

\tableofcontents


\makeatletter
\@removefromreset{section}{part}
\renewcommand\thesection{\arabic{section}}
\makeatother

\setcounter{section}{0}

\section{Introduction}

\subsection{What are Cohomological Hall algebras?}

\textit{Cohomological Hall algebras (COHAs)} are a special kind of associative algebras that arise as a convolution product on the (Borel-Moore) homology of moduli stacks of sheaves. \textit{Sheaves} can also be understood in a rather liberal way: the most well-studied examples include moduli of coherent sheaves on surfaces, Higgs sheaves, local systems or flat connections on curves, or representations of preprojective algebras of a quiver. These examples have the number $2$ in common: the resulting (abelian) categories of sheaves have all homological dimension $2$.

A uniform framework can be constructed as follows. Let $k$ be a field and let $\scrA$ be a $k$-linear abelian category and let $\calM_\scrA$ be the moduli of objects of $\scrA$, as defined by Artin and Zhang \cite{AZ_Hilbert} (see also \cite[\S7.1]{Existence_Moduli}), seen as a $k$-linear stack. One can equally form a moduli stack parametrizing extensions of objects $\calM_\scrA^{\mathsf{ext}}$, which gives rise to a correspondence
\begin{align}
	\begin{tikzcd}[ampersand replacement=\&]
		\& [-15pt] \calM_\scrA^{\mathsf{ext}} \arrow{dr}{q} \arrow{dl}[swap]{p} \\
		\calM_\scrA \times \calM_\scrA \& \& \calM_\scrA 
	\end{tikzcd}\ .
\end{align}
As understood by Dyckerhoff and Kapranov \cite{DK}, this convolution diagram is the shadow of a more complicated object, a simplicial stack satisfying a property known as the \textit{$2$-Segal property} (see \S\ref{sec:2_Segal} for a review). The category of $2$-Segal stacks is equivalent to the ($\infty$-)category of associative algebras in correspondences $\Corr^\times(\St_k)$. In particular, any lax monoidal functor
\begin{align}
	\sfH_\ast \colon \Corr^\times(\St_k) \longrightarrow \Modd_R 
\end{align}
gives a way of converting a $2$-Segal into an associative algebra (here $R$ denotes any coefficient commutative ring). Concretely, this procedure sends $\calM_\scrA$ with the above convolution diagram to $\sfH_\ast(\calM_\scrA)$ equipped with the product
\begin{align}
	\star_{\mathsf{Hall}} \coloneqq p_\ast q^! \colon \HBM(\calM_\scrA) \otimes \HBM(\calM_\scrA) \longrightarrow \HBM(\calM_\scrA) \ . 
\end{align}
For the purposes of this introduction, we will refer to such a lax monoidal functor as a \textit{homological realization functor}.

One major source of problems lies in finding such homological realizations. Indeed, one needs at least to restrict the attention to a reasonable subcategory of $\St_k$, and also impose conditions on the morphisms of correspondences:
\begin{enumerate}\itemsep=0.2cm
	\item \label{item:1-introduction} one needs the map $p$ to be lci to perform the Gysin pullback $p^!$;
	
	\item \label{item:2-introduction} one needs the map $q$ to be proper to perform the pushforward $q_\ast$.
\end{enumerate}
After \cite{Porta_Sala_Hall} was written, it quickly became commonly accepted that the best way of solving problem \eqref{item:1-introduction} is to consider $\calM_\scrA$ as a derived stack. Problem \eqref{item:2-introduction} is related to the properness of Quot schemes, and in many standard examples this is easy to deal with. The motivic framework of Khan \cite{Khan_VFC} allows for a neat construction of such homological realizations, although it has to be tweaked in order to apply to exact setting it is needed for this construction.

It is perhaps less common knowledge that $\calM_\scrA$ typically admits \textit{several derived enhancements}. One large part of the current paper grew out of the following observation:
\begin{quote}
	If $\scrC$ is a stable $k$-linear $\infty$-category equipped with a $t$-structure $\tau = (\scrC_{\geqslant 0},\scrC_{\leqslant 0})$ such that $\scrA \simeq \scrC^\heartsuit$, then Toën-Vaquié's moduli of objects $\calM_\scrC$ \cite{TV_Moduli} canonically contains an open substack $\bfCoh(\scrC,\tau)$ providing a derived enhancement for $\calM_\scrA$.
\end{quote}
Although \textit{literally} wrong, the above slogan becomes true as soon as one imposes some niceness conditions ($\scrC$ should be of finite type, and the $t$-structure $\tau$ should satisfy \textit{openness}).
The whole Part~\ref{part:foundations} is devoted to identify the precise conditions under which the Borel-Moore homology of $\bfCoh(\scrC,\tau)$ supports the structure of a COHA (which should then be thought as associated to the \textit{pair} $(\scrC,\tau)$).

\medskip

Of course, this perspective \textit{begs} the following question: is it possible to understand precisely how the Hall algebra depends on the $t$-structure? In general, answering this question is challenging. A solution in the case of the $t$-structures arising in the derived McKay correspondence is provided in \cite{DPSSV_COHA}. We refer to the introduction of that paper for a more detailed discussion of this problem and of the possible techniques to address it.

\subsection{Why study COHAs?}

This question has many possible answers, depending on whether one is more interested in the algebraic, the geometric or the physical aspects of the theory. From a philosophical point of view, COHAs should be understood as the ``largest'' algebras of (positive) Hecke operators acting on the homology of a moduli space $M_\scrA^{\mathsf{st}}$ of \textit{stable} objects of $\scrA$. This allows to consider $\sfH_\ast(M_\scrA^{\mathsf{st}})$ as a representation of a very large algebra, which in turn makes $\sfH_\ast(M_\scrA^{\mathsf{st}})$ into a rather rigid object from the algebraic perspective. Famous successful applications of this philosophy are provided by e.g. Nakajima's works \cite{Nakajima1994, Nakajima_Heisenberg,Nakajima_book, Nakajima_Quiver}. 

On the other hand, COHAs are so big that it is typically extremely challenging to get a grasp of their internal structure or of their representation theory. In recent years, many advances have been made. We briefly review them, distinguishing two important classes of cases: noncommutative ones (in the guise of the COHAs of preprojective algebras of quivers) and commutative ones (in the guise of COHAs of algebraic surfaces).

\begin{exampleintroduction}[Preprojective COHAs]
	Let $\calQ$ be a quiver, $\Pi_\calQ$ be its preprojective algebra and denote by $\bfRep(\Pi_\calQ)$ the (derived) moduli stack of finite-dimensional representations of $\Pi_\calQ$. Applying the general construction mentioned before in this setting, we obtain the \textit{cohomological Hall algebra (COHA) of the preprojective algebra} of $\calQ$ -- also referred to as the \textit{preprojective} COHA of $\calQ$, or the \textit{2d} COHA of $\calQ$.
	This algebra was first introduced by Schiffmann and Vasserot in the case where $\calQ$ is the one-loop quiver \cite{SV_Cherednik}. A \textit{K-theoretic} version of this construction—referred to as the \textit{K-theoretic Hall algebra} (KHA)—was introduced earlier by the same authors in \cite{SV_Elliptic}; in this setting, the underlying vector space is $G_0^T(\bfRep(\Pi_\calQ))$, the ($T$-equivariant) Grothendieck group of coherent sheaves on $\bfRep(\Pi_\calQ)$. Both constructions were later extended to arbitrary quivers and to any oriented Borel–Moore homology theory in \cite{YZ_COHA}.
	
	Preprojective COHAs of quivers yield the ``largest'' algebras of positive Nakajima operators acting on the equivariant cohomologies of Nakajima quiver varieties. Furthermore, they can be identified with the positive nilpotent parts $\Y_\calQ^{\mathsf{MO},+}$ of the \textit{Maulik-Okounkov Yangian} $\Y_\calQ^{\mathsf{MO}}$ of the quiver $\calQ$ \cite{MO_Yangian}, as demonstrated in \cite{BD_Yangian, SV_MO} (see also \cite{SV_Yangians}).
	Alternatively, it can be said that the preprojective COHA of $\calQ$ provides a \textit{geometric interpretation} (in terms of explicit Hecke operators) of $\Y_\calQ^{\mathsf{MO},+}$. In \cite{VV_KHA} a similar relationship in G-theory is established, when $\calQ$ is a finite or affine ADE quiver, and conjectured for arbitrary quivers: in this setting, Yangians are replaced by \textit{quantum loop algebras}.
\end{exampleintroduction}

\begin{exampleintroduction}[Surface COHAs]
	Let $S$ be a smooth (quasi-)projective surface. The construction of the COHA and KHA\footnote{We talk about KHA for the Hall algebra obtained using the $G$-theory as choice of a homological realization functor.} for coherent sheaves on $S$ has been carried out in \cite{KV_Hall}, generalizing earlier work in \cite{Minets, SS_Higgs} (that dealt with the case of the cotangent bundle of a smooth curve), and in \cite{Zhao_KHA} for the case of zero-dimensional sheaves on $S$. A \textit{categorification} of these constructions was given in \cite{Porta_Sala_Hall}, where the second and third-named authors endowed the stable $\infty$-category of complexes with coherent cohomology on suitable derived enhancements of moduli stacks of coherent sheaves on $S$ with an $\E_1$-monoidal structure. In \textit{loc.\ cit.}, a similar construction was also provided for the moduli stacks arising in the non-abelian Hodge and Riemann–Hilbert correspondences for smooth projective complex curves. The categorical Hall algebras of \cite{Porta_Sala_Hall} provide the right framework to define a categorification of \textit{BPS invariants}, called \textit{quasi-BPS categories}, currently investigated by Pădurariu and Toda (see \cite{PT_BPS_K3, PT_BPS_Higgs} and references therein).
	
	The algebraic structure of the subalgebra of zero-dimensional sheaves on $S$ has nowadays been thoroughly understood \cite{MMSV}, building on the ideas that lead to a proof of the $P=W$ conjecture \cite{HMMS}.
\end{exampleintroduction}

The above examples highlight an important \textit{asymmetry} between the quiver and the surface case. Indeed, in the first case, the picture is richer thanks to the presence of \textit{Yangians}, Hopf algebras of operators that contain COHAs as positive nilpotent halves. The Maulik-Okounkov Yangian is a ``deformation'' of the enveloping algebra $\mathsf U(\mathfrak g^{\mathsf{MO}}_\calQ[u])$ of the algebra of currents of Maulik-Okounkov Lie algebra $\mathfrak g^{\mathsf{MO}}_\calQ$ (or, to be more precise, of the BPS Lie algebra $\mathfrak g^{\mathsf{MO},T}_\calQ$, see \cite[Corollary 1.2]{BD_Yangian}). There are therefore two aspects that emerge from this picture, that clarify what is the ``quantum nature'' of the preprojective COHA:
\begin{enumerate}\itemsep=0.2cm
	\item \label{item:1-introduction-2} the existence of a canonical filtration on the preprojective COHA, whose associated graded is isomorphic (as a module) to the universal enveloping algebra of the positive nilpotent part of a Lie algebra (this is what we refer to as a \textit{deformation});
	
	\item \label{item:2-introduction-2} the existence of a \textit{larger} algebra of operators, coming with a triangular decomposition, with respect to which the COHA is a positive half.
\end{enumerate}
Both these properties are completely unclear on the geometric side, and it could be argued that one of the central challenges in the theory of surface COHAs is to construct a quantum group associated to the surface itself.
This is already challenging for cotangent bundles of smooth projective curves -- see e.g.\ \cite[\S8]{Schiffmann_Kac_polynomial}.

For the subalgebra of zero-dimensional coherent sheaves, the explicit description in terms of generators and relations obtained in \cite{MMSV} allowed to construct an isomorphism with the Yangian of the \textit{BPS Lie algebra of $S$} (in the sense of \cite{DHSM_BPS-II}), when $S$ has trivial canonical bundle. A first step toward a more complete understanding of the COHA of one-dimensional coherent sheaves on $S$ was taken in \cite{DPSSV_COHA}, co-authored by the authors of the present paper together with Schiffmann and Vasserot. In \textit{loc.\ cit.}, we considered the minimal resolution $\pi\colon S \to \C^2/G$ of a Kleinian singularity and introduced a \textit{nilpotent} COHA of coherent sheaves on $S$ \textit{set-theoretically supported} on the exceptional divisor $\pi^{-1}(0)$. Via the derived McKay correspondence (as used in \cite{DPS_McKay}), we related this nilpotent COHA of $S$ to the \textit{nilpotent} preprojective COHA of the McKay quiver $\calQ$ of $G$, thereby obtaining a description of the former in terms of the Yangian of $\calQ$.

Another key limitation in the current theory of cohomological Hall algebras is that they are generally expected to geometrically realize only positive nilpotent parts of full Yangians. A central open problem is how to realize the full Yangian structure starting from a COHA. In the classical Hall algebra setting for hereditary categories, this issue has been resolved algebraically by introducing a bialgebra structure and taking its (reduced) Drinfeld double (see, e.g., \cite{Schiffmann_Lectures} and references therein). In the cohomological setting--especially beyond the quiver case--it remains unclear how to define a coproduct that would endow a COHA with a bialgebra structure. A similar issue arises in the context of categorical Hall algebras, where notions of ``categorical'' coproduct and Drinfeld double are still lacking.

The present paper proposes a general framework for defining both cohomological (in the motivic sense à la A.~Khan) and categorical Hall algebras, and outlines a strategy for addressing the limitation \eqref{item:2-introduction-2} described above in the existing theory of COHAs and CatHAs. Below, we review our main results.

\subsection{COHAs and CatHAs: existence (with an application to Bridgeland stability conditions)}

Fix a field $k$ of characteristic zero. Let $\scrC$ be a compactly generated, $k$-linear stable $\infty$-category and let $\tau = (\scrC_{\geqslant 0},\scrC_{\leqslant 0})$ be a $t$-structure on $\scrC$.

Out of this data one can construct a moduli of objects $\calM_\scrC$, following Toën-Vaquié \cite{TV_Moduli}. This is a derived stack parametrizing families of objects in $\scrC$ satisfying a certain \textit{non-commutative analog} of the proper support condition, known as \textit{pseudo-perfectness}. We refer the reader to \S\ref{subsec:integral-transforms} for a precise definition of this notion and to \S\ref{subsec:Toen_Vaquie} for a review of Toën-Vaquié's theorem. The $t$-structure allows to define a substack $\bfCohps(\scrC,\tau) \subset \calM_\scrC$, parametrizing \textit{$\tau$-flat families} of pseudo-perfect objects in $\scrC$. We refer to \S\ref{subsec:flatness_and_openness_of_t_structures} for the precise definition of this notion.

It is formal to extend this definition to a full $2$-Segal derived stack $\calS_\bullet \bfCohps(\scrC,\tau)$. In the main body of the paper we will systematically use this simplicial notation, but in the introduction we use the more intuitive notation
\begin{align}
	\bfCohpsext(\scrC,\tau) \coloneqq \calS_2 \bfCohps(\scrC,\tau) \ . 
\end{align}
We denote by
\begin{align}
	\partial_0 \times \partial_2 \colon \bfCohpsext(\scrC,\tau) \longrightarrow \bfCohps(\scrC,\tau) \times \bfCohps(\scrC,\tau) 
\end{align}
and 
\begin{align}
	\partial_1 \colon \bfCohpsext(\scrC,\tau) \longrightarrow \bfCohps(\scrC,\tau)
\end{align}
the maps sending an extension to the outer terms (in reverse order) and to the middle term, respectively.

As anticipated in the previous section, in order to extract a COHA out of this data, we need a lax monoidal homological realization functor which has pullback for $\partial_0 \times \partial_2$ and pushforward for $\partial_1$. \textit{In principle}, one would like to either use the categorical $6$-functors formalism $\catCohb(-)$ as in \cite{Porta_Sala_Hall}, or the motivic framework of Khan \cite{Khan_VFC}. However, there are several technical difficulties that need to be addressed. The first, is that the theory of Khan, as written, does not yield a functor, essentially because the Gysin operation requires a twist. One needs to tweak a little his definition of motivic Borel-Moore homology to bypass this issue (see Definition~\ref{def:abstract_BM_homology} and Remark~\ref{rem:explicit_bigrading}). A more serious problem is related to the fact that the stack $\bfCohps(\scrC,\tau)$ is typically \textit{not quasi-compact}. This leads to consider the Borel-Moore homology groups as topological vector spaces, equipped with a topology induced by quasi-compact exhaustions. A final difficulty is the fact that Khan's theory provides the pushforward for proper morphisms, but in many cases of geometric interest $\partial_1$ only satisfies the valuative criterion of properness, without being quasi-compact (it is an example of a \textit{locally rpas morphism}, see Definition~\ref{def:modified_classes_of_morphisms}). This final issue is solved via a \textit{renormalization procedure}. We refer the reader to \S\ref{sec:homological-invariants} for a thorough discussion of these ideas, and especially to Theorem~\ref{thm:functoriality_of_genuine_BM_homology_admissible} for the final output, in the form of a lax monoidal homological realization functor.

Having taken care of the \textit{linearization procedure}, we still need to identify the good conditions on the pair $(\scrC,\tau)$ guaranteeing that we can apply our homological realization to the $2$-Segal derived stack $\calS_\bullet \bfCohps(\scrC,\tau)$. The following is our first main result, which should be thought of as a generalization and categorification of the construction of COHAs of 2CY categories in \cite[\S5]{DHSM_BPS}.
\begin{theoremintroduction}[{Theorem~\ref{thm:COHA-Coh} and Corollary~\ref{cor:induced_COHA-Serre}}]\label{thm:COHA-Coh-introduction}
	Assume that
	\begin{itemize}\itemsep0.2cm
		\item $\scrC$ is of finite type;
		\item $\sfS_\scrC^![2]$ is $t$-exact, where $\sfS_\scrC^!$ is the left Serre functor\footnote{The notion of Serre functors is introduced in Definition~\ref{def:Serre_functors}.} of $\scrC$;
		\item $\tau$ satisfies Assumption~\ref{assumption:t_structure_filtered_colimits} and it universally satisfies openness of flatness in the sense of Definition~\ref{def:openness-flatness};
		\item the map $\partial_1 \colon \bfCohpsext(\scrC,\tau) \longrightarrow \bfCohps(\scrC,\tau)$ is locally rpas.
	\end{itemize}
	Then $\catCohb_{\mathsf{pro}}( \bfCohps(\scrC, \tau) )$ has the structure of an $\E_1$-monoidal stable pro-$\infty$-category, whose underlying tensor product is given by the composition
	\begin{align}\label{eq:multiplication-introduction}
		\begin{tikzcd}[ampersand replacement=\&]
			\catCohb_{\mathsf{pro}}( \bfCohps(\scrC, \tau) ) \times \catCohb_{\mathsf{pro}}( \bfCohps(\scrC, \tau) ) \ar{r}{\boxtimes} \& \catCohb_{\mathsf{pro}}( \bfCohps(\scrC, \tau) \times \bfCohps(\scrC, \tau) ) \\
			{} \ar{r}{(\partial_1)_\ast \circ (\partial_2\times \partial_0)^\ast} \& \catCohb_{\mathsf{pro}}( \bfCohps(\scrC, \tau) )
		\end{tikzcd}\ .
	\end{align}
	Similarly, the topological abelian groups
	\begin{align}
		G_0( \bfCohps(\scrC, \tau) )\quad \text{and} \quad \HBM_\ast( \bfCohps(\scrC, \tau) )
	\end{align}
	have the structures of unital associative algebras.
\end{theoremintroduction}		

\begin{remarkintroduction}[Motivic formalisms]
	In the main body, we prove a more precise version of the above theorem, where $G_0$ and $\HBM_\ast$ is replaced by motivic Borel-Moore homology groups with more general coefficients. Here motivic is understood in an axiomatic way (see Definition~\ref{def:motivic_formalism}), which allows for an even larger coverage of examples.
\end{remarkintroduction}

\begin{remarkintroduction}[Geometric subalgebras]
	Similar results hold if we replace $\bfCohps(\scrC, \tau)$ by an open substack $\bfT$ of it such that such that for every field $\kappa$ the abelian category $\catCoh_\bfT(\scrC_\kappa,\tau_\kappa)$ is closed under extensions in $\catCohps^\heartsuit(\scrC_\kappa,\tau_\kappa)$ and the corresponding map $\partial_1 \colon \bfCohext_\bfT(\scrC,\tau) \longrightarrow \bfT$ is locally rpas. Here, we denote by $\catCoh_\bfT(\scrC_\kappa,\tau_\kappa)$ the full subcategory of $\catCohps^\heartsuit(\scrC_\kappa,\tau_\kappa)$ spanned by the objects that belong to the image of $\bfT(\kappa)$.
\end{remarkintroduction}

\begin{remarkintroduction}[Locally rpas]
	At first glance, it might seem disappointing that we do not provide any intrinsic condition in the pair $(\scrC,\tau)$ guaranteeing that $\partial_1$ is locally rpas. However, it is possible to find strong enough conditions on $(\scrC,\tau)$ guaranteeing this property: see for instance Proposition~\ref{prop:TV_properness} for a result of Toën-Vaquié in this direction that requires $\scrC$ to be smooth and proper, or \cite[Proposition~3.3.6]{Lampetti_Good_moduli} for an improved version. Nevertheless, we will also consider cases (arising from tiltings) where the locally rpas condition is satisfied but the sufficient conditions of the above criteria are not verified. For this reason, we decided to keep the abstract condition on the map $\partial_1$ in the formulation of the above theorem, to maximize its applicability.
\end{remarkintroduction}

Theorem~\ref{thm:COHA-Coh-introduction} enables the definition of COHAs and CatHAs associated with Bridgeland stability conditions. More precisely, let $X$ be a smooth projective complex variety and let $\scrD$ be a $\C$-linear strong (in the sense of \cite[Definition~3.5]{BLMNPS_Stability}) semiorthogonal component of $\catPerf(X)$ of finite cohomological amplitude (in the sense of \cite[Definition~3.7]{BLMNPS_Stability}) equipped with a Serre functor $\sfS_\scrD\simeq [2]$. Examples of such categories $\scrD$ include the category $\catPerf(S)$, with $S$ a K3 surface, or the Kuznetsov component $\mathcal{K}u(X)$, with $X$ either a Fano 3fold of Picard rank one (different from the complete intersection of a quadric and a cubic in $\PP^5$), a smooth cubic 4fold in $\PP^5$, or a Gushel-Mukai variety.

Let $\sigma$ be a stability condition on $\scrD$ with respect to a finite rank free abelian group $\Lambda$ (in the sense of \cite[Definition~21.15]{BLMNPS_Stability}). In particular, the underlying $t$-structure universally satisfies openness of flatness (cf.\ Proposition~20.8 of \textit{loc.\ cit.}). Denote by $\bfT$ either the moduli stack $\bfCohps(\scrC, \tau)$, where $\tau$ is the induced $t$-structure on $\scrC\coloneqq \Ind(\scrD)$ whose heart is $\Ind(\calP_\sigma((0,1]))$, where $\calP_\sigma$ is the slicing associated to $\sigma$, or the moduli stack $\bfCohps^{\sigma\textrm{-}\mathsf{ss}, \mu}(\scrC, \tau)$ of $\sigma$-semistable objects on $\scrC$ of fixed slope $\mu$, or the moduli stack
\begin{align}
	\bigsqcup_{\mathbf{v}\in \Z\calS}\, \bfCohps^{\sigma\textrm{-}\mathsf{ss}, \mu}(\scrC, \tau; \mathbf{v})\ ,
\end{align}
where $\bfCohps^{\sigma\textrm{-}\mathsf{ss}, \mu}(\scrD, \tau; \mathbf{v})$ denotes the moduli stack of $\sigma$-semistable objects on $\scrD$ with slope $\mu$ and with Mukai vector $\mathbf{v}\in \Lambda$. Here, $\calS$ is a set of Mukai vectors spanning a sublattice of $\Lambda$.
\begin{corollaryintroduction}[{Corollary~\ref{cor:COHA-stability-condition}}]\label{cor:COHA-stability-condition-introduction}
	$\catCohb_{\mathsf{pro}}( \bfT )$ has the structure of an $\E_1$-monoidal stable pro-$\infty$-category, whose underlying tensor product is induced by \eqref{eq:multiplication}. 
	
	Similarly, the topological abelian groups
	\begin{align}
		G_0( \bfT )\quad \text{and} \quad \HBM_\ast( \bfT )
	\end{align}
	have the structures of unital associative algebras.
\end{corollaryintroduction}
Note that the above results allow us to explore a potential definition of quasi-BPS categories, in the sense of Pădurariu and Toda, for arbitrary Bridgeland stability conditions on 2-Calabi–Yau categories.

\subsection{Representations of COHAs arising from torsion pairs}

We now go back to the \textit{limitations of the COHAs}, and explain how our methods allow to construct larger algebras of operators. To understand where the limitation comes from, recall that in the works of Nakajima \cite{Nakajima_Heisenberg, Nakajima_book} and Grojnowski \cite{Grojnowski_Hilbert} on the action of the \textit{Heisenberg algebra} on the cohomology of the Hilbert scheme of points $\mathsf{Hilb}(S)$ of a smooth projective surface $S$, the \textit{nilpotent} Nakajima Hecke operators are divided into two classes: \textit{positive} and \textit{negative}. Positive operators can be recovered from the COHA of the surface (cf.\ \cite[Proposition~7.8]{MMSV}). Negative operators arise, \textit{in principle}, by reading the Hall convolution backwards: instead of performing $\partial_{1,\ast} \circ (\partial_0 \times \partial_2)^!$ one would like to rather perform $(\partial_0 \times \partial_2)_\ast \circ \partial_1^!$ -- except that this is impossible in $2$-dimensional situation. The situation is however definitely not hopeless: Nakajima's construction was later generalized for COHAs, as well as in $K$-theory when $S = \C^2$, in \cite{SV_Elliptic,SV_Cherednik}. A more explicit approach, closer in spirit to Nakajima’s original construction, was developed by Neguţ in \cite{Negut_Shuffle, Negut_Categorification}.

In the case of $\mathsf{Hilb}(S)$, the key geometric property that accounts for the existence of a relative $2$-Segal derived stacks is the \textit{two-sided Hecke pattern property} for the stack $\bfCoh_{0\textrm{-}\!\dim}(S)$ of $0$-dimensional sheaves on $S$. This property provides the structural foundation needed to construct both left and right actions of the COHA of $0$-dimensional sheaves on $S$. It was formalized in \cite[\S5.1]{KV_Hall} (see also \cite[\S6.1]{MMSV} and Definition~\ref{def:Hecke-patterns} in the current paper).

It is worth noting that constructing two-sided Hecke patterns is highly nontrivial. The only other known example for the stack $\bfCoh_{0\textrm{-}\!\dim}(S)$ involves the moduli space of Gieseker-stable sheaves on $S$ with fixed rank and first Chern class (see \cite{Negut_Shuffle, Negut_Categorification}). To build left and right representations of other COHAs, one needs to relax this condition.

The really new contribution of the present work is a systematic method to construct both left and right representations of COHAs without relying on two-sided Hecke patterns. The input of our method is the given of a pair $(\scrC,\tau)$ satisfying the assumptions of Theorem~\ref{thm:COHA-Coh-introduction}, \textit{plus a torsion pair} $\upsilon = (\scrT,\scrF)$ on the heart $\scrC^\heartsuit$. This enables to construct two substacks $\bfCoh_{\scrT}(\scrC,\tau)$ and $\bfCoh_{\scrF}(\scrC,\tau)$ of $\bfCohps(\scrC,\tau)$. The general idea is that, if favorable conditions are imposed, then the Borel-Moore homology of the torsion substack $\bfCoh_{\scrT}(\scrC,\tau)$ inherits an algebra structure and it acts \textit{both on the left and on the right} on the Borel-Moore homology of the torsion-free substack $\bfCoh_{\scrF}(\scrC,\tau)$. The left action encodes Nakajima's positive operators, and the right action encodes negative ones.

It is possible to describe these actions via an extension of the notion of $2$-Segal derived stacks, called \textit{relative $2$-Segal derived stacks} (that encode the notion of representation in correspondences). This notion is thoroughly reviewed in \S\ref{sec:relative_Segal}, and more generally Part~\ref{part:Segal} is devoted to establish abstract criteria useful to construct such relative $2$-Segal derived stacks. In a sense, the whole paper is a consequence of the simplicial formalism developed in \S\ref{sec:1-cosk}. However, the exact details are too complicated to be properly addressed in the introduction, so we rather focus on the applications of these ideas.

The precise implementation of the ideas described above passes through a systematic use of tilted hearts. In practice, our \textit{right} action of $\bfCoh_\scrT(\scrC,\tau)$ on $\bfCoh_\scrF(\scrC,\tau)$ is constructed as a \textit{left} action of $\bfCoh_{\scrT}(\scrC,\tau_\upsilon)$ on $\bfCoh_{\scrF[1]}(\scrC,\tau)$, where $\tau_\upsilon$ denotes the tilt of $\tau$ with respect to the torsion pair $\upsilon = (\scrT,\scrF)$ (see \S\ref{subsubsec:tiltings} for a review of tiltings).

In the theorem below, we will denote by $\calS_1^\ell\bfFlagCoh^{(1),\dagger}_{\scrT,\scrF}(\scrC,\tau)$ the simplicial level one of the relative $2$-Segal stack encoding the left action, and by $\calS_1^r\bfFlagCoh^{(1),\dagger}_{\scrF[1],\scrT}(\scrC,\tau_\upsilon)$ the one encoding the right action. We can state our second main result as follows:
\begin{theoremintroduction}[{Theorem~\ref{thm:left-right-action}}]\label{thm:introduction-left-right-action}
	Let $\scrC$ be a compactly generated, stable $k$-linear $\infty$-category equipped with a $t$-structure $\tau = (\scrC_{\geqslant 0}, \scrC_{\leqslant 0})$ satisfying the assumptions of Theorem~\ref{thm:COHA-Coh-introduction}. Let $\upsilon = (\scrT, \scrF)$ be a torsion pair on $\scrC^\heartsuit$ such that
	\begin{enumerate}\itemsep0.2cm
		\item both $\scrT$ and $\scrC^\heartsuit$ are compactly generated,
		\item the inclusion $\scrT \hookrightarrow \scrC^\heartsuit$ preserves compact objects,
		\item $\sfS_\scrC^![2]$ is $t$-exact with respect to the tilted $t$-structure $\tau_\upsilon$,
		\item $\upsilon = (\scrT, \scrF)$ is open in the sense of Definition~\ref{def:open_torsion_pair},
		\item the map
			\begin{align}
				\partial_1 \colon \calS_2\bfCoh_{\scrT}(\scrC,\tau) \longrightarrow \bfCoh_{\scrT}(\scrC,\tau)
			\end{align}
		is locally rpas\footnote{This notion is introduced in Definition~\ref{def:modified_classes_of_morphisms}.}, and
	\item both maps
	\begin{align}
		\varpi_0 \colon \calS_1^\ell\bfFlagCoh^{(1),\dagger}_{\scrT,\scrF}(\scrC,\tau) \longrightarrow \bfCoh_{\scrF}(\scrC,\tau)
	\end{align}
	and
	\begin{align}
		\varpi_1 \colon \calS_1^r\bfFlagCoh^{(1),\dagger}_{\scrF[1],\scrT}(\scrC,\tau_\upsilon) \longrightarrow \bfCoh_{\scrF[1]}(\scrC,\tau)
	\end{align}
	are locally rpas.
	\end{enumerate}
	Then, 
	\begin{itemize}\itemsep0.2cm
		\item $\catCohb_{\mathsf{pro}}(\bfCoh_{\scrT}(\scrC,\tau))$ inherits the structure of a $\E_1$-monoidal pro-$\infty$-category, and
		
		\item $\catCohb_{\mathsf{pro}}(\bfCoh_{\scrF}(\scrC,\tau))$ has both the structure of a categorical left and of a categorical right module over $\catCohb_{\mathsf{pro}}(\bfCoh_{\scrT}(\scrC,\tau))$.
	\end{itemize}
	Similarly:
	\begin{itemize}\itemsep0.2cm
		\item the topological vector spaces
		\begin{align}
			G_0( \bfCoh_{\scrT}(\scrC,\tau) )\quad \text{and} \quad \HBM_\ast( \bfCoh_{\scrT}(\scrC,\tau) )
		\end{align}
		have the structures of unital associative algebras, and
		
		\item the topological vector spaces
		\begin{align}
			G_0( \bfCoh_{\scrF}(\scrC,\tau) )\quad \text{and} \quad \HBM_\ast( \bfCoh_{\scrF}(\scrC,\tau) )			
		\end{align}
		have both a left and a right module structure on$G_0( \bfCoh_{\scrT}(\scrC,\tau) )$ and $\HBM_\ast( \bfCoh_{\scrT}(\scrC,\tau) )$, respectively.
	\end{itemize}
\end{theoremintroduction}

\begin{remarkintroduction}
	As for Theorem~\ref{thm:COHA-Coh-introduction}, in the main body we prove statements for more general motivic oriented homology theories.
\end{remarkintroduction}

Thanks to the left and right module structures, we are able to give the following definition.
\begin{definitionintroduction}[{Definition~\ref{def:Yangian-torsion}}]
	The \textit{categorified quantum loop algebra} $\scrH_{(\scrT, \scrF)}$ of the pair $(\scrT, \scrF)$ is the monoidal subcategory of the monoidal $\infty$-category of endofunctors $\End( \catCohb_{\mathsf{pro}}( \bfCoh_{\scrF}(\scrC,\tau) ) )$ generated by the images of the two monoidal functors
	\begin{align}
		a_\ell&\colon \catCohb_{\mathsf{pro}}( \bfCoh_{\scrT}(\scrC,\tau) )\longrightarrow\End( \catCohb_{\mathsf{pro}}( \bfCoh_{\scrF}(\scrC,\tau) ) )\ ,\\[2pt]
		a_r&\colon \catCohb_{\mathsf{pro}}( \bfCoh_{\scrT}(\scrC,\tau) )\longrightarrow\End( \catCohb_{\mathsf{pro}}( \bfCoh_{\scrF}(\scrC,\tau) ) )\ ,
	\end{align}
	corresponding to the two module structures of $\catCohb_{\mathsf{pro}}( \bfCoh_{\scrF}(\scrC,\tau) )$.
	
	The \textit{quantum loop algebra} $\scrU_{(\scrT, \scrF)}$ of the pair $(\scrT, \scrF)$ is the subalgebra of $\End( G_0( \bfCoh_{\scrF}(\scrC,\tau) ) )$ generated by the images of the two maps of associative algebras
	\begin{align}
		a_\ell&\colon G_0( \bfCoh_{\scrT}(\scrC,\tau) )\longrightarrow\End( G_0( \bfCoh_{\scrF}(\scrC,\tau) ) )\ ,\\[2pt]
		a_r&\colon G_0( \bfCoh_{\scrT}(\scrC,\tau) )\longrightarrow\End( G_0( \bfCoh_{\scrF}(\scrC,\tau) ) )\ ,
	\end{align}
	corresponding to the two module structures of $G_0( \bfCoh_{\scrF}(\scrC,\tau) )$. Similarly, we define the \textit{Yangian} $\scrY_{(\scrT, \scrF)}$ of the pair $(\scrT, \scrF)$.
\end{definitionintroduction}
The above definition can be given also for a pair $(\bfT, \bfF)$ of geometric derived stacks locally of finite presentations over a field $k$, together with maps  $\bfT\to \bfCoh_{\scrT}(\scrC, \tau)$ and $\bfF\to \bfCoh_\scrF(\scrC, \tau)$ of derived stacks, for which $\bfF$ is a \textit{left Hecke pattern for $\bfT$ with respect to the $t$-structure $\tau$} and 
\begin{itemize}\itemsep0.2cm
	\item either a \textit{right Hecke pattern for $\bfT$ with respect to the tilted $t$-structure $\tau_\upsilon$}, or
	
	\item a \textit{right Hecke pattern for $\bfT$ with respect to the tilted $t$-structure $\tau$} (in the latter case, we talk about a \textit{two-sided Hecke pattern for $\bfT$ with respect to the $t$-structure $\tau$}).
\end{itemize}
See Definition~\ref{def:Hecke-patterns} for details. Our definition generalizes, by considering pairs of $t$-structures, that in \cite[\S5.1]{KV_Hall} (see also \cite[Definition~6.1 and Remark~6.2]{MMSV}). 

Let $x$ be a class in the Borel-Moore homology (resp.\ K-theory, or $\catCohb_{\mathsf{pro}}$) of $\bfCoh_{\scrT}(\scrC,\tau)$ and denote by $a_\ell(x)$ and $a_r(x)$ the corresponding (categorical) operators induced by the left and right actions, respectively. It is important to stress that in general the commutator $[a_\ell(x), a_r(x)]$ does not vanish\footnote{This is something already known in the literature if one uses certain families of Nakajima types operators as Negu\c{t} \cite{Negut_Shuffle} in K-theory or DeHority \cite{DeHority_KM} in cohomology: we recover their results thanks to Theorem~\ref{thm_intro:geometric_action} below.}. Nevertheless, we address from this abstract point of view the question of when two classes $x, y $ induce commuting operators. We find a criterion of geometric origin (see Corollary~\ref{cor:support_commutation_criterion}), that is applied in the geometric situation described in the next section.

\subsection{Representations of surface COHAs}

Let $S$ be a smooth projective irreducible complex surface. We apply our framework in two examples. The first concerns the (standard) torsion pair $(\catCoh_{\mathsf{tor}}(S), \catCoh_{\mathsf{t.f.}}(S))$ of the standard heart of $\catCoh(S)$ formed by torsion and torsion-free sheaves, respectively. The standard $t$-structure of $\catPerf(S)$ and this torsion pair satisfy the assumptions of Theorem~\ref{thm:introduction-left-right-action}.
Moreover, the action of the COHA/KHA of $\catCoh_{\mathsf{tor}}(S)$ preserves the rank of the torsion-free sheaves. Thus, our third main result is the following:

\begin{theoremintroduction}[{cf.\ Theorem~\ref{thm:action-torsion} and Remark~\ref{rem:action-torsion}}] \label{thm_intro:geometric_action}
	The stable pro-$\infty$-category $\catCohb_{\mathsf{pro}}( \bfCoh_{\mathsf{tor}}(S) )$ has a $\mathbb E_1$-monoidal structure. Moreover, the stable pro-$\infty$-category $\catCohb_{\mathsf{pro}}( \bfCoh_{\mathsf{t.f.}}(S; r) )$ has the structure of a left and a right categorical module over $\catCohb_{\mathsf{pro}}( \bfCoh_{\mathsf{tor}}(S) )$.
	
	Similar statements hold at the level of motivic Borel-Moore homologies (including usual Borel-Moore homology and $G_0$-theory).
\end{theoremintroduction}

We show in Corollary \ref{cor:commutators_geometric_setting} the vanishing of the commutators between (categorical operators) induced by classes in Borel-Moore homologies (resp.\ K-theory, $\catCohb_{\mathsf{pro}}$) of the moduli stacks $\bfCoh(Z_1)$ and $\bfCoh(Z_2)$ corresponding to two disjoint closed subschemes $Z_1$ and $Z_2$ of $S$.

The first part of the Theorem~\ref{thm_intro:geometric_action} has been already proved by \cite{Zhao_KHA} in the K-theoretical case for zero-dimensional sheaves and in general by \cite{KV_Hall} in the cohomological and K-theoretical case, while in \cite{Porta_Sala_Hall} in the categorified case. Moreover, the left action has been already constructed in \cite{KV_Hall} in the cohomological and K-theoretical case. The two main advances of the theorem is the categorification of the left action and above all the definition of the right action. From a technical point of view, the existence of the action is guaranteed by a general restriction mechanism established in Part~\ref{part:Segal}. Furthermore, a version of Theorem~\ref{thm_intro:geometric_action} holds also by replacing $\bfCoh_{\mathsf{tor}}(S)$ with the derived stack $\bfCoh_{0\textrm{-}\!\dim}(S)$ of zero-dimensional sheaves on $S$.

We recover -- from the viewpoint of COHAs and KHAs -- both Negu\c{t}'s construction \cite{Negut_Shuffle} of the action of the elliptic Hall algebra of $S$ on the K-theory of moduli spaces of Gieseker-stable sheaves on $S$ and DeHority's construction \cite{DeHority_KM} of the action of affinizations of Lorentzian Kac--Moody algebras on the cohomology of moduli spaces of rank one Gieseker-stable sheaves on a K3 surface. Both Negu\c{t}'s and DeHority's approach are based on the use of explicit operators given by Hecke correspondences: from one side this allows them to compute explicitly the relations between them, but on the other side this forces them to consider only certain operators since it is crucial for them to understand the geometry of Hecke correspondences. This limitation does not appear if one uses directly COHAs. Thus, our result provides a new approach to extend their results by considering bigger algebras than those obtained by them.

For the second example we consider, we first tilt the standard $t$-structure of $\catCoh(S)$ by the torsion pair whose torsion part is given by zero-dimensional sheaves. Let $\tau_\scrA$ be the corresponding tilted $t$-structure and let $\scrA$ be its heart. We call its objects \textit{perverse coherent sheaves on $S$}. We consider the torsion pair $(\scrA_{\mathsf{tor}}, \scrA_{\mathsf{t.f.}})$ of $\scrA$ whose torsion part $\scrA_{\mathsf{tor}}$ consists of rank zero complexes of $\scrA$. This torsion pair is ``dual'' to the first one via the equivalences $\D(-)\colon \scrA_{\mathsf{tor}}\xrightarrow{\sim} \catCoh_{\mathsf{tor}}(S)$ and $\D(-)[-1]\colon \scrA_{\mathsf{t.f.}} \xrightarrow{\sim} \catCoh_{\mathsf{t.f.}}(S)$. Thus, we obtain:

\begin{theoremintroduction}[{Theorem~\ref{thm:swap}}]
	\hfill
	\begin{enumerate}[leftmargin=0.6cm]
		\item The equivalence
		\begin{align}	
			\D(-) \colon \bfCoh_{\mathsf{tor}}(S, \tau_\scrA\op)\xrightarrow{\sim} \bfCoh_{\mathsf{tor}}(S)
		\end{align}
		induces an equivalence of $\mathbb E_1$-monoidal stable pro-$\infty$-categories:
		\begin{align}
			\Gamma \colon \catCohb_{\mathsf{pro}}( \bfCoh_{\mathsf{tor}}(S, \tau_\scrA\op) )\xrightarrow{\sim}\catCohb_{\mathsf{pro}}( \bfCoh_{\mathsf{tor}}(S) )\ .
		\end{align}
		
		\item The equivalences
		\begin{align}
			\D(-)[-1]\colon  \bfCoh_{\mathsf{t.f.}}(S, \tau_\scrA\op; r) \xrightarrow{\sim} \bfCoh_{\mathsf{t.f.}}(S; r)\ \ \text{and}\ \  \D(-) \colon \bfCoh_{\mathsf{tor}}(S, \tau_\scrA\op)\xrightarrow{\sim} \bfCoh_{\mathsf{tor}}(S)
		\end{align}
		induce an equivalence of left and right categorical modules over $\catCohb_{\mathsf{pro}}( \bfCoh_{\mathsf{tor}}(S) )$:
		\begin{align}
			\Psi\colon \catCohb_{\mathsf{pro}}( \bfCoh_{\mathsf{t.f.}}(S, \tau_\scrA\op; r) ) \xrightarrow{\sim} \catCohb_{\mathsf{pro}}( \bfCoh_{\mathsf{t.f.}}(S; r) )\ .
		\end{align}
	\end{enumerate}	
		
	Similar statements hold at the level of motivic Borel-Moore homology (including usual Borel-Moore homology and $G_0$-theory).	
\end{theoremintroduction}

We use the $t$-structure $\tau_\scrA$ to construct left and right representations of the categorified and cohomological Hall algebras of $\bfCoh_{0\textrm{-}\!\dim}(S)$ via \textit{Pandharipande-Thomas stable pairs} on $S$ \cite{PT_BPS, PT_Curve}. These objects are also referred to as \textit{$t$-stable Bradlow pairs} on $S$, for $t\in \R$ large enough, in \cite[\S5.1]{BLM_Virasoro}. Recall that a stable pair is a pair consisting of a pure one-dimensional sheaf $\calF$ on $S$ and a section $s\colon \scrO_S \to \calF$ with zero-dimensional cokernel. It is easy to see that this datum is equivalent to that of an exact triangle $\scrO_S\to \calF\to E$ where we ask $\calF[1]\in \scrA_{\mathsf{tor}}$ and $E\in \scrA_{\mathsf{t.f.}}$ (see Proposition~\ref{prop:stable_pairs_different_formulations}). While the action of the previous theorem is essentially induced, at the level of 2-Segal spaces, by a simultaneous restriction of the action induced by the multiplication of $\bfCoh(S,\tau_\scrA)$ and of its tilting, in the case of stable pairs the fundamental mechanism underlying the action is slightly more involved: we obtain it as a restriction of an abstract action in correspondence of the derived stack of perfect complexes $\bfPerf(S)$ on the derived stack $\bfFlagPerf^{(2)}(S)$ of flags of length $2$ (see \S\ref{sec:flag_action} for the precise definition and construction of this action). Our third main result reads:

\begin{theoremintroduction}[{cf.\ Corollaries~\ref{cor:action-zero-stable-pairs} and \ref{cor:action-zero-PT-stable-pairs}}]
	Let $\bfP(S)$ be the derived moduli stack of Pandharipande-Thomas stable pairs and by $\calP(S)$ its classical truncation. Then, the pro-$\infty$-category $\catCohb_{\mathsf{pro}}( \bfP(S) )$ has the structure of a left and right categorical module over the $\mathbb E_1$-monoidal $\infty$-category $\catCohb_{\mathsf{pro}}( \bfCoh_{0\textrm{-}\!\dim}(S) )$.
	
	A similar result holds at the level of motivic Borel-Moore homology. In particular,
	\begin{align}
		G_0( \calP(S) )\quad \text{and} \quad \HBM_\ast( \calP(S) )			
	\end{align}
	are left and right modules of $G_0( \bfCoh_{0\textrm{-}\!\dim}(S) )$ and $\HBM_\ast( \bfCoh_{0\textrm{-}\!\dim}(S) )$, respectively.	
\end{theoremintroduction}

Note that in the local surface case, Toda constructed a right categorical module structure on $\catCohb_{\mathsf{pro}}$ of Pandharipande-Thomas moduli spaces of stable pairs over the categorical Hall algebra of zero-dimensional sheaves (cf.\ \cite[\S4]{Toda_Hall_Categorical_DT}). In this case, there is no left categorical module structure because of a wall-crossing phenomenon which does not appear in our two-dimensional case.

As shown in \cite[Appendix~B]{PT_BPS}, moduli spaces of stable pairs are ``dual'' to relative Hilbert schemes of points. More precisely, let $\sfH_{1\textrm{-}\mathsf{pure}}(S)$ be the Hilbert scheme parametrizing pure one-dimensional subschemes $C \subset S$. Let $\calC\subset S\times \sfH_{1\textrm{-}\mathsf{pure}}(S)$ be the universal curve and consider the (underived) relative \textit{Hilbert scheme} $\sfHilb(\calC/\sfH_{1\textrm{-}\mathsf{pure}}(S))$ of $\sfH_{1\textrm{-}\mathsf{pure}}(S)$-flat families of zero-dimensional quotients of $\scrO_{\calC}$. Let $\bfP_\scrB(S)$ be the derived enhancement of $\sfHilb(\calC/\sfH_{1\textrm{-}\mathsf{pure}}(S))$ introduced in \S\ref{subsec:duality_and_representations} (to be more precise, see Remark~\ref{rem:relative_Hilbert}). 
By using the duality, we are able to prove:
\begin{theoremintroduction}[{cf.\ Theorem~\ref{thm:hilbert-action} and Corollary~\ref{cor:hilbert-action}}]
	The pro-$\infty$-category $\catCohb_{\mathsf{pro}}( \bfP_\scrB(S) )$ has the structure of a left and right categorical module over the $\mathbb E_1$-monoidal $\infty$-category $\catCohb_{\mathsf{pro}}( \bfCoh_{0\textrm{-}\!\dim}(S) )$.
	
	A similar result holds at the level of motivic Borel-Moore homology. In particular,
	\begin{align}
		G_0( \sfHilb(\calC/\sfH_{1\textrm{-}\mathsf{pure}}(S) )\quad \text{and} \quad \HBM_\ast( \sfHilb(\calC/\sfH_{1\textrm{-}\mathsf{pure}}(S) )			
	\end{align}
	are left and right modules of $G_0( \bfCoh_{0\textrm{-}\!\dim}(S) )$ and $\HBM_\ast( \bfCoh_{0\textrm{-}\!\dim}(S) )$, respectively.	
\end{theoremintroduction}
One can wonder if the above result holds also for $\bfCoh_{\mathsf{tor}}(S)$. We are able to lift only the left action at the level of stable pairs and only the right action at the level of relative Hilbert schemes, while for the other action, we found a \textit{no-go result} coming from a geometric constraint (cf.\ Corollary~\ref{cor:relativecotangentcomplex-stable-pairs}). 

Finally, the results above can be extended to more general stable pairs, for which $\scrO_S$ is replaced by any locally free sheaf $\calV$ of finite rank.

\subsection{COHAs and Yangians of quivers}

Our framework applies also to the quiver case: we recover the known construction of the action of the two-dimensional cohomological Hall algebra of a quiver on the cohomology of Nakajima quiver varieties associated to the same quiver, but also we construct new actions on the cohomology of other quiver varieties. Moreover, we provide a categorification of these constructions.

\begin{remarkintroduction}
	Note that in the present paper, we consider only two-dimensional COHAs of quivers. The construction of representations of one- and three-dimensional Kontsevich--Soibelman COHAs of quivers is discussed in \cite{Soibelman_COHA} (see also \cite{Franzen_COHA}). A categorification of such representations may be achieved via matrix factorizations; see \cite{Padurariu_MF} for a categorification of Kontsevich--Soibelman COHAs of quivers using this approach.
\end{remarkintroduction}

Let $\calQ$ be a quiver, let $w\in \N^{\calQ_0}$, and let $\calQ^w$ be its corresponding \textit{Crawley-Boevey} quiver (cf.\ Definition~\ref{def:CB-quiver}). Denote by $\Pi_{\calQ^w}$ the \textit{preprojective algebra} of $\calQ$ and by $\Pi_{\calQ^w}$ the \textit{derived} preprojective algebra of $\calQ$ (see Definition~\ref{def:deformed-preprojective-algebra}). 

Set $\scrC_{\calQ^w}\coloneqq \Pi_{\calQ^w}\Mod$. Then, the heart of the standard $t$-structure of $\scrC_{\calQ^w}$ is the abelian category $\Modd(\Pi_{\calQ^w})$ of representations of $\Pi_{\calQ^w}$. Let $\scrT\coloneqq \Modd(\Pi_\calQ)$ be the category of representations of $\Pi_\calQ$: it can be canonically realized as a full subcategory of $\Modd(\Pi_{\calQ^w})$, which is a torsion part of a torsion pair\footnote{We warmly thank Olivier Schiffmann for suggesting us to consider this torsion pair on the quiver setting.} $(\scrT, \scrF\coloneqq \scrT^\perp)$ of $\Modd(\Pi_{\calQ^w})$, which is open. The finite-dimensional representations belonging to $\scrF$ are exactly those who are $\infty$-co-generated in the sense of \cite[Page~261]{CB_Moment}. The corresponding moduli stack $\bfCohps^{\mathsf{s}}(\scrC_{\calQ^w}, \tau_{\mathsf{std}})$ is an open substack of $\bfCoh_{\scrF}(\scrC_{\calQ^w}, \tau_{\mathsf{std}})$. 

The moduli stack $\bfCohps(\scrC_{\calQ^w}, \tau_{\mathsf{std}})$ decomposes with respect to the dimension of the finite-dimensional representations of $\Pi_{\calQ^w}$ into open and closed substacks. In particular, $\bfCohps(\Pi_\calQ)\coloneqq \bfCoh_\scrT(\scrC_{\calQ^w}, \tau_{\mathsf{std}})$ is the open and closed substack of $\bfCohps(\scrC_{\calQ^w}, \tau_{\mathsf{std}})$ defined by the condition that the vector space at the vertex $\infty$ is zero. Denote by $\bfCohps^{\mathsf{s}}(\Pi_{\calQ^w}; w_\infty)$ the moduli stack of finite-dimensional representations of $\Pi_{\calQ^w}$ belonging to $\scrF$, for which the dimension of the vector space at the vertex $\infty$ is $w_\infty$. Note that the classical truncation of $\bfCohps^{\mathsf{s}}(\Pi_{\calQ^w}; 1)$ admits a fine moduli space, which is the Nakajima quiver variety $\calM_{\calQ, \theta}(w)$ of $\theta$-stable framed representations with framed dimension vector $w$. Here, $\theta\coloneqq (1, \ldots, 1)$.

By applying the framework described in the previous section, we get the following.
\begin{theoremintroduction}[{Theorems~\ref{thm:action-quiver} and \ref{thm:action-Nakajima-quiver-variety}}]
	The stable pro-$\infty$-category $\catCohb_{\mathsf{pro}}( \bfRep(\Pi_\calQ) )$ has a $\E_1$-monoidal structure. Moreover, the stable $\infty$-pro-category $\catCohb_{\mathsf{pro}}( \bfRep_{\scrF}(\Pi_{\calQ^w})  )$ has the structure of a left categorical module over $\catCohb_{\mathsf{pro}}( \bfRep(\Pi_\calQ) )$. 
	
	Similarly, the topological vector spaces
	\begin{align}
		G_0( \bfRep(\Pi_\calQ) )\quad \text{and} \quad \HBM_\ast( \bfRep(\Pi_\calQ) )
	\end{align}
	have the structures of unital associative algebras, and the topological vector spaces
	\begin{align}
		G_0( \bfRep_{\scrF}(\Pi_{\calQ^w}) )\quad \text{and} \quad \HBM_\ast( \bfRep_{\scrF}(\Pi_{\calQ^w}) )			
	\end{align}
	have the structures of a left $G_0( \bfRep(\Pi_\calQ) )$-module and $\HBM_\ast( \bfRep(\Pi_\calQ) )$-module, respectively.
	
	For each dimension vector $w_\infty$, the same result holds for $\bfRep_{\scrF}(\Pi_{\calQ^w}; w_\infty)$ instead of $\bfRep_{\scrF}(\Pi_{\calQ^w})$. Moreover, the same results hold equivariantly with respect to the torus $T$ introduced in \cite[\S3.3]{SV_generators}.
	
	Fix $w_\infty=1$. The stable $\infty$-pro-category 
	\begin{align}
		\bigoplus_{d\in \Z} \catCohb_{\mathsf{pro}}( \calM_{\calQ, \theta}(v, w) )_d
	\end{align}
	has the structure of a left categorical module over $\catCohb_{\mathsf{pro}}( \bfRep(\Pi_\calQ) )$. Here, $\catCohb_{\mathsf{pro}}( \calM_{\calQ, \theta}(v, w) )_d$ is the weight $d$ part of $\catCohb_{\mathsf{pro}}( \bfRep_{\scrF}(\Pi_{\calQ^w}; v, 1) )$.
	
	A similar statement holds at the level of motivic Borel-Moore homology and after replacing $\bfRep(\Pi_\calQ)$ with either $\Lambda^0_\calQ$ or $\Lambda^1_\calQ$. Furthermore, the same results hold equivariantly with respect to the torus $T$ introduced in \cite[\S3.3]{SV_generators}.
\end{theoremintroduction}
The theorem recovers the constructions of the preprojective COHA and KHA of a quiver \cite{SV_Cherednik, SV_Elliptic, SV_generators, YZ_COHA} and its categorification \cite{VV_KHA}, and the left action of the preprojective COHA of a quiver on the cohomology of Nakajima quiver varieties (cf.\ \cite{SV_Cherednik, SV_generators}).

\begin{warning*}
	The right action is not obtained since the second condition of assumption (2) of Theorem~\ref{thm:introduction-left-right-action} is not satisfied by the tilted heart. By extending the above construction to the nilpotent preprojective COHA of $\calQ$, one should obtain left and right actions, as e.g. done in \cite{SV_Cherednik}.
\end{warning*}

\subsection{Further directions}

The constructions presented in this paper open new directions for investigation, some of which are currently being pursued by us, while others are left for future work. We highlight three such directions below.

\subsubsection*{COHAs and CatHAs over the space of Bridgeland stability conditions}

Corollary~\ref{cor:COHA-stability-condition-introduction} enables the definition of COHAs and CatHAs for a fixed stability condition on a 2-Calabi–Yau category (such as the category of perfect complexes on a K3 surface or a Kuznetsov component). A natural next step is to define sheaves of associative algebras and sheaves of $\E_1$-monoidal categories over the space of stability conditions and to study their geometric properties. Similarly, one could investigate how quasi-BPS categories, in the sense of Pădurariu and Toda, vary as Bridgeland stability conditions vary. These questions will be addressed in future work.

\subsubsection{Hecke operators for one-dimensional sheaves on a smooth surface}

In \cite{GKV_Langlands}, the authors studied a generalization of Hecke operators from vector bundles on curves to surfaces, introducing Hecke operators for one-dimensional sheaves acting on spaces of functions on moduli stacks of vector bundles on a smooth surface over a finite field, satisfying certain conditions. Thanks to the framework developed in Part~\ref{part:foundations}, the authors, in collaboration with Y.~Zhao, are able to extend this construction to the $K$-theory of a suitable generalization of these moduli stacks defined over $\C$. This is part of the ongoing work \cite{DPSZ_GKV}.

\subsubsection{Variations of the stability parameter for Bradlow pairs, $P=C$ conjecture, and COHAs}

As explained in~\cite[\S5.1]{BLM_Virasoro}, the semistability condition for Bradlow pairs on a smooth projective complex surface~$S$ depends on a real parameter $t > 0$. For sufficiently large~$t$, one recovers Pandharipande–Thomas stable pairs on~$S$, while for very small~$t$, the semistability of a Bradlow pair $(\calF, s\colon \scrO_S \to \calF)$ becomes equivalent to the semistability of~$\calF$ (cf.\ \cite[Proposition5.4]{BLM_Virasoro}). It would be very interesting to explore the construction of representations of the COHA and CatHA of zero-dimensional sheaves on $S$ via moduli stacks of semistable Bradlow pairs for varying values of~$t$. In particular, it would be worthwhile to investigate the algebraic structure that emerges in the regime of very small~$t$. This framework could potentially connect to a representation-theoretic approach to the \textit{P = C} conjecture for moduli spaces of semistable one-dimensional sheaves on Del Pezzo surfaces~\cite{KPS_PC}, in a spirit akin to the proof of the \textit{P = W} conjecture in~\cite{HMMS}. This will be investigated in the future.

\subsection*{Outline}

The paper is divided into three parts.

Part~\ref{part:Segal} provides an overview of the technical machinery of 2-Segal spaces and their representations. Part~\ref{part:foundations} develops a general framework for constructing COHAs, CatHAs, and their representations using the formalism of 2-Segal spaces and representations. In particular, in \S\ref{sec:COHA-t-structure}, we define COHAs and CatHAs associated to a finite-type stable $\infty$-category~$\scrC$ equipped with a $t$-structure~$\tau$ satisfying certain assumptions. Fixing a torsion pair $(\scrT, \scrF)$ in the heart of $\scrC^\heartsuit$, \S\ref{sec:COHA_representations_torsion_pairs} applies this framework to construct a COHA and a CatHA associated to~$\scrT$, along with a representation associated to~$\scrF$. Finally, in Part~\ref{part:application}, we apply this general framework to explicit geometric examples. We consider two torsion pairs naturally associated to a smooth projective complex surface~$S$. The first torsion pair has as its torsion part the usual torsion sheaves on~$S$, and the corresponding COHA, CatHA, and their representation are studied in~\S\ref{sec:action-torsion}. The second torsion pair has as its torsion part the zero-dimensional sheaves on~$S$; the associated COHA, CatHA, and representation are investigated in~\S\ref{sec:action-torsion-perverse}. These two torsion pairs are ``dual'' to each other, and their relation to Hall algebras is also explored in the same section. \S\ref{sec:action-stable-pairs} is devoted to the construction of representations of COHAs and CatHAs via stable pairs. Finally, in~\S\ref{sec:COHA-quiver-Yangian}, we apply our construction in the quiver setting, recovering the left action of the two-dimensional COHA of a quiver on the cohomology of Nakajima quiver varieties as in~\cite{SV_generators}, and we provide a categorification of this action.

\subsection*{Notation}

For a smooth projective complex surface $S$, let $K_0(S)$ be the \textit{Grothendieck group} of $S$ and let $N(S)$ be the \textit{numerical Grothendieck group} of $S$, where the latter is defined by:
\begin{align}
	N(S)\coloneqq K_0(S)/\equiv \ .
\end{align}
Here, $F_1, F_2 \in K_0(S)$ satisfy $F_1 \equiv F_2$ if $\ch(F_1) = \ch(F_2)$. Then $N(S)$ is a finitely generated free abelian group.

We denote by $\mathsf{NS}(S)$ the \textit{Neron-Severi group} of $S$. For a coherent sheaf $E$ on $S$ whose support has dimension less than or equal to one, we denote by $\ell(E)\in \mathsf{NS}(S)$ the fundamental one cycle of $E$.

Let $\catCoh_{\leqslant 1}(S) \subset \catCoh(S)$ be the subcategory of sheaves $E$ with $\dim \mathsf{Supp}(E) \leq 1$. We define the subgroup $N_{\leqslant 1}(S) \subset N(S)$ to be
\begin{align}
	N_{\leqslant 1}(S)\coloneqq \mathsf{Im}(K_0(\catCoh_{\leqslant 1}(S)) \longrightarrow N(S) ) \ .
\end{align}
Note that we have an isomorphism $N_{\leqslant 1}(S) \simeq \mathsf{NS}(S)\oplus \Z$ sending $E$ to the pair $(\ell(E), \chi(E))$. We shall identify an element $v\in N_{\leqslant 1}(S)$ with $(\beta, n) \in \mathsf{NS}(S) \oplus \Z$ by the above isomorphism.

For any $E\in \catPerf(S)$, we set $E^\vee \coloneqq \R\calH om_{\catDb(S)}(E, \scrO_S)$; while for any coherent sheaf $\calE$ on $S$, we set $\calE^\ast \coloneqq \calH om_S(\calE, \scrO_S) \simeq \calH^0(\calE^\vee)$, which is the usual \textit{dual sheaf}.

We use the attribute ``geometric'' instead of ``algebraic'' or ``Artin'' for a (derived) stack, following the convention in \cite{Porta_Sala_Hall}.

\subsection*{Acknowledgments}

The paper was developed during the Research in Pair ``2227p: Representation theory of two-dimensional categorical Hall Algebras of curves'' at the Mathematisches Forschungsinstitut Oberwolfach by the last two-named authors and during a research visit of the second-named author at Department of Mathematics of the University of Pisa under the 2021 Visiting Fellow program of the University of Pisa. We thank both institutions for providing us an exceptional research environment. 

We would like to thank Olivier Schiffmann for various stimulating discussions and Matthew B. Young for explaining his paper \cite{Young_Segal} to us. Finally, we thank the anonymous referees for their comments and suggestions.

The last-named author acknowledges the MIUR Excellence Department Project awarded to the Department of Mathematics, University of Pisa, CUP I57G22000700001. Finally, he is a member of GNSAGA of INDAM.


\makeatletter
\@addtoreset{section}{part}
\numberwithin{section}{part}
\makeatother

\newpage
\part{The combinatorics of the Hall product}\label{part:Segal}

There are two conceptual ingredients behind the existence of Hall algebras: the first one is the notion of $2$-Segal object \cite{DK,GKT_2-Segal}, which constructs the Hall product at a \textit{nonlinear} level (i.e., as an algebra in correspondences in derived stacks); the second one is a good theory of homological realizations, that allows to linearize the result produced by the $2$-Segal formalism. 

\medskip

In this part we describe the theory of $2$-Segal spaces and their representations, collecting some scattered literature around the $2$-Segal condition that appeared in the decade that followed the introduction of this notion in \cite{DK, GKT_2-Segal}.

\section{The $2$-Segal condition}\label{sec:2_Segal}

Let $\scrC$ be a presentable $\infty$-category. Via the Yoneda embedding, we have a canonical equivalence
\begin{align}
	\Fun(\bfDelta\op, \scrC) \simeq \FunR(\PSh(\bfDelta)\op, \scrC) \ . 
\end{align}
In particular, given any simplicial presheaf $K \in \PSh(\bfDelta)$ and any $F \in \Fun(\bfDelta\op, \scrC)$, the notation
\begin{align}
	F(K) \coloneqq \lim_{([n],\alpha)\in\bfDelta_{/K}} F([n]) 
\end{align}
is well defined. If $\scrC$ is only assumed to have finite limits, the same notation is well defined provided that we restrict to finite simplicial presheaves $K$.

\begin{notation}
	Let $I$ be a totally ordered finite poset (we can always identify $I$ with $[n]$ where $n \coloneqq \vert I\vert-1$, but it is useful to not fix labels for the elements of $I$). We denote by $\Delta^I$ the corresponding representable simplicial set. A subset $J \subset I$ inherits a total order and therefore defines a subsimplicial set $\Delta^J \subset \Delta^I$. More generally, given a collection of subsets $\calT$ of $I$, we denote by $\Delta^\calT$ the simplicial set
	\begin{align}
		\Delta^\calT \coloneqq \bigcup_{J \in \calT} \Delta^J \subset \Delta^I \ . 
	\end{align}
	Often, we will commit an abuse of notation and write $\calT$ instead of $\Delta^\calT$.
\end{notation}

\begin{example}\label{eg:1_dimensional_subdivision}
	Consider the subdivision of an interval
	\begin{align}
		\dynkin[Coxeter, labels={1, 2, n, n+1}, edge length=1cm] A{}
	\end{align}
	If the subdivision consists of $n$ segments, this singles out a collection of subsets of $[n]$, that we denote
	\begin{align}
		\mathscr{J}_n \coloneqq \big\{ \{0,1\}, \{1,2\}, \{2,3\} , \cdots, \{n-1,n\} \big\} \ . 
	\end{align}
	With the above convention, given a simplicial object $F \colon \bfDelta\op \to \calC$ where $\calC$ has finite limits, we obtain
	\begin{align}
		F(\mathscr{J}_n) \simeq F(\{0,1\}) \times_{F(\{1\})} F(\{1,2\}) \times_{F(\{2\})} \cdots \times_{F(\{n-1\})} F(\{n-1,n\}) \ . 
	\end{align}
\end{example}

\begin{example}
	Let $P$ be a $n$-gon in $\R^2$, that is the convex envelope of $n$ points in general position. Any subdivision of $P$ into triangles gives rise to a collection of subsets $\calT$ of $[n]$; namely a triangle appearing in the chosen subdivision is uniquely determined by three vertices of $P$, which in turn determine a subset of $[n]$. For instance, consider the following triangulation of the regular pentagon:
	\begin{center}
		\begin{tikzpicture}
			\node[draw=none,minimum size=2cm,regular polygon,regular polygon sides=5] (a) {};
			
			\foreach \n [count=\nu from 0, remember=\n as \lastn, evaluate={\nu+\lastn}] in {5,1,2,3,4} 
			{
				\node[circle,inner sep=0,minimum size=2pt,fill=black,label=\nu*(360/4):$\nu$] at (a.corner \n) {} ;
			}
			
			\draw (a.corner 1) -- (a.corner 2) ;
			\draw (a.corner 2) -- (a.corner 3) ;
			\draw (a.corner 3) -- (a.corner 4) ;
			\draw (a.corner 4) -- (a.corner 5) ;
			\draw (a.corner 5) -- (a.corner 1) ;
			\draw (a.corner 1) -- (a.corner 3) ;
			\draw (a.corner 3) -- (a.corner 5) ;
			
			\draw[->] (-2.5,0) -- (-1.5,0) ; 
			
			\draw (-4,0)++(162:5pt) +(90:1cm) -- +(162:1cm) -- +(234:1cm) -- cycle ;
			\draw (-4,0)++(306:2pt) +(234:1cm) -- +(306:1cm) -- +(18:1cm) -- cycle ;
			\draw (-4,0)++(72:5pt) +(18:1cm) -- +(90:1cm) -- +(234:1cm) -- cycle ;
		\end{tikzpicture}
	\end{center}
	Each piece determines an isomorphic copy of the finite poset $[2] \in \bfDelta\op$, and these copies are glued along the common edges to produce the pentagon. This yields a well specified map
	\begin{align}
		F([5]) \longrightarrow F([2]) \times_{F([1])} F([2]) \times_{F([1])} F([2]) \ , 
	\end{align}
	which can be written in an unambiguous way as
	\begin{align}
		F([5]) \longrightarrow F(\{0,3,4\}) \times_{F(\{0,3\})} F(\{0,1,3\}) \times_{F(\{1,3\})} F(\{1,2,3\}) \ . 
	\end{align}
	If we let
	\begin{align}
		\calT \coloneqq \big\{\{0,3,4\}, \{0,1,3\}, \{1,2,3\}\big\} \ , 
	\end{align}
	then the above map is canonically identified with the restriction
	\begin{align}
		F([5]) \longrightarrow F(\calT) \ . 
	\end{align}
\end{example}

\begin{defin}
	Let $\scrC$ be an $\infty$-category with finite limits. A functor $F \colon \bfDelta\op \to \scrC$
	\begin{enumerate}\itemsep=0.2cm
		\item is \textit{pointed} if $F([0])$ is a terminal object in $\scrC$;
		
		\item \textit{satisfies the $1$-Segal condition} if for every $n \geqslant 2$ the morphism
		\begin{align}
			F([n]) \longrightarrow F(\mathscr{J}_n) 
		\end{align}
		is an equivalence in $\scrC$;
		
		\item \textit{satisfies the $2$-Segal condition} if for every $n \geqslant 3$ and every triangulation $\calT$ of a labeled $n$-agon in $\R^2$, the induced map
		\begin{align}
			F([n]) \longrightarrow F(\calT) 
		\end{align}
		is an equivalence in $\scrC$.
	\end{enumerate}
	For $k \in \{1,2\}$ we let $k\textrm{-}\sfSegal(\scrC)$ (resp.\ $k\textrm{-} \sfSegal_\ast(\scrC))$ be the full subcategory of $\Fun(\bfDelta\op, \scrC)$ spanned by functors satisfying the $k$-Segal condition (resp.\ pointed functors satisfying the $k$-Segal condition). \hfill $\oslash$
\end{defin}

\begin{defin}\label{def:Segal_morphisms}
	Let $\scrC$ be an $\infty$-category with finite limits. A morphism $f \colon F \to G$ in $\Fun(\bfDelta\op,\scrC)$ is said to be
	\begin{enumerate}\itemsep=0.2cm
		\item \textit{pointed} if $F([0]) \to G([0])$ is an equivalence in $\scrC$;
		
		\item \textit{$1$-Segal} if for every $n \geqslant 2$ the canonical map
		\begin{align}
			F([n]) \longrightarrow F(\mathscr{J}_n) \times_{G(\mathscr{J}_n)} G([n]) 
		\end{align}
		is an equivalence;
		
		\item \textit{$2$-Segal} if for every $n \geqslant 3$ and every triangulation $\calT$ of a labeled $n$-agon, the canonical map
		\begin{align}
			F([n]) \longrightarrow F(\calT) \times_{G(\calT)} G([n]) 
		\end{align}
		is an equivalence. \hfill $\oslash$
	\end{enumerate}
\end{defin}

\begin{proposition}\label{prop:relative_Segal_condition}
	Let $\scrC$ be an $\infty$-category with finite limits and let $k \in \{1,2\}$.
	Then:
	\begin{enumerate}\itemsep=0.2cm
		\item $k$-Segal morphisms are stable under compositions;
		\item a simplicial object $F$ is $k$-Segal if and only if the map to the terminal simplicial object $F \to \ast$ is $k$-Segal;
		\item if a morphism $f \colon F \to G$ is $1$-Segal, then it is $2$-Segal as well.
	\end{enumerate}
	In particular, any $1$-Segal object is also $2$-Segal.
\end{proposition}

\begin{proof}
	See \cite[Proposition~2.3.3]{DK} and \cite[Proposition~2.4]{Young_Segal}.
\end{proof}

\begin{definition}
	Let $\scrC$ be an $\infty$-category with finite limits. A correspondence $\sigma$ in $\Fun(\bfDelta\op, \scrC)$
	\begin{align}
		\begin{tikzcd}[column sep=small, ampersand replacement=\&]
			\& H \arrow{dr}{t} \arrow{dl}[swap]{s} \\
			F \& \& G
		\end{tikzcd}
	\end{align}
	is said to be a \textit{$2$-Segal correspondence} if both $F$ and $G$ are $2$-Segal, $t$ is $1$-Segal and for every $n \geqslant 0$ the square
	\begin{align}
		\begin{tikzcd}[ampersand replacement=\&]
			F([n]) \arrow{d} \& H([n]) \arrow{l} \arrow{d} \\
			F(\{0,n\}) \& H(\{0,n\}) \arrow{l}
		\end{tikzcd}
	\end{align}
	induced by $s$ is a pullback.
\end{definition}

It follows from Proposition~\ref{prop:relative_Segal_condition} that in the above situation, the simplicial object $H$ is also automatically $2$-Segal. We write $2\textrm{-}\sfSegal_\ast^{\leftrightarrow}(\scrC)$ for the wide subcategory of $\Corr(2\textrm{-}\sfSegal_\ast(\scrC))$ spanned by $2$-Segal correspondences. The following theorem has been proven in full generality in \cite{Godicke}.
\begin{theorem}[G{\"o}dicke]\label{thm:Godicke}
	Let $\scrC$ be an $\infty$-category with finite limits. There exists an equivalence of $\infty$-categories
	\begin{align}
		2\textrm{-}\sfSegal_\ast^{\leftrightarrow}(\scrC) \simeq \mathsf{Mon}_{\E_1}(\Corr^\times(\scrC)) \ . 
	\end{align}
\end{theorem}

\begin{proof}
	This is exactly the statement of \cite[Theorem 1.1]{Godicke}, see also Corollary 5.3 in \textit{loc.\ cit.}
\end{proof}

\begin{corollary}\label{cor:from_2_Segal_to_monoids}
	Write $2\textrm{-}\sfSegal_\ast^{(1)}(\scrC)$ for the wide subcategory of $2\textrm{-}\sfSegal_\ast(\scrC)$ spanned by $1$-Segal morphisms. Then the equivalence of Theorem~\ref{thm:Godicke} induces a well defined functor
	\begin{align}
		2\textrm{-}\sfSegal_\ast^{(1)}(\scrC) \longrightarrow \mathsf{Mon}_{\E_1}(\Corr^\times(\scrC)) \ . 
	\end{align}
\end{corollary}

\begin{proof}
	Simply notice that $2\textrm{-}\sfSegal_\ast^{(1)}(\scrC)$ can be embedded as a (non full) subcategory of $2\textrm{-}\sfSegal^{\leftrightarrow}_\ast(\scrC)$.
\end{proof}

The main example of a $2$-Segal object which is not $1$-Segal is given by the Waldhausen construction of a stable $\infty$-category.
Since it is the main responsible for the associativity of the Hall product, let us briefly recall it here.

\begin{notation}
	For every positive integer $n \geqslant 0$, we let
	\begin{align}
		\sfT_n \coloneqq \Fun([1], [n]) \ . 
	\end{align}
	A morphism $[1] \to [n]$ is tantamount to give two integers $(i,j)$ satisfying $0 \leqslant i < j \leqslant n$. Moreover $\sfT_n$ naturally inherits the structure of a (not totally ordered) poset in a natural way, and it can be represented as follows:
	\begin{align}
		\begin{tikzcd}[ampersand replacement=\&]
			(0,0) \arrow{r} \& (0,1) \arrow{r} \arrow{d} \& (0,2) \arrow{r} \arrow{d} \& \cdots \arrow{r} \& (0,n) \arrow{d} \\
			\& (1,1) \arrow{r} \& (1,2) \arrow{r} \arrow{d} \& \cdots \arrow{r} \& (1,n) \arrow{d} \\
			\& \& (2,2) \arrow{r} \& \cdots \arrow{r} \& (2,n) \arrow{d} \\
			\& \& \& \ddots \& \vdots \arrow{d} \\
			\& \& \& \& (n,n)
		\end{tikzcd} \ .
	\end{align}
	This construction depends cosimplicially in $n$, i.e.\ it gives rise to a functor
	\begin{align}
		\sfT \colon \mathbf \Delta \longrightarrow \Cat \ . 
	\end{align}
\end{notation}

\begin{defin}[Waldhausen construction]\label{def:Waldhausen}
	Let $\calE$ be a stable $\infty$-category.
	We define $\calS_n \calE$ as the full subcategory of $\Fun(\sfT_n, \calE)$ spanned by those functors $F$ satisfying the following pairs of assumptions:
	\begin{enumerate}\itemsep=0.2cm
		\item for every $0 \leqslant i \leqslant n$, $F(i,i) \simeq 0$;
		
		\item for every pair of positive integers $0 \leqslant i < j < n$ the square
		\begin{align}
			\begin{tikzcd}[ampersand replacement=\&]
				F(i,j) \arrow{r} \arrow{d} \& F(i+1,j) \arrow{d} \\
				F(i,j+1) \arrow{r} \& F(i+1,j+1)
			\end{tikzcd} 
		\end{align}
		is a pullback square in $\calE$. \hfill $\oslash$
	\end{enumerate}
\end{defin}
Since $\sfT_n$ depends cosimplicially in $n$, the construction $\Fun(\sfT_n, \calE)$ depends simplicially in $n$. It is then straightforward to check that the Waldhausen construction defines a simplicial object
\begin{align}
	\calS_\bullet \calE \colon \mathbf \Delta\op \longrightarrow \Cat_\infty \ . 
\end{align}

\begin{proposition}\label{prop:Waldhausen_2_Segal}
	Let $\calE$ be a stable $\infty$-category.
	The Waldhausen construction $\calS_\bullet \calE$ is a $2$-Segal object with values in $\Cat_\infty$.
\end{proposition}

\begin{proof}
	See \cite[Theorem~7.3.3]{DK}.
\end{proof}

\section{Relative Segal objects}\label{sec:relative_Segal}

As Corollary~\ref{cor:from_2_Segal_to_monoids} shows, every $2$-Segal objects gives rise to a monoid in the category of correspondences. In this paper we will be mostly interested by modules over monoids in correspondences. This notion can be formulated similarly to the Segal condition, working with \textit{relative simplicial objects}.

\begin{notation}
	We let
	\begin{align}
		\sfrel \coloneqq \{\frakm \to \fraka\} \ . 
	\end{align}
	As an abstract category, it is equivalent to $\Delta^1$.
\end{notation}

\begin{definition}
	Let $\scrC$ be an $\infty$-category. A \textit{relative simplicial object with values in $\scrC$} is a functor $F^\star \colon \bfDelta\op \times \sfrel \to \scrC$.
\end{definition}

\begin{notation}\label{notation:components_m_a}
	Given a relative simplicial object $F^\star$, we write
	\begin{align}
		F^{\fraka} \coloneqq F \vert_{\bfDelta\op \times \{\fraka\}} \ , \qquad F^{\frakm} \coloneqq F \vert_{\bfDelta\op \times \{\frakm\}} \ . 
	\end{align}
	Compatibly, given $[n] \in \bfDelta\op$, we write
	\begin{align}
		F^{\fraka}([n]) \coloneqq F([n],\fraka) \quad \text{and} \quad F^{\frakm}([n]) \coloneqq F([n],\frakm) \ . 
	\end{align}
	We refer to $F^{\frakm}$ as the \textit{underlying module object} and to $F^{\fraka}$ as the \textit{underlying algebra object}.
\end{notation}

\begin{definition}\label{def:relative_1_Segal}
	Let $\scrC$ be an $\infty$-category with finite limits.
	\begin{enumerate}\itemsep=0.2cm
		\item A relative simplicial object $F^\star$ is said to be \textit{pointed} if $F^{\fraka}([0])$ is a terminal object in $\scrC$.
		
		\item A relative simplicial object $F^\star$ is said to be \textit{relative left $1$-Segal} if the restriction $F^{\fraka}$ is $1$-Segal and for every $n \geqslant 0$ the canonical maps make the square
		\begin{align}
			\begin{tikzcd}[ampersand replacement=\&]
				F^{\frakm}([n]) \arrow{r} \arrow{d} \& F^{\frakm}(\{n\}) \arrow{d} \\
				F^{\fraka}([n]) \arrow{r} \& F^{\fraka}(\{n\})
			\end{tikzcd} 
		\end{align}
		into a pullback in $\scrC$.
		
		\item A relative simplicial object $F^\star$ is said to be \textit{relative $2$-Segal} if $F^{\fraka}$ is $2$-Segal, $F^{\frakm}$ is $1$-Segal and for every $n \geqslant 3$ and every $0 \leqslant i < j \leqslant n$ the canonical maps make the square
		\begin{align}
			\begin{tikzcd}[ampersand replacement=\&]
				F^{\frakm}([n]) \arrow{r} \arrow{d} \& F^{\frakm}(\{0,1,\ldots,i,j,j+1,\ldots,n\}) \arrow{d} \\
				F^{\fraka}(\{i,\ldots,j\}) \arrow{r} \& F^{\fraka}(\{i,j\})
			\end{tikzcd}
		\end{align}		
		into a pullback in $\scrC$.
	\end{enumerate}
	We let $1\textrm{-}\sfSegal^{\ell\textrm{-}\sfrel}_\ast(\scrC)$ (resp.\ $2\textrm{-}\sfSegal^{\sfrel}(\scrC)$) be the subcategory of $\Fun(\bfDelta\op \times \Delta^1, \scrC)$ spanned by pointed relative left $1$-Segal (resp.\ relative $2$-Segal) objects.
\end{definition}

\begin{variant}
	There is an obvious \textit{right} variant of the relative $1$-Segal condition, where instead of using the map $[n] \to \{n\}$ in $\bfDelta\op$ one uses $[n] \to \{0\}$. We denote the corresponding $\infty$-category by $1\textrm{-}\sfSegal^{r\textrm{-}\sfrel}_\ast(\scrC)$.
\end{variant}

\begin{warning}
	The reader should be cautious with the terminology: while a morphism $f \colon F \to G$ is the same thing as a relative simplicial object (see Definition~\ref{def:relative_1_Segal}), saying that $f$ is $1$-Segal is not the same as saying that it is relative left $1$-Segal. However, if $G$ is itself $1$-Segal, then if $f_\bullet$ is relative left (or right) $1$-Segal, it follows that $X_\bullet$ is $1$-Segal (see \cite[Proposition~2.1]{Young_Segal}), and therefore that the morphism $f_\bullet$ is $1$-Segal as well (since any morphism between $1$-Segal spaces is automatically $1$-Segal).
\end{warning}

\begin{definition} \label{def:1_Segal_morphism}
	Let $\scrC$ be an $\infty$-category with finite limits. We say that a morphism $f^\star \colon F^\star \to G^\star$ of relative simplicial objects in $\scrC$ is
	\begin{enumerate}\itemsep=0.2cm
		\item \textit{relative left (resp.\ right) $1$-Segal} if the morphism $f^{\fraka} \colon F^{\fraka} \to G^{\fraka}$ is $1$-Segal and for every $n \geqslant 0$ the square
		\begin{align}
			\begin{tikzcd}[ampersand replacement=\&]
				F^{\frakm}([n]) \arrow{r} \arrow{d} \& F^{\fraka}([n]) \times_{F^{\fraka}([0])} F^{\frakm}([0]) \arrow{d} \\
				G^{\frakm}([n]) \arrow{r} \& G^{\fraka}([n]) \times_{F^{\fraka}([0])} F^{\frakm}([0])
			\end{tikzcd}
		\end{align}
		is a pullback, where the $\frakm$-component of the horizontal maps are induced by $\Delta^0 \xrightarrow{n} \Delta^n$ (resp.\  by $\Delta^0 \xrightarrow{0} \Delta^n$);
		
		\item \textit{relative $2$-Segal} if the morphism $f^{\fraka} \colon F^{\fraka} \to G^{\fraka}$ is $1$-Segal and for every $n \geqslant 1$ the square
		\begin{align}
			\begin{tikzcd}[ampersand replacement=\&]
				F^{\frakm}([n]) \arrow{r} \arrow{d} \& F^{\fraka}(\mathscr{J}_{n-1}) \times_{F^{\fraka}([0])} F^{\frakm}(\{n-1,n\}) \arrow{d} \\
				G^{\frakm}([n]) \arrow{r} \& G^{\fraka}(\mathscr{J}_{n-1}) \times_{G^{\fraka}([0])} G^{\frakm}(\{n-1,n\})
			\end{tikzcd}
		\end{align}
		is a pullback.
	\end{enumerate}
	We write $2\textrm{-}\sfSegal_\ast^{2\sfrel}(\scrC)$ for the (non full) subcategory of $2 \textrm{-} \sfSegal_\ast(\scrC)$ spanned by relative $2$-Segal morphisms.
\end{definition}

The following two propositions generalize \cite[Propositions~2.1 and 2.4]{Young_Segal} to morphisms. The proofs are straightforward generalizations of the proofs in \textit{loc.\ cit.}, and they are left to the interested reader.
\begin{proposition} \label{prop:checking_relative_1_Segal}
	Let $f^\star \colon F^\star \to G^\star$ be a morphism in $\Fun(\bfDelta\op \times \sfrel, \scrC)$. Assume that the morphism $f^{\fraka} \colon F^{\fraka} \to G^{\fraka}$ is $1$-Segal. Then $f^\star$ is relative left $1$-Segal if and only if the morphism $f^{\frakm} \colon F^{\frakm} \to G^{\frakm}$ is $1$-Segal and the square
	\begin{align}
		\begin{tikzcd}[ampersand replacement=\&]
			F^{\frakm}([1]) \arrow{r}{s_1 \times \partial^{\frakm}_1} \arrow{d} \& F^{\fraka}([1]) \times_{F^{\fraka}([0])} F^{\frakm}([0]) \arrow{d} \\
			G^{\frakm}([1]) \arrow{r}{s_1 \times \partial^{\frakm}_1} \& G^{\fraka}([1]) \times_{G^{\fraka}([0])} G^{\frakm}([0])
		\end{tikzcd} 
	\end{align}
	is a pullback.
\end{proposition}

\begin{proposition} \label{prop:1_Segal_implies_2_Segal}
	Let $f^\star \colon F^\star \to G^\star$ be a morphism in $\Fun(\bfDelta\op \times \sfrel, \scrC)$. If $f^\star$ is relative left or right $1$-Segal, then it is relative $2$-Segal.
\end{proposition}

Finally, the following theorem explains how to extract a module in correspondences out of a relative $2$-Segal object.
\begin{theorem}[G{\"o}dicke]\label{thm:Godicke_II}
	There exists a well defined $\infty$-functor
	\begin{align}
		2\textrm{-}\sfSegal_\ast^{2\sfrel}(\scrC) \longrightarrow \mathsf{LMod}(\Corr^\times(\scrC)) 
	\end{align}
	compatible with the functor $1\textrm{-}\sfSegal_\ast^{\ell\textrm{-}\sfrel}(\scrC) \to \mathsf{LMod}(\scrC)$ of \cite[Proposition 4.2.2.9]{Lurie_HA}.
\end{theorem}

\begin{proof}
	This is a special case of \cite[Corollary 5.4]{Godicke}.
\end{proof}

\section{Self and flag representations}

Any associative (unital) monoid $A$ acts on himself by multiplication on the left and on the right. It is easy to give a purely simplicial description of this phenomenon. The following construction generalizes \cite[Example~4.2.2.4]{Lurie_HA} and \cite[Proposition~2.5]{Young_Segal}, and it plays an important role in constructing the representations of COHAs and CatHAs later on in this paper.

\begin{construction}\label{construction:self_action}
	Let $m \geqslant 1$ be an integer and define the function $a_m \colon \sfrel \to \Z$ by setting
	\begin{align}
		a_m(\ast) \coloneqq \begin{cases}
			m - 1 & \text{if } \ast = \mathfrak m \ ,\\ -1 & \text{if } \ast = \mathfrak a \ .
		\end{cases}
	\end{align}
	Consider now the functors
	\begin{align}
		\tensor*[^\ell]{\mathsf{act}}{^{(m)}}, \ \tensor*[^r]{\mathsf{act}}{^{(m)}} \colon \bfDelta\op \times \sfrel \longrightarrow \bfDelta\op  
	\end{align}
	defined by
	\begin{align}
		\tensor*[^\ell]{\mathsf{act}}{^{(m)}}( [n], \ast ) \coloneqq [n] \star [a_m(\ast)] \quad \text{and} \quad \tensor*[^r]{\mathsf{act}}{^{(m)}}([n],\ast) \coloneqq [a_m(\ast)] \star [n] \ , 
	\end{align}
	where $\star$ denotes the join operation.
	
	\medskip
	
	Given an $\infty$-category $\scrC$ with finite limits and a simplicial object $A \colon \bfDelta\op \to \scrC$, we write
	\begin{align}
		\tensor*[^\ell]{A}{^{(m)}} \coloneqq A \circ \tensor*[^\ell]{\mathsf{act}}{^{(m)}} \quad \text{and} \quad \tensor*[^r]{A}{^{(m)}} \coloneqq A \circ \tensor*[^r]{\mathsf{act}}{^{(m)}} \ , 
	\end{align}
	and we refer to them as the \textit{left} and \textit{right $m$-flags action objects}. Observe that by definition we have $\tensor*[^\ell]{A}{^{(m)}}([0],0) = A([m]) = \tensor*[^r]{A}{^{(m)}}([0],0)$.
\end{construction}

\begin{proposition}\label{prop:flag_action}
	Let $m \geqslant 1$ be an integer and let $A \colon \bfDelta\op \to \scrC$ be a simplicial object. If $A$ is (unital) $1$-Segal, then $\tensor*[^\ell]{A}{^{(m)}}$ and $\tensor*[^r]{A}{^{(m)}}$ are (unital) relative left and right $1$-Segal spaces, respectively. If $A$ is (unital) $2$-Segal, then both $\tensor*[^\ell]{A}{^{(m)}}$ and $\tensor*[^r]{A}{^{(m)}}$ are (unital) relative $2$-Segal spaces.
\end{proposition}

\begin{proof}
	When $m = 1$ this coincides exactly with \cite[Example~4.2.2.4]{Lurie_HA} (in the $1$-Segal case) and with \cite[Proposition~2.5]{Young_Segal} (in the $2$-Segal case). The proof of the general case is dealt with in an identical way.
\end{proof}

\begin{example}
	Let us spell out this construction in the special case of the \textit{Waldhausen construction} of a stable $\infty$-category $\calE$ (see Definition~\ref{def:Waldhausen}).
	Set therefore $A \coloneqq \calS_\bullet \calE$. We unravel below the meaning of $\tensor*[^\ell]{A}{^{(m)}}$ for $m = 1, 2$:
	\begin{itemize}\itemsep=0.2cm
		\item When $m = 1$, we have
		\begin{align}
			\tensor*[^\ell]{A}{^{(1)}}([0]) \simeq A([1]) \simeq \calE \quad \text{and} \quad \tensor*[^\ell]{A}{^{(1)}}([1]) \simeq A([2]) \simeq \calS_2\calE
		\end{align}
		so $\calE$ is the recipient of the action, and the action is given by the correspondence
		\begin{align}
			\begin{tikzcd}[ampersand replacement=\&]
				\calS_2\calE \arrow{r}{\ev_{02}} \arrow{d}[swap]{\ev_{01} \times \ev_{12}} \& \calE \\
				\calE \times \calE \& \phantom{\calE} 
			\end{tikzcd}\ .
		\end{align}
		Pictorially, we can represent this as the extension
		\begin{align}
			\begin{tikzcd}[ampersand replacement=\&]
				0 \arrow{r} \& M_{01} \arrow{r} \arrow{d} \& M_{02} \arrow{d} \\
				\& 0 \arrow{r} \& A_{12} \arrow{d} \\
				\& \& 0
			\end{tikzcd}\ ,
		\end{align}
		where $A_{12}$ suggests that it is the ``algebra element'', while $M_{01}$ and $M_{02}$ are the ``module elements''.
		
		\item When $m = 2$, we have
		\begin{align}
			\tensor*[^\ell]{A}{^{(2)}}([0]) \simeq A([2]) \simeq \calS_2\calE \qquad \text{and} \qquad \tensor*[^\ell]{A}{^{(2)}}([1]) \simeq \calS_3\calE \ , 
		\end{align}
		so this time $\calS_2 \calE$ is the recipient of the action, and the action is given by the correspondence
		\begin{align}
			\begin{tikzcd}[ampersand replacement=\&]
				\calS_3 \calE \arrow{r}{\ev_{013}} \arrow{d}[swap]{\ev_{012} \times \ev_{23}} \& \calS_2 \calE \\
				\calS_2 \calE \times \calE \& \phantom{\calS_2 \calE}
			\end{tikzcd} \ .
		\end{align}
		More graphically, we can represent this as the flag
		\begin{align}
			\begin{tikzcd}[ampersand replacement=\&]
				0 \arrow{r} \& M_{01} \arrow{r} \arrow{d} \& M_{02} \arrow{r} \arrow{d} \& M_{03} \arrow{d} \\
				\& 0 \arrow{r} \& M_{12} \arrow{r} \arrow{d} \& M_{13} \arrow{d} \\
				\& \& 0 \arrow{r} \& A_{23} \arrow{d} \\
				\& \& \& 0
			\end{tikzcd} 
		\end{align}
		We can heuristically describe this situation saying that the algebra element $A_{23}$ acts on the extension $M_{01} \to M_{02} \to M_{12}$ producing the extension $M_{01} \to M_{03} \to M_{13}$.
	\end{itemize}
\end{example}

When $m \geqslant 2$ the above construction can be slightly refined as we shall explain below in Construction~\ref{construction:Vflags_generalized}. First, we fix the following notation:
\begin{notation}
	We let $\bfDelta_+$ be the augmented simplicial category. Its objects are possibly empty finite linearly ordered sets. Conventionally, we denote by $[-1]$ the empty poset, which becomes the initial object of $\bfDelta_+$. Moreover, for every $[n] \in \bfDelta_+$ we set
	\begin{align}
		[n] \star [-2] \coloneqq [-1] \quad \text{and} \quad [-2] \star [n] \coloneqq [-1] \ . 
	\end{align}
	Given an integer $m \geqslant 1$ we let
	\begin{align}
		\sigma_m \colon \bfDelta\op \times \sfrel \longrightarrow \bfDelta\op_+ 
	\end{align}
	by setting
	\begin{align}
		\sigma_m([n],\ast) \coloneqq [a_m(\ast)-1] \ . 
	\end{align}
	Notice that when $\ast = \mathfrak a$, $a_{m}(\mathfrak a) - 1 = -2$, and therefore by convention we have $\sigma_m([n],\mathfrak a) = [-1]$ for every $[n] \in \bfDelta\op$.
\end{notation}

\begin{remark}\label{rem:Vflags_generalized}
	Let $\scrC$ be an $\infty$-category with finite limits and let $A \colon \bfDelta\op \to \scrC$ be a simplicial object. We extend $A$ to an augmented simplicial object $A_+ \colon \bfDelta\op_+\to \scrC$ by setting $A_+([-1]) \coloneqq \bf1$, where $\bf1$ denotes the final object of $\scrC$. Then the composite $A_+ \circ \sigma_m$ is canonically identified with the morphism $A([m-1]) \to \bf1$, where both the source and the target are seen as constant simplicial objects. Since $\bf1$ is the final object of $\scrC$, it is immediate to see that $A_+ \circ \sigma_m$ is a relative (left $1$-Segal and hence) $2$-Segal space.
\end{remark}

\begin{construction}\label{construction:Vflags_generalized}
	Let $m \geqslant 2$ be an integer.
	The boundary maps
	\begin{align}
		\partial_{m-1} \colon [m-1] \longrightarrow [m-2] \quad \text{and} \quad \partial_0 \colon [m-1] \longrightarrow [m-2] 
	\end{align}
	which avoid $m-1\in [m-1]$ and $0\in [m-1]$ respectively, induce natural transformations
	\begin{align}
		\tensor*[^\ell]{\delta}{_{m-1}} \colon \tensor*[^\ell]{\mathsf{act}}{^{(m)}} \longrightarrow \sigma_{m-1} \quad \text{and} \quad \tensor*[^r]{\delta}{_{0}} \colon \tensor*[^r]{\mathsf{act}}{^{(m)}} \longrightarrow \sigma_{m-1} \ . 
	\end{align}
	Let $\scrC$ be $\infty$-category with finite limits and let $A \colon \bfDelta\op \to \scrC$ be a simplicial object.
	Denoting by $\bf1$ the final object of $\scrC$, let $\sfV \colon \mathbf1 \to A([m-2])$ be a morphism in $\scrC$.
	Reviewing $\id_{\bf1} \colon \bf1 \to \bf1$ as a constant relative simplicial object in $\scrC$, we review $\sfV$ as a morphism
	\begin{align}
		\begin{tikzcd}[ampersand replacement=\&]
			\bf1 \arrow{d}{\id_{\bf1}} \arrow{r}{\sfV} \& A([m-2]) \arrow{d} \\
			\bf1 \arrow[equal]{r} \& \bf1 
		\end{tikzcd} \ ,
	\end{align}
	where the columns are seen as constant simplicial objects in $\scrC$. Furthermore, Remark~\ref{rem:Vflags_generalized} canonically identifies the right column with $A_+ \circ \sigma_{m-1}$. We therefore review $\sfV$ as a morphism $\id_{\bf1} \to A_+ \circ \sigma_{m-1}$. At this point, we define $\tensor*[^\ell]{A}{^{(m),\dagger}_\sfV}$ and $\tensor*[^r]{A}{^{(m),\dagger}_\sfV}$ as the fiber products
	\begin{align}
		\begin{tikzcd}[ampersand replacement=\&]
			\tensor*[^\ell]{A}{^{(m),\dagger}_\sfV} \arrow{r} \arrow{d} \& \id_{\bf1} \arrow{d}{\sfV} \\
			\tensor*[^\ell]{A}{^{(m)}}\arrow{r}{\tensor*[^\ell]{\delta}{_0}} \& A_+ \circ \sigma_{m-1}
		\end{tikzcd} \quad \text{and} \quad \begin{tikzcd}[ampersand replacement=\&]
			\tensor*[^r]{A}{^{(m),\dagger}_\sfV} \arrow{r} \arrow{d} \& \id_{\bf1} \arrow{d}{\sfV} \\
			\tensor*[^r]{A}{^{(m)}} \arrow{r}{\tensor*[^r]{\delta}{_{m-1}}} \& A_+ \circ \sigma_{m-1} 
		\end{tikzcd}\ ,
	\end{align}
	taken in $\Fun(\bfDelta\op \times \sfrel, \scrC)$.
\end{construction}

\begin{remark}\label{rem:empty-case}
	The previous construction makes equally sense for $m = 1$, but since $A_+([1-2]) = A_+([-1]) = \bf1$, there is only one possible choice for the morphism $\sfV$, and the resulting objects $\tensor*[^\ell]{A}{^{(1),\dagger}_\sfV}$ and $\tensor*[^r]{A}{^{(1),\dagger}_\sfV}$ coincide with $\tensor*[^\ell]{A}{^{(1)}}$ and $\tensor*[^r]{A}{^{(1)}}$, respectively.
\end{remark}

In virtue of the previous remark, the following result subsumes Proposition~\ref{prop:flag_action}.
\begin{proposition}\label{prop:Vflags_generalized_2-Segal}
	Let $m \geqslant 2$ be an integer. Let $\scrC$ be an $\infty$-category with finite limits. Let $A \colon \bfDelta\op \to \scrC$ be a simplicial object and let $\sfV \colon \mathbf{1} \to A([m-2])$ be any morphism. If $A$ is (unital) $1$-Segal, then both $\tensor*[^\ell]{A}{^{(m),\dagger}_\sfV}$ and $\tensor*[^r]{A}{^{(m),\dagger}_\sfV}$ are (unital) relative left and right $1$-Segal spaces, respectively. If $A$ is (unital) $2$-Segal, then both $\tensor*[^\ell]{A}{^{(m),\dagger}_\sfV}$ and $\tensor*[^r]{A}{^{(m),\dagger}_\sfV}$ are (unital) $2$-Segal spaces.
\end{proposition}

\begin{proof}
	Observe that both $\id_{\bf1} \colon \bf1 \to \bf1$ and $A_+ \circ \sigma_m$ are both unital relative left and right $1$-Segal spaces. Recall moreover that the collection of (unital) relative left and right $1$-Segal and (unital) relative $2$-Segal spaces are closed under limits in $\Fun(\bfDelta\op \times \Delta^1, \scrC)$. Thus, the conclusion follows from Proposition~\ref{prop:flag_action} and the very definition of $\tensor*[^\ell]{A}{^{(m),\dagger}_\sfV}$ and $\tensor*[^r]{A}{^{(m),\dagger}_\sfV}$.
\end{proof}

\section{$1$-Coskeletality}\label{sec:1-cosk}

In practice, it is often redundant to check all the Segal conditions. Indeed, if the simplicial objects in consideration are sufficiently finite, it is sufficient to check only finitely many Segal conditions. In this subsection we make this statement precise, thanks to the notion of $1$-coskeletality.

\begin{notation}
	Let $j_1 \colon \bfDelta\op_{\leqslant 1} \hookrightarrow \bfDelta$ be the canonical inclusion.
	We write
	\begin{align}
		\mathsf{cosk}_1 \coloneqq j_{1,\ast} \circ j_1^\ast \colon \Fun(\bfDelta\op, \scrC) \longrightarrow \Fun(\bfDelta\op, \scrC) \ , 
	\end{align}
	and we refer to it as the $1$-coskeleton functor. Here $j_1^\ast$ denotes the restriction along $j_1$ and $j_{1,\ast}$ denotes the right Kan extension along $j_1$. Observe that there is a canonical natural transformation $\id \to \mathsf{cosk}_1$, corresponding to the unit of the adjunction $j_1^\ast \dashv j_{1,\ast}$.
\end{notation}

\begin{definition} \label{def:1_coskeletal_morphism}
	We say that a morphism $f \colon F \to G$ in $\Fun(\bfDelta\op, \scrC)$ is \textit{$1$-coskeletal} if the square
	\begin{align}
		\begin{tikzcd}[column sep = large, ampersand replacement=\&]
			F \arrow{r}{f} \arrow{d} \& G \arrow{d} \\
			\mathsf{cosk}_1(F) \arrow{r}{\mathsf{cosk}_1(f)} \& \mathsf{cosk}_1(G)
		\end{tikzcd} 
	\end{align}
	is a pullback.
\end{definition}

\begin{lemma}\label{lem:skeleton_coskeleton}
	Let $f \colon F \to G$ be a $1$-coskeletal morphism. Then for every morphism $K \to H$ in $\PSh(\bfDelta)$, the square
	\begin{align}
		\begin{tikzcd}[ampersand replacement=\&]
			F(K) \arrow{r} \arrow{d} \& F(K) \times_{G(K)} G(H) \arrow{d} \\
			F(\mathsf{sk}_1(H)) \arrow{r} \& F( \mathsf{sk}_1(K) ) \times_{G( \mathsf{sk}_1(K) )} G( \mathsf{sk}_1(H) )
		\end{tikzcd} 
	\end{align}
	is a pullback.
\end{lemma}

\begin{proof}
	Unraveling the definition of $1$-coskeletal morphism and using the fact that the evaluations $\ev_{K}$ and $\ev_{H}$ commute with limits, we readily reduce to the case where $G \to \mathsf{cosk}_1(G)$ is an equivalence. In this case, the map $F \to \mathsf{cosk}_1(F)$ is an equivalence as well, and therefore we are reduced to check that for every $M \in \Fun(\bfDelta\op,\scrC)$ and every $L \in \PSh(\bfDelta)$ the map
	\begin{align}
		\ev_{L}( \mathsf{cosk}_1(M) ) \longrightarrow \ev_{\mathsf{sk}_1(L)}( \mathsf{cosk}_1(M) ) 
	\end{align}
	is an equivalence. Composing with the functors $\Hom_{\scrC}(T,-)$ as $T$ varies on a set of generators for $\scrC$, we are immediately reduced to the case $\scrC = \mathsf{Spc}$, where the latter is the $\infty$-category of spaces. In this case, the conclusion follows since the adjunction $\mathsf{sk}_1 \dashv \mathsf{cosk}_1$ holds.
\end{proof}

\begin{proposition} \label{prop:1_Segal_0_coskeletal}
	Let $f \colon F \to G$ be a morphism in $\Fun(\bfDelta\op, \scrC)$. If $f$ is $0$-coskeletal, then it is $1$-Segal. Moreover, $f$ is pointed $1$-Segal if and only if $f$ is an equivalence.
\end{proposition}

\begin{proof}
	Since $f$ is $0$-coskeletal, the square
	\begin{align}
		\begin{tikzcd}[ampersand replacement=\&]
			F \arrow{r} \arrow{d}{f} \& \mathsf{cosk}_0(F) \arrow{d}{\mathsf{cosk}_0(f)} \\
			G \arrow{r} \& \mathsf{cosk}_0(G)
		\end{tikzcd} 
	\end{align}
	is a pullback.
	Now, we can canonically identify $\mathsf{cosk}_0(F)$ with the \v{C}ech nerve of $F([0]) \to \bf1$, where $\bf1$ is the final object in $\scrC$.
	Thus, the proof of \cite[Proposition~6.1.2.11]{HTT} shows that $\mathsf{cosk}_0(F)$ is $1$-Segal. For the same reason, $\mathsf{cosk}_0(G)$ is $1$-Segal, and therefore map $\mathsf{cosk}_0(f)$ is $1$-Segal as well. Since the above square is a pullback, we deduce that $f$ is $1$-Segal as well, proving the first half of the statement. For the second half, saying that $f$ is pointed is equivalent to say that the induced map $f_0 \colon X_0 \to Y_0$ is an equivalence. Since $\mathsf{cosk}_0(f)$ is the right Kan extension along $\bfDelta\op_{\leqslant 0} \to \bfDelta\op$ of the morphism $f_0$, the conclusion follows.
\end{proof}

\begin{proposition}\label{prop:1_Segal_1_coskeletal}
	Let $f \colon F \to G$ be a morphism of simplicial objects. If $f$ is $1$-coskeletal, then it is $1$-Segal if and only if the canonical map
	\begin{align}
		F([2]) \longrightarrow F(\mathscr{J}_2) \times_{G(\mathscr{J}_2)} G([2]) 
	\end{align}
	is an equivalence. Furthermore, it is pointed $1$-Segal if and only if in addition to this condition the map $X([0]) \to Y([0])$ is an equivalence.
\end{proposition}

\begin{proof}
	The pointed case is an obvious consequence of the unpointed one, so we only deal with the latter. Let $j \colon \bfDelta\op_{\leqslant 1} \hookrightarrow \bfDelta\op$ be the natural inclusion. Fix an integer $n \geqslant 1$. Since $f$ is $1$-coskeletal, Lemma~\ref{lem:skeleton_coskeleton} provides us with the following pulllback square:
	\begin{align}
		\begin{tikzcd}[ampersand replacement=\&]
			F([n]) \arrow{r} \arrow{d} \& X(\mathscr{J}_n) \times_{G( \mathscr{J}_n )} G( [n] ) \arrow{d} \\
			F( \mathsf{sk}_1(\Delta^n) ) \arrow{r} \& F( \mathscr{J}_n ) \times_{G( \mathscr{J}_n)} G( \mathsf{sk}_1(\Delta^n) ) 
		\end{tikzcd} \ .
	\end{align}	
	It is therefore enough to prove that the condition of the statement implies that the bottom row is an equivalence. We proceed by induction on $n$. When $n = 1$, $\mathsf{sk}_1(\Delta^1) = \Delta^1 = \mathscr{J}_1$, and therefore the statement is trivial. When $n = 2$, the assertion is true by assumption. Let therefore $n \geqslant 3$. For $j = n+1, n, n-1, \ldots, 0$ we define inductively simplicial sets $\mathscr{J}_n(j)$ and $\mathsf{sk}_1(\Delta^n)(j)$ together with morphisms
	\begin{align}
		\alpha_j \colon \mathscr{J}_n(j) \longrightarrow \mathsf{sk}_1(\Delta^n)(j) \quad \text{and} \quad \beta_j \colon \mathsf{sk}_1(\Delta^n)(j) \to \mathsf{sk}_1(\Delta^n)(j-1) \ , 
	\end{align}
	as follows:
	\begin{enumerate}\itemsep=0.2cm
		\item for $j = n+1$ we set $\mathsf{sk}_1(\Delta^n)(n+1) \coloneqq \mathsf{sk}_1(\Delta^{n-1})$, $\mathscr{J}_n(n+1) \coloneqq \mathscr{J}_{n-1}$ and we let $\alpha_{n+1}$ be the canonical map; furthermore we conventionally set $\mathsf{sk}_1(\Delta^n)(n+2) = \emptyset$ so that $\beta_{n+1}$ is uniquely determined;
		
		\item for $j = n$ we set $\mathscr{J}_n(n) \coloneqq \mathscr{J}_n$, and we define $\mathsf{sk}_1(\Delta^n)(n)$, $\alpha_{n}$ and $\beta_n$ via the following diagram:
		\begin{align}
			\begin{tikzcd}[ampersand replacement=\&]
				\Delta^0 \arrow{r}{n-1} \arrow{d}{0} \& \mathscr{J}_{n-1} \arrow{d} \arrow{r}{\alpha_{n+1}} \& \mathsf{sk}_1(\Delta^{n-1}) \arrow{d}{\beta_{n}} \\
				\Delta^1 \arrow{r} \& \mathscr{J}_n(n) \arrow{r}{\alpha_{n}} \& \mathsf{sk}_1(\Delta^n)(n) 
			\end{tikzcd} \ ,
		\end{align}
		where we require the square on the right to be a pushout;
		
		\item given $\mathscr{J}_n(j)$, $\mathsf{sk}_1(\Delta^n)(j)$, $\alpha_j$ and $\beta_{j+1}$ for $j \leqslant n$, we define the same data at level $j-1$ by forming the following diagram:
		\begin{align}\label{eq:attaching_cells}
			\begin{tikzcd}[ampersand replacement = \&]
				\Delta^0 \amalg \Delta^0 \arrow{r}{(j-1, n)} \arrow{d}{(0,1)} \& \mathscr{J}_n(j) \arrow{d} \arrow{r}{\alpha_j} \& \mathsf{sk}_1(\Delta^n)(j) \arrow{d}{\beta_j} \\
				\Delta^1 \arrow{r} \& \mathscr{J}_n(j-1) \arrow{r}{\alpha_{j-1}} \& \mathsf{sk}_1(\Delta^n)(j-1) 
			\end{tikzcd}\ ,
		\end{align}
		where both squares are required to be pushouts.
	\end{enumerate}
	Observe that by construction $\mathsf{sk}_1(\Delta^n)(0) = \mathsf{sk}_1(\Delta^n)$ and that the canonical map $\mathscr{J}_n \to \mathsf{sk}_1(\Delta^n)$ can be factored as the composition
	\begin{align}
		\mathscr{J}^n = \mathscr{J}^n(n) \xrightarrow{\alpha_{n}} \mathsf{sk}_1(\Delta^n)(n) \xrightarrow{\beta_{n-1}} \mathsf{sk}_1(\Delta^n)(n-1) \xrightarrow{\beta_{n-2}} \cdots 	\xrightarrow{\beta_{1}} \mathsf{sk}_1(\Delta^n)(0) = \mathsf{sk}_1(\Delta^n) \ . 
	\end{align}
	This reduces us to check that the map
	\begin{align}
		F( \mathsf{sk}_1(\Delta^n)(n) ) \longrightarrow F(\mathscr{J}^n(n)) \times_{G( \mathscr{J}_n(n))} G( \mathsf{sk}_1(\Delta^n)(n) )
	\end{align}
	induced by $\alpha_n$ and the maps
	\begin{align}
		F( \mathsf{sk}_1(\Delta^n)(j-1)) \longrightarrow F( \mathsf{sk}_1(\Delta^n)(j) ) \times_{G( \mathsf{sk}_1(\Delta^n)(j))} G( \mathsf{sk}_1(\Delta^n)(j-1) ) \ , 
	\end{align}
	induced respectively by the $\beta_j$'s are equivalences.
	
	\smallskip
	
	For $\alpha_n$, it is enough to observe that the relevant map fits in the following pullback square
	\begin{align} \label{eq:1_Segal_on_2_cell_base}
		\begin{tikzcd}[ampersand replacement = \&]
			F(\mathsf{sk}_1(\Delta^{n})(n))\arrow{r} \arrow{d} \& F( \mathscr{J}_n(n)) \times_{G(\mathscr{J}_n(n))} G( \mathsf{sk}_1(\Delta^1)(n) ) \arrow{d} \\
			F( \mathsf{sk}_1(\Delta^{n-1}) ) \arrow{r} \& G (\mathscr{J}_{n-1} ) \times_{G( \mathscr{J}_{n-1})} G( \mathsf{sk}_1(\Delta^{n-1}) ) 
		\end{tikzcd}\ ,
	\end{align}
	and that the inductive hypothesis of level $n-1$ guarantees that the bottom row is an equivalence.
	
	\smallskip
	
	We now deal with the morphisms $\beta_j$ for $j \leqslant n-1$. Observe that the outer square in the diagram \eqref{eq:attaching_cells} can be alternatively be decomposed in the following ladder:
	\begin{align}
		\begin{tikzcd}[ampersand replacement=\&]
			\Delta^0 \amalg \Delta^0 \arrow{r}{(0,2)} \arrow{d} \& \Lambda^2_1 \arrow{r}{(j-1,j,n)} \arrow{d} \& \mathsf{sk}_1(\Delta^n)(j) \arrow{d}{\beta_j} \\
			\Delta^1 \arrow{r} \& \mathsf{sk}_1(\Delta^2) \arrow{r} \& \mathsf{sk}_1(\Delta^n)(j-1) 
		\end{tikzcd} \ .
	\end{align}
	Since the left and the outer squares are pushout, the same goes for the right one. Thus, the induced square
	\begin{align}
		\begin{tikzcd}[ampersand replacement = \&]
			F( \mathsf{sk}_1(\Delta^n)(j-1) ) \arrow{r} \arrow{d} \& F( \mathsf{sk}_1(\Delta^n)(j) ) \times_{G( \mathsf{sk}_1(\Delta^n)(j)) } G( \mathsf{sk}_1(\Delta^n)(j-1)) \arrow{d} \\
			F( \mathsf{sk}_1(\Delta^2) ) \arrow{r} \& F( \Lambda^2_1 ) \times_{G( \Lambda^2_1 )} G( \mathsf{sk}_1(\Delta^2) )
		\end{tikzcd}
	\end{align}
	is a pullback, and the bottom row is an equivalence. The conclusion follows.
\end{proof}

\section{Induced algebras and modules}

When $A$ is a monoid (in sets) and $B \subset A$ is a subset such that $1_A \in B$ and which is closed under multiplication in $A$, then $B$ inherits canonically the structure of a submonoid. In this section we provide a lift of this statement working with the Segal combinatorics.

As usual, we fix an $\infty$-category $\scrC$ with finite limits. Let $s \colon F^{\frakm} \to F^{\fraka}$ be a relative $2$-Segal object with values in $\scrC$. Given morphisms $A \to F^{\fraka}([1])$ and $M \to F^{\frakm}([0])$ in $\scrC$, we are going to spell out some sufficient conditions that guarantee that the simplicial structure of $F^\star$ induces a relative $2$-Segal space $F^\star_{A,M}$ satisfying
\begin{align}
	F^{\fraka}_{A,M}([1]) = A \quad \text{and} \quad F^{\frakm}_{A,M}([0]) = M \ .
\end{align}

\begin{notation}\label{notation:category_I}
	Let $\I$ be the (non full) subcategory of $\bfDelta\op \times \sfrel$ depicted as follows:
	\begin{align}
		\left\{ \begin{tikzcd}[ampersand replacement=\&]
			{} \& ([0], \frakm) \arrow{d} \\
			([1],\fraka) \& ([0], \fraka) \arrow{l}[swap]{\sigma_0}
		\end{tikzcd} \right\} \ . 
	\end{align}
	We let
	\begin{align}
		\begin{tikzcd}[ampersand replacement=\&]
			\bfDelta\op_{\leqslant 0} \times \sfrel \arrow{r}{i_0} \& \I \arrow{r}{i_1} \& \bfDelta\op_{\leqslant 1} \times \sfrel \arrow{r}{j_1} \& \bfDelta\op \times \sfrel
		\end{tikzcd} 
	\end{align}
	be the natural inclusions.
	We set
	\begin{align}
		\ev_{\I} \coloneqq j_1 \circ i_1 \quad \text{and} \quad j_0 \coloneqq \ev_{\I} \circ i_0 = j_1 \circ i_1 \circ i_0 \ .
	\end{align}
	Compatibly with Notation~\ref{notation:components_m_a}, we let $\I^{\fraka}$ and $\I^{\frakm}$ be the full subcategories of $\I$ spanned, respectively, by the objects $\{([1],\fraka), ([0],\fraka)\}$ and $([0],\frakm)$. For $\star \in \{\fraka, \frakm\}$, we write $i_0^\star$, $j_0^\star$, $i_1^\star$, $j_1^\star$ and $\ev_{\I}^\star$ for the restriction of the above functors to $\bfDelta\op_{\leqslant 0} \times \{\star\}$, $\I^\star$, $\bfDelta\op_{\leqslant 1} \times \{\star\}$ and $\bfDelta\op \times \{\star\}$, respectively.
\end{notation}

\begin{notation}
	Given an object $\calA = (\pi \colon M \to A_0, s \colon A_0 \to A_1) \in \Fun(\I, \scrT)$, we set
	\begin{align}
		\calA^{\fraka} \coloneqq (s \colon A_0 \to A_1) \quad \text{and} \quad \calA_0 \coloneqq (\pi \colon M \to A_0) \ . 
	\end{align}
\end{notation}

\begin{lemma}\label{lem:I_coskeletality}
	Let $\calA = (\pi \colon M \to A_0, s \colon A_0 \to A_1) \in \Fun(\I, \scrT)$.
	\begin{enumerate}\itemsep=0.2cm
		\item \label{item:I_coskeletality-1} There is a canonical equivalence
		\begin{align}
			\ev_{\I,\ast}^{\fraka}( \calA^{\fraka} ) \simeq \ev_{\I,\ast}(\calA)^{\fraka} \coloneqq \ev_{\I,\ast}(\calA) \vert_{\bfDelta\op \times \{\fraka\}} \ .
		\end{align}
		In particular, the $\fraka$-component of $\ev_{\I,\ast}(\calA)$ is $1$-coskeletal.
		
		\item \label{item:I_coskeletality-2} The square
		\begin{align} \label{eq:I_coskeletality}
			\begin{tikzcd}[ampersand replacement=\&]
				\ev_{\I,\ast}(\calA)^{\frakm} \arrow{r} \arrow{d} \&j_{0,\ast}^{\frakm}(M) \arrow{d}{j_{0,\ast}(\pi)} \\
				\ev_{\I,\ast}(\calA)^{\fraka} \arrow{r} \& j_{0,\ast}^{\fraka}(A_0)
			\end{tikzcd}
		\end{align}
		is a pullback in $\Fun(\bfDelta\op, \scrT)$. In particular, the structural morphism $\ev_{\I,\ast}(\calA)^\mathfrak{m} \to \ev_{\I,\ast}(\calA)^\mathfrak{a}$ between the $\frakm$- and the $\fraka$-components of $\ev_{\I,\ast}(\calA)$ is $0$-coskeletal, and hence $1$-Segal.
	\end{enumerate}	
\end{lemma}

\begin{proof}
	The first part of statement~\eqref{item:I_coskeletality-1} follows directly from the fact that there are no morphisms from objects of $\I^{\fraka}$ to objects of $\I^{\frakm}$. The second half follows immediately from the functoriality of right Kan extensions, as $\ev_{\I}^{\fraka}$ factors through $j_1^{\fraka} \colon \bfDelta\op_{\leqslant 1} \times \{\fraka\} \hookrightarrow \bfDelta\op \times \{\fraka\}$.
	
	We now prove statement~\eqref{item:I_coskeletality-2}. By definition, $\ev_{\I}^{\frakm}$ is just the inclusion $\bfDelta\op_{\leqslant 0} \times \{\frakm\} \to \bfDelta\op \times \{\frakm\}$, while $j_0^{\fraka}$ is the inclusion $\bfDelta\op_{\leqslant 0} \times \{\fraka\} \to \bfDelta\op \times \{\fraka\}$. It follows that both $\ev_{\I,\ast}^{\frakm}(M)$ and $j_{0,\ast}^{\fraka}(A_0)$ are $0$-coskeletal, and so the same goes for the morphism between them. Thus, given that the square \eqref{eq:I_coskeletality} is a pullback, it follows that the structural morphism $\ev_{\I,\ast}(\calA)^\mathfrak{m} \to \ev_{\I,\ast}(\calA)^\mathfrak{a}$ is $0$-coskeletal, and Proposition~\ref{prop:1_Segal_0_coskeletal} guarantees that it is $1$-Segal. So, we only have to check that the square \eqref{eq:I_coskeletality} is a pullback square.
	
	We start with the following observation. Let $k \in \{0,1\}$ and consider the functor
	\begin{align}
		j_{k,\ast} \colon \Fun\big(\bfDelta\op_{\leqslant k} \times \sfrel, \scrC\big) \longrightarrow \Fun\big(\bfDelta\op \times \sfrel, \scrC\big) 
	\end{align}
	given by right Kan extension along $j_k$.
	
	Thanks to the canonical equivalences
	\begin{align}
		\Fun\big(\bfDelta\op_{\leqslant k} \times \sfrel, \scrC\big) \simeq \Fun\big(\bfDelta\op_{\leqslant k}, \Fun(\sfrel,\scrC)\big) \quad \text{and} \quad \Fun\big(\bfDelta\op \times \sfrel, \scrC\big) \simeq \Fun\big(\bfDelta\op, \Fun(\sfrel,\scrC)\big) 
	\end{align}
	and to the fact that limits in $\Fun(\sfrel,\scrC)$ are computed level-wise, we see that given a functor $F^\star \colon \bfDelta\op_{\leqslant k} \times \sfrel \to \scrC$, one has canonical equivalences
	\begin{align}
		j_{k,\ast}(F^\star)^{\fraka} \simeq j_{k,\ast}^{\fraka}(F^{\fraka}) \quad \text{and} \quad j_{k,\ast}(F^\star)^{\frakm} \simeq j_{k,\ast}^{\frakm}(F^{\frakm}) \ . 
	\end{align}
	It follows that the square \eqref{eq:I_coskeletality} is obtained by applying $j_{1,\ast}$ to the square
	\begin{align}\label{eq:I_coskeletality_2}
		\begin{tikzcd}[ampersand replacement=\&]
			i_{1,\ast}(\calA)^{\frakm} \arrow{d} \arrow{r} \& i_\ast^{\frakm}(M) \arrow{d} \\
			i_{1,\ast}(\calA)^{\fraka} \arrow{r} \& i_\ast^{\fraka}(A_0) 
		\end{tikzcd}\ .
	\end{align}
	Since $j_{1,\ast}$ commutes with limits, it is therefore enough to check that this is a pullback square.
	
	A direct computation shows that $i_{1,\ast}(\calA)$ is canonically identified with the diagram
	\begin{align}
		\begin{tikzcd}[ampersand replacement=\&]
			A_1 \times M \times M \arrow[shift left = 4pt]{r} \arrow[shift right = 4pt]{r} \arrow[leftarrow]{r} \arrow{d} \& M \arrow{d} \\
			A_1 \times A_0 \times A_0 \arrow[shift left = 4pt]{r} \arrow[shift right = 4pt]{r} \arrow[leftarrow]{r} \& A_0
		\end{tikzcd} \ , 
	\end{align}
	and consequently $i_{1,\ast}(\calA)^{\frakm} \colon \bfDelta\op_{\leqslant 1} \to \scrT$ is identified with the top horizontal part of the above diagram, while $i_{1,\ast}(\calA)^{\fraka}$ is identified with the bottom horizontal part.
	On the other hand, recalling that $\calA_0 = (\pi \colon M \to A)$, we can identify $i_\ast(\calA_0)$ with the diagram
	\begin{align}
		\begin{tikzcd}[ampersand replacement=\&]
			M \times M \arrow[shift left = 4pt]{r} \arrow[shift right = 4pt]{r} \arrow[leftarrow]{r} \arrow{d} \& M \arrow{d} \\
			A_0 \times A_0 \arrow[shift left = 4pt]{r} \arrow[shift right = 4pt]{r} \arrow[leftarrow]{r} \& A_0
		\end{tikzcd} \ , 
	\end{align}
	and consequently $i_{0,\ast}^{\frakm}(M) \simeq i_{0,\ast}(\calA_0)^{\frakm}$ is the top horizontal part of the diagram, while $i_{0,\ast}^{\fraka}(A_0) \simeq i_{0,\ast}(\calA_0)^{\fraka}$ is the bottom horizontal part. Thus, in order to check that the square \eqref{eq:I_coskeletality_2} is a pullback, we have to check that the two squares
	\begin{align}
		\begin{tikzcd}[ampersand replacement=\&]
			M \arrow{d} \arrow[equal]{r} \& M \arrow{d} \\
			A_0 \arrow[equal]{r} \& A_0
		\end{tikzcd} \quad \text{and} \quad \begin{tikzcd}[ampersand replacement=\&]
			A_1 \times M \times M \arrow{r} \arrow{d} \& M \times M \arrow{d} \\
			A_1 \times A_0 \times A_0 \arrow{r} \& A_0 \times A_0
		\end{tikzcd} 
	\end{align}
	are pullbacks. Since the horizontal arrows of the right square are the canonical projections, both statements are obvious, and the conclusion follows.
\end{proof}

\begin{definition}\label{def:Hecke_datum}
	Let $F^\star \in \Fun(\bfDelta\op \times \Delta^1, \scrT)$ be a relative simplicial object. A \textit{Hecke datum $\rho$ for $F^\star$} is the given of an object $\calA = (s \colon A_0 \to A_1, \pi \colon M \to A_0)$ in $\Fun(\I, \scrT)$ and a morphism $\bff \coloneqq f^\star \colon \calA \to \ev_{\I}^\star(F^\star)$.
\end{definition}	
We can represent a Hecke datum $\rho = (\calA, \bff)$ as the datum of the following commutative diagram:
\begin{align}
	\begin{tikzcd}[ampersand replacement=\&]
		M \arrow{r} \arrow{d}{f_0^{\frakm}} \& A_0 \arrow{r} \arrow{d}{f_0^{\fraka}} \& A_1 \arrow{d}{f_1^{\mathfrak{a}}} \\
		F^{\frakm}([0]) \arrow{r} \& F^{\fraka}([0]) \arrow{r} \& F^{\fraka}([1])
	\end{tikzcd}
\end{align}

\begin{construction} \label{construction:Hecke_pattern}
	Let $F^\star \in \Fun(\bfDelta\op \times \sfrel, \scrT)$ be a relative simplicial object and let $\rho = (\calA, \bff)$ be a Hecke datum for $F^\star$. We define the relative simplicial object $F^\star_\rho$ as the fiber product
	\begin{align}
		\begin{tikzcd}[column sep = large, ampersand replacement=\&]
			F^\star_\rho \arrow{r} \arrow{d} \& F^\star \arrow{d} \\
			\ev_{\I,\ast}(\calA) \arrow{r}{\ev_{\I,\ast}(\bff)} \& \ev_{\I,\ast} \ev_{\I}^\ast(F^\star) 
		\end{tikzcd}\ .
	\end{align}
\end{construction}

\begin{proposition} \label{prop:Hecke_pattern}
	In the setting of Construction~\ref{construction:Hecke_pattern}, assume that the square
	\begin{align}\label{eq:Hecke_I}
		\begin{tikzcd}[ampersand replacement = \&]
			F^{\fraka}_\rho([2]) \arrow{d} \arrow{r} \& F^{\fraka}([2]) \arrow{d}{\partial_0 \times \partial_2} \\
			A_1 \times_{A_0} A_1 \arrow{r} \& F^{\fraka}([1]) \times_{F^{\fraka}([0])} F^{\fraka}([1])
		\end{tikzcd}
	\end{align}
	is a pullback. Then the morphism between the $\fraka$-components
	\begin{align}
		F^{\fraka}_\rho \longrightarrow X^{\fraka} 
	\end{align}
	is $1$-Segal. If in addition the square
	\begin{align}
		\begin{tikzcd}[column sep = large, ampersand replacement=\&]
			F^{\frakm}_\rho([1]) \arrow{r}{\partial_1} \arrow{d} \& A_1 \times_{A_0} M  \arrow{d} \\
			F^{\frakm}([1]) \arrow{r}{\partial_1} \& F^{\fraka}([1]) \times_{F^{\fraka}([0])} X^{\frakm}([0])
		\end{tikzcd}
	\end{align}
	is a pullback, then the morphism $F^\star_\rho \to F^\star$ is relative left $1$-Segal. Furthermore, replacing $\partial_0$ by $\partial_1$ in the above square and asking it to be a pullback implies that $F^\star_\rho \to F^\star$ is relative right $1$-Segal.
\end{proposition}

\begin{proof}
	Applying Lemma~\ref{lem:I_coskeletality}--\eqref{item:I_coskeletality-1}, we see that the morphism $F^{\fraka}_\rho \to F^{\fraka}$ is $1$-coskeletal. Using Proposition~\ref{prop:1_Segal_1_coskeletal} we reduce the check that the square \eqref{eq:Hecke_I} is a pullback, which holds by assumption. This proves the first half of the proposition. For the second half, observe our assumption together with Proposition~\ref{prop:checking_relative_1_Segal} shows that it is enough to prove that the map $F^{\frakm}_\rho \to F^{\frakm}$ is $1$-Segal. Set
	\begin{align}
		G^{\frakm} \coloneqq F^{\fraka}_\rho \times_{F^{\fraka}} F^{\frakm} \ , 
	\end{align}
	so that the former map factors as
	\begin{align}
		\begin{tikzcd}[column sep = small, ampersand replacement=\&]
			F^{\frakm}_\rho \arrow{r}{\varphi} \& G^{\frakm} \arrow{r}{\psi} \& F^{\frakm} 
		\end{tikzcd} \ .
	\end{align}
	The map $\psi \colon G^{\frakm} \to F^{\frakm}$ is by definition the base change of the map $F^{\fraka}_\rho \to F^{\fraka}$ we just proved to be $1$-Segal. Thus, $\psi$ is $1$-Segal itself. To complete the proof, it is therefore enough to prove that $\varphi$ is $1$-Segal. For this, define $M^{\frakm}$ by asking the right bottom square in the following commutative diagram
	\begin{align}
		\begin{tikzcd}[ampersand replacement=\&]
			F^{\frakm}_\rho \arrow{r}{\varphi} \arrow{d} \& G^{\frakm} \arrow{r} \arrow{d} \& F^{\frakm} \arrow{d} \\
			\ev_{\I,\ast}(\calA)^{\frakm} \arrow{r} \arrow{dr} \& M^{\frakm} \arrow{r} \arrow{d} \& \ev_{\I,\ast} \big( \ev_{\I}^\ast(F^\star) \big)^{\frakm} \arrow{d} \\
			\& \ev_{\I,\ast}( \ev_{\I}^\ast(\calA))^{\fraka} \arrow{r} \& \ev_{\I,\ast}(\ev_{\I}^\ast(F^\star))^{\fraka}
		\end{tikzcd}
	\end{align}
	to be a pullback. Lemma~\ref{lem:I_coskeletality}--\eqref{item:I_coskeletality-2} implies that the bottom diagonal map and the rightmost lower vertical map are $1$-Segal. Thus, the map $M^{\frakm} \to \ev_{\I,\ast}( \ev_{\I}^\ast(\calA))^{\fraka}$ is $1$-Segal. It automatically follows that $\ev_{\I,\ast}(\calA)^{\frakm} \to M^{\frakm}$ is $1$-Segal. To complete the proof, it is therefore enough to check that the upper left square is a pullback. Observe that the upper outer rectangle and the bottom lower right square are pullbacks by definition. So, the transitivity property of pullback squares reduces us to check that the right outer vertical rectangle is a pullback. However, it can alternatively be decomposed as
	\begin{align}
		\begin{tikzcd}[ampersand replacement=\&]
			G^{\frakm} \arrow{r} \arrow{d} \& F^{\frakm} \arrow{d} \\
			F^{\fraka}_\rho \arrow{r} \arrow{d} \& F^{\fraka} \arrow{d} \\
			\ev_{\I,\ast}( \ev_{\I}^\ast(\calA))^{\fraka} \arrow{r} \& \ev_{\I,\ast}(\ev_{\I}^\ast(F^\star))^{\fraka}
		\end{tikzcd}  \ ,
	\end{align}
	and now both squares are pullback by definition. The conclusion follows.
\end{proof}

\begin{corollary} \label{cor:Hecke_pattern}
	In the setting of Proposition~\ref{prop:Hecke_pattern}, if in addition $F^\star$ is a relative $2$-Segal space, the same goes for $F^\star_\rho$.
\end{corollary}

\begin{proof}
	This is an immediate consequence of Propositions~\ref{prop:1_Segal_implies_2_Segal} and \ref{prop:Hecke_pattern}.
\end{proof}

\section{The $\Lambda$-graded Segal condition}\label{sec:Lambda_graded}

When dealing with cohomological and categorical Hall algebras, it is often important to keep track of an additional grading, that geometrically comes from decomposing the relevant moduli stack of sheaves into connected components corresponding to topological invariants. In favorable situations, the set of connected components form an abelian monoid, and there is a compatibility between the Hall product and the group structure. In this section, we introduce a general framework to keep track of this compatibility at the $2$-Segal level.

\medskip

Let $(\Lambda^{\fraka},+,0)$ be an abelian monoid and let $\Lambda^{\frakm}$ be a $\Lambda^{\fraka}$-set. We can encode the group and the representation structure into a $1$-Segal relative simplicial set, which, committing a slight abuse of notation, we simply denote by
\begin{align}
	\Lambda^\star \colon \bfDelta\op \times \sfrel \longrightarrow \bfSet \ . 
\end{align}
Write
\begin{align}
	\pi \colon \bfLambda^{\star, \opp} \longrightarrow \bfDelta\op \times \sfrel 
\end{align}
for the associated cocartesian fibration. We also let
\begin{align}
	\pi^{\fraka} \colon \bfLambda^{\fraka,\mathsf{op}} \to \bfDelta\op \times \{\fraka\} 
\end{align}
be the pullback of $\pi$ along $\bfDelta\op \times \{\fraka\} \hookrightarrow \bfDelta\op \times \sfrel$. An object in $\bfLambda^{\star,\opp}$ is a $4$-tuple $([n],\star,\underline{\bfv}, w)$, where:
\begin{itemize}\itemsep=0.2cm
	\item $([n],\star) \in \bfDelta\op \times \sfrel$;
	
	\item $(\underline{\bfv},w) \in (\Lambda^{\fraka})^n \times \Lambda^{\frakm}$.
\end{itemize}
To better keep track of indexes, we enumerate the components of $\underline{\bfv}$ via the following convention:
\begin{align}
	\underline{\bfv} = (\bfv_{0,1}, \bfv_{1,2}, \ldots, \bfv_{n-1,n}) \in (\Lambda^{\fraka})^n \ . 
\end{align}
Given a subset $S = \{i_0 < i_1 < \cdots < i_m\} \subset [n]$, we write
\begin{align}
	\underline{\bfv}_S = \Big( \sum_{i_k \leqslant j < i_{k+1}} \bfv_{j,j+1} \Big)_{k = 0, 1, \ldots, m} \in (\Lambda^{\fraka})^{m} \ . 
\end{align}

\begin{definition}
	Let $\scrC$ be an $\infty$-category with finite limits. A \textit{$\Lambda^{\fraka}$-graded (resp.\ $\Lambda^\star$-graded) simplicial object with values in $\scrC$} is the given of a functor $F \colon \bfLambda^{\fraka,\opp} \to \scrC$ (resp.\ $F \colon \bfLambda^{\star,\opp} \to \scrC$).
\end{definition}

\begin{definition}\label{def:Lambda_graded_2_Segal}
	Let $\scrC$ be an $\infty$-category with finite limits. A functor $F \colon \bfLambda^{\fraka,\opp} \to \scrC$
	\begin{enumerate}\itemsep=0.2cm
		\item is \textit{pointed} if $F([0],0)$ is a terminal object in $\scrC$;
		
		\item \textit{satisfies the $1$-Segal condition} if for every $n \geqslant 2$ and every $\underline{\bfv} \in \Lambda^n$ the morphism
		\begin{align}
			F([n],\underline{\bfv}) \longrightarrow F( \{0,1\}, \bfv_{0,1} ) \times_{F(\{1\})} \cdots \times_{F(\{n-1\})} F(\{n-1,n\},\bfv_{n-1,n}) 
		\end{align}
		is an equivalence;
		
		\item \textit{satisfies the $2$-Segal condition} if for every $n \geqslant 3$, every $0 \leqslant i < j \leqslant n$ and every $\underline{\bfv} \in \Lambda^n$, the square
		\begin{align}
			\begin{tikzcd}[ampersand replacement=\&]
				F([n], \underline{\bfv}) \arrow{r} \arrow{d} \& F(\{0,\ldots,i,j,\cdots,n\}, \underline{\bfv}_{0,\cdots,i,j,\ldots n}) \arrow{d} \\
				F(\{i,i+1,\ldots, j\}, \underline{\bfv}_{i,i+1,\ldots,j}) \arrow{r} \& F(\{i,j\}, \underline{\bfv}_{i,j})
			\end{tikzcd} 
		\end{align}
		is a pullback in $\scrC$.
	\end{enumerate}
	For $k \in \{1,2\}$ we let $k_{\Lambda^{\fraka}}\textrm{-}\sfSegal(\scrC)$ (resp.\ $k_{\Lambda^{\fraka}}\textrm{-} \sfSegal_\ast(\scrC))$ be the full subcategory of $\Fun(\bfLambda^{\fraka,\opp}, \scrC)$ spanned by functors satisfying the $k$-Segal condition (resp.\ pointed functors satisfying the $k$-Segal condition).
\end{definition}

We leave to the reader to adapt the notion of relative $1$- and $2$-Segal object to the $\Lambda^\star$-graded setting (see Definition~\ref{def:relative_1_Segal}), as well as the analogous notion for morphisms of $\Lambda^{\fraka}$-graded and $\Lambda^\star$-graded objects (see Definitions~\ref{def:Segal_morphisms} and \ref{def:1_Segal_morphism}).

\medskip

In order to formulate an analogue of Gödicke's Theorems~\ref{thm:Godicke} \& \ref{thm:Godicke_II}, we need a brief digression on Day's convolution product.
\begin{recollection}[Day's convolution]\label{recollection:Day_convolution}
	Let $\scrC^\otimes$ be a symmetric monoidal $\infty$-category with coproducts. Reviewing $\Lambda$ as a symmetric monoidal discrete category, we can endow $\Fun(\Lambda, \scrC)$ with a symmetric monoidal structure known as Day's convolution. Given $A, B \colon \Lambda \to \scrC$, their Day's convolution is informally defined by the rule
	\begin{align}
		(A \otimes_\Lambda B)(\bfv) \coloneqq \coprod_{\bfv_1 + \bfv_2 = \bfv} A(\bfv_1) \otimes B(\bfv_2) \ . 
	\end{align}
	We refer to \cite[\S2.2.6]{Lurie_HA} for a thorough review of the theory in the $\infty$-categorical setting, and to \cite{Porta_Teyssier_Day} for an approach that does not rely as much on the theory of $\infty$-operads. The universal property of Day's convolution described in \cite[Theorem~2.2.6.2]{Lurie_HA} implies the following two properties:
	\begin{enumerate}\itemsep=0.2cm
		\item If $F \colon \scrC^\otimes \to \scrD^\otimes$ is a lax symmetric monoidal functor, the induced functor given by composition with $F$:
		\begin{align}
			\Fun(\Lambda, \scrC) \longrightarrow \Fun(\Lambda, \scrD) 
		\end{align}
		naturally inherits a lax symmetric monoidal structure.
		
		\item If $\alpha \colon \Lambda \to \Lambda'$ is a morphism of abelian monoids, then left Kan extension along $\alpha$
		\begin{align}
			\alpha_! \colon \Fun(\Lambda, \scrC) \longrightarrow \Fun(\Lambda', \scrC) 
		\end{align}
		naturally inherits a natural lax symmetric monoidal structure.
		In particular, when $\Lambda' = \ast$, we deduce that the functor sending $F \colon \Lambda \to \scrC$ to
		\begin{align}
			\bigoplus_{\bfv \in \Lambda} F(\bfv) \in \scrC 
		\end{align}
		has a natural lax symmetric monoidal structure.
	\end{enumerate}
\end{recollection}

In practice, we are interested in taking $\scrC$ to be an $\infty$-topos (typically, the $\infty$-topos of derived stacks) with cartesian monoidal structure. Day's convolution equips $\Fun(\Lambda^{\fraka}, \scrC)$ with a symmetric monoidal structure that we denote $\times_\Lambda$, which progagates to $\Corr(\Fun(\Lambda^{\fraka}, \scrC))$. Similarly, $\Fun(\Lambda^{\frakm}, \scrC)$ acquires a categorical action of $\Fun(\Lambda^{\fraka}, \scrC)$, which propagates to a categorical action of $\Corr(\Fun(\Lambda^{\fraka}, \scrC))$ on $\Corr(\Fun(\Lambda^{\frakm}, \scrC))$. If $F \in \Fun(\bfLambda\op, \scrC)$ is a $\Lambda^{\fraka}$-graded simplicial object, we have natural maps
\begin{align}
	\partial_0 \times \partial_2 \colon F_2(\bfv_{0,1}, \bfv_{1,2}) \longrightarrow F_1(\bfv_{1,2}) \times F_1(\bfv_{0,1}) 
\end{align}
as well as
\begin{align}
	\partial_1 \colon F_2(\bfv_{0,1},\bfv_{1,2}) \longrightarrow F_1(\bfv_{0,1} + \bfv_{1,2}) \ , 
\end{align}
which induce for every $\bfv \in \Lambda$ a correspondence
\begin{align}\label{eq:Lambda_graded_multiplication_in_correspondences}
	\begin{tikzcd}[ampersand replacement=\&]
		{} \& \displaystyle \coprod_{\bfv_{0,1} + \bfv_{1,2} = \bfv} F_2(\bfv_{0,1},\bfv_{1,2}) \arrow{dl}[swap]{\partial_0 \times \partial_2} \arrow{dr}{\partial_1} \\
		(F_1 \times_\Lambda F_1)(\bfv) \& \& F_1(\bfv)
	\end{tikzcd} \ .
\end{align}
The same argument given in \cite{Godicke} yields (we leave the details to the motivated reader):
\begin{proposition}\label{prop:Lambda_graded_2_Segal}
	Let $\scrC$ be an $\infty$-topos. There exists a well defined $\infty$-functor
	\begin{align}
		\TwoSeggraded{\Lambda^{\fraka}}_\ast^{(1)}(\scrC) \longrightarrow \Alg_{\E_1}(\Corr^{\times_\Lambda}(\Fun(\Lambda^{\fraka}, \scrC))) \ , 
	\end{align}
	that sends a $\Lambda$-graded $2$-Segal object $F$ to $F_1 \coloneqq F([1],-) \colon \Lambda \to \calE$, equipped with the multiplication given by the correspondence~\eqref{eq:Lambda_graded_multiplication_in_correspondences}. Similarly, there exists an $\infty$-functor
	\begin{align}
		\TwoSeggraded{\Lambda^{\star}}_\ast^{2\sfrel}(\scrC) \longrightarrow \mathsf{LMod}( \Corr^{\times_\Lambda}( \Fun(\Lambda^{\fraka}, \scrC) ) ; \Corr^{\times_\Lambda}( \Fun(\Lambda^{\frakm}, \scrC) ) )
	\end{align}
	lifting Theorem~\ref{thm:Godicke_II} to the $\Lambda^\star$-graded setting.
\end{proposition}

\section{Induced $\Lambda$-gradings}

We complete this combinatorial part describing a $\Lambda$-graded analogue of Corollary~\ref{cor:Hecke_pattern}. In this section we assume that $\scrC$ is an $\infty$-topos and we fix an abelian monoid $(\Lambda^{\fraka}, +, 0)$ and a $\Lambda^{\fraka}$-set $\Lambda^{\frakm}$. We denote as in the previous section the projection
\begin{align}
	\pi \colon \bfLambda^{\star,\opp} \longrightarrow \bfDelta\op \times \sfrel \ . 
\end{align}
It gives rise to an adjunction $\pi_! \dashv \pi^\ast$:
\begin{align}
	\pi_! \colon \Fun( \bfLambda^{\star,\opp}, \scrC) \leftrightarrows \Fun(\bfDelta\op \times \sfrel, \scrC) \colon \pi^\ast \ . 
\end{align}

\begin{lemma}\label{lem:Lambda_graded_2_Segal}
	Let $k \in \{1,2\}$. Assume that $\scrC$ is an $\infty$-topos.Then, both $\pi_!$ and $\pi^\ast$ preserve the $\Lambda^\star$-graded $k$-Segal conditions.
\end{lemma}

\begin{proof}
	The statement concerning $\pi^\ast$ is trivial. After unraveling the definitions, we see that it is enough to know that in $\scrC$ coproducts commute with pullbacks. However, in an $\infty$-topos coproducts commute with all ($0$-)connected limits, whence the conclusion.
\end{proof}

\begin{corollary}\label{cor:induced_Lambda_grading}
	Assume that $\scrC$ is an $\infty$-topos. Let $F \colon \bfLambda^{\star,\opp} \to \scrC$ be a $\Lambda^\star$-graded $2$-Segal object.
	Let $G \colon \bfDelta\op \times \sfrel \to \scrC$ be a relative $2$-Segal object and let $f \colon G \to \pi_!(F)$ be a morphism. Define
	\begin{align}
		G_F \coloneqq \pi^\ast(G) \times_{\pi^\ast \pi_!(F)} F \in \Fun(\bfLambda^{\star,\opp}, \scrC) \ . 
	\end{align}
	Then $G_F$ satisfies the $\Lambda^\star$-graded $2$-Segal condition.
\end{corollary}

\begin{proof}
	It follows directly from Lemma~\ref{lem:Lambda_graded_2_Segal} and the stability of the $\Lambda^\star$-graded $2$-Segal condition under limits.
\end{proof}

\begin{definition}
	In the setting of Corollary~\ref{cor:induced_Lambda_grading}, we refer to $G_F$ as the $\Lambda^\star$-grading on $G$ induced by $(F,f)$.
\end{definition}

\newpage
\part{COHAs, CatHAs, and their representations: Foundations}\label{part:foundations}

In this part, we introduce \textit{categorical and cohomological (in a motivic sense) Hall algebras} (COHAs and CatHAs, in the following) associated to a stable $\infty$-category $\scrC$ equipped with a $t$-structure $\tau$ satisfying some conditions. Moreover, given a torsion pair $(\scrT, \scrF)$ on the heart $\scrC^\heartsuit$ of the $t$-structure $\tau$ satisfying some assumptions, we define a COHA and CatHA associated to $\scrT$, together with a (categorical) representation of it associated to $\scrF$.

To introduce COHAs and CatHAs we make use of a good theory of ``realizations'', that allows to linearize the result produced by the $2$-Segal formalism discussed in Part~\ref{part:Segal}. In line with the nowadays classical literature on the subject \cite{SV_Cherednik,SV_Yangians,SV_generators,SS_Higgs,YZ_COHA,KV_Hall,Porta_Sala_Hall,PT_BPS,PT_Curve}, we consider in this paper two large families of realizations:
\begin{itemize}\itemsep=0.2cm
	\item \textit{motivic realizations};
	\item \textit{categorical realizations}.
\end{itemize}
By motivic realizations we mean all the realizations that can be constructed within A.~Khan's axiomatic framework \cite{Khan_VFC,Khan_Voevodsky_criterion,Khan_Weaves}. By categorical realizations we essentially mean the derived ($\infty$-)category of bounded coherent sheaves, for which the relevant foundational theory has been fully developed in \cite{Porta_Sala_Hall}, heavily building on \cite{Gaitsgory_Rozenblyum_Study_I,Gaitsgory_Rozenblyum_Study_II}.

\section{Homological invariants of derived stacks}\label{sec:homological-invariants}

In this section we define (motivic) Borel-Moore homology for geometric derived stacks, and establish its basic functorialities. We use in an essential way the motivic theory developed by A.\ Khan and its collaborators. Nevertheless, we need to slightly modify and extend the definitions of \cite{Khan_VFC} in order for the statement of Theorems~\ref{thm:functoriality_of_BM_homology} and \ref{thm:functoriality_of_genuine_BM_homology_admissible} to make sense. Notably, our tweaks have two goals:
\begin{enumerate}\itemsep=0.2cm
	\item combine \cite[Theorems~3.12, 3.13 \& Remark 3.7]{Khan_VFC} in a \textit{single} functoriality statement, see Theorem~\ref{thm:functoriality_of_BM_homology};
	
	\item deal with geometric stacks that are \textit{not necessarily quasi-compact}, see Theorem~\ref{thm:functoriality_of_genuine_BM_homology_admissible}.
\end{enumerate}
We thus provide a general framework for motivic homology theories, which we specialize to concrete settings in Examples~\ref{ex:genuine_motivic_formalism} and \ref{ex:topological_formalism}.

\subsection{Motivic formalisms in the quasi-compact setting} \label{subsec:motivic_formalism}

We briefly review Khan's theory of intersection theory on derived stacks \cite{Khan_VFC,Khan_Voevodsky_criterion,Khan_Ravi_Generalized,Khan_Virtual_localization}. We work over a field $k$ of characteristic zero. We let $\dASp_k$ denote the $\infty$-category of derived algebraic spaces almost of finite type over $k$ (in particular, all objects in $\dASp_k$ are assumed to be quasi-compact). To set up the theory, we mostly follow \cite{Khan_Voevodsky_criterion}.

\begin{definition}\label{def:pre_motivic_formalism}
	A \textit{pre-motivic formalism} is a functor
	\begin{align}
		\bfD^\ast \colon \dASp_k\op \longrightarrow \CAlg(\PrLotimes)
	\end{align}
	with values in presentably symmetric monoidal $\infty$-categories satisfying the following assumptions:
	\begin{enumerate}\itemsep=0.2cm
		\item \label{item:pre_motivic_formalism-1} for every smooth morphism $f \colon T \to S$, the inverse image functor $f^\ast \colon \bfD^\ast(T) \to \bfD^\ast(S)$ admits a left adjoint, denoted $f_\sharp$, and called the $\sharp$-image.
		
		\item \label{item:pre_motivic_formalism-2} The $\sharp$-image satisfies both the projection and the base-change formulas.
		
		\item \label{item:pre_motivic_formalism-3} For any finite family $(S_\alpha)_\alpha$ in $\dASp_k$, the canonical map
		\begin{align}
			\bfD^\ast\big( \coprod_\alpha S_\alpha \big) \longrightarrow \prod_\alpha \bfD^\ast(S_\alpha) 
		\end{align}
		is an equivalence.
	\end{enumerate}
\end{definition}

\begin{notation}
	Given a pre-motivic formalism $\bfD^\ast$, we denote for every $S \in \dASp_k$ by $\mathbf 1_S$ the tensor unit of $\bfD^\ast(S)$.
\end{notation}

\begin{definition}
	Let $\bfD^\ast$ be a pre-motivic formalism. Let $S \in \dASp_k$ and let $\calE$ be a finite locally free sheaf on $S$. Set $E \coloneqq \V_S(\calE) \coloneqq \Spec_S(\Sym_{\scrO_S}(\calE))$ and write $\pi \colon E \to S$ for the projection and $s \colon S \to E$ for the zero section. The \textit{Thom twist} is the endofunctor of $\bfD^\ast(S)$ given by
	\begin{align}
		\langle \calE \rangle \coloneqq \pi_\sharp \circ s_\ast \colon \bfD^\ast(S) \longrightarrow \bfD^\ast(S) \ . 
	\end{align}
	We also write $\mathsf{Th}_S(\calE) \coloneqq \mathbf 1_S\langle \calE \rangle$, and we note that $\langle \calE \rangle \simeq (-) \otimes \mathsf{Th}_S(\calE)$.
\end{definition}

\begin{notation}
	When $\calE = \scrO_S^n$, we write $\langle n \rangle$ instead of $\langle \scrO_S^n \rangle$. For $\calF \in \bfD^\ast(S)$, we also set
	\begin{align}
		\calF(n) \coloneqq \calF\langle n \rangle[-2n] \ , 
	\end{align}
	and we refer to $\calF(n)$ as the \textit{Tate twist}. See \cite[Notation~1.34 \& Remark~1.35]{Khan_Voevodsky_criterion}.
\end{notation}

\begin{definition}\label{def:motivic_formalism}
	Let $\bfD^\ast$ be a pre-motivic formalism. We say that $\bfD^\ast$ is a \textit{motivic formalism} if it satisfies the following three conditions:
	\begin{enumerate}\itemsep=0.2cm
		\item \textit{Homotopy invariance.} For every $S \in \dASp_k$ and every vector bundle $p \colon E \to S$, the unit map
		\begin{align}
			\id_{\bfD^\ast(S)} \longrightarrow p_\ast p^\ast 
		\end{align}
		is an equivalence.
		
		\item \textit{Localization.} For every open-closed pair
		\begin{align}
			\begin{tikzcd}[ampersand replacement=\&]
				U \arrow[hook]{r}{j} \& S \& Z \arrow[hook']{l}[swap]{i}
			\end{tikzcd}
		\end{align}
		in $\dASp$, the functor $i_\ast$ is fully faithful with essential image spanned by the objects in the kernel of $j^\ast$.
		
		\item \textit{Thom isomorphisms.} For every $S \in \dASp_k$ and every finite locally free sheaf $\calE$ on $S$, the Thom twist $\langle \calE \rangle$ is a self-equivalence of $\bfD^\ast(S)$.
	\end{enumerate}
\end{definition}

\begin{remark}
	When $\bfD^\ast$ has Thom isomorphisms, the motivic J-homomorphism of \cite[\S16.2]{Bachmann_Hoyois_Norms} provides us with a canonical natural transformation of functors with values in $\mathsf{Mon}_{\E_\infty}^{\mathsf{gp}}(\Spc)$ 
	\begin{align}
		\langle - \rangle \colon K(-) \longrightarrow \Pic(\bfD^\ast(-)) \ , 
	\end{align}
	where $K(S)$ denotes the Thomason-Trobaugh $K$-theory space of $S$. In particular, for every $S \in \dASp_k$ and every virtual vector bundle $\calE \in K_0(S)$, there is an associated self-equivalence $\langle \calE \rangle$. See also \cite[Remark~1.32]{Khan_Voevodsky_criterion}.
\end{remark}

It is shown in \cite[Theorems~2.24 \& 2.34]{Khan_Voevodsky_criterion} that out of a motivic formalism $\bfD^\ast$ the $\ast$-direct image satisfies proper base-change and that exceptional operations can also be defined in a unique way. In turn, the universal property of the $(\infty,2)$-category of correspondences (see \cite[\S7.3]{Gaitsgory_Rozenblyum_Study_I}) produces a six-functors formalism, that is a lax-monoidal functor
\begin{align}
	\bfD^\ast_! \colon \Corr(\dASp_k)_{\mathsf{all}, \mathsf{ft}} \longrightarrow \CAlg(\PrLotimes) \ . 
\end{align}
Alternatively, one can obtain this extension combining \cite{Khan_Voevodsky_criterion} with \cite[Proposition~A.5.10]{Mann_p_adic_six_functors} (see also \cite[Theorem~4.6]{Scholze_Six_functors}). In the special case where $\bfD^\ast$ is Voevodsky's motivic formalism (see Example~\ref{ex:genuine_motivic_formalism} below), this construction had already been performed in Khan's thesis \cite{Khan_Motivic_localization}. If $\bfD^\ast$ is a motivic formalism satisfying étale descent, one obtains via \cite[Theorem~4.3]{Khan_Weaves} a six-functors formalism
\begin{align}
	\bfD^\ast_! \colon \Corr(\dGeomqcqs_k)_{\mathsf{all}, \mathsf{ft}}^{\mathsf{proper}} \longrightarrow \CAlg(\PrLotimes) \ ,
\end{align}
defined on all qcqs geometric (e.g.\ qcqs Artin) derived stacks locally almost of finite type.

\begin{remark}\label{rem:Scholze_six_operations}
	This last extension result is delicate, in that it is not a straightforward consequence of the construction techniques of \cite[Part~III]{Gaitsgory_Rozenblyum_Study_I}. It can be obtained out of the $(\infty,2)$-categorical techniques of \cite{Liu_Zheng}, and it has initially announced in \cite[Theorem~A.5]{Khan_VFC}. However, the alternative proof (independent of \cite{Liu_Zheng}) contained a gap. This has been fixed in \cite[Theorem~4.3]{Khan_Weaves}. An alternative proof can be obtained via the universal $!$-descent topology, as sketched in \cite[Theorem~5.19, Propositions~6.6 and 6.19]{Scholze_Six_functors} (see also \cite[Appendix~5]{Mann_p_adic_six_functors}).
\end{remark}

\begin{example}[Genuine motivic formalism]\label{ex:genuine_motivic_formalism}
	The genuine motivic formalism is the functor first introduced in \cite{Morel_Voevodsky}
	\begin{align}
		\motD^\ast \coloneqq \sfSH^\ast \colon (\dSchqcqs_k)\op \longrightarrow \CAlg(\PrL) 
	\end{align}
	assigning to $S \in \dSchqcqs_k$ its Voevodsky's stable homotopy category $\sfSH(S)$ (see \cite[Definition~2.38]{Robalo_Bridge} for the precise construction). It follows from \cite{Khan_Motivic_localization} that this is a motivic formalism in the sense of Definition~\ref{def:motivic_formalism}. See also \cite[Example~2.5]{Khan_Voevodsky_criterion}. It follows from \cite{Khan_Motivic_localization} that the resulting six-functors formalism factors through underived schemes $\mathsf{Sch}_k^{\mathsf{qcqs}}$ (in fact, the same is true for every motivic formalism, as a consequence of the homotopy invariance property).
	
	The motivic formalism $\sfSH^\ast$ does not satisfy étale descent. The extension to geometric derived stacks is therefore obstructed, and for this reason it is important to consider mild variations on this construction.
	Let $\calA \in \CAlg(\sfSH(\Spec(k)))$ be an $\E_\infty$-motivic ring spectrum and set
	\begin{align}
		\motD^\ast_\calA(S) \coloneqq \Modd_{p_S^\ast(\calA)}( \sfSH^\ast(S) ) \ . 
	\end{align}
	It is straightforward to check that this is again a motivic formalism. When $\calA = \mathbb S_k$ is the motivic sphere spectrum, this construction gives back $\sfSH^\ast$. Taking $\calA = \sfH R$ the Eilenberg-Maclane spectrum associated to a commutative ring of characteristic zero, the resulting motivic formalism
	\begin{align}
		\sfDM_R^\ast \coloneqq \motD_{\sfH R} ^\ast
	\end{align}
	satisfies étale descent, and it gives therefore rise to an extended six-functors formalism
	\begin{align}
		\sfDM^\ast_{R,!} \colon \Corr(\dGeomqcqs)_{\mathsf{all},\mathsf{ft}}^{\mathsf{proper}} \longrightarrow \CAlg(\PrL) \ . 
	\end{align}
\end{example}

\begin{remark}[Orientations]\label{rem:orientations}
	Recall that $\mathsf{MGL} \in \sfSH(\Spec(k))$ is the motivic spectrum of algebraic cobordism \cite{Panin_Universality_MGL}, and that it classifies orientations of motivic spectra \cite[Definition~2.1.2 \& Proposition~2.2.6]{Deglise_Orientations} (see also \cite{Vezzosi_Brown_Peterson}). In particular, given a $\E_\infty$-motivic ring spectrum $\calA$ equipped with a morphism $\mathsf{MGL} \to \calA$ of $\E_\infty$-motivic ring spectra, \cite[Proposition~1.4.7]{Deglise_Fasel_Quadratic} and \cite[Remark~2.1.4]{Deglise_Fasel_Borel} yield for every $S \in \dSchqcqs_k$ and every virtual vector bundle $\calE \in \mathsf K^0(S)$ a canonical identification
	\begin{align}
		\langle \calE \rangle \simeq \langle \mathsf{rank}(\calE) \rangle 
	\end{align}
	between Thom twists in $\motD_{\calA}(S)$. Recall also that $\sfH R$ is canonically oriented (see for instance \cite{Hoyois_MGL}). By descent, it follows that Thom twists are trivialized in $\sfDM^\ast_R(S)$ for any $S \in \dGeomqcqs_k$.
\end{remark}

\begin{remark}[Universality of $\sfSH$]\label{rem:universality}
	The construction of $\sfSH$ given in \cite{Robalo_Bridge} endows it with a universal property, which makes it universal among motivic formalisms. In particular, given any motivic formalism $\bfD^\ast$ there exists a unique symmetric monoidal transformation
	\begin{align}
		\sfSH^\ast \longrightarrow \bfD^\ast 
	\end{align}
	compatible with $\sharp$-direct images along smooth morphisms, inverse images, tensor products and arbitrary Thom  twists.
	See \cite[Remark~2.14]{Khan_Voevodsky_criterion}. Thanks to \cite[Theorems~2.24 \& 2.34]{Khan_Voevodsky_criterion} and \cite[Proposition~A.5.10]{Mann_p_adic_six_functors} this propagates to a natural transformation
	\begin{align}
		\sfSH^\ast_! \longrightarrow \bfD^\ast_! 
	\end{align}
	between the associated six-functors formalisms defined on $\dASp_k$.
\end{remark}

\begin{example}[Topological formalism]\label{ex:topological_formalism}
	Fix $k = \C$ and recall from \cite[\S3]{Porta_Yu_Representability} the existence of an analytification functor
	\begin{align}
		(-)\an \colon \dSchqcqs_\C \longrightarrow \Top \ , 
	\end{align}
	that sends $X \in \dSchqcqs_\C$ to the underlying topological space of its analytification. Fix a (derived) commutative ring $R$ and consider the functor
	\begin{align}
		\Shhyp(-;R)^\otimes \colon \Top\op \longrightarrow \CAlg(\PrLotimes_{\mathsf{st}}) \ , 
	\end{align}
	that sends a topological space $X$ to the stable $\infty$-category $\Shhyp(X;R)$ of hypercomplete sheaves with coefficients in $\Modd_R$, equipped with its naturally induced symmetric monoidal structure (see \cite[Recollection~1.15 \& Notation~1.16]{HPT} for some background). We define
	\begin{align}
		\BettiD_R^\ast \coloneqq \Shhyp((-)\an;R) \colon (\dSchqcqs_\C)\op \longrightarrow \CAlg(\PrLotimes_{\mathsf{st}}) \ , 
	\end{align}
	and we refer to it as the \textit{topological formalism with coefficients in $R$}.
	
	\medskip
	
	The $\sharp$-pushforward introduced in \cite[Proposition 2.5 \& Notation 2.7]{HPT}, together with its base change property \cite[Corollary~2.8]{HPT} readily imply that $\BettiD_R^\ast$ is a pre-motivic formalism. It satisfies homotopy invariance as a consequence of \cite[Theorem~2.13]{HPT}. Localization is standard in this case (see e.g.\ \cite[Recollection~B.1.6]{HPT_Beyond} and \cite[Proposition~A.8.15]{Lurie_HA}). Finally, it has Thom isomorphisms as a consequence of \cite[Proposition~7.7 \& Corollary~7.9]{Volpe_Six_operations}. In particular, \cite[Corollary~7.10]{Volpe_Six_operations} guarantees that it is $\otimes$-invertible, with $\otimes$-inverse given by $s^!(\mathbf 1_E)$. Thus, $\BettiD_R^\ast$ is a motivic formalism.
\end{example}

\begin{remark}[Topological Tate twist]\label{rem:topological_Tate_twist}
	Using \cite[Corollary~7.10]{Volpe_Six_operations}, we find
	\begin{align}
		\mathbf 1_S(1) \coloneqq \mathbf 1_S \langle \scrO_S \rangle[-2] \ , 
	\end{align}
	and $\mathbf 1_S\langle \scrO_S \rangle$ is the dual of $s^!(\mathbf 1_{\A^1_\C \times S})$, where $s$ denotes the zero section of $\A^1_\C \times S$. Since $\A^1_\C \times S$ is the trivial vector bundle over $S$ of \textit{real} dimension $2$, we obtain
	\begin{align}
		s^!(\mathbf 1_{\A^1_\C \times S}) \simeq \mathbf 1_S[-2] \ , 
	\end{align}
	whence $\mathbf 1_S(1) \simeq \mathbf 1_S$. In other words, the Tate twist is trivial in the topological setting.
\end{remark}

\begin{remark}\label{rem:comparison_motivic_to_topological}
	It follows from the universality of $\sfSH^\ast$ (see Remark~\ref{rem:universality}) that there exists a unique transformation of motivic formalisms
	\begin{align}
		\sfDM^\ast_R \longrightarrow \BettiD \ . 
	\end{align}
	Combining \cite[Theorems~2.24 \& 2.34]{Khan_Voevodsky_criterion}, \cite[Proposition~A.5.10]{Mann_p_adic_six_functors} and \cite[Theorem~4.3]{Khan_Weaves}, we obtain a well defined transformation of six-functors formalisms
	\begin{align}
		\sfDM^\ast_{R,!} \longrightarrow \BettiD^\ast_! 
	\end{align}
	defined on $\dGeomqcqs_\C$. Notice that this transformation gives rise in particular to a symmetric monoidal functor
	\begin{align}
		\sfDM^\ast_{R}(X) \longrightarrow \BettiD^\ast(X) 
	\end{align}
	for every qcqs geometric derived stack $X$, which therefore sends $\sfH R \in \sfDM^\ast_{R,!}(X)$ to the tensor unit $\mathbf 1_X$ of $\BettiD^\ast(X)$. 
\end{remark}

\begin{rem}[Topological Thom twists]\label{rem:topological_orientations}
	Let $X \in \dGeomqcqs_k$ and let $\calE \in K^0(X)$ be a virtual vector bundle. Then, it follows from Remark~\ref{rem:orientations} that there is a canonical trivialization $\langle \calE \rangle \simeq \langle \mathsf{rank}(\calE) \rangle$ inside $\sfDM^\ast_R(X)$. Since the transformation $\sfDM^\ast_R \to \BettiD_R^\ast$ arising from the universality of $\sfSH^\ast$ (see Remark~\ref{rem:universality}) commutes with Thom twists, we conclude that there is a canonical equivalence
	\begin{align}
		\langle \calE \rangle \simeq \langle \mathsf{rank}(\calE) \rangle 
	\end{align}
	inside $\BettiD^\ast_R(X)$ as well. This can be obtained more directly, observing that complex vector bundles are always canonically oriented.\hfill $\triangle$
\end{rem}

\subsection{Borel-Moore homology in the quasi-compact setting}

We fix a motivic formalism $\bfD^\ast$ and consider the associated six-functors formalism
\begin{align}
	\bfD^\ast_! \colon \Corr(\dGeomqcqs_k)_{\mathsf{all}, \mathsf{ft}}^{\mathsf{proper}} \longrightarrow \CAlg(\PrLotimes) \ . 
\end{align}
We can now introduce bivariant Borel-Moore homology:
\begin{definition}[Bivariant Borel-Moore homology]\label{def:abstract_naive_BM_homology}
	Let
	\begin{align}
		\begin{tikzcd}[column sep=small, ampersand replacement=\&]
			X \arrow{rr}{f} \arrow{dr}[swap]{a} \& \& Y \arrow{dl}{b} \\
			{} \& S
		\end{tikzcd} 
	\end{align}
	be a commutative triangle in $\dGeomqcqs$. Let $\calA \in \bfD^\ast(S)$. We define the \textit{naive bivariant $\bfD$-Borel-Moore chains of $f$ with coefficients in $\calA$} as the spectrum
	\begin{align}
		\CBMDnaive(f/S;\calA) \coloneqq \R\Hom_{\bfD^\ast(S)}\big( \mathbf 1_S, a_{\ast} f^! b^\ast(\calA) \big) \in \Sp \ . 
	\end{align}
	We also define the \textit{naive bivariant $\bfD$-Borel-Moore homology groups of $f$ with coefficients in $\calA$} as the graded abelian group
	\begin{align}
		\HBMDnaive(f/S;\calA) \coloneqq \sfH_\ast\big( \CBMDnaive(f/S;\calA) \big) \ .
	\end{align}
\end{definition}

\begin{notation}\label{notation:homology_and_cohomology}
	In the setting of the above definition, we consider the following two special cases:
	\begin{itemize}\itemsep=0.2cm
		\item when $b = \id_S$ and $f = a$, we set
		\begin{align}
			\CBMDnaive(X/S;\calA) \coloneqq \CBMDnaive(a/S;\calA) \quad \text{and} \quad \HBMDnaive(X/S;\calA) \coloneqq \HBMDnaive(a/S;\calA)
		\end{align}
		and we respectively refer to them as the \textit{naive relative $\bfD$-Borel-Moore homology chains of $X$ with coefficients in $\calA$} and the \textit{naive relative $\bfD$-Borel-Moore homology groups of $X$ with coefficients in $\calA$}.
		
		\item when $f = \id_X$ and $a = b$, we set
		\begin{align}
			\CDnaive(X/S;\calA) \coloneqq \CBMDnaive(\id_X/S;\calA) \quad \text{and} \quad \HDnaive(X/S;\calA) \coloneqq \HBMDnaive(\id_X/S;\calA)
		\end{align}
		and we respectively refer to them as the \textit{naive relative $\bfD$-cochains of $X$ with coefficients in $\calA$} and the \textit{naive relative $\bfD$-cohomology groups of $X$ with coefficients in $\calA$}.
	\end{itemize}
\end{notation}

\begin{remark}
	Let $R$ be a commutative ring and assume that $\bfD$ is $R$-linear. Then $\CBMDnaive(f/S;\calA)$ becomes canonically a complex of $R$-modules, and $\HBMDnaive(f/S;\calA)$ becomes a graded $R$-module.
\end{remark}

\begin{example}\label{ex:motivic_examples}
	Consider the genuine motivic formalism $\sfDM_R^\ast$ of Example~\ref{ex:genuine_motivic_formalism}.
	\begin{enumerate}\itemsep=0.2cm
		\item \label{item:motivic_examples-1} \textit{Chow groups.}
		Take $\calA = \sfH R$ the motivic Eilenberg-MacLane $\E_\infty$-ring spectrum. In this case we simply write
		\begin{align}
			\Hmotnaive(X/S;R) \coloneqq \Hmotnaive(X/S;\sfH R) \ , 
		\end{align}
		and we refer to it as the \textit{naive rational motivic homology}. When $S = \Spec(k)$ and $X$ is a quasi-projective scheme over $k$, $\Hmotnaive(X/S;\Q)$ (as well as its twists) can be explicitly computed by Bloch's cycle complex, see \cite[Example~2.10]{Khan_VFC}; in particular one has the relation
		\begin{align}
			\tensor*[^{\naive}]{\sfH}{^{\mathsf{mot}}_{2n}}(X/\Spec(\C); \Q(n)) \simeq \mathsf A_n(X)_\Q \ , 
		\end{align}
		where $\mathsf A_n(X)$ denotes the group of $n$-dimensional algebraic cycles on $X$ up to rational equivalence. This extends to Artin stacks locally of finite type over $\Spec(k)$ admitting a stratification by quotient stacks, yielding an identification with the rationalization of Kresch's groups \cite{Kresch_Chow_groups}.
		
		\item \label{item:motivic_examples-2} \textit{G-theory.} Take $\calA = \mathsf{KGL}^{\mathsf{et}}$, the étale hypersheafification of the algebraic $K$-theory spectrum. Then $\calA$ carries a canonical $\E_\infty$-ring structure, and the associated Borel-Moore homology theory coincides with algebraic $\sfG$-theory. See \cite[Example~2.13]{Khan_VFC} or \cite[Corollary~3.3.7]{Jin_Algebraic_G_theory}. 
	\end{enumerate}
\end{example}

\begin{ex}\label{ex:topological_example}
	Working in the topological formalism of Example~\ref{ex:topological_formalism}, this construction recovers in particular the construction of Borel-Moore homology for f-Artin stacks given in \cite[\S3.2]{KV_Hall}. In particular, the extension of the formalism of constructible sheaves in Equation~(3.2.1) in \textit{loc.\ cit.} is covered by the extension theorem \cite[Theorem~4.3]{Khan_Weaves}, and \cite[Proposition~A.5.10]{Mann_p_adic_six_functors} guarantees that this extension is compatible with all the six operations. \hfill $\triangle$
\end{ex}

\begin{variant}[Twisted bivariant Borel-Moore homology]
	Given $f \colon X \to Y$ in $\dGeomqcqs_{/S}$ as in Definition~\ref{def:abstract_naive_BM_homology}, and given $\calA\in \bfD^\ast(S)$ and $\calL \in \Pic(\bfD^\ast(X))$, we define the \textit{$\calL$-twisted version of bivariant Borel-Moore chains} by
	\begin{align}
		\CBMDnaive(f/S;\calA)\langle \calL \rangle \coloneqq \R\Hom_{\bfD^\ast(X)}\big( \calL, f^! b^\ast(\calA) \big) \ .
	\end{align}
	Similarly, we set
	\begin{align}
		\HBMDnaive(f/S;\calA)\langle \calL \rangle \coloneqq \sfH_\ast\big( \CBMDnaive(f/S;\calA) \langle \calL \rangle \big) \ . 
	\end{align}
	Notice that every class $\calE \in K^0(X)$ gives rise to such an object, via the Thom twist construction $\calL \coloneqq \mathbf 1_X\langle \calE \rangle$. In particular, for every $n \in \Z$ we have a well defined object $\CBMDnaive(f/S;\calA)\langle n \rangle$.\hfill $\oslash$
\end{variant}

\begin{remark}
	Notice that $\mathbf 1_X \simeq a^\ast(\mathbf 1_S)$. Thus, the adjunction $a^\ast \dashv a_\ast$ shows that the above twisted version recovers Definition~\ref{def:abstract_naive_BM_homology} when $\calL = \mathbf 1_X$.
\end{remark}

\begin{definition}\label{def:abstract_BM_homology}
	In the same context of Definition~\ref{def:abstract_naive_BM_homology}, fix an abelian subgroup $\Gamma \subseteq \Pic(\bfD^\ast(S))$.
	We define the \textit{relative $(\bfD,\Gamma)$-Borel-Moore chains with coefficients in $\calA$} as the $\Gamma$-graded spectrum
	\begin{align}
		\CBMDGamma(f/S;\calA) \coloneqq \bigoplus_{\calL \in \Gamma} \CBMDnaive(f/S;\calA)\langle \calL \rangle \in \Sp^{\Gamma}\ .
	\end{align}
	We also define the \textit{relative $(\bfD,\Gamma)$-Borel-Moore homology groups with coefficients in $\calA$ twisted by $\calL$} as the $(\Z\times\Gamma)$-graded abelian group
	\begin{align}
		\HBMDGamma_\ast(f/S;\calA) \coloneqq \bigoplus_{\calL \in \Gamma} \HBMDnaive(f/S; \calA)\langle \calL \rangle \ .
	\end{align}
	When $\Gamma = \Z\langle1\rangle$ is the abelian subgroup generated by the Thom twist $\langle 1 \rangle$, we simply write $\CBMD(f/S;\calA)$ (resp.\ $\HBMD(f/S;\calA)$) instead of $\CBMDGamma(f/S;\calA)$ (resp.\ $\HBMDGamma(f/S;\calA)$).
\end{definition}

\begin{remark}
	Note that Remark~\ref{rem:comparison_motivic_to_topological} yields a well defined morphism
	\begin{align}
		\Cmot(X/S;\sfH R) \longrightarrow \CBM(X/S; R) \ . 
	\end{align}
\end{remark}

\begin{remark}\label{rem:explicit_bigrading}
	Taking $\Gamma = \Z\langle 1 \rangle$ and focusing on Borel-Moore homology, we have more explicitly:
	\begin{align}
		\HBMD_i(X/S;\calA) & \simeq \bigoplus_{n \in \Z} \tensor*[^{\naive}]{\sfH}{^{\bfD}_{2n+i}}(X/S;\calA)(n) \\
		& \simeq \bigoplus_{n \in \Z} \sfH_{2n+i}\big( \R\Hom_{\bfD^\ast(S)}( \mathbf 1_S(n), a_\ast a^!(\calA) ) \big) \ .
	\end{align}
	Assume on the other hand that the Tate twist admits a square root $\mathbf 1_X(1/2)$ and define the \textit{half} Thom twist setting
	\begin{align}
		(-)\langle 1/2 \rangle \coloneqq (-)[1](1/2) \ . 
	\end{align}
	Let $\Gamma = \Z\langle 1/2 \rangle$ be the subgroup generated by $\langle 1/2 \rangle$. Then
	\begin{align}
		\HBMDGamma_i(X/S;\calA) \simeq \bigoplus_{n \in \Z} \tensor*[^{\naive}]{\sfH}{^{\bfD}_{n+i}}(X/S;\calA)(n/2) \ . 
	\end{align}
	This information can be assembled together in the following table. The black rows correspond to the choice $\Gamma = \Z \langle 1 \rangle$, while the blue ones can only be defined in presence of a square root of the Tate twist, so the table as a whole corresponds to the choice $\Gamma = \Z\langle 1/2 \rangle$:
	\[
	\begin{tblr}{
			cells = {c,m},
			hline{2} = {},
			vlines = {},
			row{4,6,8,10} = {fg={blue}}
		}
		& \sfH_1( \CBMDnaiveshort\langle\ast\rangle ) & \sfH_0( \CBMDnaiveshort\langle\ast\rangle ) & \sfH_{-1}( \CBMDnaiveshort\langle\ast\rangle ) & \sfH_{-2}( \CBMDnaiveshort\langle\ast\rangle ) \\
		\vdots & \vdots & \vdots & \vdots & \vdots \\
		\CBMDnaiveshort\langle2\rangle & \HBMDnaiveshort{5}(2) & \HBMDnaiveshort{4}(2) & \HBMDnaiveshort{3}(2) & \HBMDnaiveshort{2}(2) \\
		\CBMDnaiveshort\langle3/2\rangle & \HBMDnaiveshort{4}(3/2) & \HBMDnaiveshort{3}(3/2) & \HBMDnaiveshort{2}(3/2) & \HBMDnaiveshort{1}(3/2) \\
		\CBMDnaiveshort\langle1\rangle & \HBMDnaiveshort{3}(1) & \HBMDnaiveshort{2}(1) & \HBMDnaiveshort{1}(1) & \HBMDnaiveshort{0}(1) \\
		\CBMDnaiveshort\langle1/2\rangle & \HBMDnaiveshort{2}(1/2) & \HBMDnaiveshort{1}(1/2) & \HBMDnaiveshort{0}(1/2) & \HBMDnaiveshort{-1}(1/2) \\
		\CBMDnaiveshort\langle0\rangle & \HBMDnaiveshort{1}(0) & \HBMDnaiveshort{0}(0) & \HBMDnaiveshort{-1}(0) & \HBMDnaiveshort{-2}(0)\\
		\CBMDnaiveshort\langle-1/2\rangle & \HBMDnaiveshort{0}(-1/2) & \HBMDnaiveshort{-1}(-1/2) & \HBMDnaiveshort{-2}(-1/2) & \HBMDnaiveshort{-3}(-1/2) \\
		\CBMDnaiveshort\langle-1\rangle & \HBMDnaiveshort{-1}(-1) & \HBMDnaiveshort{-2}(-1) & \HBMDnaiveshort{-3}(-1) & \HBMDnaiveshort{-4}(-1) \\
		\CBMDnaiveshort\langle-3/2\rangle & \HBMDnaiveshort{-2}(-3/2) & \HBMDnaiveshort{-3}(-3/2) & \HBMDnaiveshort{-4}(-3/2) & \HBMDnaiveshort{-5}(-3/2) \\
		\CBMDnaiveshort\langle-2\rangle & \HBMDnaiveshort{-3}(-2) & \HBMDnaiveshort{-4}(-2) & \HBMDnaiveshort{-5}(-2) & \HBMDnaiveshort{-6}(-2) \\
		\vdots & \vdots & \vdots & \vdots & \vdots \\
	\end{tblr} 
	\]
\end{remark}

\begin{notation}
	Following Example~\ref{ex:motivic_examples}--\eqref{item:motivic_examples-2}, we denote by $G(-)$ the algebraic G-theory group corresponding to the motivic Borel-Moore homology group for $\calA = \mathsf{KGL}^{\mathsf{et}}$, and we set $G_0(-)\coloneqq \pi_0 G(-)$. While, we denote by $\HBM_\ast(-)$, the \textit{usual} Borel-Moore homology group, i.e., following Example~\ref{ex:topological_example} the motivic Borel-Moore homology group for $\calA=\Q$ and $\Gamma = \Z\langle 1/2 \rangle$.
\end{notation}

\subsection{Operations for Borel-Moore homology}\label{subsec:operations_BM_homology}

We now discuss operations in extended Borel-Moore homology. These operations are constructed out of the six-functors formalism in a rather standard way. We will therefore not give the details of the constructions, but limit ourselves to give precise references.

\medskip

Fix a motivic formalism $\bfD^\ast$, a base scheme $S \in \dSchqcqs_k$. Fix as well a coefficient ring $\calA \in \CAlg(\bfD^\ast(S))$ and an abelian subgroup $\Gamma \subseteq \Pic(\bfD^\ast(S))$. Before starting to overview the operations, let me introduce two important definitions. The first one concerns orientations; to formulate it, recall first from Remark~\ref{rem:universality} that there is a natural symmetric monoidal transformation $\sfSH^\ast \to \bfD^\ast$. We denote by $\mathbf{MGL}^{\bfD}_S$ the image of the algebraic cobordism spectrum $\mathbf{MGL}_S \in \mathbf{SH}(S)$ of Remark~\ref{rem:orientations} inside $\bfD^\ast(S)$.

\begin{definition}\label{def:oriented}
	We say that an object $\calA\in \CAlg(\bfD^\ast(S))$ is \textit{oriented} if there exists a morphism
	\begin{align}
		\mathbf{MGL}^{\bfD}_S \longrightarrow \calA 
	\end{align}
	inside $\CAlg(\bfD^\ast(S))$. The given of such a morphism is referred to as an \textit{orientation on $\calA$}.
\end{definition}

\begin{remark}
	It follows from Remark~\ref{rem:orientations} that if $\calA$ is oriented then there are trivializations
	\begin{align}
		\calA\langle \calE \rangle \simeq \calA\langle \mathsf{rank}(\calE) \rangle
	\end{align}
	for every $\calE \in K^0(S)$, that canonically depend on the choice of the orientation on $\calA$. See also Remark~\ref{rem:topological_orientations}.
\end{remark}

\begin{definition}\label{def:closed_under_Thom_twists}
	We say that an abelian subgroup $\Gamma \subset \Pic(\bfD^\ast(S))$ is \textit{closed under Thom twists} if for every $n \in \Z$ and every $\calL \in \Gamma$, one has $\calL \langle n \rangle \in \Gamma$.
\end{definition}

We are now ready to list the operations. Consider a commutative triangle
\begin{align}
	\begin{tikzcd}[column sep=small,ampersand replacement=\&]
		X \arrow{rr}{f} \arrow{dr}[swap]{a_X} \& \& Y \arrow{dl}{a_Y} \\
		{} \& S
	\end{tikzcd} 
\end{align}
in $\dGeomqcqs_S$. We also assume $\calA$ to be oriented and $\Gamma$ to be closed under Thom twists.

\begin{construction}[Operations]\label{constr:operations}
	\hfill
	\begin{enumerate}\itemsep0.2cm
		\item \label{constr:operations-proper-pushforward} \textit{Proper pushforward:} if $f$ is representable by proper algebraic spaces, \cite[Theorem~2.34-(ii)]{Khan_Voevodsky_criterion} supplies a canonical equivalence $f_! \simeq f_\ast$ between functors from $\bfD^\ast(X)$ to $\bfD^\ast(Y)$. The counit of the adjunction $f_\ast \simeq f_! \dashv f^!$ induces then a transformation
		\begin{align}
			f_\ast \colon \CBMDnaive(X / S; \calA) \longrightarrow \CBMDnaive(Y / S; \calA) \ .  
		\end{align}
		Replacing $\calA$ by $\calA \otimes \calL$ and summing over $\calL \in \Gamma$, we obtain a well defined morphism
		\begin{align}
			f_\ast \colon \CBMDGamma(X/S;\calA) \longrightarrow \CBMDGamma(Y/S;\calA) \ , 
		\end{align}
		which we refer to as the proper pushforward. We denote in the same way the morphism at the level of (naive) Borel-Moore homology groups
		\begin{align}
			f_\ast \colon \HBMDGamma_\ast(X / S; \calA) \longrightarrow \HBMDGamma_\ast(Y/S;\calA) \ . 
		\end{align}
		
		\item \label{constr:operations-smooth-pullback} \textit{smooth pullback:} if $f$ is smooth, then \cite[Theorem~A.13]{Khan_VFC} supplies a canonical natural equivalence
		\begin{align}
			\mathsf{pur}_f \colon f^! \longrightarrow f^\ast\langle \LL_{X/Y} \rangle \ . 
		\end{align}
		Using this equivalence, we find a morphism
		\begin{align}
			f^! \colon \CBMDnaive(Y/S;\calA) \longrightarrow \CBMDnaive(X/S;\calA)\langle \LL_{X/Y} \rangle \ , 
		\end{align}
		which we refer to as the smooth pullback. Since $\calA$ is oriented and $\Gamma$ is closed under Thom twists, we can trivialize the Thom twist $\langle \LL_{X/Y}\rangle$, therefore obtaining a morphism
		\begin{align}
			f^! \colon \CBMDGamma(Y/S;\calA) \longrightarrow \CBMDGamma(X/S;\calA) \ . 
		\end{align}
		We denote in the same way the morphism at the level of (naive) Borel-Moore homology groups
		\begin{align}
			f^! \colon \HBMDGamma(Y/S;\calA) \longrightarrow \HBMDGamma(X/S;\calA) \ . 
		\end{align}
		
		\item \label{constr:operations-Gysin-pullback} \textit{Gysin pullback:} if $f$ is derived lci, then the derived deformation to the normal cone of \cite{Khan_Rydh_Virtual_cartier} paired with the specialization map of \cite[\S3.1]{Khan_VFC} allows to define a morphism
		\begin{align}
			f^! \colon \CBMDnaive(Y/S;\calA) \longrightarrow \CBMDnaive(X/S;\calA)\langle \LL_{X/Y} \rangle \ , 
		\end{align}
		which we refer to as the Gysin pullback. As above, since $\calA$ is oriented and $\Gamma$ is closed under Thom twists we can trivialize the Thom twist $\langle \LL_{X/Y}\rangle$, therefore obtaining a morphism
		\begin{align}
			f^! \colon \CBMDGamma(Y/S;\calA) \longrightarrow \CBMDGamma(X/S;\calA) \ . 
		\end{align}
		We denote in the same way the morphism at the level of (naive) Borel-Moore homology groups
		\begin{align}
			f^! \colon \HBMDGamma(Y/S;\calA) \longrightarrow \HBMDGamma(X/S;\calA) \ . 
		\end{align}
		
		\item \label{constr:operations-cap-product} \textit{Cap products}: let $a \colon X \to S$ be a morphism in $\dGeomqcqs_k$. The cap product is the map
		\begin{align}
			\cap \colon \CDnaive(X/S;\calA) \otimes \CBMDnaive(X/S;\calA) \longrightarrow \CBMDnaive(X/S;\calA) 
		\end{align}
		that corresponds under the adjunction $a^\ast \dashv a_\ast$ and $a_! \dashv a^!$ to the composition
		\begin{align}
			\begin{tikzcd}[ampersand replacement=\&]
				a_!( a^\ast a_\ast f^! b^\ast(\calA) \otimes a^\ast b_\ast b^!(\calA) )  \ar{r}{a_!(\varepsilon^\ast_a \otimes \varepsilon^\ast_b)} \& a_!( f^! b^\ast(\calA) \otimes f^\ast b^!(\calA) ) \ar{r}{\mathrm{Ex}^{!\ast}_\otimes(f)} \& a_!f^!( b^\ast(\calA) \otimes b^!(\calA)) \\
				{} \ar{r}{\mathrm{Ex}^{!\ast}_\otimes(b)} \& a_! f^! b^!( \calA \otimes \calA ) \simeq a_! a^!(\calA \otimes \calA) \ar{r}{\varepsilon^!_a} \& \calA \otimes \calA \ar{r}{m} \& \calA
			\end{tikzcd}\ .
		\end{align}
		Let now $\calL_1, \calL_2 \in \Pic(\bfD^\ast(S))$ and set $\calL \coloneqq \calL_1 \otimes \calL_2$. Tensoring the above chain of morphism by $\calL$ and using the projection formula, we obtain a well defined morphism
		\begin{align}
			\cap \colon \CDnaive(X/S;\calA)\langle \calL_1 \rangle \otimes \CBMDnaive(X/S;\calA) \langle \calL_2 \rangle \longrightarrow \CBMDnaive(X/S;\calA) \langle \calL \rangle \ . 
		\end{align}
		Further summing over $\calL_1, \calL_2\in \Gamma$ yields a well defined morphism
		\begin{align}
			\cap \colon \mathsf C_{\bfD,\Gamma}^\bullet(X/S;\calA) \otimes \CBMDGamma(X/S;\calA) \longrightarrow \CBMDGamma(X/S;\calA) 
		\end{align}
		and similarly at the level of (naive) Borel-Moore homology groups.
		
		\item \label{constr:operations-exterior-product} \textit{Exterior products:} let
		\begin{align}
			\begin{tikzcd}[ampersand replacement=\&]
				X \times_S Y \rar{p_1} \dar{p_2} \arrow{dr}{p} \& Y \dar{b} \\
				X \rar{a} \& S
			\end{tikzcd}
		\end{align}
		be a pullback square in $\dGeomqcqs_k$. The multiplication on $\calA$ induces for $\calE_1 \in \mathsf K^0(X)$ and $\calE_2 \in \mathsf K^0(Y)$ a well defined morphism
		\begin{align}
			\boxtimes \colon \CBMDnaive(X/S;\calA)\langle \calE_1 \rangle \otimes \CBMDnaive(Y/S;\calA)\langle \calE_2 \rangle \longrightarrow \CBMDnaive(X \times_S Y/S; \calA) \langle \calE_1 \boxtimes \calE_2 \rangle \ . 
		\end{align}
		Since $\calA$ is oriented and $\Gamma$ is closed under Thom twists, this induces a well defined morphism
		\begin{align}
			\boxtimes \colon \CBMDGamma(X/S;\calA) \otimes \CBMDGamma(Y/S;\calA) \longrightarrow \CBMDGamma(X \times_S Y/S;\calA) \ . 
		\end{align}
		We denote in the same way the morphism at the level of (naive) Borel-Moore homology groups
		\begin{align}
			\boxtimes \colon \HBMDGamma(X/S;\calA) \otimes \HBMDGamma(Y/S;\calA) \longrightarrow \HBMDGamma(X \times_S Y/S;\calA) \ . 
		\end{align}
	\end{enumerate}
\end{construction}

\begin{remark}
	The Borel-Moore homology defined in \cite{Khan_VFC} is what we refer to as the naive Borel-Moore homology. As the above formulas show, the naive version cannot become a functor with respect to the Gysin pullback operation, because of the appearance of the Thom twist. This is solved summing over all integral twists, and using the given orientation on $\calA$ to trivialize $\langle \LL_{X/Y} \rangle$.
\end{remark}

\begin{definition}\label{def:vector-bundle-stack}
	Let $Y \to S$ be a morphism in $\dGeomqcqs_k$ and let $\calE \in \catPerf(Y)$ be a perfect complex of tor-amplitude $[-k,1]$ for some integer $k \geqslant -1$.  We shall call the derived stack 
	\begin{align}
		\V_Y(\calE[-1]) \coloneqq \Spec_Y(\Sym_{\scrO_Y}(\calE[-1]))
	\end{align}
	the \textit{derived stack of co-sections of $\calE[-1]$} or the \textit{vector bundle stack} associated to $\calE$.
\end{definition}
We denote by $\pi \colon X \to Y$ the canonical projection. Then, $\pi$ is representable by smooth $(k+1)$-Artin stacks.

\begin{remark}\label{rem:Gysin_vector_bundle}
	Let $Y \to S$ be a morphism in $\dGeomqcqs_k$ and let $\calE \in \catPerf(Y)$ be a perfect complex of tor-amplitude $[-k,1]$ for some integer $k \geqslant -1$. Let $X=\V_Y(\calE[-1])$ be the vector bundle stack associated to $\calE$. Then, Construction~\ref{constr:operations}--\eqref{constr:operations-smooth-pullback} induces a well defined pullback
	\begin{align}
		\pi^! \colon \CBMDnaive(Y;\calA) \longrightarrow \CBMDnaive(X;\calA)\langle - \pi^\ast \calE[-1] \rangle \simeq \CBMDnaive(X;\calA) \langle \pi^\ast \calE \rangle \ . 
	\end{align}
	It follows from the homotopy-invariance assumption (in the form of \cite[Proposition~A.10]{Khan_VFC}) and the invertibility of $\mathsf{pur}_\pi$ that the morphism $\pi^!$ is an equivalence. In particular, if $\calA$ is oriented and $\Gamma$ is closed under Thom twists, we obtain an equivalence
	\begin{align}
		\pi^! \colon \CBMDGamma(Y/S;\calA) \simeq \CBMDGamma(X/S;\calA) \ ,
	\end{align}
	that does not respect the natural extra $\Gamma$-grading. A similar isomorphism holds at the level of (naive) Borel-Moore homology groups.
\end{remark}

\begin{remark}[Derived invariance \& deformation to the normal cone]\label{rem:BM-reduced-stack}
	The homotopy invariance axiom implies that the six-functors formalism
	\begin{align}
		\bfD^\ast_! \colon \Corr(\dSchqcqs_k)_{\mathsf{all}} \longrightarrow \CAlg(\PrL) 
	\end{align}
	factors through underived schemes (see \cite{Khan_Motivic_localization}). Equivalently, for every derived scheme $S$, the inclusion of the truncation $j \colon \trunc{S} \to S$ induces an equivalence
	\begin{align}
		j^\ast \colon \sfSH(S) \longrightarrow \sfSH(\trunc{S}) \ , 
	\end{align}
	which is compatible with the formation of all the six operations. As a consequence, Borel-Moore homology is completely insensitive to the derived structure, and so the reader might wonder how the derived lci condition is actually relevant here. The reason is that the derived deformation to the normal bundle of \cite{Khan_Rydh_Virtual_cartier} takes a derived lci morphism $f \colon X \to Y$ to a deformation $\scrD_f \to \A^1_k$, which away from $0 \in \A^1_k$ coincides with $Y \times \A^1_k$, and whose fiber at $0$ is given by $\V_X(\LL_{X/Y}[-1])$. In particular, the classical scheme $\trunc{{\scrD_f}}$ still remembers the derived structure of the morphism $f$, and therefore so does the Gysin pullback.
\end{remark}

These constructions satisfy the usual compatibilities:
\begin{theorem}[Khan]\label{thm:functoriality_of_BM_homology}
	Let $S \in \dGeomqcqs$, $\calA \in \CAlg(\bfD^\ast(S))$ and fix an abelian subroup $\Gamma \subset \Pic(\bfD^\ast(S))$. Assume that $\calA$ is oriented and that $\Gamma$ is closed under Thom twists.
	Then the construction
	\begin{align}
		\HBMDGamma_0(-/S;\calA) \colon \Corr^\times(\dGeomqcqs_S)_{\mathsf{lci},\mathsf{rpas}} \longrightarrow \Modd_\C^\heartsuit 
	\end{align}
	that sends $X \to S$ to $\HBMDGamma_0(X/S;\calA)$ (disregarding the natural extra grading) and a correspondence
	\begin{align}
		\begin{tikzcd}[column sep=small,ampersand replacement=\&]
			{} \& Z \arrow{dr}{p} \arrow{dl}[swap]{f} \\
			X \& \& Y
		\end{tikzcd} 
	\end{align}
	where $f$ is lci and $p$ is representable by proper algebraic spaces to the composite
	\begin{align}
		p_\ast \circ f^! \colon \HBMDGamma_0(X/S; \calA) \longrightarrow \HBMDGamma_0(Y/S;\calA) 
	\end{align}
	defines a lax monoidal functor.
	
	Moreover, let $(\bfD', \calA', \Gamma')$ be a second choice of a motivic formalism, an oriented ring of coefficients and an an abelian subgroup stable under Thom twists. Let $s \colon \bfD \to \bfD'$ be a symmetric monoidal natural transformation inducing an inclusion $s(\Gamma) \subseteq \Gamma'$ and let $\phi \colon s(\calA) \to \calA$ be a morphism in $\CAlg(\bfD'(S))$. Then the pair $(s,\phi)$ induces a lax symmetric monoidal transformation
	\begin{align}
		\HBMDGamma_0(-/S;\calA) \longrightarrow \mathsf{H}^{\bfD',\Gamma'}_0(-/S;\calA') \ . 
	\end{align}
\end{theorem}

\begin{proof}
	The functoriality follows combining \cite[Theorems~3.12, 3.13 \& Remark 3.7]{Khan_VFC}, while the lax monoidal structure comes from the construction of exterior products.
\end{proof}

\subsection{Beyond the quasi-compact case I}\label{subsec:beyond_qc}

Fix a motivic formalism $\bfD^\ast$, and consider the associated (extended) six-functors formalism
\begin{align}
	\bfD^\ast_! \colon \Corr( \dGeomqcqs_k )_{\mathsf{ft},\mathsf{all}} \longrightarrow \CAlg(\PrLotimes) \ ,
\end{align}
that we discussed in \S~\ref{subsec:motivic_formalism}. Consider the functor
\begin{align}
	\Pro(-) \colon \PrLotimes \longrightarrow \Cat_\infty \ , 
\end{align}
that sends an $\infty$-category $\scrC$ to the associated large $\infty$-category of pro-objects $\Pro(\scrC)$. It follows from \cite[Corollary~4.8.1.14]{Lurie_HA} that the above functor carry a canonical lax symmetric monoidal structure, and therefore it extends to a well defined functor
\begin{align}
	\Pro(-) \colon \CAlg(\PrLotimes) \longrightarrow \CAlg(\Cat_\infty) \ . 
\end{align}
We write $\proD^\ast_!$ for the induced functor
\begin{align}
	\proD^\ast_! \colon \Corr(\dGeomqcqs_k)_{\mathsf{ft},\mathsf{all}} \longrightarrow \CAlg(\Cat_\infty) \ . 
\end{align}

\begin{remark}
	At the functoriality level, this extension simply sends the functors $f^\ast \dashv f_\ast$ and $f_! \dashv f^!$ to the adjoint pairs $\Pro(f^\ast) \dashv \Pro(f_\ast)$ and $\Pro(f_!) \dashv \Pro(f^!)$.
\end{remark}

We now apply \cite[Proposition~5.16]{Scholze_Six_functors} (equivalently, \cite[Proposition~A.5.16]{Mann_p_adic_six_functors}) to extend $\proD^\ast_!$ to a six-functors formalism
\begin{align}
	\proD^\ast_! \colon \Corr\big( \PSh(\dGeomqcqs_k) \big)_{\mathsf{ft},\mathsf{all}} \longrightarrow \CAlg(\Cat_\infty) \ , 
\end{align}
where now $\mathsf{ft}$ denotes the collection of morphism that are \textit{representable} by qcqs geometric derived stacks of finite type. By \cite[Theorem~4.20]{Scholze_Six_functors}, there exists a minimal class of morphisms $\overline{\mathsf{ft}}$ containing $\mathsf{ft}$ and which is stable under disjoint unions, local on the target, local on the source, and tame in the sense of \cite[Definition~4.18]{Scholze_Six_functors}, for which there is a uniquely defined further extension
\begin{align}
	\proD^\ast_! \colon \Corr\big( \PSh(\dGeomqcqs_k) \big)_{\overline{\mathsf{ft}},\mathsf{all}} \longrightarrow \CAlg(\Cat_\infty) \ . 
\end{align}

\begin{lemma}
	Let $f \colon X \to Y$ be a morphism in $\PSh(\dGeomqcqs_k)$. Assume that $Y \in \dGeomqcqs_k$ and that there exists an exhaustion by open Zariski substacks $\{X_\alpha\}_{\alpha \in A}$ of $X$ such that:
	\begin{enumerate}\itemsep=0.2cm
		\item for each $\alpha \in A$, $X_\alpha$ is qcqs;
		
		\item for each $\alpha \in A$, the map $X_\alpha \to Y$ is of finite type.
	\end{enumerate}
	Then $f$ belongs to $\overline{\mathsf{ft}}$.
\end{lemma}

\begin{proof}
	Since open Zariski covers are of universal $!$-descent in the sense of \cite[Definition 4.14]{Scholze_Six_functors}, the conclusion directly follows from the fact that $\overline{\mathsf{ft}}$ is tame.
\end{proof}

Thanks to the above lemma, we can restrict the above six-functors formalism to
\begin{align}
	\proD^\ast_! \colon \Corr\big( \Ind_{\mathsf{open}}(\dGeomqcqs_k) )_{\ind\textrm{-}\mathsf{ft},\mathsf{all}} \longrightarrow \CAlg(\Cat_\infty) \ ,  
\end{align}
where the notations of Appendix~\ref{appendix:ind_objects} are in use. We apply Proposition~\ref{prop:indization} taking $\scrT \coloneqq \dSt_k$, $P$ the property ``being geometric'' and $Q$ the property ``being an open immersion''. It is shown in \cite[Lemma~A.1]{Porta_Sala_Hall} that every geometric derived stack $X \in \dSt_k$ is $(P,Q)$-admissible, and that any quasi-compact open exhaustion is $(P,Q)$-admissible. We then obtain a fully faithful embedding
\begin{align}\label{eq:indization}
	\dGeom_k^{\mathsf{qs}} \hookrightarrow \Ind_{\mathsf{open}}(\dGeomqcqs_k) \ .
\end{align}
Notice that this embedding takes the class of morphisms \textit{locally of finite type} (denoted $\mathsf{lft}$) to the class $\ind\textrm{-}\mathsf{ft}$. We can therefore further restrict the previous six-functor formalism to
\begin{align}\label{eq:proD}
	\proD^\ast_! \colon \Corr\big( \dGeom^{\mathsf{qs}}_k \big)_{\mathsf{lft},\mathsf{all}} \longrightarrow \CAlg(\Pro(\Cat_\infty)) \ . 
\end{align}

\begin{remark}\label{rem:pro_six_operations}
	It is equally possible to run the above extension procedure using $\bfD^\ast_!$ instead of $\proD^\ast_!$. The use of the pro-version has as a concrete consequence a change of functoriality for morphisms that are not quasi-compact. For instance, if $\scrX$ is a quasi-separated (but not quasi-compact) geometric stack and $p \colon \scrX \to \Spec(k)$ denotes the canonical morphism, then
	\begin{align}
		\tensor*[^{\mathsf{pro}}]{p}{_\ast} \colon \proD(\scrX) \longrightarrow \proD(\Spec(k)) \simeq \Pro(\bfD^\ast(k)) 
	\end{align}
	can be computed explicitly as follows: fix a quasi-compact exhaustion $\{\scrU_i\}_{i \in I}$ of $\scrX$.
	Then given $\calF \in \proD(\scrX)$, one has
	\begin{align}
		\tensor*[^{\mathsf{pro}}]{p}{_\ast}(\calF) \simeq \flim_{i \in I\op} p_{i,\ast}(\calF\vert _{\scrU_i}) \in \Pro(\bfD^\ast(\Spec(k))) \ ,
	\end{align}
	whereas
	\begin{align}
		p_\ast(\calF) \simeq \lim_{i \in I\op} p_{i,\ast}( \calF_{\scrU_i}) \in \bfD^\ast(\Spec(k)) \ . 
	\end{align}
	Notice that this formula applies in particular for objects $\calF \in \bfD^\ast(\scrX)$. On the other hand, the functor
	\begin{align}
		\tensor*[^{\mathsf{pro}}]{p}{^!} \colon \Pro(\bfD^\ast(\Spec(k))) \longrightarrow \proD(\scrX) 
	\end{align}
	can be explicitly described saying that each restriction to $\scrU_i$ is computed as
	\begin{align}
		\tensor*[^{\mathsf{pro}}]{p}{^!}(\calF) \vert _{\scrU_i} \simeq \Pro(p^!_i)(\calF) \ . 
	\end{align}
	In particular, $\tensor*[^{\mathsf{pro}}]{p}{^!}$ takes $\bfD^\ast(\Spec(k))$ to $\bfD^\ast(\scrX)$.
\end{remark}

\begin{remark}
	For the purposes of this paper, the use of pro-objects to encode the quasi-compact topology on Borel-Moore groups is more than enough. Nevertheless, one could equally replace pro-objects with \textit{condensed ones}. The above construction would then go through yielding an extended six-functors formalism
	\begin{align}
		\tensor*[^{\mathsf{cond}}]{\bfD}{^\ast_!} \colon \Corr(\dGeomqs_k)_{\mathsf{lft},\mathsf{all}} \longrightarrow \CAlg(\Pro(\Cat_\infty)) \ . 
	\end{align}
	In fact, this formalism would take values in $\CAlg(\Pro(\PrL))$, and the associated pushforward to the point would take values in $\mathsf{Cond}(\bfD^\ast(k))$. The condensed variant exhibits a better behavior with respect to the $t$-structure, but we will not need this in the rest of the paper.
\end{remark}

We can now mimic the definition of Borel-Moore homology given in \S\ref{subsec:motivic_formalism}. For this, we continue to work with our fixed motivic formalism $\bfD^!_\ast$ and its associated extended six-functors formalism
\begin{align}
	\bfD^!_\ast \colon \Corr(\dGeomqs_k)_{\mathsf{lft}, \mathsf{all}} \longrightarrow \CAlg(\Cat_\infty) \ .
\end{align}

\begin{definition}
	Let $S \in \dGeomqcqs_k$ be a qcqs geometric derived stack and fix $\calA \in \bfD^\ast(S)$ and an abelian subgroup $\Gamma \subset \bfD^\ast(S)$. Let $a \colon \scrX \to S$ be a morphism in $\dGeomqs_k$ which is in $\mathsf{lft}$. We define:
	\begin{itemize}\itemsep=0.2cm
		\item the \textit{naive relative $\bfD$-Borel-Moore chains with coefficients in $\calA$} as the pro-spectrum
		\begin{align}
			\CBMDnaive(\scrX / S; \calA) \coloneqq \R\Hom_{\proD(S)}(\mathbf 1_S, \tensor*[^{\mathsf{pro}}]{a}{_\ast} \tensor*[^{\mathsf{pro}}]{a}{^!}(\calA)) \in \Pro(\Sp) \ .  
		\end{align}
		
		\item the \textit{extended relative $(\bfD,\Gamma)$-Borel-Moore chains with coefficients in $\calA$} as the $\Gamma$-graded pro-spectrum
		\begin{align}
			\CBMDGammaext(\scrX / S; \calA) \coloneqq \bigoplus_{\calL \in \Gamma} \R\Hom_{\proD(X)}( \tensor*[^{\mathsf{pro}}]{a}{^\ast}(\calL), \tensor*[^{\mathsf{pro}}]{a}{^!}(\calA) ) \in \Pro(\Sp^{\Gamma}) \ . 
		\end{align}
	\end{itemize}
	When $\Gamma = \Z\langle 1 \rangle$ is generated by the Thom twist $\langle 1 \rangle$, we drop the superscript and simply write $\CBMD(\scrX / S; \calA)$.
\end{definition}

In order to extract the Borel-Moore homology groups in this setting, we need the following terminology.
\begin{notation}
	Let $\scrC$ be a stable $\infty$-category equipped with a $t$-structure $\tau = (\scrC_{\geqslant 0}, \scrC_{\leqslant 0})$. The $\infty$-category of pro-objects $\Pro(\scrC)$ admits an induced $t$-structure $\tensor[^{\mathsf{pro}}]{\tau}{} \coloneqq (\Pro(\scrC_{\geqslant 0}), \Pro(\scrC_{\leqslant 0}))$, whose heart coincides with $\Pro(\scrC^\heartsuit)$. We therefore have associated homology functors
	\begin{align}
		\sfH_i \colon \Pro(\scrC) \longrightarrow \Pro(\scrC^\heartsuit) 
	\end{align}
	that send a pro-object $\scrF = \flim_{j \in J} F_j$ to
	\begin{align}
		\sfH_i(\scrF) \coloneqq \flim_{j \in J} \sfH_i(F_j) \ . 
	\end{align}
\end{notation}

\begin{warning}
	Assume that $\scrC$ is presentable, so that we have at our disposal realization functors
	\begin{align}
		\lim \colon \Pro(\scrC) \longrightarrow \scrC \qquad \text{and} \qquad \Pro(\scrC^\heartsuit) \longrightarrow \scrC^\heartsuit \ . 
	\end{align}
	The square
	\begin{align}
		\begin{tikzcd}[ampersand replacement=\&]
			\Pro(\scrC) \arrow{r}{\sfH_i} \arrow{d}{\lim} \& \Pro(\scrC^\heartsuit) \arrow{d}{\lim} \\
			\scrC \arrow{r}{\sfH_i} \& \scrC^\heartsuit
		\end{tikzcd} 
	\end{align}
	typically does \textit{not} commute. For instance, when $\scrC = \Sp$ (or when $\scrC$ is the derived category of abelian groups) and the pro-object is indexed by $\N$, then lack of commutativity is controlled by Milnor's short exact sequences
	\begin{align}
		0 \longrightarrow {\lim_n}^1 \sfH_{i+1}(\scrF_n) \longrightarrow \sfH_i( \lim_n \scrF_n ) \longrightarrow \lim_n \sfH_i(\scrF_n) \longrightarrow 0 \ . 
	\end{align}
\end{warning}

\begin{definition}
	Let $S \in \dGeomqcqs_k$ be a qcqs geometric derived stack and fix $\calA \in \bfD^\ast(S)$ and an abelian subgroup $\Gamma \subset \Pic(\bfD^\ast(S))$. Let $a \colon \scrX \to S$ be a morphism in $\dGeomqs$ which is nil-lft. We define:
	\begin{itemize}\itemsep=0.2cm
		\item the \textit{naive relative $\bfD$-Borel-Moore homology groups with coefficients in $\calA$} as the graded pro-abelian group
		\begin{align}
			\HBMDnaive(\scrX / S; \calA) \coloneqq \sfH_\ast\big( \CBMDnaive(\scrX / S; \calA)  \big) \in \Pro(\mathsf{Ab}) \ . 
		\end{align}
		
		\item the \textit{extended relative $(\bfD,\Gamma)$-Borel-Moore homology groups with coefficients in $\calA$} as the $(\Z \times \Gamma)$-graded pro-abelian group
		\begin{align}
			\HBMDGammaext( \scrX / S; \calA ) \coloneqq \bigoplus_{\calL \in \Gamma} \sfH_\ast\big( \CBMD(\scrX / S; \calA)\langle \calL \rangle \big) \in \Pro(\mathsf{Ab}^{\Z \times \Gamma}) \ . 
		\end{align}
	\end{itemize}
	When $\Gamma = \Z\langle 1 \rangle$ is generated by the Thom twist $\langle 1 \rangle$, we drop the superscript and simply write $\HBMD_\ast(\scrX / S; \calA)$.
\end{definition}

\begin{remark}
	We leave it to the reader to spell out the definition of cohomology and bivariant homology in the above context, following Notation~\ref{notation:homology_and_cohomology}.
\end{remark}

\begin{remark}[Underlying topology I]\label{rem:quasi-compact_topology}
	Let $\scrF = \flim_{j \in J} F_j \in \Pro(\mathsf{Ab})$ be a pro-abelian group. The realization
	\begin{align}
		\lim_{j \in J} F_j \in \mathsf{Ab} 
	\end{align}
	canonically inherits the structure of a \textit{topological abelian group}, that can be concretely described as follows: each $F_j$ is endowed with the \textit{discrete} topology, and the limit is equipped with the induced limit topology. In particular, this procedure allows to see each $\bfD$-Borel-Moore homology group $\HBMDGamma_\ast(\scrX / S; \calA)$ as a topological abelian group. Notice that the pro-structure is induced by the pro-pushforward operation, as in Remark~\ref{rem:pro_six_operations}. In particular, when $\scrX$ is qcqs, $\HBMDGamma_\ast(\scrX / S; \calA)$ is equipped with the discrete topology. On the other hand, when $\scrX$ is only quasi-separated (or admissible), then
	\begin{align}
		\HBMDGammaext(\scrX / S; \calA) \simeq \lim_{i \in I\op} \HBMDGamma_\ast(\scrU_i / S; \calA) \ , 
	\end{align}
	where $\{\scrU_i\}_{i \in I}$ is a quasi-compact (or admissible) open exhaustion of $\scrX$ and each $\HBMDGamma_\ast(\scrU_i / S; \calA)$ is equipped with the discrete topology. For this reason, we refer to this as the \textit{quasi-compact topology on $\HBMDGamma_\ast(\scrX / S; \calA)$}.
\end{remark}

It is essentially straightforward to upgrade the operations to the extended setting.

\begin{remark}[Continuity of Borel-Moore operations]\label{rem:extension-functoriality}
	The operations discussed above are naturally compatible with the quasi-compact topology discussed in \S\ref{subsec:beyond_qc}.
	More specifically, fix a qcqs geometric derived stack $S$ and a morphism $f \colon \scrX \to \scrY$ of quasi-separated geometric derived stacks over $S$ locally of finite type. Then:
	\begin{enumerate}\itemsep=0.2cm
		\item \label{rem:extension-functoriality-proper-pushforward} The construction of the proper pushforward relies on the identification $f_\ast \simeq f_!$ for $f$ proper. Assume that $f$ is representable by proper algebraic spaces. For any qcqs open exhaustion $\{\scrU_i\}_{i\in I}$ of $\scrY$ the pullback $\{\scrU_i \times_\scrY \scrX\}_{i \in I}$ is a qcqs open exhaustion of $\scrX$. In particular, this induces a morphism of $\Gamma$-graded pro-abelian groups
		\begin{align}
			f_\ast \colon \HBMDGammaext(\scrX / S; \calA) \longrightarrow \HBMDGammaext(\scrY / S; \calA) \ . 
		\end{align}
		In line with Remark~\ref{rem:quasi-compact_topology}, we can say that $f_\ast$ is continuous for the quasi-compact topology.
		
		\item \label{rem:extension-functoriality-smooth-pullback} The smooth pullback relies on the purity equivalence $\mathsf{pur}_f \colon f^! \simeq f^\ast\langle \LL_f \rangle$ that holds when $f$ is smooth and $Y$ is qcqs. Assume that $f$ is representable by quasi-compact geometric stacks. Then if $\{\scrU_i\}_{i \in I}$ is a qcqs open exhaustion of $\scrY$, the pullback $\{\scrU_i \times_\scrY \scrX\}_{i\in I}$ is an admissible open exhaustion of $\scrX$. This implies that the smooth pullback extends to a morphism of graded pro-abelian groups
		\begin{align}
			f^\ast \colon \HBMDnaive(\scrY / S ; \calA) \longrightarrow \HBMDnaive(\scrX / S; \calA)\langle \LL_{\scrX / \scrY} \rangle \ . 
		\end{align}
		When $\calA$ is oriented, this induces a morphism of pro-abelian groups
		\begin{align}
			f^\ast \colon \HBMDGammaext(\scrY / S ; \calA) \longrightarrow \HBMDGammaext(\scrX / S; \calA) \ , 
		\end{align}
		which does not respect the extra $\Gamma$-graduation. As before, this can be summarized saying that both functors $f^\ast$ are continuous for the quasi-compact topology.
		
		\item \label{rem:extension-functoriality-Gysin-pullback} The Gysin pullback relies on the deformation to the normal cone and on the specialization map. This carries over for morphisms $f \colon \scrX \to \scrY$ that are representable by quasi-compact and derived lci geometric derived stacks: reasoning as in the previous point, we obtain a morphism of pro-objects (hence continuous for the quasi-compact topology). Thus, we obtain continuous morphisms
		\begin{align}
			f^! \colon \HBMDnaive(\scrY / S; \calA) \longrightarrow \HBMDnaive(\scrX / S; \calA) \langle \LL_{\scrX / \scrY} \rangle 
		\end{align}
		and
		\begin{align}
			f^! \colon \HBMDGammaext(\scrY / S; \calA) \longrightarrow \HBMDGammaext(\scrX / S; \calA) \ .
		\end{align}
		
		\item \label{rem:extension-functoriality-exterior-product} The exterior product only uses the formalism of six-operations, and therefore it carries over without any change in the qs setting. Indeed, notice that if $\{\scrU_i\}_{i \in I}$ and $\{\calV_j\}_{j \in J}$ are qcqs open exhaustions of $\scrX$ and $\scrY$ respectively, then $\{\scrU_i \times \calV_j\}_{(i,j) \in I \times J}$ is a qcqs open exhaustion of $\scrX \times \scrY$.
		From here, one immediately deduces a morphism of pro-objects
		\begin{align}
			\boxtimes \colon \HBMDGammaext(\scrX / S; \calA) \otimes \HBMDGammaext(\scrY / S; \calA) \longrightarrow \HBMDGammaext(\scrX \times_S \scrY / S; \calA) \ . 
		\end{align}
		This morphism induces in particular a morphism of pro-abelian groups
		\begin{align}
			\boxtimes \colon \HBMDGammaextzero(\scrX / S; \calA) \otimes \HBMDGammaextzero(\scrY / S; \calA) \longrightarrow \HBMDGammaextzero(\scrX \times_S \scrY / S; \calA) 
		\end{align}
		that does not respect the extra $\Gamma$-grading. Notice that the tensor product on the left is taken inside $\Pro(\mathsf{Ab})$; in particular its realization as a topological abelian group should rather be thought as a \textit{completed tensor product}. When we wish to forget the pro-structure and only remember the induced pro-discrete topology, we will therefore write
		\begin{align}
			\widehat{\boxtimes} \colon \HBMDGammaextzero(\scrX / S; \calA) \widehat{\otimes} \HBMDGammaextzero(\scrY / S; \calA) \longrightarrow \HBMDGammaextzero(\scrX \times_S \scrY / S; \calA) 
		\end{align}
		for the induced continuous morphism.
	\end{enumerate}
\end{remark}

In the following theorem we denote by $\mathsf{qc.lci}$ the collection of morphisms in $\dGeomqs_S$ that are representable by quasi-compact and derived lci geometric derived stacks.

\begin{theorem}\label{thm:functoriality_of_extended_BM_homology_admissible}
	Let $S \in \dGeomqcqs_k$, $\calA \in \CAlg(\bfD^\ast(S))$ and let $\Gamma \subseteq \Pic(\bfD^\ast(S))$ be an abelian subgroup. Assume that $\calA$ is oriented and that $\Gamma$ is closed under Thom twists. Then the construction
	\begin{align}
		\HBMDGammaextzero(-/S; \calA) \colon \Corr^\times(\dGeomqs_S)_{\mathsf{qc.lci},\mathsf{rpas}} \longrightarrow \Pro(\Modd_R^\heartsuit) 
	\end{align}
	that sends $\scrX \to S$ to the pro-object $\HBMDGamma_0(\scrX/S;\calA)$ (disregarding the extra $\Gamma$-grading) and a correspondence
	\begin{align}
		\begin{tikzcd}[column sep=small,ampersand replacement=\&]
			{} \& \scrZ \arrow{dr}{p} \arrow{dl}[swap]{f} \\
			\scrX \& \& \scrY
		\end{tikzcd} 
	\end{align}
	where $f$ is representable by lci geometric derived stacks and $p$ is representable by proper algebraic spaces to the composite
	\begin{align}
		p_\ast \circ f^! \colon \HBMDGammaextzero(\scrX / S; \calA) \longrightarrow \HBMDGammaextzero(\scrY / S; \calA)
	\end{align}
	defines a lax-monoidal functor. Moreover, if $(\bfD', \calA', \Gamma')$ is a second motivic formalism with a choice of an oriented ring of coefficients $\calA'$ and an abelian subgroup $\Gamma'$ closed under Thom twist, then a morphism
	\begin{align}
		(s,\phi) \colon (\bfD, \calA, \Gamma) \longrightarrow (\bfD', \calA', \Gamma') 
	\end{align}
	as in the second half of Theorem~\ref{thm:functoriality_of_BM_homology} induces a lax symmetric monoidal transformation
	\begin{align}
		\HBMDGammaextzero(-/S;\calA) \longrightarrow \tensor*[^{\mathsf{ext}}]{\sfH}{^{\bfD', \Gamma'}_0}(-/S;\calA') \ . 
	\end{align}
\end{theorem}

\subsection{Beyond the quasi-compact case II: renormalization}\label{subsec:beyond_qc_II}

The extended Borel-Moore homology groups defined in \S\ref{subsec:beyond_qc} do not always yield the desired result.
The origin of the problem can be traced to the following basic example:

\begin{example}\label{eg:infinitely_many_points}
	Let
	\begin{align}
		X \coloneqq \bigsqcup_{\Z} \Spec(k) \ . 
	\end{align}
	This is a scheme which is not quasi-compact. Its extended Borel-Moore homology groups are computed as
	\begin{align}
		\HBMDext(X) \simeq \prod_\Z R \ , 
	\end{align}
	This is a simplification of the situation that typically occurs when defining the Hall product, for which we would need a pushforward
	\begin{align}
		p_\ast \colon \HBMDext(X) \longrightarrow \HBMDext(\Spec(k)) \simeq R \ . 
	\end{align}
	However, since $X$ is an infinite disjoint union, the structural map $X \to \Spec(k)$ is not proper, and therefore the pushforward is not well defined. This can be solved redefining (or \textit{renormalizing}) Borel-Moore homology in order to force the direct sum to appear in this kind of situations.
\end{example}

\begin{definition}
	We say that a geometric derived stack $\scrX$ is \textit{connected} if whenever $\scrU \to \scrX$ is both an open Zariski and a closed immersion then either $\scrU = \emptyset$ or $\scrU = \scrX$. We denote by $\dGeomqsconn_k$ the full subcategory of $\dGeom_k$ spanned by quasi-separated and connected stacks locally almost of finite presentation.
\end{definition}

Let $S \in \dGeomqcqs_k$, $\calA \in \CAlg(\bfD^\ast(S))$ and let $\Gamma \subseteq \Pic(\bfD^\ast(S))$ be an abelian subgroup. Assume that $\calA$ is oriented and that $\Gamma$ is closed under Thom twists. Then Theorem~\ref{thm:functoriality_of_extended_BM_homology_admissible} supplies by restriction a functor
\begin{align}
	\HBMDGammaextzero(-/S; \calA) \colon (\dGeomqsconn_S)_{\mathsf{rpas}} \longrightarrow \Pro(\Modd_R^\heartsuit) \ .
\end{align}
We set
\begin{align}
	\Pro^\sqcup(\Modd_R^\heartsuit) \coloneqq \PSh^\sqcup(\Pro(\Modd_R^\heartsuit)) \ , 
\end{align}
where we are using the notation from Appendix~\ref{appendix:ind_objects}. Then we consider the functor
\begin{align}
	\HBMDGammaextzero(-/S; \calA) \colon (\dGeomqsconn_S)_{\mathsf{rpas}} \longrightarrow \Pro^\sqcup(\Modd_R^\heartsuit) \ .
\end{align}
induced by the above $\HBMDGammaextzero(-/S;\calA)$ with the natural inclusion $\Pro(\Modd_R^\heartsuit) \hookrightarrow \Pro^\sqcup(\Modd_R^\heartsuit)$. We \textit{left} Kan extend this composition along the natural inclusion
\begin{align}
	(\dGeomqsconn_S)_{\mathsf{rpas}} \hookrightarrow (\dGeomqs_S)_{\mathsf{rpas}} \ . 
\end{align}
We denote the output of this operation by
\begin{align}
	\HBMDGamma_0(-/S;\calA) \colon (\dGeom_S)_{\mathsf{rpas}} \longrightarrow \Pro^\sqcup(\Modd_R^\heartsuit) \ .
\end{align}
A simple inspection reveals that
\begin{align}\label{eq:genuine_BM}
	\HBMDGamma_0(\scrX/S; \calA) = \fbigoplus_{\alpha \in \pi_0(\scrX)} \HBMDGammaextzero(\scrX_\alpha/S;\calA) \ ,
\end{align}
where the direct sum ranges over all connected components $\scrX_\alpha$ of $\scrX$. Here $\fbigoplus$ denotes the formal coproduct taken in $\Pro^\sqcup(\Modd_R^\heartsuit)$. We refer to $\HBMDGamma_0(\scrX/S;\calA)$ simply as the \textit{Borel-Moore homology of $\scrX$ with coefficients in $\calA$}. It is bigraded (by $\pi_0(\scrX)$ and by $\Gamma$), but the operations that we discuss below do not respect this bigrading.

\begin{remark}[Underlying topology II]\label{rem:underlying_topology_II}
	This is the continuation of Remark~\ref{rem:quasi-compact_topology}. Write $\mathsf{TopAb}$ for the category of topological abelian groups, which is a complete and cocomplete category. We have in particular a realization functor
	\begin{align}
		\Pro^\sqcup(\Modd_R^\heartsuit) \longrightarrow \PSh^\sqcup(\mathsf{TopAb}) \longrightarrow \mathsf{TopAb} \ . 
	\end{align}
	Via this functor, $\HBMDGamma_0(\scrX/S;\calA)$ is realized to the topological abelian group
	\begin{align}
		\bigoplus_{\alpha \in \pi_0(\scrX)} \HBMDGammaextzero(\scrX_\alpha/S;\calA) \ , 
	\end{align}
	where each $\HBMDGammaextzero(\scrX_\alpha/S;\calA)$ is topologized as in Remark~\ref{rem:quasi-compact_topology}, and the coproduct is given the induced colimit topology. In applications, we will systematically think $\HBMDGammaextzero(\scrX/S;\calA)$ as a topological abelian group in this way, forgetting the more refined $\Pro^\sqcup$-structure.
\end{remark}

\begin{remark}[Difference with the extended groups]\label{rem:difference_with_extended}
	Notice that the formula \eqref{eq:genuine_BM} does not hold for extended Borel-Moore homology groups, even after forgetting to topological abelian groups. Indeed, one has instead
	\begin{align}
		\HBMDGammaextzero(\scrX/S;\calA) \simeq \prod_{\alpha \in \pi_0(\scrX)} \HBMDGammaextzero(\scrX_\alpha/S;\calA) \ , 
	\end{align}
	as Example~\ref{eg:infinitely_many_points} shows.
\end{remark}

The operations for $\HBMDGammaextzero(-/S;\calA)$ discussed in Remark~\ref{rem:extension-functoriality} carry over to $\HBMDGamma_0(-/S;\calA)$ with minimal modifications.
In fact, the renormalization procedure discussed above, allows to extend the class of morphisms for which the proper pushforward is defined.
In order to properly discuss them, we need to introduce two special classes of morphisms.

\begin{defin}\label{def:modified_classes_of_morphisms}
	Let $f \colon \scrX \to \scrY$ be a morphism of quasi-separated geometric derived stacks.
	We say that $f$ is:
	\begin{enumerate}\itemsep=0.2cm
		\item \textit{locally rpas} ($\mathsf{lrpas}$) if for every connected component $\scrX_0 \subset \scrX$, the composite map $\scrX_0 \to \scrY$ is representable by proper algebraic spaces;
		
		\item \textit{finitely connected} ($\mathsf{fconn}$) if for every connected component $\scrY_0 \subset \scrY$, the geometric derived stack
		\begin{align}
			f^{-1}(\scrY_0) \simeq \scrX \times_\scrY \scrY_0 
		\end{align}
		has finitely many connected components;
		
		\item \textit{universally finitely connected} ($\mathsf{ufconn}$) if it stays finitely connected after base change along any map $\scrZ \to \scrY$ in $\dGeomqs_k$. \hfill $\oslash$
	\end{enumerate}
\end{defin}

Before continuing the discussion of the functoriality of the Borel-Moore homology discussed above, let us record the following properties:

\begin{lemma}\label{lem:properties_locallly_rpas}
	The class $\mathsf{lrpas}$ of locally rpas morphisms between quasi-separated geometric derived stacks is closed under pullbacks and compositions.
\end{lemma}

\begin{proof}
	We start by discussing closure under pullbacks. Let
	\begin{align}
		\begin{tikzcd}[ampersand replacement=\&]
			\scrW \arrow{r} \arrow{d}{g} \& \scrX \arrow{d}{f} \\
			\scrZ \arrow{r} \& \scrY
		\end{tikzcd} 
	\end{align}
	be a pullback square in $\dGeomqs_k$. Assume first that $f$ is locally rpas. Let
	\begin{align}
		\scrX \simeq \bigsqcup_{\alpha \in \pi_0(\scrX)} \scrX_\alpha 
	\end{align}
	be the decomposition of $X$ into connected components. Since $\dSt_k$ is an $\infty$-topos, colimits are universal and therefore
	\begin{align}
		\scrW \simeq \bigsqcup_{\alpha \in \pi_0(\scrX)} \scrW \times_\scrZ \scrX_\alpha \ . 
	\end{align}
	Each morphism $\scrW \times_\scrZ \scrX_\alpha \to \scrZ$ is then rpas by assumption. Although $\scrW \times_\scrZ \scrX_\alpha$ need not to be connected (nor quasi-compact), the inclusion of each connected component $\scrW_{\alpha,0} \subset \scrW \times_\scrZ \scrX_\alpha$ is in particular a closed immersion, so that the composite morphism $\scrW_{\alpha,0} \to \scrZ$ is again rpas.
	
	\medskip
	
	Let now
	\begin{align}
		\begin{tikzcd}[ampersand replacement=\&]
			\scrX \arrow{r}{f} \& \scrY \arrow{r}{g} \& \scrZ
		\end{tikzcd} 
	\end{align}
	be two composable morphisms in $\dGeomqs_k$. Assume that both $f$ and $g$ are locally rpas. Let $\scrX_0 \subset \scrX$ be a connected component. Then the induced map $\scrX_0 \to \scrY$ is rpas, and it factors through a connected component $\scrY_0 \subset \scrY$. Thus, the restriction of $g\circ f$ to $\scrX_0$ can be written as the composition
	\begin{align}
		\scrX_0 \longrightarrow \scrY_0 \longrightarrow \scrZ \ , 
	\end{align}
	where now both maps are rpas. It follows that the same goes for the composite, so that $g\circ f$ is locally rpas.
\end{proof}

\begin{lemma}\label{prop:properties_fconn}
	\hfill
	\begin{enumerate}\itemsep=0.2cm
		\item \label{item:properties_fconn-1} The class $\mathsf{fconn}$ of finitely connected morphisms between quasi-separated derived stacks is closed under compositions.
		
		\item \label{item:properties_fconn-2} If $f \colon \scrX \to \scrY$ is quasi-compact and each connected component of $\scrY$ is quasi-compact, then $f$ is finitely connected.
		
		\item \label{item:properties_fconn-3} For $\scrX \in \dGeomqs_k$ and $\calF \in \catPerf(X)$, the projection
		\[ \pi \colon \V_\scrX(\calF) \longrightarrow \scrX \]
		is universally finitely connected.
	\end{enumerate}
\end{lemma}

\begin{proof}
	Let
	\begin{align}
		\begin{tikzcd}[ampersand replacement=\&]
			X \arrow{r}{f} \& Y \arrow{r}{g} \& Z
		\end{tikzcd} 
	\end{align}
	be two composable morphisms in $\dGeomqs_k$, and assume that both are finitely connected. Let $Z_0 \subset Z$ be a connected component. Then $Z_0 \times_Z Y$ consists of finitely many connected components $Y_1, Y_2, \cdots, Y_n$. We then have
	\begin{align}
		Z_0 \times_Z X \simeq \bigsqcup_{i = 1}^n Y_i \times_Y X \ , 
	\end{align}
	and each $Y_i \times_Y X$ has only finitely many connected components by assumption. Therefore, the same goes for $Z_0 \times_Z X$, which proves that the composite $gf$ is finitely connected. This proves \eqref{item:properties_fconn-1}.
	
	\medskip
	
	Statement \eqref{item:properties_fconn-2} follows directly from the definitions. We now prove statement \eqref{item:properties_fconn-3}. First notice that if $f \colon \scrY \to \scrX$ is any morphism, then
	\begin{align}
		\scrY \times_\scrX \V_\scrX(\calF) \simeq \V_\scrY(f^\ast(\calF)) \ . 
	\end{align}
	Thus, it suffices to argue that if $\scrX$ is connected, then so is $\V_\scrX(\calF)$. Choose a smooth atlas $\{u_i \colon X_i \to \scrX\}_{i \in I}$ with each $X_i$ nonempty, connected and affine. For each $i \in I$, we denote by $\scrX_i$ for the open substack of $\scrX$ image of $u_i$. Since each $X_i$ is connected, the same goes for $\scrX_i$.
	
	\medskip
	
	We claim that there exists a well-ordering on $I$ with the following property: for every $i \in I$,
	\begin{align}
		\scrX_i \cap \bigcup_{j < i} \scrX_j \ne \emptyset \ . 
	\end{align}
	To construct such a well-order, we proceed by transfinite induction. To begin with, we pick any $i_0 \in I$. Assume now that for some ordinal $\alpha$, the indexes $i_\kappa$ have been constructed for all $\kappa < \alpha$. Write $I_{<\alpha} \coloneqq \{i_k\}_{k < \alpha}$. Then the inductive hypothesis implies that
	\begin{align}
		\scrX_{<\alpha} \coloneqq \bigcup_{i \in I_{<\alpha}} \scrX_i 
	\end{align}
	is connected. Since $\scrX$ is connected, there must be an index $j \in I \smallsetminus I_{< \alpha}$ such that
	\begin{align}
		\scrX_j \cap \scrX_{<\alpha} \ne \emptyset \ . 
	\end{align}
	Indeed, if this was not the case, the two open substacks
	\begin{align}
		\scrX_{<\alpha} \qquad \text{and} \qquad \bigcup_{i \in I \smallsetminus I_{<\alpha}} \scrX_i 
	\end{align}
	would be two disjoint open subsets of $\scrX$ and their union would be the whole $\scrX$, thus contradicting the connectedness of $\scrX$. We can therefore define $i_\alpha$ to be any index $j$ for which $\scrX_j$ has non-empty intersection with $\scrX_{<\alpha}$. This completes the inductive step.
	
	\medskip
	
	Set $\calF_i \coloneqq u_i^\ast(\calF)$. Then
	\begin{align}
		\{ v_i \colon \V_{X_i}(\calF_i) \to \V_{\scrX}(\calF) \}_{i \in I} 
	\end{align}
	is a jointly surjective family of smooth morphisms, with the following property: if $\scrV_i$ denotes the open substack of $\V_{\scrX}(\calF)$ image of $v_i$, then for every $i \in I$ we have
	\begin{align}
		\scrV_i \cap \bigcup_{j < i} \scrV_j \neq \emptyset \ , 
	\end{align}
	where the order is the well-order constructed above. Thus, it suffices to show that each $\scrV_i$ is connected, and for this it is in fact sufficient to prove that $\V_{X_i}(\calF_i)$ is connected.
	
	\medskip
	
	In other words, we can assume without loss of generality that $\scrX = X$ is affine. At this point, the same inductive procedure on the tor-amplitude of $\calF$ described in the proof of \cite[Proposition~A.10]{Khan_VFC} allows to conclude.
\end{proof}

\begin{remark}[Operations]
	Let $f \colon \scrX \to \scrY$ be a morphism between quasi-separated geometric derived stacks.
	\begin{enumerate}\itemsep=0.2cm
		\item Assume that $f$ is locally rpas.
		Then the proper pushforward of Remark~\ref{rem:extension-functoriality}-\eqref{rem:extension-functoriality-proper-pushforward} induces a well defined continuous morphism
		\begin{align}
			f_\ast \colon \HBMDGamma_0(\scrX/S;\calA) \longrightarrow \HBMDGamma_0(\scrY/S;\calA) \ . 
		\end{align}
		This simply follows applying the proper pushforward of Remark~\ref{rem:extension-functoriality}-\eqref{rem:extension-functoriality-proper-pushforward} to each connected component of $\scrX$, and using the formula \eqref{eq:genuine_BM}.
		
		\item Assume that $f$ is derived lci, quasi-compact and finitely connected.
		Then the Gysin pullback of Remark~\ref{rem:extension-functoriality}-\eqref{rem:extension-functoriality-Gysin-pullback} induces a well defined continuous morphism
		\begin{align}
			f^! \colon \HBMDGamma_0(\scrY/S;\calA) \longrightarrow \HBMGamma_0(\scrX/S;\calA) \ . 
		\end{align}
		This simply follows applying the Gysin pullback of Remark~\ref{rem:extension-functoriality}-\eqref{rem:extension-functoriality-Gysin-pullback} to each connected component of $\scrY$, and using the formula \eqref{eq:genuine_BM}.
		
		\item Exterior products morphisms are obtained gluing together those of Remark~\ref{rem:extension-functoriality}-\eqref{rem:extension-functoriality-exterior-product} on each connected component, using the fact that, by design, the tensor products in $\Pro^\sqcup(\Modd^\heartsuit_R)$ commute with formal direct sums in both variables.
	\end{enumerate}
\end{remark}

In this setting, we can finally state the final version of the functoriality for Borel-Moore homology.

\begin{theorem}\label{thm:functoriality_of_genuine_BM_homology_admissible}
	Let $S \in \dGeomqcqs_k$, $\calA \in \CAlg(\bfD^\ast(S))$ and let $\Gamma \subseteq \Pic(\bfD^\ast(S))$ be an abelian subgroup. Assume that $\calA$ is oriented and that $\Gamma$ is closed under Thom twists. Then the construction
	\begin{align}
		\HBMDGamma_0(-/S; \calA) \colon \Corr^\times(\dGeomqs_S)_{\mathsf{qc.lci}\:\cap\:\mathsf{ufconn},\:\mathsf{lrpas}} \longrightarrow \Pro^\sqcup(\Modd_R^\heartsuit) 
	\end{align}
	that sends $\scrX \to S$ to the pro-object $\HBMDGamma_0(\scrX/S;\calA)$ (disregarding the extra $\Gamma$-grading) and a correspondence
	\begin{align}
		\begin{tikzcd}[column sep=small,ampersand replacement=\&]
			{} \& \scrZ \arrow{dr}{p} \arrow{dl}[swap]{f} \\
			\scrX \& \& \scrY
		\end{tikzcd} 
	\end{align}
	where $f$ is finitely connected and representable by quasi-compact lci geometric derived stacks, and $p$ is locally rpas to the composite
	\begin{align}
		p_\ast \circ f^! \colon \HBMDGamma_0(\scrX / S; \calA) \longrightarrow \HBMDGamma_0(\scrY / S; \calA)
	\end{align}
	defines a lax-monoidal functor. Moreover, if $(\bfD', \calA', \Gamma')$ is a second motivic formalism with a choice of an oriented ring of coefficients $\calA'$ and an abelian subgroup $\Gamma'$ closed under Thom twist, then a morphism
	\begin{align}
		(s,\phi) \colon (\bfD, \calA, \Gamma) \longrightarrow (\bfD', \calA', \Gamma') 
	\end{align}
	as in the second half of Theorem~\ref{thm:functoriality_of_BM_homology} induces a lax symmetric monoidal transformation
	\begin{align}
		\HBMDGamma_0(-/S;\calA) \longrightarrow \sfH^{\bfD', \Gamma'}_0(-/S;\calA') \ . 
	\end{align}
\end{theorem}

\begin{remark}
	Notice that the $\infty$-category
	\begin{align}
		\Corr^\times(\dGeomqs_S)_{\mathsf{qc.lci}\:\cap\:\mathsf{ufconn},\:\mathsf{lrpas}} 
	\end{align}
	is well defined thanks to the fact that the classes $\mathsf{lrpas}$ and $\mathsf{qc.lci} \: \cap \: \mathsf{ufconn}$ are closed under composition and pullback. For the former, this is the content of Lemma~\ref{lem:properties_locallly_rpas}. For the latter, stability under pullbacks holds by definition, while stability under composition follows directly from Lemma~\ref{prop:properties_fconn}--\eqref{item:properties_fconn-1}.
\end{remark}

\begin{remark}\label{rem:Gaitsgory_renormalization}
	The renormalization procedure considered above is akin to the \textit{renormalized de Rham pushforward} of Gaitsgory and Drinfeld \cite[\S9]{Gaitsgory_Drinfeld}. Their definition crucially relies on the fact that the exceptional inverse image $p^!$ is a functor between dualizable objects in $\PrL$, in order to define the renormalized pushforward as the dual of $p^!$. We stress that we do not know, currently, whether $\bfD(\scrX)$ is dualizable in the main examples that we consider in the current paper and in the related project \cite{DPSSV_COHA} (both for the choices of $\bfD$ and of $\scrX$). For this reason, we decided to adopt this more hands-on approach to define Borel-Moore homology. The drawback of our approach is that it yields the full functoriality only at the level of homology \textit{groups}. This limitation is already present in the quasi-compact case, and in any case we do not need any more refined functoriality in this paper. Still, we observe that, if one can solve the $\infty$-functoriality issue in the quasi-compact case, the Gaitsgory-Drinfeld approach to renormalization is likely to yield Theorem~\ref{thm:functoriality_of_genuine_BM_homology_admissible} at the chain level as well.
\end{remark}

\subsection{Borel-Moore homology for $\Lambda$-graded stacks}\label{subsec:Lambda-graded}

We finally consider a last extension of Borel-Moore homology. In constructing the COHAs, one typically works with a stack $\scrX$ that can be presented as a disjoint union
\begin{align}
	\scrX = \coprod_{\bfv \in \Lambda} \scrX(\bfv) \ ,
\end{align}
where in addition $\Lambda$ has a monoid structure, and there is a compatibility between this monoid structure and the Hall multiplication. In order to make this compatibility precise, we now discuss a final extension of Theorem~\ref{thm:functoriality_of_genuine_BM_homology_admissible}.

\medskip

Fix an abelian monoid $(\Lambda,+)$. We define the category of \textit{$\Lambda$-graded derived stacks} as:
\begin{align}
	\LambdadSt_k \coloneqq \Fun(\Lambda, \dSt_k) \ , 
\end{align}
and we consider it with the symmetric monoidal structure induced by $\Lambda$ via Day's convolution (see Recollection~\ref{recollection:Day_convolution}). This symmetric monoidal structure propagates to the $\infty$-category of correspondences $\Corr(\LambdadSt_k)$.

\begin{definition}
	Let $P$ be a property of derived stacks (resp.\ of morphisms of derived stacks). A $\Lambda$-graded derived stack $F$ (resp.\ a morphism $F \to G$) in $\Fun(\Lambda, \dSt_k)$ is said to have the property $P$ if for every $\bfv \in \Lambda$ the derived stack $F(\bfv)$ (resp.\ the morphism $F(\bfv) \to G(\bfv)$) has the property $P$.
\end{definition}

We denote by $\LambdadGeomqs_k$ the full subcategory of $\LambdadSt_k$ spanned by admissible indgeometric stacks. Concretely,
\begin{align}
	\LambdadGeomqs_k \coloneqq \Fun(\Lambda, \dGeomqs_k) \ . 
\end{align}
In the same way, we have a well defined symmetric monoidal $\infty$-category of correspondences
\begin{align}
	\Corr^\times(\LambdadGeomqs_k)_{\mathsf{qc.lci}\:\cap\:\mathsf{ufconn},\mathsf{lrpas}} \ . 
\end{align}
Fix now a coefficient ring $R$ of characteristic zero. Once again, Day's convolution endows the category
\begin{align}
	\Lambda\textrm{-} \Pro^\sqcup( \Modd_R^\heartsuit ) \coloneqq \Fun(\Lambda, \Pro^\sqcup(\Modd_R^\heartsuit)) 
\end{align}
with a symmetric monoidal structure. Combining the formal properties of Day's convolution recalled in Recollection~\ref{recollection:Day_convolution} with Theorem~\ref{thm:functoriality_of_genuine_BM_homology_admissible}, we obtain:
\begin{theorem}\label{thm:BM_Lambda_graded_functoriality}
	Let $\bfD^\ast$ be a motivic formalism. Let $S \in \dGeomqcqs$, $\calA \in \CAlg(\bfD^\ast(S))$ and let $\Gamma \subseteq \Pic(\bfD^\ast(S))$ be an abelian subgroup. Assume that $\calA$ is oriented and that $\Gamma$ is closed under Thom twists. Then, the construction
	\begin{align}
		\HBMDGamma_0(-/S; \calA) \colon \Corr(\Fun(\Lambda,\dGeomqs_S))_{\mathsf{qc.lci}\:\cap\:\mathsf{ufconn},\mathsf{lrpas}} \longrightarrow \Lambda\textrm{-}\Pro^\sqcup(\Modd_R^\heartsuit) 
	\end{align}
	that sends $\scrX \to S$ to
	\begin{align}
		\HBMDGamma_0(\scrX/S;\calA) \coloneqq \bigoplus_{\bfv \in \Lambda} \HBMDGamma_0(\scrX(\bfv)/S;\calA)
	\end{align}
	and whose functoriality is given as in Theorem~\ref{thm:functoriality_of_genuine_BM_homology_admissible} defines a lax symmetric monoidal functor.
\end{theorem}

\begin{remark}
	Note that a similar result holds at the categorical level by using the framework developed in \cite[Appendix~A]{Porta_Sala_Hall}.
\end{remark}

\section{Moduli stack of flat objects}\label{sec:moduli-stack}

In this section, we introduce various moduli stacks of objects and discuss their geometricity under suitable assumptions. The most relevant  stack for us will be the moduli stack of objects on a stable $\infty$-category which are \textit{flat} with respect to a $t$-structure. In the geometric setting, this derived stack is a natural derived enhancement of the usual classical stack of flat families of properly supported coherent sheaves on a smooth quasi-projective complex variety.

\subsection{Brief review of presentable $\infty$-categories}

We denote by $\PrL$ the $\infty$-category of presentable $\infty$-categories with functors that are left adjoints. We also denote by $\PrLomega$ the faithful (but non-full) subcategory of $\PrL$ of compactly generated presentable $\infty$-categories and functors that preserve compact objects. We refer to \cite[\S5]{HTT} for the details of the theory.

\medskip

There is a symmetric monoidal $\infty$-category structure on $\PrL$, whose underlying tensor product can be characterized in terms of bicocontinuous functors:

\begin{definition}
	Let $\scrC, \scrD$ and $\scrE$ be presentable $\infty$-categories. We say that a bifunctor
	\begin{align}
		F \colon \scrC \times \scrD \longrightarrow \scrE 
	\end{align}
	is \textit{bicocontinuous} if it is cocontinuous in both variables, that is if for every $c \in \scrC$ and $d \in \scrD$ both functors
	\begin{align}
		F(-,d) \colon \scrC \longrightarrow \scrE \quad \text{and} \quad F(c,-) \colon \scrD \longrightarrow \scrE 
	\end{align}
	commute with colimits. We denote by $\Fun^{\sfL \times \sfL}( \scrC \times \scrD, \scrE )$ the full subcategory of $\Fun(\scrC \times \scrD, \scrE)$ spanned by bicocontinuous functors.
\end{definition}

We can summarize the main results of \cite[\S4.8.1]{Lurie_HA} as follows:
\begin{theorem}[Lurie]\label{thm:main_properties_tensor_product_of_categories}
	\hfill
	\begin{enumerate}\itemsep=0.2cm
		\item \label{thm:main_properties_tensor_product_of_categories:corepresentable} For any pair of presentable $\infty$-categories $\scrC$ and $\scrD$, the functor
		\begin{align}
			\Fun^{\sfL \times \sfL}( \scrC \times \scrD, -) \colon \PrL \longrightarrow \Cat_\infty 
		\end{align}
		is corepresentable. We denote by $\scrC \otimes \scrD$ the corepresentative.
		
		\item \label{thm:main_properties_tensor_product_of_categories:monoidal_structure} The tensor product defined at the previous point endows $\PrL$ with a closed symmetric monoidal structure. Moreover, for every pair of presentable $\infty$-categories $\scrC$ and $\scrD$ the $\infty$-category of cocontinuous functors $\FunL(\scrC, \scrD)$ is itself presentable, and the natural evaluation map
		\begin{align}
			\ev_\scrC \colon \scrC \otimes \FunL(\scrC, \scrD) \longrightarrow \scrD 
		\end{align}
		exhibits $\FunL(\scrC,\scrD)$ as an exponential of $\scrD$ for $\scrC$.
		
		\item \label{thm:main_properties_tensor_product_of_categories:monoidal_structure:compact_generation} If $\scrC$ and $\scrD$ are compactly generated, then so is $\scrC \otimes \scrD$. In particular, $\PrLomega$ acquires a symmetric monoidal structure and the natural colimit-preserving inclusion functor
		\begin{align}
			\PrLomega \hookrightarrow \PrL 
		\end{align}
		acquires a symmetric monoidal structure.
		
		\item \label{thm:main_properties_tensor_product_of_categories:functor_formula} For every pair of presentable $\infty$-categories $\scrC$ and $\scrD$ the $\infty$-category $\FunR(\scrC\op, \scrD)$ is presentable and there is a canonical equivalence
		\begin{align}
			\scrC \otimes \scrD \simeq \FunR(\scrC\op, \scrD) \ . 
		\end{align}
		
		\item \label{thm:main_properties_tensor_product_of_categories:adjointness} Let $F \colon \scrC \leftrightarrows \scrD \colon G$ be a pair of adjoint functors. If both are in $\PrL$, then for every presentable $\infty$-category $\scrE$ the pair of functors $F \otimes \id_\scrE$ and $G \otimes \id_\scrE$ are again adjoint.
	\end{enumerate}
\end{theorem}

\begin{proof}
	Statement~\eqref{thm:main_properties_tensor_product_of_categories:corepresentable} follows from \cite[Propositions~4.8.1.3 \&~4.8.1.15]{Lurie_HA}, although to find the exact formulation of the universal property of $\scrC \otimes \scrD$ given above one has to unravel Notation 4.8.1.2-(iv) in \textit{loc.\ cit.} Statement \eqref{thm:main_properties_tensor_product_of_categories:monoidal_structure} follows from \cite[Proposition~4.8.1.15 \& Remark~4.8.1.18]{Lurie_HA}. Statement \eqref{thm:main_properties_tensor_product_of_categories:monoidal_structure:compact_generation} is proven in \cite[Proposition 5.3.2.11]{Lurie_HA}.
	Finally, point~\eqref{thm:main_properties_tensor_product_of_categories:functor_formula} is the content of \cite[Lemma~4.8.1.16 \& Proposition~4.8.1.17]{Lurie_HA}. As for \eqref{thm:main_properties_tensor_product_of_categories:functor_formula}, it simply follows from the functoriality of the construction: namely, the adjunction $F \dashv G$ gives rise to unit and counit transformations
	\begin{align}
		\eta \colon \id_\scrC \longrightarrow G \circ F \quad \text{and} \quad \varepsilon \colon F \circ G \longrightarrow \id_\scrD 
	\end{align}
	satisfying the triangular identities. The functoriality of $(-) \otimes \scrE$ shows then that $\eta \otimes \id_\scrE$ and $\varepsilon \otimes \id_\scrE$ still satisfy the triangular identities, whence the conclusion.
\end{proof}

We refer to objects in $\CAlg(\PrL)$ as \textit{presentably symmetric monoidal $\infty$-categories}. Given such an object $\scrV^\otimes \in \CAlg(\PrL)$, we set
\begin{align}
	\PrL_\scrV \coloneqq \Modd_{\scrV}(\PrL) \ . 
\end{align}
As a consequence of \cite[Proposition~4.8.2.18]{Lurie_HA}, if $\scrV$ is stable, every $\scrV$-module $\scrC$ is automatically stable. When $k$ is a derived commutative ring, we simply set
\begin{align}
	\PrL_k \coloneqq \PrL_{\Modd_k} \ . 
\end{align}

\begin{variant}
	Since $\Modd_k$ is compactly generated and rigid, we can equally introduce the notation
	\begin{align}
		\PrLomega_k \coloneqq \Modd_{\Modd_k}(\PrLomega) \ . 
	\end{align}
	See \cite[Proposition~2.4]{Binda_Porta_GAGA} for a justification of this definition.
\end{variant}

We now review the basic finiteness conditions for objects in $\PrLomega_k$. To begin with, recall that every object in $\PrLomega_k$ is dualizable inside $\PrL_k$, see \cite[Proposition~D.7.2.3]{Lurie_SAG}. We write $\ev_{\scrC}$ and $\mathsf{coev}_\scrC$ for the evaluation and coevaluation morphisms.

\begin{definition}
	Let $\scrC \in \PrLomega_k$. We say that $\scrC$ is:
	\begin{enumerate}\itemsep=0.2cm
		\item \textit{proper} if the evaluation map $\ev_\scrC$ preserves compact objects;
		
		\item \textit{smooth} if the coevaluation map $\mathsf{coev}_\scrC$ preserves compact objects;
		
		\item \textit{of finite type} if it is a compact object in $\PrLomega_k$.
	\end{enumerate}
\end{definition}

The first two definitions are motivated by the following:
\begin{example}\label{ex:geometric_smooth_and_proper}
	Let $X$ be a derived scheme locally almost of finite type over a field of characteristic zero $k$. Then, $X$ is proper (resp.\ smooth) if and only if $\catQCoh(X) \in \PrLomega_k$ is proper (resp.\ smooth) in the above sense. For properness, this is a particular case of \cite[Proposition~11.1.4.3]{Lurie_SAG}. For smoothness, this follows combining \cite[Theorem~11.3.6.1]{Lurie_SAG}.
\end{example}

Most of the following results have been obtained in \cite{TV_Moduli}.
\begin{proposition}\label{prop:basics_smooth_and_proper}
	Let $\scrC \in \PrLomega_k$. Then:
	\begin{enumerate}\itemsep=0.2cm
		\item \label{prop:basics_smooth_and_proper-1} $\scrC$ is proper if and only if for every pair of compact objects $M, N \in \scrC$, $\Hom_\scrC(M,N) \in \Modd_k$ is perfect; this is further equivalent to ask that $\ev_\scrC \colon \scrC \otimes_k \scrC^\vee \to \Modd_k$ admits a $k$-linear right adjoint $\ev_\scrC^r$;
		
		\item \label{prop:basics_smooth_and_proper-2} $\scrC$ is smooth if and only if $\ev_\scrC \colon \scrC \otimes_k \scrC^\vee \to \Modd_k$ admits a left adjoint $\ev_\scrC^\ell$.
		
		\item \label{prop:basics_smooth_and_proper-3} if $\scrC$ is of finite type, it is also smooth;
		
		\item \label{prop:basics_smooth_and_proper-4} if $\scrC$ is smooth and proper, it is also of finite type;
		
		\item \label{prop:basics_smooth_and_proper-5} if $\scrC$ is smooth, then it admits a single compact generator;
		
		\item \label{prop:basics_smooth_and_proper-6} Assume that $\scrC$ is smooth (resp.\ proper) and let $E$ be a single compact generator for $\scrC$.
		Then an object $F \in \scrC$ is compact if (resp.\ only if) $\Hom_\scrC(E,F) \in \Modd_k$ is perfect.
	\end{enumerate}
\end{proposition}

\begin{proof}
	Statement~\eqref{prop:basics_smooth_and_proper-1} follows directly from the definitions. For Statement~\eqref{prop:basics_smooth_and_proper-2}, observe that $\coev_{\scrC}$ preserves compact objects if and only if it admits a $k$-linear right adjoint. Unraveling the definitions, we see that \cite[Lemma~2.13]{Christ_Relative_CY_perverse_schobers} shows that this is equivalent to ask that $\ev_{\scrC}$ admits a left adjoint. For Statement~\eqref{prop:basics_smooth_and_proper-3}, observe first that if $\scrC$ is of finite type then it must admit a single compact generator. Then, the statement follows combining \cite[Lemma~2.11 \& Proposition~2.14]{TV_Moduli}. We now turn to point \eqref{prop:basics_smooth_and_proper-4}: if $\scrC$ is smooth and proper, it is dualizable as an object in $\PrLomega_k$. Thus, we can rewrite
	\begin{align}
		\Map_{\PrLomega_k}(\scrC,-) \simeq \big( \scrC^\vee \otimes (-) \big)^\simeq \ , 
	\end{align}
	and therefore $\scrC$ is a compact object in $\PrLomega_k$. Statement~\eqref{prop:basics_smooth_and_proper-5} is the content of \cite[Proposition~11.3.2.4]{Lurie_SAG}. Statement \eqref{prop:basics_smooth_and_proper-6} follows from \eqref{prop:basics_smooth_and_proper-5} and \cite[Propositions~4.6.4.4 \& 4.6.4.12]{Lurie_HA}.
\end{proof}

\begin{remark}\label{rem:finiteness_Efimov}
	In contrast to Example~\ref{ex:geometric_smooth_and_proper}, Lunts showed in \cite[Theorem~6.3]{Lunts_Categorical_resolution} that the stable $\infty$-category $\IndCoh(X) = \Ind(\catCohb(X))$ of indcoherent sheaves on $X$ is always smooth. Kontsevich conjectured that it should in fact be of finite type, and this was more recently proven by Efimov \cite{Efimov_Finiteness}.
\end{remark}

\begin{proposition}\label{prop:new_finite_type_from_old_ones}
	\hfill
	\begin{enumerate}\itemsep=0.2cm
		\item \label{item:new_finite_type_from_old_ones-1} If $\scrC, \scrD \in \PrLomega_k$ are of finite type, then so is $\scrC \otimes_k \scrD$.
		
		\item \label{item:new_finite_type_from_old_ones-2} Let $I$ be a compact object in $\Cat_\infty$ and let $\scrC \in \PrLomega_k$ be of finite type.
		Then $\Fun(I,\scrC)$ is again in $\PrLomega_k$ and it is of finite type.
		
		\item \label{item:new_finite_type_from_old_ones-3} Let $\scrC_\bullet \colon A \to \PrLR_k$ be a diagram such that $\scrC_a$ is compactly generated for every $a\in A$.
		Set
		\begin{align}
			\scrC \coloneqq \lim_{a \in A} \scrC_a \ , 
		\end{align}
		the limit being computed in $\PrL$. Then $\scrC$ is compactly generated. Furthermore, if $\scrC_a$ is of finite type for every $a \in \scrC$ and $A$ is a compact $\infty$-category, then $\scrC$ is of finite type as well.
	\end{enumerate}
\end{proposition}

\begin{proof}
	We first prove Statement~\eqref{item:new_finite_type_from_old_ones-1}. Proposition~\ref{prop:basics_smooth_and_proper}--\eqref{prop:basics_smooth_and_proper-3} \& \eqref{prop:basics_smooth_and_proper-5} allow to choose compact generators $E$ and $E'$ for $\scrC$ and $\scrD$, respectively.
	Set
	\begin{align}
		R \coloneqq \Hom_\scrC(E,E) \qquad \text{and} \qquad R' \coloneqq \Hom_\scrD(E',E') \ , 
	\end{align}
	so that $\scrC \simeq \Modd_R$ and $\scrD \simeq \Modd_{R'}$. Using \cite[Corollary~2.12]{TV_Moduli}, we see that $R$ and $R'$ are compact in $\mathsf{Alg}_{\E_1}(\Modd_k)$ and that it is enough to prove that the same holds for $R \otimes_k R'$. However, $R$ and $R'$ are retract of finite objects in $\mathsf{Alg}_{\E_1}(\Modd_k)$, and therefore the same goes for their tensor product.
	
	For Statement~\eqref{item:new_finite_type_from_old_ones-2}, we refer to \cite[Lemma~7.1.9]{Porta_Teyssier_Exodromy} and for Statement~\eqref{item:new_finite_type_from_old_ones-3} to \cite[Lemma~17.3.3]{Porta_Teyssier_Stokes}.
\end{proof}

Let $\scrC, \scrD \in \PrLomega_k$. In particular, they are dualizable in $\PrL_k$ and evaluations $\ev_\scrC$, $\ev_\scrD$ and coevaluations $\coev_\scrC$, $\coev_\scrD$ induce an adjunction
\begin{align}
	(\scrC \otimes_k -) \dashv (\scrC^\vee \otimes_k -)  \ . 
\end{align}
In particular, we obtain canonical equivalences
\begin{align}
	&\begin{tikzcd}[ampersand replacement=\&]
		\Theta_{\scrC,\scrD} \colon \FunL_k( \scrC \otimes \scrD^\vee, \Modd_k ) \ar{r}{\sim} \& \FunL_k(\scrC, \scrD)
	\end{tikzcd}\ ,
	\\
	&\begin{tikzcd}[ampersand replacement=\&]
		\Xi_{\scrC,\scrD} \colon \FunL_k( \Modd_k, \scrC^\vee \otimes \scrD ) \ar{r}{\sim} \& \FunL_k(\scrC, \scrD)
	\end{tikzcd}\ .
\end{align}
We denote by $\Theta_{\scrC,\scrD}^{-1}$ and $\Xi_{\scrC,\scrD}^{-1}$ their inverses. When $\scrD = \scrC$ we write $\Theta_\scrC$ and $\Xi_\scrC$ instead of $\Theta_{\scrC, \scrC}$ and $\Xi_{\scrC,\scrC}$.

\begin{example}
	Let $\scrD = \scrC$. Then $\Theta_{\scrC}^{-1}(\id_\scrC) = \ev_\scrC$ and $\Xi_\scrC^{-1}(\id_\scrC) = \coev_\scrC$.
\end{example}

\begin{definition}
	Let $f \colon \scrC \to \scrD$ be in $\PrLomega_k$. We say that:
	\begin{enumerate}\itemsep=0.2cm
		\item $f$ is \textit{left dualizable} if $\Theta_{\scrC,\scrD}^{-1}(f)$ admits a (automatically $k$-linear) left adjoint.
		In this case, we define its \textit{left dual} as
		\begin{align}
			f^! \coloneqq \Xi_{\scrD,\scrC}\big( \tau \circ \mathsf{l.adj}( \Theta_{\scrC,\scrD}^{-1}(f) ) \big) \ , 
		\end{align}
		where $\tau \colon \scrC \otimes \scrD^\vee \simeq \scrD^\vee \otimes \scrC$ is the symmetry of the tensor product and $\mathsf{l.adj}$ denotes the left adjoint.
		
		\item $f$ is \textit{right dualizable} if $\Theta_{\scrC,\scrD}^{-1}(f)$ admits a $k$-linear right adjoint (that is, the right adjoint should be itself a functor in $\PrL_k$). In this case, we define its \textit{right dual} as
		\begin{align}
			f^\ast \coloneqq \Xi_{\scrD,\scrC}\big( \tau \circ \mathsf{r.adj}( \Theta_{\scrC,\scrD}^{-1}(f) ) \big) \ . 
		\end{align}
	\end{enumerate}
\end{definition}

This characterization of smoothness and properness given in Proposition~\ref{prop:basics_smooth_and_proper} allows to introduce the following notion.
\begin{definition}\label{def:Serre_functors}
	Let $\scrC \in \PrLomega_k$.
	\begin{enumerate}\itemsep=0.2cm
		\item If $\scrC$ is smooth, the \textit{left Serre functor of $\scrC$} is the left dual of the identity of $\scrC$. We often write $\sfS_\scrC^!$ instead of $\id_\scrC^!$.
		
		\item If $\scrC$ is proper, the \textit{right Serre functor of $\scrC$} is the right dual of the identity of $\scrC$. We often write $\sfS_\scrC$ instead of $\id_\scrC^\ast$.
	\end{enumerate}
\end{definition}

We have the following:
\begin{proposition}\label{prop:basics_Serre_functors}
	Let $\scrC \in \PrLomega_k$.
	\begin{enumerate}\itemsep=0.2cm
		\item \label{prop:basics_Serre_functors-1} If $\scrC$ is proper, then there exists an equivalence of functors
		\begin{align}
			\Hom_\scrC(-_2, -_1)^\vee \simeq \Hom(-_1,\sfS_\scrC(-_2)) \colon (\scrC^\omega)\op \times \scrC^\omega \to \Modd_k \ . 
		\end{align}
		
		\item \label{prop:basics_Serre_functors-2} If $\scrC$ is smooth, then there exists an equivalence of functors
		\begin{align}
			\Hom_\scrC(-_2, -_1)^\vee \simeq \Hom(\sfS_\scrC^!(-_1), -_2) \colon (\scrC^{\mathsf{ps}})\op \times \scrC^\omega \to \Modd_k \ . 
		\end{align}
		
		\item \label{prop:basics_Serre_functors-3} If $\scrC$ is smooth and proper, then $\sfS_\scrC^!$ and $\sfS_\scrC$ are mutually inverse (and in particular they are auto-equivalences of $\scrC$).
	\end{enumerate}
\end{proposition}

\begin{proof}
	Statements~\eqref{prop:basics_Serre_functors-1} and \eqref{prop:basics_Serre_functors-3} are proven in \cite[Lemma~2.20]{Christ_Relative_CY_perverse_schobers}, and Statement~\eqref{prop:basics_Serre_functors-2} is proven in \cite[Lemma~2.23]{Christ_Relative_CY_perverse_schobers} (see also \cite[Corollary~2.5]{Brav_Dyckerhoff_Relative_CY_II}).
\end{proof}

\begin{example}\label{ex:smooth_Serre_functor_geometric_example}
	Let $X$ be a smooth separated algebraic variety. In this case, $\catQCoh(X) \simeq \IndCoh(X)$. Notice that the auto-duality $(-)^\vee \colon \catPerf(X)\op \simeq \catPerf(X)$ induces an identification $\catQCoh(X)^\vee \simeq \catQCoh(X)$, so that \cite[Theorem~4.7]{BenZvi_Francis_Nadler_Integral_transforms} supplies a canonical identification
	\begin{align}
		\catQCoh(X) \otimes_k \catQCoh(X)^\vee \simeq \catQCoh(X \times X) \ . 
	\end{align}
	Observe that the spans
	\begin{align}
		\coev_X \colon \begin{tikzcd}[ampersand replacement=\&]
			X \arrow{r}{\Delta_X} \arrow{d}{\pi_X} \& X \times X \\
			\Spec(k)
		\end{tikzcd} \quad \text{and} \quad \ev_X \colon \begin{tikzcd}[ampersand replacement=\&]
			X \arrow{r}{\pi_X} \arrow{d}{\Delta_X} \& \Spec(k) \\
			X \times X
		\end{tikzcd}
	\end{align}
	exhibit $X$ as a dualizable object in $\Corr(\dSch_k)$. Applying $\IndCoh$, we therefore find that $\catQCoh(X)$ is dualizable, with evaluation map given by
	\begin{align}
		\pi_{X,\ast} \circ \Delta_X^! \colon \catQCoh(X) \otimes \catQCoh(X)^\vee \longrightarrow \Modd_k \ . 
	\end{align}
	The left adjoint to $\ev_{\catQCoh(X)}$ is then $\Delta_{X,\ast} \circ \pi_X^\ast$. From here, one obtains a canonical identification
	\begin{align}
		\sfS_{\catQCoh(X)}^!(F) \simeq F \otimes_k \omega_X^{-1} \ , 
	\end{align}
	where $\omega_X \coloneqq \pi_X^!(k)$ is the dualizing complex of $X$. Since $X$ is smooth, it can further be shown that $\omega_X \simeq \det(\LL_X)[\mathsf{dim}(X)]$.
\end{example}

\subsection{Families of objects}

Let $k$ be a field of characteristic zero and let $A$ be a derived commutative $k$-algebra. Let $\scrC \in \PrLomega_k$. The moduli of objects of $\scrC$ is a derived stack that parametrizes objects in $\scrC$. It should therefore send a test derived affine scheme $S = \Spec(A)$ to an \textit{$S$-family of objects in $\scrC$}:
\begin{definition}
	Let $\scrC \in \PrLomega_k$ and let $S \in \dAff_k$. An \textit{$S$-family of objects in $\scrC$} is an object in $\scrC_S \coloneqq \scrC \otimes_k \catQCoh(S)$.
\end{definition}

Before delving into a more careful analysis, let us consider the following example:
\begin{example}\label{ex:moduli_of_sheaves_on_varieties}
	Let $X$ be a (possibly derived) scheme over $k$ and set $\scrC \coloneqq \catQCoh(X)$. We further assume that $X$ is quasi-compact and separated. In this case, we have
	\begin{align}
		\catQCoh(X) \simeq \Ind(\catPerf(X)) \ . 
	\end{align}
	In particular, $\catQCoh(X)$ is compactly generated and \cite[Theorem~4.7]{BenZvi_Francis_Nadler_Integral_transforms} supplies a canonical equivalence
	\begin{align}
		\catQCoh(X) \otimes_k \catQCoh(S) \simeq \catQCoh(X \times S) 
	\end{align}
	for every $S \in \dAff_k$. In other words, an $S$-family of objects in $\catQCoh(X)$ is a quasi-coherent sheaf on $X \times S$, matching the usual intuition from algebraic geometry.
\end{example}

The \textit{large moduli of objects of $\scrC$} is the functor
\begin{align}
	\widehat{\calM}_\scrC \colon \dAff_k\op \longrightarrow \mathsf{Spc} 
\end{align}
defined by the rule
\begin{align}
	\widehat{\calM}_\scrC(S) \coloneqq \scrC_S^\simeq \ , 
\end{align}
where $(-)^\simeq$ denotes the maximal $\infty$-groupoid. Since $\scrC$ is compactly generated, Theorem~\ref{thm:main_properties_tensor_product_of_categories}--\eqref{thm:main_properties_tensor_product_of_categories:functor_formula} allows to rewrite
\begin{align}
	\scrC_S = \scrC \otimes_k \catQCoh(S) \simeq \FunR_k(\scrC\op, \catQCoh(S)) \simeq \Fun_k((\scrC^\omega)\op, \catQCoh(S)) \ . 
\end{align}
This formula immediately shows that $\widehat{\calM}_\scrC$ satisfies étale hyperdescent, but as Example~\ref{ex:moduli_of_sheaves_on_varieties} shows this has no chances of being a geometric derived stack. In classical algebraic geometry, the definition of the geometric stack of coherent sheaves \cite[\href{https://stacks.math.columbia.edu/tag/08KA}{\S08KA}]{stacks-project} imposes two restrictions: proper support with respect to the base and flatness. In this \textit{non-commutative} setting, we can formulate analogous conditions, using the language of \textit{integral transforms} and of \textit{$t$-structures}.

\subsection{Integral transforms}\label{subsec:integral-transforms}

Given an $S$-family $M \in \scrC_S$, \cite[Proposition~4.8.1.17]{Lurie_HA} allows to review it as a limit preserving functor $\scrC\op \to \catQCoh(S)$. We denote the restriction of this functor to $(\scrC^\omega)\op$ by
\begin{align}
	\Phi_M \colon (\scrC^\omega)\op \longrightarrow \catQCoh(S) \ , 
\end{align}
and we refer to it as the \textit{integral transform of $M$}.
\begin{example}
	In the setting of Example~\ref{ex:moduli_of_sheaves_on_varieties}, given $\calF \in \catQCoh(X \times S)$, its integral transform is the functor
	\begin{align}
		\Phi_\calF \colon \catPerf(X \times S) \longrightarrow \catQCoh(S) 
	\end{align}
	given by
	\begin{align}
		\Phi_\calF(\calG) \coloneqq \pr_{S,\ast}\big( \calHom_{X \times S}(\calG, \calF) \big) \in \catQCoh(S) \ , 
	\end{align}
	where $\calHom_{X \times S}$ denotes the internal hom in $\catQCoh(X \times S)$ and $\pr_S \colon X \times S \to S$ is the natural projection.
\end{example}

\begin{definition}
	Let $\scrC \in \PrLomega_k$ and let $S \in \dAff_k$. An $S$-family $M \in \scrC_S$ of objects in $\scrC$ is said to be \textit{properly supported} if its integral transform $\Phi_M$ takes values in $\catCohb(S)$.
\end{definition}

This terminology is justified by the following result.
\begin{theorem}
	Let $X$ be a quasi-compact and separated derived scheme locally almost of finite presentation over $\Spec(k)$. Let $S \in \dAff_k$ be a derived affine locally almost of finite presentation and let $\calF \in \catQCoh(X \times S)$. Then, the transform $\Phi_\calF$ takes values in $\catCohb(S)$ if and only if $\calF$ is bounded coherent and its (geometrically defined) support is proper relative to $S$.
\end{theorem}

\begin{proof}
	This is a special case of \cite[Theorem~3.0.2]{BZNP_Integral_transform}.
\end{proof}

Notice that the integral transform not only controls the support of $\calF$, but it equally forces some algebraic finiteness. In the same vein, we introduce the notion of pseudo-perfect family, which is the main finiteness condition on families of objects considered in \cite{TV_Moduli}.
\begin{definition}
	Let $\scrC\in \PrLomega_k$ and let $S \in \dAff_k$. An $S$-family $M \in \scrC_S$ of objects in $\scrC$ is said to be \textit{pseudo-perfect} if its integral transform $\Phi_M$ takes values in $\catPerf(S)$. We denote by $\scrC_S^{\mathsf{ps}}$ the full subcategory of $\scrC_S$ spanned by pseudo-perfect objects.
\end{definition}

With this terminology, Proposition~\ref{prop:basics_smooth_and_proper}--\eqref{prop:basics_smooth_and_proper-6} becomes:
\begin{lemma}\label{lem:pseudo_perfect_and_smoothness}
	Let $A$ be a derived commutative $k$-algebra and let $\scrC \in \PrLomega_A$.
	\begin{enumerate}\itemsep=0.2cm
		\item If $\scrC$ is proper, then any compact object is pseudo-perfect.
		
		\item If $\scrC$ is smooth, then any pseudo-perfect object is compact.
	\end{enumerate}
\end{lemma}

Let $\scrC \in \PrLomega_k$.
We denote by
\begin{align}
	\bfPerfps(\scrC) \colon \dAff_k\op \longrightarrow \Spc
\end{align}
the functor sending $S$ to the $\infty$-groupoid of pseudo-perfect families of objects in $\scrC$. We refer to $\bfPerfps(\scrC)$ as the \textit{moduli of objects of $\scrC$}.

\begin{example}\label{ex:geometric_moduli_of_objects}
	Assume that $X$ is smooth. Then, Lemma~\ref{lem:pseudo_perfect_and_smoothness} guarantees that if $\calF \in \catQCoh(X \times S)$ is pseudo-perfect it is also perfect. It follows from this fact and the above theorem that the underived points of $\bfPerfps(\catQCoh(X))$ correspond to families of properly supported perfect complexes on $X$. On the other hand, if $X$ is proper and underived, then every perfect complex on $X$ is pseudo-perfect but the converse typically does not hold: for instance, when $X = \Spec(k[t]/(t^2))$ and $S = \Spec(k)$, then $\calF = k$ is pseudo-perfect but not perfect.
\end{example}

\subsection{Toën-Vaquié's theorem}\label{subsec:Toen_Vaquie}

We can now state Toën-Vaquié's celebrated geometricity result \cite[Theorem~0.2]{TV_Moduli}.
\begin{theorem}[Toën-Vaquié]\label{thm:Toen_Vaquie_moduli}
	Let $\scrC \in \PrLomega_k$ be of finite type. Then $\bfPerfps(\scrC)$ is a locally geometric derived stack locally of finite type. Furthermore, the tangent complex at a point $x \colon S \to \bfPerfps(\scrC)$ corresponding to a pseudo-perfect family of objects $M_x \in \scrC_S$ is given by
	\begin{align}
		x^\ast \mathbb T_{\bfPerfps(\scrC)} \simeq \mathsf{Hom}_\scrC(M_x, M_x)[1] \ . 
	\end{align}
\end{theorem}

\begin{example}
	Let $X$ be a smooth variety over $\Spec(k)$. Then $\scrC \coloneqq \catQCoh(X) \simeq \IndCoh(X)$ is of finite type, and therefore its moduli of objects $\bfPerfps(\scrC)$ is locally geometric. It follows from Example~\ref{ex:geometric_moduli_of_objects} that $\bfPerfps(\scrC)$ is identified with the moduli stack
	\begin{align}
		\bfMap_{\mathsf{prop}}(X, \bfPerf_k) 
	\end{align}
	parametrizing \textit{properly supported} families of perfect complexes on $X$.
\end{example}

\begin{example}\label{ex:Waldhausen_finite_type}
	Let $\scrC\in \PrLomega_k$ be of finite type and let $n \geqslant 0$ be an integer. Applying Proposition~\ref{prop:new_finite_type_from_old_ones}--\eqref{item:new_finite_type_from_old_ones-2} with $I = \Delta^{n-1}$ implies that the Waldhausen construction $\calS_n \scrC$ is again of finite type (see also \cite[Lemma~4.5]{Porta_Sala_Hall}). In particular, the moduli of objects $\calS_n\bfPerfps(\scrC)\coloneqq\bfPerfps(\calS_n \scrC)$ is locally geometric and locally of finite type. Notice that when $n = 2$, the stack $\calS_2\bfPerfps(\scrC)$ is the moduli stack parametrizing \textit{extensions} of objects in $\scrC$.
\end{example}

The moduli of object is a very large stack. Consider for instance the following example:
\begin{example}
	Let $X$ be a smooth and proper variety over $\Spec(k)$ and set $\scrC \coloneqq \catQCoh(X)$. Observe that there is a canonical morphism
	\begin{align}
		j_X \colon X \longrightarrow \bfPerfps(\scrC) \simeq \bfMap(X,\bfPerf_k) 
	\end{align}
	corresponding to $\Delta_\ast \scrO_X \in \catPerf(X \times X)$. Informally speaking, this map can be thought as sending a point $x \in X$ to the skyscraper $k_x \in \catPerf(X)$. The argument given in \cite[Proposition~5.6]{Toen_Vaquie_Algebraisation} shows that: (i) $j_X$ factors through the open substack $\bfPerfps(\scrC)^{\mathsf{simp}}$ of simple objects; (ii) after composing with the projection to the coarse moduli space of $\bfPerfps(\scrC)^{\mathsf{simp}}$, $j_X$ becomes a Zariski open immersion. It follows that for every other smooth and proper variety $Y$ which is derived equivalent to $X$, one can review $Y$ inside $\bfPerfps(\scrC)$.
\end{example}

A standard technique involves using a $t$-structure on $\scrC$ to cut out smaller substacks of $\bfPerfps(\scrC)$, as we are going to discuss now.

\subsection{Generalities on $t$-structures}

We will not review the theory of $t$-structures (for which we refer the reader to \cite[\S1.2.1]{Lurie_HA}), but focus on a couple of desirable properties of $t$-structures that will be needed later on. 

Let $k$ be a derived commutative ring and let $\scrC \in \PrLomega_k$. We fix a $t$-structure $\tau = (\scrC_{\geqslant 0}, \scrC_{\leqslant 0})$ on $\scrC$ which satisfies the following:

\begin{assumption}\label{assumption:t_structure_filtered_colimits}
	The $t$-structure $\tau$ is compatible with filtered colimits, right complete and either of the following is satisfied:
	\begin{enumerate}[label=(A.\arabic*)]\itemsep=0.2cm
		\item \label{item-1-assumption:t_structure_filtered_colimits} the $t$-structure $\tau$ is left complete;
		
		\item \label{item-2-assumption:t_structure_filtered_colimits} for every faithfully flat morphism $f \colon A \to B$ in $\CAlg_k$, the functor
		\begin{align}
			f^\ast \colon \scrC \otimes_k A \longrightarrow \scrC \otimes_k B 
		\end{align}
		is conservative.
	\end{enumerate}
\end{assumption}

\begin{rem}\label{rem:A1_implies_A2}
	Assumption~\ref{item-1-assumption:t_structure_filtered_colimits} implies Assumption~\ref{item-2-assumption:t_structure_filtered_colimits}, as \cite[Lemma~D.6.4.5 and Proposition~D.6.4.6]{Lurie_SAG} show. \hfill $\triangle$
\end{rem}
We write $\scrC^{\heartsuit}$ for the heart of the $t$-structure $\tau$, and for $M \in \scrC$ we denote by $\pi_n(M) \in \scrC^\heartsuit$ its $n$-th homotopy group. In case of ambiguity, we use the heaver notations $\scrC^{\heartsuit_\tau}$ and $\tensor*[^{\tau}]{\pi}{_n}(M)$.

\subsubsection{Induced $t$-structures}\label{subsubsec:induced_t_structures}

Let $A \in \dCAlg_k$. We write
\begin{align}
	\scrC_A \coloneqq \scrC \otimes_k \Modd_A \ . 
\end{align}
Denote by $f \colon k \to A$ the structural morphism. Then forgetful functor $f_\ast \colon \Modd_A \to \Modd_k$ commutes with colimits and therefore it induces a forgetful functor
\begin{align}
	f_\ast \otimes \id_\scrC \colon \scrC_A \longrightarrow \scrC \ , 
\end{align}
that we simply denote $f_\ast$. This basic construction enjoys several nice properties:
\begin{proposition}\label{prop:projection_formula}
	Let $f \colon k \to A$ be a morphism of derived commutative rings, seen as an object in $\dCAlg_k$. Let $\scrC \in \PrLomega_k$. Then:
	\begin{enumerate}\itemsep=0.2cm
		\item \label{item:conservativity_forgetful} the forgetful functor $f_\ast \colon \scrC_A \to \scrC$ is conservative;
		
		\item \label{item:projection_formula-1} For any $G, F \in \scrC$ and $M \in \Modd_k$, if $G \in \scrC^\omega$ the canonical map
		\begin{align}
			M \otimes_k \Hom_{\scrC}(G,F) \longrightarrow \Hom_{\scrC}(G, M \otimes_k F) 
		\end{align}
		is an equivalence.
		
		\item  \label{item:projection_formula-2} For any $M \in \Modd_A$ and $F \in \scrC$, the canonical map
		\begin{align}\label{eq:projection_formula}
			f_\ast(M) \otimes_k F \longrightarrow f_\ast( M \otimes_A f^\ast(F) )
		\end{align}
		is an equivalence.
	\end{enumerate}
\end{proposition}

\begin{proof}
	Under the equivalence
	\begin{align}
		\scrC_A \simeq \FunR(\Modd_A\op, \scrC) \ , 
	\end{align}
	$f_\ast$ is identified with evaluation at $A$. Since $A$ is a compact stable generator for $\Modd_A$, it follows that $f_\ast$ is conservative.
	This proves \eqref{item:conservativity_forgetful}. For \eqref{item:projection_formula-1}, since both source and target of the map commute with colimits in $M$, one is readily reduced to the case $M = A$, which is obvious. The same argument applies to \eqref{item:projection_formula-2}, modulo the observation that we already made that $f_\ast \colon \scrC_A\to \scrC$ commutes with all colimits.
\end{proof}

Write $(\scrC_A)_{\geqslant 0}$ (resp.\ $(\scrC_A)_{\leqslant 0}$) for the full subcategory of $\scrC_A$ spanned by the objects $M$ such that $\sfU_A(M) \in \scrC$ belongs to $\scrC_{\geqslant 0}$ (resp.\ to $\scrC_{\leqslant 0}$). Since $A$ is itself connective, it is straightforward to check that
\begin{align}
	\tau_A \coloneqq \big( (\scrC_A)_{\geqslant 0}, (\scrC_A)_{\leqslant 0} \big) 
\end{align}
is a $t$-structure on $\scrC_A$ which satisfies again Assumption~\ref{assumption:t_structure_filtered_colimits}. We refer to $\tau_A$ as the \textit{induced $t$-structure}. In what follows, unless stated explicitly, we will consider $\scrC_A$ equipped with the induced $t$-structure $\tau_A$. We use similar notations when working with $\dAff_k$ instead of $\dCAlg_k$ (in particular, if $S \in \dAff_k$, we write $\tau_S$ for the induced $t$-structure on $\scrC_S$).

\begin{corollary}\label{cor:exactness}
	The functor $f_\ast$ is $t$-exact. Moreover, its left adjoint
	\begin{align}
		f^\ast \colon \scrC \longrightarrow \scrC_B \ .
	\end{align}
	is right $t$-exact.
	
	If $f$ is a faithfully flat morphism in $\CAlg$, the functor $f^\ast$ is $t$-exact. If furthermore the $t$-structure $\tau$ is left complete, then $f^\ast$ is conservative as well.
\end{corollary}

\begin{proof}
	The first part of the proof is obvious from the construction of induced $t$-structures.
	
	From the first part of the statement, $f_\ast$ is $t$-exact, while $f^\ast$ is right $t$-exact, while by Proposition~\ref{prop:projection_formula}--\eqref{item:conservativity_forgetful} $f_\ast$ is conservative.
	Therefore, to prove that $f^\ast$ is $t$-exact it is enough to prove that the monad $T_B \coloneqq f_\ast \circ f^\ast$ is $t$-exact.
	In other words, we have to prove that the endofunctor
	\begin{align}
		B \otimes_A - \colon \scrC \longrightarrow \scrC
	\end{align}
	is $t$-exact.
	Since $B$ is flat as $A$-module, Lazard's theorem (see \cite[Theorem~7.2.2.15]{Lurie_HA}) implies that we can represent $B$ as filtered colimit of free $A$-modules of finite rank. Since the functor $A \otimes_A - \simeq \id_\scrC$ is $t$-exact and the $t$-structure on $\scrC$ is compatible with filtered colimits, the conclusion follows.
	The fact that $f^\ast$ is conservative is true either by hypothesis (if Assumption A.2 is satisfied), or as a consequence of \cite[Lemma~D.6.4.5 and Proposition~D.6.4.6]{Lurie_SAG}, if Assumption A.1 is satisfied. This proves the last part of the proposition.
\end{proof}

\subsubsection{Almost perfect objects}

Having a $t$-structure, one can introduce the following generalization of almost perfect modules.
\begin{definition}
	Let $\scrC \in \PrLomega_A$ and let $\tau$ be a $t$-structure satisfying Assumption~\ref{assumption:t_structure_filtered_colimits}. We say that an object $F \in \scrC$ is
	\begin{enumerate}\itemsep=0.2cm
		\item \textit{$\tau$-finitely $n$-presented} if $\tau_{\leqslant n}(F)$ is compact in $\scrC_{\leqslant n}$;
		
		\item \textit{$\tau$-almost perfect} if it is $\tau$-finitely $n$-presented for every $n \in \mathbb Z$.
	\end{enumerate}
	We write $\catAPerf(\scrC,\tau)$ for the full subcategory of $\scrC$ spanned by $\tau$-almost perfect objects.
\end{definition}

\begin{rem}
	Since $\tau$ is right complete, every $\tau$-finitely $n$-presented object is eventually coconnective.\hfill $\triangle$
\end{rem}

In particular, when $\scrC = \catQCoh(X)$ for some derived scheme $X$ and $\tau$ denotes its standard $t$-structure, then $\catAPerf(\scrC,\tau)$ canonically coincides with $\catAPerf(X)$. Almost perfect objects enjoy the following standard list of properties.
\begin{proposition} \label{prop:basic_aperf}
	The following assertions hold:
	\begin{enumerate}\itemsep=0.2cm
		\item \label{prop:basic_aperf:compact_implies_aperf} Every compact object of $\scrC$ is almost perfect.
		
		\item \label{prop:basic_aperf:stability} Almost perfect objects are closed under shifts, finite colimits and retracts in $\scrC$.
		
		\item \label{prop:basic_aperf:homotopy_groups_compact} Let $F \in \catAPerf(\scrC,\tau)$. Then, $\pi_n(F)$ is a compact object in $\scrC^\heartsuit$ for every integer $n \in \Z$.
		
		\item \label{prop:basic_aperf:homotopy_groups_almost_perfect} Let $F \in \scrC$ and assume that $F$ is eventually coconnective.
		If for every $n \in \Z$, the homotopy group $\pi_n(F)$ is an almost perfect object of $\scrC$, then $F \in \catAPerf(\scrC,\tau)$.
	\end{enumerate}
\end{proposition}

\begin{proof}
	Left as an exercise for the reader (see also \cite[\S~C.6.4]{Lurie_SAG}).
\end{proof}

As observed in \cite[Remark~0.17]{Neeman_single_generator}, almost perfect object do not fully depend on the choice of the $t$-structure but rather only on their equivalence class for the relation of being \textit{relatively bounded}, in the following sense:
\begin{definition}\label{def:bounded-structures}
	Let $\tau^{(1)} = (\scrC_{\geqslant 0}^{(1)}, \scrC_{\leqslant 0}^{(1)})$ and $\tau^{(2)} = (\scrC_{\geqslant 0}^{(2)}, \scrC_{\leqslant 0}^{(2)})$ be two $t$-structures on $\scrC$. We say that $\tau_{(2)}$ is \textit{right $\tau_{(1)}$-bounded} (resp.\ \textit{left $\tau_{(1)}$-bounded}) if there exists an integer $n \in \Z$ such that
	\begin{align}
		\scrC_{\geqslant 0}^{(2)} \subseteq \scrC_{\geqslant n}^{(1)} \qquad \big(\text{resp.\ } \scrC_{\leqslant 0}^{(2)} \subseteq \scrC_{\leqslant n}^{(1)}\, \big) \ . 
	\end{align}
\end{definition}

\begin{rem}\label{rem:left_right_boundedness}
	Equivalently, $\tau_{(2)}$ is right $\tau_{(1)}$-bounded if the identity $\id_{\scrC} \colon (\scrC, \tau_{(2)}) \to (\scrC, \tau_{(1)})$ is right $t$-exact up to a shift. In particular, $\tau_{(2)}$ is right $\tau_{(1)}$-bounded if and only if $\tau_{(1)}$ is left $\tau_{(2)}$-bounded.\hfill$\triangle$
\end{rem}

Now, we summarize how the properties of a $t$-structure discussed in the previous sections are inherited to a relatively bounded $t$-structure.
\begin{proposition}\label{prop:two_t_structures_properties}
	Let $\tau_{(1)}$ and $\tau_{(2)}$ be two $t$-structures on $\scrC$. If $\tau_{(2)}$ is right $\tau_{(1)}$-bounded, then
	\begin{align}
		\catAPerf(\scrC, \tau_{(2)}) \subseteq \catAPerf(\scrC, \tau_{(1)}) \ . 
	\end{align}
	If $\tau_{(2)}$ is both left and right $\tau_{(1)}$-bounded then
	\begin{align}
		\catAPerf(\scrC, \tau_{(1)}) = \catAPerf(\scrC, \tau_{(2)}) \ . 
	\end{align}
\end{proposition}

\begin{proof}
	Thanks to Remark~\ref{rem:left_right_boundedness}, both statements follow from \cite[Proposition~C.6.6.13-(1)]{Lurie_SAG} applied to $\id_\scrC$.
\end{proof}

We finish this section by mentioning the following result provides a relation between $t$-structures and Ind-completions.
\begin{proposition}[{\cite[Proposition~2.13]{AGH_t-structure}, \cite[Lemma~C.2.4.3]{Lurie_SAG}}]\label{prop:t-structure-ind-completion}
	Let $\scrD$ be a small stable $\infty$-category with a $t$-structure $(\scrD^{\leqslant 0}, \scrD^{\geqslant 0})$. Then, $\Ind(\scrC)$  inherits a $t$-structure $(\Ind(\scrC)^{\leqslant 0}, \Ind(\scrD)^{\geqslant 0})$, where $\Ind(\scrD)^{\leqslant 0}$ is the essential image of the fully faithful functor $\Ind(\scrD^{\leqslant 0})\to \Ind(\scrC)$, and $\Ind(\scrD)^{\geqslant 0}$ is defined similarly. The $t$-structure $(\Ind(\scrD)^{\leqslant 0}, \Ind(\scrD)^{\geqslant 0})$ is compatible with filtered colimits and such that the inclusion functor $\scrD \to \Ind(\scrD)$ is $t$-exact. Moreover, if the $t$-structure on $\scrD$ is bounded below, then $\Ind(\scrD)$ is right complete.
\end{proposition}
We will call the $t$-structure on $\Ind(\scrC)$ constructed in the above proposition the \textit{induced} $t$-structure. 

\subsubsection{Torsion pairs and tiltings}\label{subsubsec:tiltings}

A $t$-structure can be tweaked by means of a \textit{torsion pair} on its heart via a procedure known as \textit{tilting}. This is a great source of interesting examples that we will be concerned with later on.
\begin{definition}\label{def:torsion_pair}
	Let $\scrA$ be an abelian category. A \textit{torsion pair} on $\scrA$ is a pair $\upsilon = (\scrT, \scrF)$ of full subcategories of $\scrA$ such that:
	\begin{enumerate} \itemsep=0.2cm
		\item for every $T \in \scrT$ and $F \in \scrF$, one has $\Hom_{\scrA}(T, F) = 0$;
		
		\item every $E \in \scrA$ fits into an exact sequence
		\begin{align}
			0\longrightarrow T\longrightarrow E\longrightarrow F\longrightarrow 0 
		\end{align}
		where $T \in \scrT$ and $F \in \scrF$.
	\end{enumerate}
	In this case, we refer to $\scrT$ as the \textit{torsion part} of $\upsilon$ and to $\scrF$ as the \textit{torsion-free part} of $\upsilon$.
\end{definition}

The proof of next lemma consists of standard category theory and we leave it to the reader.
\begin{lemma}\label{lem:torsion_pairs_extensions_mono_epi}
	Let $\scrA$ be an abelian category and let $\upsilon = (\scrT, \scrF)$ be a torsion pair on $\scrA$.
	Then:
	\begin{enumerate}\itemsep=0.2cm
		\item both $\scrT$ and $\scrF$ are closed under extensions;
		
		\item if $T \to T'$ is an epimorphism in $\scrA$ and $T \in \scrT$, then $T' \in \scrT$;
		
		\item \label{item:torsion_pairs_extensions_mono_epi-3} if $F' \to F$ is a monomorphism in $\scrA$ and $F \in \scrF$, then $F' \in \scrF$.
	\end{enumerate}
\end{lemma}

\begin{construction}[Tilting]
	Let $\upsilon = (\scrT, \scrF)$ be a torsion pair on $\scrC^\heartsuit$.
	Define
	\begin{align}
		\tensor*[^{\upsilon}]{\scrC}{_{\geqslant 0}} \coloneqq \Big\{ F \in \scrC_{\geqslant 0} \mid \pi_0(F) \in \scrT \Big\} \quad \text{and} \quad \tensor*[^{\upsilon}]{\scrC}{_{\leqslant 0}} \coloneqq \Big\{ F \in \scrC_{\leqslant 1} \mid \pi_1(F) \in \scrF \Big\} \ . 
	\end{align}
	It is easy to see that the pair $\tau_\upsilon \coloneqq (\tensor*[^{\upsilon}]{\scrC}{_{\geqslant 0}}, \tensor*[^{\upsilon}]{\scrC}{_{\leqslant 0}})$ is again a $t$-structure on $\scrC$. We refer to $\tau_\upsilon$ as the \textit{$\upsilon$-tilted $t$-structure obtained from $\tau$}. We denote by $\tensor*[^{\upsilon}]{\scrC}{^\heartsuit}$ the corresponding heart and by $\tensor*[^{\upsilon}]{\pi}{_i}$ the corresponding homotopy groups.
\end{construction}

\begin{remark}\label{rem:assumption_A_for_tiltings}
	Assume that $\scrC^\heartsuit$ and $\scrT$ are compactly generated, and that the inclusion $\scrT \hookrightarrow \scrC^\heartsuit$ commutes with compact objects. This implies that $\scrF$ is closed under filtered colimits in $\scrC^\heartsuit$. Thus, if in addition the original $t$-structure $\tau$ is compatible with filtered colimits, it follows from the above explicit description of $\tensor*[^{\upsilon}]{\scrC}{_{\leqslant 0}}$ that the tilted $t$-structure $\tau_\upsilon$ is also compatible with filtered colimits.
\end{remark}

\begin{remark}[Relative boundedness]\label{rem:relative_boundedness}
	Let $\tau^{(1)} = (\scrC^{(1)}_{\geqslant 0}, \scrC^{(1)}_{\leqslant 0})$ and $\tau^{(2)} = (\scrC^{(2)}_{\geqslant 0}, \scrC^{(2)}_{\leqslant 0})$ be two $t$-structures on $\scrC$ and write $\scrC^{\heartsuit_1}$ and $\scrC^{\heartsuit_2}$ for the respective hearts. It follows from \cite[Lemma~1.1.2]{Polishchuk_t-structures} that if the inclusions
	\begin{align}\label{eq:relative_boundedness}
		\scrC^{(1)}_{\geqslant 1} \subseteq \scrC^{(2)}_{\geqslant 0} \subseteq \scrC^{(1)}_{\geqslant 0}
	\end{align}
	hold, then
	\begin{align}
		\upsilon \coloneqq \big(\scrC^{\heartsuit_1} \cap \scrC^{\heartsuit_2}, \scrC^{\heartsuit_1} \cap \scrC^{\heartsuit_2}[1] \big) 
	\end{align}
	defines a torsion pair on $\scrC^{\heartsuit_1}$ and $\tau^{(2)} = \tau^{(1)}_\upsilon$. The vice-versa being obvious, it follows that $\tau^{(2)}$ is obtained as a tilting from $\tau^{(1)}$ if and only if the inclusions \eqref{eq:relative_boundedness} hold.
\end{remark}
The following follows straightforwardly from the above remark and Proposition~\ref{prop:t-structure-ind-completion}.
\begin{corollary}\label{cor:completion-torsion-pair}
	Let $\scrD$ be a stable $\infty$-category endowed with a $t$-structure $\sigma$ and let $\upsilon=(\scrT, \scrF)$ be a torsion pair of $\scrD^\heartsuit$. Then, the pair $(\Ind(\scrT), \Ind(\scrF))$ is the canonical torsion pair on $\Ind(\scrD)^\heartsuit$ induced by $\upsilon$, where the heart is with respect to the canonical $t$-structure induced by $\sigma$.  
\end{corollary}

We collect the following result which will be useful later on. First, we recall the notion of Serre subcategory.
\begin{definition}
	Let $\scrA$ be an abelian category. A full subcategory $\scrA'$ of $\scrA$ is a \textit{Serre subcategory}, if for any short exact sequence
	\begin{align}
		0\longrightarrow E'\longrightarrow E \longrightarrow E'' \longrightarrow 0
	\end{align}
	in $\scrA$, we have $E\in \scrA'$ if and only if $E', E''\in \scrA'$. 
\end{definition}

\begin{lemma}\label{lem:tor_Serre}
	Let $\upsilon = (\scrT, \scrF)$ be a torsion pair on $\scrC^\heartsuit$. If $\scrT$ is a Serre subcategory of $\scrC^\heartsuit$, then it is a Serre subcategory of $\scrC^{\heartsuit_\upsilon}$ as well.
\end{lemma}

\begin{proof}
	Since $\scrT$ is the torsion-free part of a torsion pair on $\tensor*[^{\upsilon}]{\scrC}{^\heartsuit}$ by Remark~\ref{rem:relative_boundedness}, Lemma~\ref{lem:torsion_pairs_extensions_mono_epi} shows that it is closed under extensions and subobjects. It is then enough to prove that it is closed under quotients. To prove this, let
	\begin{align}\label{eq:torsion}
		0 \longrightarrow T_1 \longrightarrow T \longrightarrow T_2 \longrightarrow 0 
	\end{align}
	be a short exact sequence in $\tiltedheart$ and assume that $T \in \scrT$. Then, $T_1 \in \scrT$ as well, and passing to the associated long exact sequence with respect to the $t$-structure $\tau$ we find
	\begin{align}
		0 \longrightarrow \calH^{-1}_\tau(T_2) \longrightarrow \calH^{0}_\tau(T_1) \longrightarrow \calH^{0}_\tau(T) \longrightarrow \calH^{0}_\tau(T_2) \longrightarrow 0 \ . 
	\end{align}
	Since $T_2 \in \tiltedheart$ by assumption, $\calH_\tau^{-1}(T_2) \in \scrF$. On the other hand, $T_1 \simeq \calH_\tau^0(T_1)$ is torsion. 
	Therefore the injectivity of the map $\calH_\tau^{-1}(T_2) \to T_1$ together with the assumption that $\scrT$ is a Serre subcategory of $\scrC^\heartsuit$ implies that $\calH_\tau^{-1}(T_2)$ belongs to $\scrT$ as well. It follows that $\calH_\tau^{-1}(T_2) = 0$ and thus, $T_2 \in \scrT$.
	\end{proof}

\subsection{Flatness and openness of $t$-structures}\label{subsec:flatness_and_openness_of_t_structures}

We keep fixing a derived commutative ring $k$ and $\scrC \in \PrLomega_k$. We assume that $\scrC$ is of finite type and we fix a $t$-structure $\tau = (\scrC_{\geqslant 0}, \scrC_{\leqslant 0})$ on $\scrC$ satisfying Assumption~\ref{assumption:t_structure_filtered_colimits}. Given $S \in \dAff_k$, we consider $\scrC_S$ equipped with the induced $t$-structure $\tau_S$. Furthermore, given $N \in \catQCoh(S)$ and $M \in \scrC_S$ we write $N \otimes_S M$ for the object in $\scrC_S$ obtained by the canonical action of $\catQCoh(S)$ on $\scrC_S$.

\begin{definition}\label{def:flatness_wrt_t_structure}
	Let $S \in \dAff_k$. A family $M \in \scrC_S$ is said to be \textit{$\tau$-flat relative to $S$} if for every $N \in \catQCoh^\heartsuit(S)$, one has $N \otimes_S M \in \scrC_S^\heartsuit$.
\end{definition}

\begin{lemma} \label{lem:base_change_tor_amplitude}
	Let $f \colon A \to B$ be a morphism in $\CAlg$ and let $a \le b$ be two integers. Then the induced functor
	\begin{align}
		f^\ast \colon \scrC \longrightarrow \scrC_B 
	\end{align}
	takes $\scrC^{\tau \textrm{-} [a,b]}$ to $\scrC_B^{\tau \textrm{-} [a,b]}$.
\end{lemma}

\begin{proof}
	Recall from Proposition~\ref{prop:projection_formula}--\eqref{item:conservativity_forgetful} that the forgetful functor $f_\ast \colon \scrC_B \to \scrC$ is $t$-exact. Let $M \in \Modd_B^\heartsuit$ and let $F \in \scrC^{\tau \textrm{-} [a,b]}$. We have to check that $f^\ast(F) \otimes_B M$ is in cohomological amplitude $[a,b]$. It is enough to check that the same is true for $f_\ast( f^\ast(F) \otimes_B M )$. Using the projection formula of Proposition~\ref{prop:projection_formula}--\eqref{item:projection_formula-2}, the conclusion follows from the fact that $F$ has tor-amplitude in $[a,b]$ relative to $A$.
\end{proof}

\begin{lemma} \label{lem:flat_objects_descent}
	Let $f \colon A \to B$ be a faithfully flat morphism in $\CAlg$ and let $a \le b$ be two integers. If the $t$-structure $\tau$ is left complete, then an object $F \in \scrC$ belongs to $\scrC^{\tau \textrm{-} [a,b]}$ if and only if $f^\ast(F)$ belongs to $\scrC_B^{\tau \textrm{-} [a,b]}$.
\end{lemma}

\begin{proof}
	The ``only if'' direction follows from Lemma~\ref{lem:base_change_tor_amplitude}.
	
	Let us show the ``if'' part. Since the $t$-structure is left complete, Corollary~\ref{cor:exactness} ensures that the functor $f^\ast \colon \scrC \longrightarrow \scrC_B$ is $t$-exact and conservative. Let $M \in \Modd_A^\heartsuit$. Then $M \otimes_A F$ belongs to $\scrC^{\tau \textrm{-} [a,b]}$ if and only if $f^\ast( M \otimes_A F)$ belongs to $(\scrC_B)^{\tau \textrm{-} [a,b]}$. Since $f^\ast(M) \in \Modd_B^\heartsuit$, the conclusion follows from the equivalence $f^\ast( M \otimes_A F ) \simeq f^\ast(M) \otimes_B f^\ast(F)$.
\end{proof}

\begin{notat}
	We denote by $\bfCohps(\scrC,\tau)$ the substack of the moduli of objects $\bfPerfps(\scrC)$ corresponding to families that are both $\tau$-flat and pseudo-perfect. \hfill $\oslash$
\end{notat}
Lemma~\ref{lem:base_change_tor_amplitude} guarantees that indeed $\bfCohps(\scrC,\tau)$ defines a substack of $\bfPerfps(\scrC)$. In general, it is not geometric unless some stronger condition is imposed on the $t$-structure $\tau$. To formulate the appropriate notion of openness, we introduce the following intermediate definition.
\begin{definition}
	Let $S \in \dAff_k$ and let $M \in \scrC_S$. The \textit{flat locus functor of $M$} is the functor
	\begin{align}
		\Phi_M \colon (\dAff_k)_{/S}\op \longrightarrow \Spc 
	\end{align}
	defined by
	\begin{align}
		\Phi_M\big(T \xrightarrow{f} S\big) \coloneqq \begin{cases}
			\ast & \text{if } f^\ast(M) \in \scrC_T \text{ is } \tau\textrm{-flat relative to $T$}\ , \\
			\emptyset & \text{otherwise}\ .
		\end{cases}
	\end{align}
\end{definition}

The following series of technical lemmas are needed to establish the geometric properties of the flat locus functor, see Proposition~\ref{prop:properties_of_flat_locus} below.

\begin{lemma} \label{lem:flatness_on_pi_0}
	Let $a, b \in \Z$, $a \leqslant b$. Let $A \in \dCAlg_k$ and let $\scrC \in \PrLomega_A$. An object $F \in \scrC$ has tor-amplitude $[a,b]$ relative to $A$ if and only if $\calH^0(A) \otimes_A F \in \scrC_{\calH^0(A)}$ has tor-amplitude $[a,b]$ relative to $\calH^0(A)$.
\end{lemma}

\begin{proof}
	The ``only if'' direction follows directly from Lemma~\ref{lem:base_change_tor_amplitude}.
	
	Let us show the ``if'' direction. Denote by $f \colon A \to \calH^0(A)$ the canonical map. Assume that $f^\ast(F)$ has tor-amplitude $[a,b]$ relative to $\calH^0(A)$ and let $M \in \Modd_A^{\heartsuit}$. Since the $t$-exact functor $f_\ast$ induces an equivalence
	\begin{align}
		\Modd^\heartsuit_A \simeq \Modd^\heartsuit_{\calH^0(A)} \ , 
	\end{align}
	we can write $f_\ast(M)$ instead of $M$. Now, the projection formula of Proposition~\ref{prop:projection_formula}--\eqref{item:projection_formula-2} provides the following canonical equivalence:
	\begin{align}
		f_\ast(M) \otimes_A F \simeq f_\ast( M \otimes_{\calH^0(A)} f^\ast(F) ) \ . 
	\end{align}
	Since $f_\ast \colon \scrC_{\calH^0(A)} \to \scrC$ is $t$-exact, the conclusion immediately follows from the assumption that $f^\ast(F)$ has tor-amplitude $[a,b]$ relative to $\calH^0(A)$.
\end{proof}

\begin{lemma} \label{lem:central_devissage}
	Let $f \colon A \to B$ be a morphism in $\dCAlg_k$. Let $F \in \scrC$ and assume that $f^\ast(F)$ is $B$-flat. Let $\scrE$ be the full subcategory of $\Modd_A^\heartsuit$ spanned by those $M \in \Modd_A^\heartsuit$ such that $M \otimes_A F \in \scrC^\heartsuit$. Then $\scrE$ has the following properties:
	\begin{enumerate}\itemsep=0.2cm
		\item \label{lem:central_devissage:filtered_colimits} it is closed under filtered colimits;
		
		\item \label{lem:central_devissage:extensions} it is closed under retracts and extensions;
		
		\item \label{lem:central_devissage:image} it contains the essential image of $f_\ast \colon \Modd_B^\heartsuit \to \Modd_A^\heartsuit$.
	\end{enumerate}
\end{lemma}

\begin{proof}
	Properties \eqref{lem:central_devissage:filtered_colimits} and \eqref{lem:central_devissage:extensions}  follow directly from the definitions and we leave the details to the reader.
		
	We prove property \eqref{lem:central_devissage:image}. Let $M \in \Modd_B^\heartsuit$. Then Proposition~\ref{prop:projection_formula}--\eqref{item:projection_formula-2} yields
	\begin{align}
		f_\ast(M) \otimes_A F \simeq f_\ast( M \otimes_B f^\ast(F) ) \ . 
	\end{align}
	Since $f^\ast(F)$ is $B$-flat, we see that $M \otimes_B f^\ast(F) \in \scrC_B^\heartsuit$. Since $f_\ast$ is conservative, it follows that $f_\ast(M) \otimes_A F \in \scrC^\heartsuit$. In other words, $f_\ast(M) \in \scrE$.
\end{proof}

\begin{lemma} \label{prop:flatness_on_reduced}
	Let $A \in \dCAlg_k$ and let $I \subseteq \calH^0(A)$ be a nilpotent ideal. Let $f \colon A \to \calH^0(A) / I$ be the canonical map. Then, an object $F \in \scrC$ is $A$-flat if and only if $f^\ast(F) \in \scrC_{\calH^0(A)/I}$ is $\calH^0(A)/I$-flat. In particular, if $A$ is Noetherian, then $F \in \scrC$ is $A$-flat if and only if $A_{\mathsf{red}} \otimes_A F$ is $A_{\mathsf{red}}$-flat.
\end{lemma}

\begin{proof}
	The ``only if'' direction follows directly from Lemma~\ref{lem:base_change_tor_amplitude}.
	
	Let us show the ``if'' direction. In virtue of Lemma~\ref{lem:flatness_on_pi_0} we can replace $A$ by $\calH^0(A)$ and therefore assume that $A$ is underived to begin with. Write $B \coloneqq A / I$. Let $\scrE \subseteq \Modd_A^\heartsuit$ be the full subcategory spanned by those discrete $A$-modules $M$ such that $M \otimes_A F$ belongs to $\scrC^\heartsuit$. We claim that $\scrE = \Modd_A^\heartsuit$.
	Lemma~\ref{lem:central_devissage}--\eqref{lem:central_devissage:filtered_colimits} shows that it is enough to prove that $\scrE$ contains every finitely presented, discrete $A$-module. Fix therefore such a module $M \in (\Modd_A^\heartsuit)^{\mathsf{fp}}$. We prove by induction on $n \ge 1$ that
	\begin{align}
		(M / I^n M) \otimes_A F 
	\end{align}
	belongs to $\scrC^\heartsuit$. Since $I$ is nilpotent, there exists an integer $n \gg 0$ such that $I^n = 0$, so that the conclusion will follow.
	When $n = 1$, $M / I M$ is supported on $A / I = B$. Indeed, if we write $N$ for $M / I M$ equipped with its canonical $B$-module structure, then we have
	\begin{align}
		M / I M \simeq f_\ast(N) \ . 
	\end{align}
	Therefore, Lemma~\ref{lem:central_devissage}--\eqref{lem:central_devissage:image} implies that $M / I M \in \scrE$. This proves the basis of the induction. For the induction step, consider the short exact sequence
	\begin{align}
		0 \longrightarrow I^n M / I^{n+1} M \longrightarrow M / I^{n+1} M \longrightarrow M / I^n M \longrightarrow 0 \ . 
	\end{align}
	Lemma~\ref{lem:central_devissage}--\eqref{lem:central_devissage:extensions} shows that it is enough to prove that both $I^n M / I^{n+1} M$ and $M / I^n M$ belong to $\scrC^\heartsuit$. For the latter, it is enough to apply the inductive hypothesis. For the former, it is enough to observe that $I \cdot ( I^n M / I^{n+1} M ) = 0$. In other words, $I^n M / I^{n+1} M$ can be written as $f_\ast(N)$ for some $N \in \Modd_B$. Thus, the conclusion follows once again from Lemma~\ref{lem:central_devissage}--\eqref{lem:central_devissage:image}.
\end{proof}

\begin{proposition}\label{prop:properties_of_flat_locus}
	Let $S \in \dAff_k$ and let $M \in \scrC_S$. The derived prestack $\Phi_M$ satisfies étale hyperdescent. Furthermore, the canonical map $i_M \colon \Phi_M \to S$ is $(-1)$-truncated, nilcomplete, infinitesimally cohesive and formally étale.
\end{proposition}

\begin{proof}
	It is clear from the definition that $i_M$ is $(-1)$-truncated. For nilcompleteness, we have to prove that for every $\Spec(B) \to \Spec(A)$, the canonical map
	\begin{align}
		\Phi_M(B) \longrightarrow \lim_{n \leqslant 0} \Phi_M(\tau^{\geqslant n}(B))
	\end{align}
	is an equivalence. However, it follows from Lemma~\ref{lem:flatness_on_pi_0} that $B \otimes_A F$ is $B$-flat if and only if $\calH_\tau^0(B) \otimes_A F$ is $\calH_\tau^0(B)$-flat. Thus $\Phi_F(B) = \ast$ if and only if $\Phi_F(\calH_\tau^0(B)) = \ast$, whence the conclusion. The same reasoning applies to infinitesimal cohesiveness, using this time Proposition~\ref{prop:flatness_on_reduced}. Finally, let $\Spec(B') \to \Spec(B)$ be a square-zero extension and let
	\begin{align}
		\begin{tikzcd}[ampersand replacement=\&]
			\Spec(B') \arrow{d} \arrow{r} \& \Phi_F \arrow{d}{i_F} \\
			\Spec(B) \arrow{r} \arrow[dashed]{ur} \& \Spec(A)
		\end{tikzcd}
	\end{align}
	be a commutative diagram. This implies that $B' \otimes_A F$ is $B'$-flat. Proposition~\ref{prop:flatness_on_reduced} implies that $B \otimes_A F$ is $B$-flat as well, and therefore that there exists a unique way to solve the above lifting problem. Thus $i_F$ is formally étale.
\end{proof}

\begin{corollary}\label{cor:flat_locus}
	Let $S \in \dAff_k$ and let $M \in \scrC_S$. Then $\Phi_M$ is representable by a geometric derived stack locally almost of finite type if and only if the map $i_M \colon \Phi_M \to S$ is an open Zariski immersion. In particular, if this is the case then $\Phi_M$ is a derived scheme.
\end{corollary}

\begin{defin}\label{def:openness-flatness}
	We say that the $t$-structure $\tau$ \textit{universally satisfies openness of flatness} if for every $S \in \dAff_k$ and every pseudo-perfect family $M \in \scrC_S$, the flat locus functor $\Phi_M$ is an open subscheme of $S$. \hfill $\oslash$
\end{defin}

\begin{remark}
	Note that the above definition is a non-commutative (and more general) analog of \cite[Definition~10.4]{BLMNPS_Stability}.
\end{remark}

\begin{proposition}\label{prop:openness}
	The $t$-structure $\tau$ universally satisfies openness of flatness if and only if the structural morphism
	\begin{align}\label{eq:structural-morphism}
		\bfCohps(\scrC,\tau) \longrightarrow \bfPerfps(\scrC) 
	\end{align}
	is representable by open Zariski immersions. In particular, when this is the case, $\bfCohps(\scrC,\tau)$ is a geometric derived stack locally of finite presentation over $k$.
\end{proposition}

\begin{proof}
	Let $S \in \dAff_k$ be an affine derived scheme and let $S \to \bfPerfps(\scrC)$ be a morphism classifying a pseudo-perfect family $M \in \scrC_S$. Unraveling the definitions, we see that the square
	\begin{align}
		\begin{tikzcd}[ampersand replacement=\&]
			\Phi_M \arrow{r} \arrow{d} \& \bfCohps(\scrC,\tau) \arrow{d} \\
			S \arrow{r} \& \bfPerfps(\scrC)
		\end{tikzcd} 
	\end{align}
	is a derived pullback square. The conclusion then follows from Corollary~\ref{cor:flat_locus}.
\end{proof}

When tilting a $t$-structure $\tau$ by a torsion-pair $\upsilon$, it is desirable to know if the result $\tau_\upsilon$ still satisfies openness of flatness. This question has been addressed in \cite[Appendix~A]{AB_Moduli}, and the discussion there can be summarized as follows.

\begin{construction} \label{construction:torsion_pair_substacks}
	Given $(\scrC,\tau)$ and a torsion pair $\upsilon = (\scrT,\scrF)$ in $\scrC^\heartsuit$, we define the derived stack $\bfCoh_{\scrT}(\scrC,\tau)$ as the fiber product
	\begin{align}
		\begin{tikzcd}[ampersand replacement=\&]
			\bfCoh_{\scrT}(\scrC,\tau) \arrow{r} \arrow{d} \& \bfCohps(\scrC,\tau_\upsilon) \arrow{d} \\
			\bfCohps(\scrC,\tau) \arrow{r} \& \bfPerfps(\scrC)
		\end{tikzcd} \ .
	\end{align}
	Let $[1] \colon \bfPerfps(\scrC) \to \bfPerfps(\scrC)$ be the morphism corresponding to the shift by $1$ in $\scrC$. Then we define $\bfCoh_{\scrF}(\scrC,\tau)$ as the fiber product
	\begin{align}
		\begin{tikzcd}[ampersand replacement=\&]
			\bfCoh_{\scrF}(\scrC,\tau) \arrow{rr} \arrow{d} \& \& \bfCohps(\scrC,\tau_\upsilon) \arrow{d} \\
			\bfCohps(\scrC,\tau) \arrow{r} \& \bfPerfps(\scrC) \arrow{r}{[1]} \& \bfPerfps(\scrC)
		\end{tikzcd} \ .
	\end{align}
\end{construction}

\begin{defin}\label{def:open_torsion_pair}
	We say that a torsion pair $\upsilon = (\scrT, \scrF)$ in $\scrC^\heartsuit$ is \textit{open} if the morphisms
	\begin{align}
		\bfCoh_\scrT(\scrC,\tau) \longrightarrow \bfCohps(\scrC,\tau) \quad \text{and} \quad \bfCoh_\scrF(\scrC,\tau) \longrightarrow \bfCohps(\scrC,\tau)
	\end{align}
	are representable by Zariski open immersions.
\end{defin}

\begin{warning}\label{warning:families_of_torsion_objects}
	Let $\upsilon = (\scrT,\scrF)$ be an open torsion pair on $\scrC^\heartsuit$. Let $S \in \Aff_k$ be an \textit{un}derived affine scheme. By construction, the inclusions
	\begin{align}
		(\scrC_S)_{\geqslant 1} \subseteq \tensor*[^{\upsilon}]{(\scrC_S)}{_{\geqslant 0}} \subseteq (\scrC_S)_{\geqslant 0} 
	\end{align}
	hold, and therefore Remark~\ref{rem:relative_boundedness} implies that the induced $t$-structure $(\tau_\upsilon)_S$ is obtained as the tilting of $\tau_S$ by the torsion pair $\upsilon_S = (\scrC_S^{\heartsuit} \cap \scrC_S^{\heartsuit_\upsilon}, \scrC_S^{\heartsuit} \cap \scrC_S^{\heartsuit_\upsilon}[1])$. Since $S$ is underived, an element $M \in \scrC_S$ classifying a point of $\bfCoh_\scrT(\scrC,\tau)$ automatically belongs to the torsion part of $\upsilon_S$. However, the reader should be aware that the converse is typically false: in order for an object in the torsion part of $\upsilon_S$ to classify an $S$-point of $\bfCoh_\scrT(\scrC,\tau)$ it has to be \textit{both} $\tau_S$-flat and $(\tau_S)_\upsilon$-flat. This is always the case when $S$ is the spectrum of a field, but it is otherwise typically false.
\end{warning}

\begin{remark}
	Alternatively, an element $M \in \scrC_S$ classifies a point of $\bfCoh_\scrT(\scrC,\tau)$ if and only if it is $\tau_S$-flat and for every morphism $\Spec(K) \to S$ with $K$ a field, $M_K$ is torsion with respect to the induced torsion pair $\upsilon_K$.
\end{remark}

\begin{remark}[{Comparison with \cite[Definition~A.2]{AB_Moduli}}]
	Let $X$ be a proper scheme of finite presentation over a base field $\Spec(K)$ and take $\scrC \coloneqq \catQCoh(X)$. The above definition of the stacks $\bfCoh_\scrT(\scrC,\tau)$ and $\bfCoh_\scrF(\scrC,\tau)$ differs from \cite[Definition~A.2]{AB_Moduli}. In \textit{loc.\ cit.}, the torsion (resp.\ torsion-free) substack of $\bfCoh(\scrC,\tau)$ parametrizes the flat families of coherent sheaves whose geometric fibers are torsion (resp.\ torsion-free). As it follows from Warning~\ref{warning:families_of_torsion_objects}, the families parametrized by our $\bfCoh_\scrT(\scrC,\tau)$ satisfy this condition. In general, the converse is unclear; however if the substack of torsion (resp.\ torsion-free) sheaves in the sense of \cite[Definition~A.2]{AB_Moduli} are open (in the sense of \textit{loc.\ cit.}), then it follows from \cite[Theorem A.8]{AB_Moduli} that $\tau_\upsilon$ universally satisfies openness of flatness; this implies that our $\bfCoh_\scrT(\scrC,\tau)$ (resp.\ $\bfCoh_\scrF(\scrC,\tau)$) is open inside $\bfCoh(\scrC,\tau)$. Thus, in this situation, we have two open substacks of $\bfCoh(\scrC,\tau)$ parametrizing torsion (resp.\ torsion-free) sheaves, arising respectively from Construction~\ref{construction:torsion_pair_substacks} and from \cite[Definition~A.2]{AB_Moduli}.
	As they agree on closed points by Warning~\ref{warning:families_of_torsion_objects}, it follows that they coincide.
\end{remark}

The same argument given in \cite[Theorem A.8]{AB_Moduli} implies:
\begin{proposition}\label{prop:openness-tilted}
	Assume that $t$-structure $\tau$ universally satisfies openness of flatness and the torsion pair $\upsilon = (\scrT, \scrF)$ in $\scrC^\heartsuit$ is open. Then, $\tau_\upsilon$ universally satisfies openness of flatness.
\end{proposition}

\begin{rem}\label{rem:torsion-pair-openness}
	If we assume that the tilted $t$-structure $\tau_\upsilon$ universally satisfies openness of flatness, by Proposition~\ref{prop:openness}, the derived stack $\bfCohps(\scrC,\tau_\upsilon)$ is a geometric derived stack locally of finite presentation over $k$, which is open inside $\bfPerfps(\scrC)$. Moreover, the natural morphisms $\bfCoh_{\scrT}(\scrC,\tau) \to \bfCohps(\scrC,\tau_\upsilon)$ and $[1] \colon \bfCoh_{\scrF}(\scrC,\tau) \to \bfCohps(\scrC,\tau_\upsilon)$ are representable by open immersions. \hfill $\triangle$
\end{rem}

\subsection{Quot-schemes and properness}

We keep fixing a field of characteristic zero $k$, $\scrC \in \PrLomega_k$ and a $t$-structure $\tau = (\scrC_{\geqslant 0}, \scrC_{\leqslant 0})$ satisfying Assumption~\ref{assumption:t_structure_filtered_colimits}. We also assume that $\scrC$ is of finite type and that $\tau$ universally satisfies openness of flatness. Under these conditions, the derived stack $\bfCohps(\scrC,\tau)$ is a geometric derived stack. We define $\calS_2 \bfCohps(\scrC,\tau)$ as the stack parametrizing extensions of pseudo-perfect and $\tau$-flat families of objects in $\scrC$; formally speaking, this is the fiber product
\begin{align}
	\begin{tikzcd}[ampersand replacement=\&]
		\calS_2 \bfCohps(\scrC,\tau) \arrow{r} \arrow{d} \& \calS_2 \bfPerfps(\scrC) \arrow{d}{\partial_0 \times \partial_1 \times \partial_2} \\
		\bfCohps(\scrC,\tau)^{\times 3} \arrow{r} \& \bfPerfps(\scrC)^{\times 3}
	\end{tikzcd} \ .
\end{align}
More informally, whenever $S \in \dAff_k$ is a test scheme, an $S$-point of $\calS_2 \bfCohps(\scrC,\tau)$ is a fiber sequence
\begin{align}
	M_{01} \longrightarrow M_{02} \longrightarrow M_{12} 
\end{align}
in $\scrC_S$, where all the objects are pseudo-perfect and $\tau$-flat. In particular, if $S$ is \textit{un}derived, then $M_{ij} \in \scrC_S^\heartsuit$ and therefore the above fiber sequence is in fact a short exact sequence in the abelian category $\scrC_S^\heartsuit$.This justifies the following:
\begin{definition}
	Let $S \in \dAff_k$ and let $M \in \scrC_S$ be a $\tau$-flat and pseudo-perfect object. The \textit{derived Quot associated to $M$} is the fiber product
	\begin{align}
		\begin{tikzcd}[ampersand replacement=\&]
			\bfQuot_{M}(\scrC,\tau) \arrow{r} \arrow{d} \& \calS_2 \bfCohps(\scrC, \tau) \arrow{d}{\partial_1} \\
			S \arrow{r}{M} \& \bfCohps(\scrC,\tau) 
		\end{tikzcd}  \ .
	\end{align}
\end{definition}

In general, these Quot schemes need not to be proper, although it is of course a very desirable property that the $t$-structure might have:
\begin{definition}\label{def:t_structure_proper}
	We say that the $t$-structure $\tau$ is \textit{weakly proper} if for every $S \in \dAff_k$ and every $\tau$-flat and pseudo-perfect object $M \in \scrC_S$, the associated $\bfQuot_M(\scrC,\tau)$ satisfies the strong existence part and the uniqueness part of the valuative criterion of properness (in the sense of \cite[Definition~11.8]{BLMNPS_Stability}).
\end{definition}

\begin{remark}
	In the geometric context considered in \cite[\S11]{BLMNPS_Stability}, the authors proved that $\bfQuot_M(\scrC,\tau)$ satisfies the strong existence and the uniqueness parts of the valuative criterion (cf.\ Definition~11.8 in \textit{loc.\ cit.}) under a somehow strong assumption on $\tau$ (see Definition~6.15 in \textit{loc.\ cit.}), which is satisfied by those t-structures induced by Bridgeland's stability conditions in the sense of Definition~20.5 of \textit{loc.\ cit.} (see Corollary~20.10 in \textit{loc.\ cit.}).
\end{remark}

In the geometric setting, to check that the standard $t$-structure is weakly proper in the above sense, the classical argument of \cite[\href{https://stacks.math.columbia.edu/tag/0DM4}{Tag 0DM4}]{stacks-project} applies. This argument has been adapted to the non-commutative setting as early as in \cite[Lemma E3.3]{AZ_Hilbert}, and it has been revisited by Toën and Vaquié in \cite[Proposition 4.1]{TV_Points}. In order to state it, we need to review a couple of extra properties of $t$-structures: \textit{coherence} and \textit{noetherianity}.

\subsubsection{Coherence}

The notion of coherence for a $t$-structure is a non-commutative analogue of the more classical notion of coherence for associative rings. In general it is not true that $\catAPerf(\scrC,\tau)$ inherits a $t$-structure. This happens exactly when the $t$-structure is \textit{coherent}:
\begin{definition}
	A $t$-structure $\tau$ on $\scrC$ satisfying Assumption~\ref{assumption:t_structure_filtered_colimits} is said to be \textit{coherent} if every compact object of $\scrC^\heartsuit$ is almost perfect in $\scrC$ with respect to $\tau$.
\end{definition}

\begin{lemma}
	The following statements are equivalent:
	\begin{enumerate}\itemsep=0.2cm
		\item \label{item:1} the $t$-structure $\tau$ is coherent;
		
		\item \label{item:2} for every $F \in \catAPerf(\scrC,\tau)$ and every integer $n \in \Z$, $\pi_n(F)$ is almost perfect as an object of $\scrC$;
		
		\item \label{item:3} the $t$-structure $\tau$ induces a $t$-structure on $\catAPerf(\scrC,\tau)$;
		
		\item \label{item:4} the prestable $\infty$-category $\scrC_{\geqslant 0}$ is coherent in the sense of \cite[Definition~C.6.5.1]{Lurie_SAG}.
	\end{enumerate}
	Furthermore, if $\tau$ is coherent, then an eventually connective object $F$ is almost perfect if and only if for every integer $n \in \Z$, the homotopy group $\pi_n(F)$ is compact as object of $\scrC^{\heartsuit}$.
\end{lemma}

\begin{proof}
	The implications \eqref{item:1} $\Rightarrow$ \eqref{item:2} $\Leftrightarrow$ \eqref{item:4} are easy exercises.
	
	To prove that \eqref{item:2} implies \eqref{item:1}, let $F \in (\scrC^\heartsuit)^\omega$. Since $\scrC$ is compactly generated, it follows that there exists a compact object $G \in \scrC^\omega$ such that $F \simeq \pi_0(F)$ is a retract of $\pi_0(G)$, so that Proposition~\ref{prop:basic_aperf}--\eqref{item:2} \& \eqref{item:3} implies that $F$ is itself almost perfect with respect to $\tau$. As for the equivalence \eqref{item:1} $\Leftrightarrow$ \eqref{item:4}, one should first observe that condition (b) in \cite[Definition~C.6.5.1]{Lurie_SAG} is automatically satisfied because $\scrC$ is assumed to be compactly generated. So we have to argue that \eqref{item:1} is equivalent to condition (a) in \cite[Definition~C.6.5.1]{Lurie_SAG}. This readily follows combining the implication \eqref{item:1} $\Rightarrow$ \eqref{item:2}, Proposition~\ref{prop:basic_aperf}--\eqref{item:4} and \cite[Proposition~C.6.4.5]{Lurie_SAG}.
\end{proof}

\begin{example}
	Let $R$ be a connective $\E_1$-ring. Then the standard $t$-structure on $\mathsf{LMod}_R$ is coherent if and only if $R$ is left coherent in the sense of \cite[Definition 7.2.4.16]{Lurie_HA}. This follows from the above lemma and \cite[Proposition 7.2.4.18]{Lurie_HA}. In particular, if $X$ is a noetherian derived scheme, then the standard $t$-structure on $\catQCoh(X)$ is coherent.
\end{example}

\subsubsection{Noetherianity}

First, we recall the notion of Grothendieck abelian category.
\begin{definition}[{\cite[Definition 1.3.5.1]{Lurie_HA}}]
	Let $\scrA$ be an abelian category. We say that $\scrA$ is \textit{Grothendieck} if it is presentable and the collection of monomorphisms in $\scrA$ is closed under filtered colimits.
\end{definition}

Recall now the notion of local noetherianity for abelian categories.
\begin{definition}
	Let $\scrA$ be a Grothendieck abelian category. We say that $\scrA$ is \textit{locally noetherian} if it is compactly generated and any subobject of a compact object is itself compact.
\end{definition}
We refer the reader to \cite[\S~C.6.8]{Lurie_SAG} for a more thorough discussion of the notion of noetherianity. Notice that the notion of local noetherianity introduced in \cite[Definition C.6.8.5]{Lurie_SAG} does not coincide a priori with the one defined above. Nevertheless, the two definitions do agree, as a consequence of Corollaries~C.6.8.8 \& C.6.8.9 in \textit{loc.\ cit.}

\begin{definition}
	Let $\scrC \in \PrLomega_k$ and let $\tau$ be a $t$-structure on $\scrC$ satisfying Assumption~\ref{assumption:t_structure_filtered_colimits}.
	We say that the pair $(\scrC,\tau)$ is \textit{locally noetherian} if it is coherent and $\scrC^\heartsuit$ is locally noetherian.
\end{definition}

For the record, let us recall:
\begin{proposition}[Non-commutative Hilbert's basis theorem]
	Let $A \in \dCAlg_k$ be a derived commutative $k$-algebra and let $\scrC \in \PrLomega_A$ be a compactly generated $A$-linear presentable $\infty$-category.
	Let $\tau$ be a noetherian $t$-structure satisfying Assumption~\ref{assumption:t_structure_filtered_colimits}.
	If $f \colon A \to B$ is almost of finite presentation in $\dCAlg_k$, then the induced $t$-structure $\tau_B$ on $\scrC_B \coloneqq \scrC \otimes_A B$ is again noetherian.
\end{proposition}

\begin{proof}
	This is a special case of \cite[Proposition D.5.6.1]{Lurie_SAG}.
\end{proof}

\begin{remark}
	If $(\scrC, \tau)$ is only assumed to be coherent instead of locally noetherian, \cite[Proposition D.5.5.1]{Lurie_SAG} shows that the above statement holds provided that $f$ is in addition assumed to be quasi-finite.
\end{remark}

\subsubsection{The properness criterion}\label{subsubsec:properness}

Having introduced this terminology, we can finally state:
\begin{theorem}[Toën-Vaquié]\label{prop:TV_properness}
	Let $\scrC \in \PrLomega_k$ be smooth and proper and let $\tau$ be a $t$-structure on $\scrC$ satisfying Assumption~\ref{assumption:t_structure_filtered_colimits}. Assume that $(\scrC, \tau)$ is locally noetherian and that it universally satisfies openness of flatness. Then, $(\scrC, \tau)$ is also weakly proper.
\end{theorem}

\begin{proof}
	See \cite[Proposition~4.1]{TV_Points}.
\end{proof}

Unfortunately, for the main applications that we have in mind, the above proposition is not sufficient. The reason is that the operation of tilting often destroys noetherianity, as the following example shows.

\begin{example}[Tilting destroys noetherianity]\label{ex:non_noetherian_tilting}
	Let $A$ be a DVR with fraction field $K$. Since $A$ is noetherian, the standard $t$-structure on $\Modd_A$ is noetherian as well. Inside $\Modd_A^\heartsuit$, consider
	\begin{align}
		(\Modd_A^\heartsuit)_{\mathsf{tors}} \coloneqq \big\{M \in \Modd_A^\heartsuit \mid M \otimes_A K \simeq 0\big\} \ . 
	\end{align}
	This is the torsion part of a torsion pair $\upsilon$. The corresponding torsion-free part $(\Modd_A^\heartsuit)_{\mathsf{t.f.}}$ consists of those $M \in \Modd_A^\heartsuit$ for which multiplication by a uniformizer $\varpi$ of $A$ is injective. Observe that the object $K / A$ is torsion, while both $A$ and $K$ are torsion-free. Therefore, starting with the short exact sequence
	\begin{align}
		0 \longrightarrow A \longrightarrow K \longrightarrow K / A \longrightarrow 0 
	\end{align}
	in $\Modd_A^\heartsuit$ and rotating it, we obtain a fiber sequence
	\begin{align}
		K / A \longrightarrow A[1] \longrightarrow K[1] 
	\end{align}
	in $\Modd_A$. Since all the objects are in $\tensor*[^{\upsilon}]{\Modd}{^{\heartsuit}_{A}}$, this is a short exact sequence in the perverse heart. In particular, the map $K / A \to A[1]$ is \textit{injective} in the perverse heart. It is easy to see that $A[1]$ is compact in $\tensor*[^{\upsilon}]{\Modd}{^{\heartsuit}_{A}}$ and that $K/A$ is not. In other words, $\tensor*[^{\upsilon}]{\Modd}{^{\heartsuit}_{A}}$ is not noetherian, and therefore the tilted $t$-structure $\tau_\upsilon$ is also not noetherian.
\end{example}

We conclude with an explicit example that will be useful later on.
\begin{example}\label{ex:Kollar_husks}
	Let $X$ be a smooth and proper scheme. For an integer $m \geqslant 0$, say that $F \in \catQCoh^\heartsuit(X)$ is \textit{supported in dimension $\leqslant m$} if for every coherent subsheaf $\calF' \subset \calF$ one has
	\begin{align}
		\dim \mathsf{supp}(\calF') \leqslant m \ . 
	\end{align}
	Write $\scrT_m$ for the full subcategory of $\catQCoh^\heartsuit(X)$ spanned by sheaves supported in dimension $\leqslant m$, and write $\calF_m$ for the right orthogonal of $\scrT_m$ inside $\catQCoh^\heartsuit(X)$. It is easy to check that $\upsilon_m \coloneqq (\scrT_m, \calF_m)$ form a torsion pair on $\catQCoh^\heartsuit(X)$. It was observed in \cite[Remark~2.14]{Toda_Stable} that the heart of the corresponding tilted $t$-structure $\tensor*[^{\upsilon_m}]{\tau}{}$ is not noetherian.
\end{example}

Using the theory of Koll\'ar's \textit{husks} \cite{Kollar_Hulls_and_husks}, we can nevertheless prove that $\tensor*[^{\upsilon_m}]{\tau}{}$ is close to be weakly proper:
\begin{proposition}\label{prop:Kollar_husks}
	Fix $S \in \Aff_k$ and $F \in (\calF_m)_S$. Then, the associated $\bfQuot_F(\catQCoh(X), \tensor*[^{\upsilon_m}]{\tau}{})$ satisfies the strong existence part and the uniqueness part of the valuative criterion of properness.
\end{proposition}

\begin{proof}
	This follows combining \cite[Theorem~3.6 and Proposition~3.8]{Lin_Wang_Xia_Tilted_quot}, which rely on \cite{Kollar_Hulls_and_husks}.
\end{proof}

\section{COHAs, CatHAs, and their representations arising from $m$-flags}\label{sec:flag_action}

In this section, we construct cohomological and categorical Hall algebras and their representations via stacks of $m$-flags in greater generality (encompassing the geometric framework developed in \cite{Porta_Sala_Hall}). We shall make use of the theory of \textit{(relative) $2$-Segal spaces}, introduced in Part~\ref{part:Segal}, which we apply to the moduli stack of pseudo-perfect objects.

We work over a base ring $k$.
We fix a motivic formalism $\bfD^\ast$, an algebra of coefficients $\calA \in \CAlg(\bfD^\ast(\Spec(k)))$, and an abelian subgroup $\Gamma \subseteq \Pic(\bfD^\ast(\Spec(k)))$, and we make the following:

\begin{assumption}\label{assumption:motivic_formalism}
	The algebra of coefficients $\calA$ is oriented and the abelian subgroup $\Gamma$ is closed under Thom twists in $\Pic(\bfD^\ast(\Spec(k)))$.
\end{assumption}

\subsection{2-Segal spaces and modules via the moduli of pseudo-perfect objects} \label{subsec:m_flags}

Fix a base ring $k$ and a $k$-linear compactly generated stable $\infty$-category $\scrC\in \PrLomega_k$. 

Fix an integer $m \geqslant 1$. We set
\begin{align}
	\bfFlagPerfps^{(m)}(\scrC) \coloneqq \calS_m \bfPerfps(\scrC) \ , 
\end{align}
where the latter is introduced in Example~\ref{ex:Waldhausen_finite_type}, and we refer to this as the \textit{derived stack of $m$-flags of pseudo-perfect objects} in $\scrC$.

\begin{remark}\label{rem:flag_action_coordinates}
	Unraveling the definitions, we see that the derived stack $\bfFlagPerfps^{(m)}(\scrC)$ can informally be described as the functor sending $A \in \CAlg_k$ to diagrams of the form
	\begin{align}\label{eq:flags}
		\begin{tikzcd}[ampersand replacement=\&]
			0 \arrow{r} \& F_{0,1} \arrow{r} \arrow{d} \& F_{0,2} \arrow{r} \arrow{d} \& \cdots \arrow{r} \& F_{0,m-1} \arrow{r} \arrow{d} \& F_{0,m} \arrow{d} \\
			{} \& 0 \arrow{r} \& F_{1,2} \arrow{r} \& \cdots \arrow{r} \& F_{1,m-1} \arrow{r} \arrow{d} \& F_{1,m} \arrow{d} \\
			\& \& \& \ddots \& \vdots \arrow{d} \& \vdots \arrow{d} \\
			\& \& \& \& F_{m-2,m-1} \arrow{d} \arrow{r} \& F_{m-2,m} \arrow{d} \\
			\& \& \& \& 0 \arrow{r} \& F_{m-1,m} \arrow{d} \\
			\& \& \& \& \& 0 
		\end{tikzcd} \ ,
	\end{align}
	where every square is a pullback in $\Fun_A\big((\scrC^\omega)\op, \catPerf(A)\big)$. We refer to such a diagram as an \textit{$m$-flag of pseudo-perfect objects of $\scrC$}.
	
	When $m=2$, further unraveling the definitions shows that $\bfFlagPerfps^{(2)}(\scrC)=\calS_2 \bfPerfps(\scrC)$ parametrizes diagrams $\F$ of the form
	\begin{align}\label{eq:extensions}
		\begin{tikzcd}[ampersand replacement=\&]
			0 \arrow{r} \& F_{0,1} \arrow{r} \arrow{d} \& F_{0,2} \arrow{d} \\
			\& 0 \arrow{r} \& F_{1,2} \arrow{d} \\
			\& \& 0
		\end{tikzcd} \ , 
	\end{align}
	where the central square is asked to be a pullback.
	In other words, we can identify $\bfFlagPerfps^{(2)}(\scrC)=\calS_2 \bfPerfps(\scrC)$ with the stack $\bfPerfpsext(\scrC)$ parametrizing \textit{extensions} of pseudo-perfect objects in $\scrC$.
	With respect to this representation and to the canonical identification $\calS_1 \bfPerfps(\scrC) \simeq \bfPerfps(\scrC)$, we see that\footnote{Here and in what follows, $\partial_\ast$ denotes the $\ast$-th face map.}
	\begin{align}\label{eq:partial}
		\partial_0(\F) = F_{1,2}\ , \qquad \partial_1(\F) = F_{0,2}\ , \qquad \partial_2(\F) = F_{0,1} \ . 
	\end{align}
\end{remark}

We further fix an $(m-1)$-flag $\sfV \in \bfFlagPerfps^{(m-1)}(\scrC)$, which we represent as the following diagram:
\begin{align}
	\begin{tikzcd}[ampersand replacement=\&]
		0 \arrow{r} \& V_{0,1} \arrow{r} \arrow{d} \& V_{0,2} \arrow{r} \arrow{d} \& \cdots \arrow{r} \& V_{0,m-1} \arrow{d} \\
		{} \& 0 \arrow{r} \& V_{1,2} \arrow{r} \& \cdots \arrow{r} \& V_{1,m-1} \arrow{d} \\
		\& \& \& \ddots \& \vdots \arrow{d} \\
		\& \& \& \& V_{m-2,m-1} \arrow{d} \\
		\& \& \& \& 0
	\end{tikzcd} \ .
\end{align}
The boundary map $\partial_m \colon [m] \to [m-1]$ in $\mathbf \Delta\op$ induces a forgetful functor
\begin{align}
	\partial_m \colon \bfFlagPerfps^{(m)}(\scrC) \longrightarrow \bfFlagPerfps^{(m-1)}(\scrC) \ ,
\end{align}
which indeed forgets the last column from the right-hand-side.

\begin{definition}\label{def:V_flags}
	Let $\sfV \in \bfFlagPerfps^{(m-1)}(\scrC)$.
	We define the derived stack $\bfFlagPerfps^{(m),\dagger}(\scrC;\sfV)$ of $\sfV$-flags of length $m$ as the fiber product
	\begin{align}
		\begin{tikzcd}[ampersand replacement=\&]
			\bfFlagPerfps^{(m),\dagger}(\scrC;\sfV) \arrow{r} \arrow{d} \&  \bfFlagPerfps^{(m)}(\scrC) \arrow{d}{\partial_m}\\
			\Spec(k) \arrow{r}{x_\sfV} \& \bfFlagPerfps^{(m-1)}(\scrC)
		\end{tikzcd}\ ,
	\end{align}
	where $x_\sfV \colon \Spec(k) \to \bfFlagPerfps^{(m-1)}(\scrC)$ is the map corresponding to the flag $\sfV$.
\end{definition}

Applying the \textit{right} version of Construction~\ref{construction:Vflags_generalized} to $\calS_\bullet\bfPerfps(\scrC)$ and the $(m-1)$-flag $\sfV$, we obtain a relative $2$-Segal space that we denote as\footnote{The discrepancy between the choice of either left or right version in Construction~\ref{construction:Vflags_generalized} and the superscript $\ast$ in $\calS^\ast_\bullet(-)$ is justified from the fact that, in this way, the notation is compatible with the classical literature on (classical, cohomological, K-theoretical) Hall algebras and their representations. See \cite{Porta_Sala_Hall}.}
\begin{align} \label{eq:universal_m_flag_representation}
	u^\ell_\bullet \colon \calS^\ell_\bullet \bfFlagPerfps^{(m),\dagger}(\scrC;\sfV) \longrightarrow \calS_\bullet \bfPerfps(\scrC) \ .
\end{align}
When $m = 1$, the choice of $\sfV$ is empty (cf.\ Remark~\ref{rem:empty-case}).
In this case, we therefore simply denote the above simplicial object by
\begin{align}
	u^\ell_\bullet \colon \calS^\ell_\bullet \bfFlagPerfps^{(1)}(\scrC) \longrightarrow \calS_\bullet \bfPerfps(\scrC) \ . 
\end{align}

\begin{remark}\label{rem:flag_action_coordinates_II}	
	In light of Theorem~\ref{thm:Godicke_II}, we think of $\bfPerfps(\scrC)$ as an algebra in the $\infty$-category of correspondences $\Corr(\dSt)$ acting on the derived stack $\bfFlagPerfps^{(m)}(\scrC;\sfV)$. The action is implemented by the correspondence
	\begin{align}
		\begin{tikzcd}[ampersand replacement=\&]
			\calS_1^\ell \bfFlagPerfps^{(m),\dagger}(\scrC;\sfV) \arrow{r}{\varpi_0} \arrow{d}{u_1^\ell \times \varpi_1} \& 	\bfFlagPerfps^{(m),\dagger}(\scrC;\sfV) \\
			\bfPerfps(\scrC) \times \bfFlagPerfps^{(m),\dagger}(\scrC;\sfV) 
		\end{tikzcd} \ .
	\end{align}
	We shall describe explicitly the maps $\varpi_0, \varpi_1$, and $u_1^\ell$. Unraveling the definition, we see that $\calS_1^\ell \bfFlagPerfps^{(m),\dagger}(\scrC; \sfV)$ fits in the following fiber product:
	\begin{align}\label{eq:identifications-2}
		\begin{tikzcd}[ampersand replacement = \&]
			\calS_1^\ell \bfFlagPerfps^{(m),\dagger}(\scrC; \sfV) \arrow{r} \arrow{d} \& \calS_{m+1} \bfPerfps(\scrC) \arrow{d}{\partial_{m,m+1}} \\
			\Spec(k) \arrow{r}{x_{\sfV}} \& \calS_{m-1} \bfPerfps(\scrC) \ ,
		\end{tikzcd}
	\end{align}
	where $\partial_{m,m+1}$ is induced by the map $[m+1] \to [m-1]$ in $\mathbf \Delta\op$ that avoids $m$ and $m+1$ inside $[m+1]$, while the morphism $x_{\sfV}$ is the morphism classifying the $(m-1)$-flag $\sfV$. In other words, $\calS_1^\ell \bfFlagPerfps^{(m),\dagger}(\scrC;\sfV)$ can be informally described as the derived stack parametrizing diagrams of the form
	\begin{align}\label{eq:flag_action}
		\begin{tikzcd}[ampersand replacement=\&]
			0 \arrow{r} \& V_{0,1} \arrow{r} \arrow{d} \& V_{0,2} \arrow{r} \arrow{d} \& \cdots \arrow{r} \& V_{0,m-1} \arrow{r} \arrow{d} \& F_{0,m} \arrow{d} \arrow{r} \& F_{0,m+1} \arrow{d} \\
			{} \& 0 \arrow{r} \& V_{1,2} \arrow{r} \& \cdots \arrow{r} \& V_{1,m-1} \arrow{r} \arrow{d} \& F_{1,m} \arrow{d} \arrow{r} \& F_{1,m+1} \arrow{d} \\
			\& \& \& \ddots \& \vdots \arrow{d} \& \vdots \arrow{d} \& \vdots \arrow{d} \\
			\& \& \& \& V_{m-2,m-1} \arrow{d} \arrow{r} \& F_{m-2,m} \arrow{d} \arrow{r} \& F_{m-2,m+1} \arrow{d} \\
			\& \& \& \& 0 \arrow{r} \& F_{m-1,m} \arrow{d} \arrow{r} \& F_{m-1,m+1} \arrow{d} \\
			\& \& \& \& \& 0 \arrow{r} \& F_{m,m+1} \arrow{d} \\
			\& \& \& \& \& \& 0
		\end{tikzcd}
	\end{align}
	The morphism $u_1^\ell$ sends the diagram \eqref{eq:flag_action} to $F_{m,m+1}$, the morphism $\varpi_0$ sends the diagram \eqref{eq:flag_action} to the full $\sfV$-flag determined by the chain
	\begin{align}
		F_{0,m+1} \longrightarrow F_{1,m+1} \longrightarrow \cdots \longrightarrow F_{m-1,m+1} \ , 
	\end{align}
	while the morphism $\varpi_1$ sends the diagram \eqref{eq:flag_action} to the full $\sfV$-flag determined by the chain
	\begin{align}
		F_{0,m} \longrightarrow F_{1,m} \longrightarrow \cdots \longrightarrow F_{m-1,m} \ . 
	\end{align}
	\end{remark}

\subsection{COHAs, CatHAs, and their representations}\label{subsec:COHA-representations}

Let $\bfH$ and $\bfM$ be derived stacks and let
\begin{align}
	\jmath \colon \bfH \longrightarrow \bfPerfps(\scrC) \qquad \text{and} \qquad \imath \colon \bfM \longrightarrow \bfFlagPerfps^{(m),\dagger}(\scrC;\sfV) 
\end{align}
be two morphisms. Choose as well a morphism $0_{\bfH} \colon \Spec(k) \to \bfH$ lifting the zero morphism $\Spec(k) \to \bfPerfps(\scrC)$ classifying the zero  pseudo-perfect object.
We can arrange these data into the following diagram $\rho$:
\begin{align}
	\rho\coloneqq \begin{tikzcd}[ampersand replacement=\&]
		\bfM \arrow{r} \arrow{d}{\imath} \& \Spec(k) \arrow{r}{0_\bfH} \arrow[equal]{d} \& \bfH \arrow{d}{\jmath} \\
		\bfFlagPerfps^{(m),\dagger}(\scrC;\sfV) \arrow{r} \& \Spec(k) \arrow{r}{0} \& \bfPerfps(\scrC) 
	\end{tikzcd} \ .
\end{align}
We obtain in this way a Hecke datum for the relative simplicial derived stack
\begin{align}
	u^\ell_\bullet \colon \calS^\ell_\bullet \bfFlagPerfps^{(m),\dagger}(\scrC;\sfV) \longrightarrow \calS_\bullet \bfPerfps(\scrC) \ ,
\end{align}
in the sense of Definition~\ref{def:Hecke_datum}.
Thus, applying Construction~\ref{construction:Hecke_pattern}, we obtain a new relative simplicial derived stack
\begin{align}
	u^\ell_\bullet \colon \calS^\ell_\bullet \bfFlagPerf_{\bfH,\bfM}^{(m),\dagger}(\scrC;\sfV) \longrightarrow \calS_\bullet \bfPerf_\bfH(\scrC) \ . 
\end{align}

\begin{rem}
	Let $(\Lambda,+)$ be a monoid and assume that $\calS^\ell_\bullet \bfFlagPerfps^{(m),\dagger}(\scrC;\sfV) \to \calS_\bullet \bfPerf(\scrC)$ admits a $\Lambda^\star$-graded structure. Then Corollary~\ref{cor:induced_Lambda_grading} provides induced $\Lambda^\star$-gradings on both simplicial derived stacks $\calS^\ell_\bullet \bfFlagPerf_{\bfH,\bfM}^{(m),\dagger}(\scrC;\sfV)$ and $\calS_\bullet \bfPerf_\bfH(\scrC)$. \hfill $\triangle$
\end{rem}

\begin{remark}\label{rem:relative-2-Segal}
	Unwinding the definitions, we see that $\calS_1 \bfPerf_\bfH(\scrC) \simeq \bfH$ and that $\calS_2 \bfPerf_\bfH(\scrC)$ fits in the following pullback square
	\begin{align}\label{eq:extensions-T}
		\begin{tikzcd}[ampersand replacement=\&]
			\calS_2 \bfPerf_\bfH(\scrC) \arrow{r} \arrow{d} \& \calS_2 \bfPerfps(\scrC) \arrow{d}{\partial_0 \times \partial_1 \times \partial_2} \\
			\bfH \times \bfH \times \bfH \arrow{r} \& \bfPerfps(\scrC) \times \bfPerfps(\scrC) \times \bfPerfps(\scrC) 
		\end{tikzcd} \ ,
	\end{align}
	where the maps $\partial_0, \partial_1$, and $\partial_2$ has been introduced in Formula~\eqref{eq:partial}.
	
	Similarly, $\calS^\ell_0 \bfFlagPerf_{\bfH,\bfM}^{(m),\dagger}(\scrC;\sfV) \simeq \bfM$ and $\calS^\ell_1 \bfFlagPerf_{\bfH,\bfM}^{(m),\dagger}(\scrC;\sfV)$ fits in the following pullback square:
	\begin{align}
		\begin{tikzcd}[ampersand replacement=\&]
			\calS^\ell_1 \bfFlagPerf_{\bfH, \bfM}^{(m),\dagger}(\scrC;\sfV) \arrow{r} \arrow{d} \& \calS_{1}^\ell \bfFlagPerfps^{(m),\dagger}(\scrC;\sfV) \arrow{d}{u_1^\ell\times \varpi_0 \times \varpi_1} \\
			\bfH\times \bfM \times \bfM \arrow{r} \& \bfPerfps(\scrC) \times \bfFlagPerfps^{(m),\dagger}(\scrC;\sfV) \times \bfFlagPerfps^{(m),\dagger}(\scrC;\sfV)
		\end{tikzcd} \ ,
	\end{align}
	where the maps $\varpi_0, \varpi_1$, and $u_1^\ell$ has been introduced in Remark~\ref{rem:flag_action_coordinates_II}.
\end{remark}

At this point, Corollary~\ref{cor:Hecke_pattern} immediately implies the following.
\begin{proposition}\label{prop:2-Segal-algebra} 
	Assume that the square
	\begin{align}\label{eq:pullback-square-T}
		\begin{tikzcd}[ampersand replacement=\&]
			\calS_2 \bfPerf_\bfH(\scrC) \arrow{r} \arrow{d} \& \calS_2 \bfPerfps(\scrC) \arrow{d}{\partial_0 \times \partial_2} \\
			\bfH \times \bfH \arrow{r} \& \bfPerfps(\scrC) \times \bfPerfps(\scrC)
		\end{tikzcd}
	\end{align}
	is a pullback. 
	Then $\calS_\bullet \bfPerf_\bfH(\scrC)$ is a $2$-Segal stack. 
\end{proposition}

To realize categorical Hall algebras, we apply the framework developed in \cite[\S4.2]{Porta_Sala_Hall}, based on the theory of correspondences defined in \cite{Gaitsgory_Rozenblyum_Study_I, Gaitsgory_Rozenblyum_Study_II}. In particular, our categorification will be at the level of the pro-category $\catCohb_{\mathsf{pro}}(-)$, which is a pro-enhancement of the usual stable $\infty$-category $\catCohb(-)$ of locally cohomologically bounded complexes. On the other hand, to realize cohomological Hall algebras we apply the motivic framework discussed in \S\ref{subsec:beyond_qc}, \S\ref{subsec:operations_BM_homology}, and \S\ref{subsec:Lambda-graded}. 

In the following statement, the terminology from Definition~\ref{def:modified_classes_of_morphisms} is in use. We have the following.
\begin{corollary} \label{cor:COHA}
	Assume that $\bfH$ is a quasi-separated geometric derived stack locally almost of finite presentation over $k$ and the square \eqref{eq:pullback-square-T} is a pullback. In addition, assume that
	\begin{enumerate}[label=(\roman*)]\itemsep0.2cm
		\item \label{item:algebra-(i)} the map
		\begin{align}
			\partial_0 \times \partial_2 \colon \calS_2\bfPerf_\bfH(\scrC) \longrightarrow \bfH \times \bfH 
		\end{align}
		is quasi-compact, finitely connected and derived lci, and
		\item \label{item:algebra-(ii)} the map
		\begin{align}
			\partial_1 \colon \calS_2 \bfPerf_\bfH(\scrC) \longrightarrow \bfH 
		\end{align}
		is locally rpas.	\end{enumerate}
	Then, $\catCohb_{\mathsf{pro}}( \bfH )$ has the structure of an $\E_1$-monoidal stable pro-$\infty$-category, whose underlying tensor product is given by the composition
	\begin{align}
		\catCohb_{\mathsf{pro}}( \bfH ) \times \catCohb_{\mathsf{pro}}( \bfH ) \xrightarrow{\boxtimes} \catCohb_{\mathsf{pro}}( \bfH \times \bfH ) \xrightarrow{(\partial_1)_\ast \circ (\partial_0 \times \partial_2)^\ast} \catCohb_{\mathsf{pro}}( \bfH )\ .
	\end{align}
	
	Similarly, let $\bfD^\ast$ be a motivic formalism, and fix $\calA \in \CAlg(\bfD^\ast(\Spec(k)))$ and $\Gamma \subseteq \Pic(\bfD^\ast(\Spec(k)))$ such that Assumption~\ref{assumption:motivic_formalism} is satisfied. Then, the topological vector space $\HBMDGamma_0(\bfH;\calA)$ becomes a unital associative algebra. In particular,
	\begin{align}
		G_0( \bfH )\quad \text{and} \quad \HBM_\ast( \bfH )
	\end{align}
	become unital associative algebras.
\end{corollary}		

\begin{proof}
	First, recall from Proposition~\ref{prop:2-Segal-algebra} that $\calS_\bullet \bfPerf_\bfH(\scrC)$ is a $2$-Segal object. Therefore, Theorem~\ref{thm:Godicke} allows to review it as an algebra object in $\Corr^\times(\dSt_k)$. Our assumptions \ref{item:algebra-(i)} and \ref{item:algebra-(ii)} allows to review it as an algebra object in the smaller $\Corr^\times(\dGeomqs_S)_{\mathsf{qc.lci}\:\cap\:\mathsf{ufconn},\:\mathsf{lrpas}}$. We can therefore apply the lax monoidal functor $\HBMDGamma_0(-;\calA)$ provided by Theorem~\ref{thm:functoriality_of_genuine_BM_homology_admissible} to linearize this algebra structure. In the categorical case, we can use the functor $\catCohb_{\mathsf{pro}}(-)$ studied in \cite[\S4.2]{Porta_Sala_Hall} in place of $\HBMDGamma_0(-;\calA)$ to obtain the same conclusion.
\end{proof}

\begin{remark}\label{rem:Lambda_graded_multiplication}
	If, in addition $\calS_\bullet \bfPerf_\bfH(\scrC)$ admits a $\Lambda^{\mathfrak{a}}$-graded structure for a monoid $(\Lambda,+)$, then $\catCohb_{\mathsf{pro}}(\bfH)$, $G_0(\bfH)$ and $\HBM_\ast(\bfH)$ inherit canonical $\Lambda$-grading, and the Hall multiplication is canonically compatible with the $\Lambda$-graded structure. This simply follows running the above proof using Theorem~\ref{thm:BM_Lambda_graded_functoriality} instead of Theorem~\ref{thm:functoriality_of_genuine_BM_homology_admissible}.
\end{remark}

Now we introduce generalizations to the setting of Segal spaces of the notions of left and right Hecke patterns introduced in \cite[\S5.1]{KV_Hall} (see also \cite[Definition~6.1 and Remark~6.2]{MMSV}).
\begin{definition}
	We say that $\bfM$ is a \textit{right 1-Segal Hecke pattern for $\bfH$} if 
	\begin{align}\label{eq:diagram-1-u}
		\begin{tikzcd}[ampersand replacement=\&]
			\calS^\ell_1 \bfFlagPerf_{\bfH, \bfM}^{(m),\dagger}(\scrC;\sfV) \arrow{r} \arrow{d} \& \calS_1^\ell \bfFlagPerfps^{(m),\dagger}(\scrC;\sfV) \arrow{d}{u_1^\ell\times \varpi_1} \\
			\bfH \times \bfM  \arrow{r} \&  \bfPerfps(\scrC) \times \bfFlagPerfps^{(m),\dagger}(\scrC;\sfV)
		\end{tikzcd} \ ,
	\end{align}
	is a pullback, while $\bfM$ is a \textit{left 1-Segal Hecke pattern for $\bfH$} if 
	\begin{align}\label{eq:diagram-0-u}
		\begin{tikzcd}[ampersand replacement=\&]
			\calS^\ell_1 \bfFlagPerf_{\bfH, \bfM}^{(m),\dagger}(\scrC;\sfV) \arrow{r} \arrow{d} \& \calS_1^\ell \bfFlagPerfps^{(m),\dagger}(\scrC;\sfV) \arrow{d}{u_1^\ell\times \varpi_0} \\
			\bfH \times \bfM \arrow{r} \& \bfPerfps(\scrC) \times \bfFlagPerfps^{(m),\dagger}(\scrC;\sfV) 
		\end{tikzcd} \ ,
	\end{align}
	is a pullback.
	\end{definition}

Applying again Corollary~\ref{cor:Hecke_pattern}, we obtain the following.
\begin{proposition}\label{prop:2-Segal-representation}	
	Assume that 
	\begin{enumerate}\itemsep=0.2cm
		\item the square \eqref{eq:pullback-square-T} is a pullback, and 
		\item \label{item:2-Segal-representation-(2)} $\bfM$ is either a right or a left 1-Segal Hecke pattern for $\bfH$.
	\end{enumerate}	
	Then $\calS^\ell_\bullet \bfFlagPerf_{\bfH,\bfM}^{(m),\dagger}(\scrC;\sfV) \to \calS_\bullet \bfPerf_\bfH(\scrC)$ is a relative $2$-Segal stack.
\end{proposition}

\begin{warning}\label{warning:left_and_right_relative_2_Segal}
	A relative $2$-Segal space induces simultaneously a left and a right module structure in correspondences, accordingly to whether one uses $(u_1^\ell \times \varpi_0,\varpi_1)$ or $(u_1^\ell \times \varpi_1, \varpi_0)$. In other words, the structure of relative $2$-Segal space encodes the associativity of \textit{both} actions. Indeed, if for example we check that the square \eqref{eq:diagram-1-u} is a pullback, then the map
	\begin{align}
		(\calS_\bullet^\ell \bfFlagPerf_{\bfH,\bfM}^{(m),\dagger}(\scrC;\sfV) \to \calS_\bullet \bfPerf_{\bfH}(\scrC)) \longrightarrow (\calS_\bullet^\ell \bfFlagPerfps^{(m),\dagger}(\scrC;\sfV) \to \calS_\bullet \bfPerfps(\scrC)) 
	\end{align}
	is relative \textit{right} $1$-Segal (hence relative $2$-Segal), and the latter implies that the same map is relative \textit{left} $1$-Segal (hence relative $2$-Segal).
\end{warning}

Similarly to Corollary~\ref{cor:COHA}, we get:
\begin{corollary} \label{cor:COHA-representations}
	Assume that 
	\begin{enumerate}[label=(\roman*)]\itemsep0.2cm
		\item \label{item:COHA-representations-1} both $\bfH$ and $\bfM$ are quasi-separated geometric derived stacks locally almost of finite presentation over $k$;
		
		\item \label{item:COHA-representations-2} the assumptions in Corollary~\ref{cor:COHA} are satisfied by $\bfH$, and the condition \eqref{item:2-Segal-representation-(2)} of Proposition~\ref{prop:2-Segal-representation} is satisfied by $\bfH, \bfM$;
		
		\item \label{item:COHA-representations-3} the map
		\begin{align}
			u_1^\ell \times \varpi_1 \colon \calS^\ell_1 \bfFlagPerf_{\bfH,\bfM}^{(m),\dagger}(\scrC;\sfV) \longrightarrow \bfH\times \bfM
		\end{align}
		is quasi-compact, finitely connected and derived lci, and 
		
		\item \label{item:COHA-representations-4} the map
		\begin{align}
			\varpi_0 \colon \calS^\ell_1 \bfFlagPerf_{\bfH,\bfM}^{(m),\dagger}(\scrC;\sfV) \longrightarrow \bfM
		\end{align}
		is locally rpas.	\end{enumerate}
	Then, $\catCohb_{\mathsf{pro}}(\bfM)$ has the structure of a left categorical module over the $\E_1$-monoidal $\infty$-category $\catCohb_{\mathsf{pro}}(\bfH)$, whose underlying action is given by the composition
	\begin{align}
		\catCohb_{\mathsf{pro}}( \bfH ) \times \catCohb_{\mathsf{pro}}( \bfM ) \xrightarrow{\boxtimes} \catCohb_{\mathsf{pro}}( \bfH \times \bfM ) \xrightarrow{(\varpi_0)_\ast \circ (u_1^\ell\times \varpi_1 )^\ast} \catCohb_{\mathsf{pro}}( \bfM )\ .
	\end{align}
	
	Similarly, let $\bfD^\ast$ be a motivic formalism, and fix $\calA \in \CAlg(\bfD^\ast(\Spec(k)))$ and $\Gamma \subseteq \Pic(\bfD^\ast(\Spec(k)))$ such that Assumption~\ref{assumption:motivic_formalism} is satisfied. Then, the topological vector space $\HBMDGamma_0(\bfM;\calA)$ becomes a left module over  $\HBMDGamma_0(\bfH;\calA)$. In particular,
	\begin{align}
		G_0( \bfM )\quad \text{and} \quad \HBM_\ast( \bfM )			
	\end{align}
	are left modules over $G_0( \bfH )$ and $\HBM_\ast( \bfH )$, respectively.
\end{corollary}

\begin{proof}
	To begin with,
	\begin{align}
		 \calS^\ell_\bullet \bfFlagPerf_{\bfH,\bfM}^{(m),\dagger}(\scrC;\sfV) \longrightarrow \calS_\bullet \bfPerf_\bfH(\scrC) 
	\end{align}
	is a relative $2$-Segal stack by Proposition~\ref{prop:2-Segal-representation}. Therefore, Theorem~\ref{thm:Godicke_II} allows to review
	\begin{align}
		\calS_0\bfFlagPerf^{(m),\dagger}_{\bfH,\bfM}(\scrC;\sfV) \simeq \bfM 
	\end{align}
	as a left module over $\bfH$ in the $\infty$-category of correspondences $\Corr^\times(\dSt_k)$. Our assumptions \ref{item:COHA-representations-1}--\ref{item:COHA-representations-4} imply that $\bfM$ is a left module over $\bfH$ in the smaller $\Corr^\times(\dGeomqs_S)_{\mathsf{qc.lci}\:\cap\:\mathsf{ufconn},\:\mathsf{lrpas}}$. We can therefore apply the lax monoidal functor $\HBMDGamma_0(-;\calA)$ provided by Theorem~\ref{thm:functoriality_of_genuine_BM_homology_admissible} to linearize this action. In the categorical case, we can use the functor $\catCohb_{\mathsf{pro}}(-)$ studied in \cite[\S4.2]{Porta_Sala_Hall} in place of $\HBMDGamma_0(-;\calA)$ to obtain the same conclusion.
\end{proof}

\begin{remark}
	Running Construction~\ref{construction:self_action} instead of Construction~\ref{construction:Vflags_generalized}, we would get a relative $2$-Segal space
	\begin{align}
		u^\ell_\bullet \colon \calS_\bullet^\ell \bfFlagPerfps^{(m)}(\scrC) \longrightarrow \calS_\bullet \bfPerfps(\scrC) \ , 
	\end{align}
	and all the results of this section equally apply to this setup. For $m = 1$ there is no difference between the $\sfV$-marked version and the unmarked one.
\end{remark}

Instead of the right version, we can apply the left version of Construction~\ref{construction:Vflags_generalized} to $\calS_\bullet\bfPerfps(\scrC)$ and the $(m-1)$-flag $\sfV$, we obtain a relative $2$-Segal space that we denote as
\begin{align}
	u^r_\bullet \colon \calS^r_\bullet \bfFlagPerfps^{(m),\dagger}(\scrC;\sfV) \longrightarrow \calS_\bullet \bfPerfps(\scrC) \ .
\end{align}
Then, all the statements above hold also after replacing ``left'' with ``right'' and  vice versa (see Warning~\ref{warning:left_and_right_relative_2_Segal}).
In particular, the same proof of Corollary~\ref{cor:COHA-representations} implies:
\begin{corollary}\label{cor:COHA-representations-right}
	Assume that 
	\begin{enumerate}[label=(\roman*)]\itemsep0.2cm
		\item both $\bfH$ and $\bfM$ are quasi-separated geometric derived stacks locally almost of finite presentation over $k$;
		
		\item the assumptions in Corollary~\ref{cor:COHA} are satisfied by $\bfH$, and the condition \eqref{item:2-Segal-representation-(2)} of for the ``right'' version of Proposition~\ref{prop:2-Segal-representation} is satisfied by $\bfH, \bfM$;
		
		\item the map
		\begin{align}
			\varpi_0 \times u_1^r \colon \calS^r_1 \bfFlagPerf_{\bfM,\bfH}^{(m),\dagger}(\scrC;\sfV) \longrightarrow \bfM\times \bfH 
		\end{align}
		is quasi-compact, finitely connected and derived lci, and
		
		\item  the map
		\begin{align}
			\varpi_1 \colon \calS^r_1 \bfFlagPerf_{\bfM,\bfH}^{(m),\dagger}(\scrC;\sfV) \longrightarrow \bfM
		\end{align}
		is locally rpas.
	\end{enumerate}
	Then, $\catCohb_{\mathsf{pro}}(\bfM)$ has the structure of a right categorical module over the $\E_1$-monoidal stable pro-$\infty$-category $\catCohb_{\mathsf{pro}}(\bfH)$, whose underlying action is given by the composition
	\begin{align}
		\catCohb_{\mathsf{pro}}( \bfM ) \times \catCohb_{\mathsf{pro}}( \bfH ) \xrightarrow{\boxtimes} \catCohb_{\mathsf{pro}}( \bfH \times \bfM ) \xrightarrow{(\varpi_1)_\ast \circ (\varpi_0\times u_1^r)^\ast} \catCohb_{\mathsf{pro}}( \bfM )\ .
	\end{align}
	
	Similarly, let $\bfD^\ast$ be a motivic formalism, and fix $\calA \in \CAlg(\bfD^\ast(\Spec(k)))$ and $\Gamma \subseteq \Pic(\bfD^\ast(\Spec(k)))$ such that Assumption~\ref{assumption:motivic_formalism} is satisfied. Then, the topological vector space $\HBMDGamma_0(\bfM;\calA)$ is a right module over $\HBMDGamma_0(\bfH;\calA)$. In particular,
	\begin{align}
		G_0( \bfM )\quad \text{and} \quad \HBM_\ast( \bfM )			
	\end{align}
	are right modules over $G_0( \bfH )$ and $\HBM_\ast( \bfH )$, respectively.
\end{corollary}

\begin{remark}\label{rem:Lambda_graded_action}
	If, in addition $\calS_\bullet^\ell \bfPerf_{\bfH,\bfM}^{(m),\dagger}(\scrC;\sfV) \to \calS_\bullet \bfPerf_{\bfH}(\scrC)$ and its right variant admits a $\Lambda^{\star}$-graded structure for a monoid $(\Lambda,+)$, then the actions constructed in Corollaries~\ref{cor:COHA-representations} and \ref{cor:COHA-representations-right} canonically lift to $\Lambda$-graded representations. This stems, exactly as in Remark~\ref{rem:Lambda_graded_multiplication} from the possibility of using Theorem~\ref{thm:BM_Lambda_graded_functoriality} instead of Theorem~\ref{thm:functoriality_of_genuine_BM_homology_admissible} in the proof of the above two corollaries.
\end{remark}

\subsubsection{Induced properties}

Choose maps $\alpha \colon \bfH' \to \bfH$ and $\mu \colon \bfM' \to \bfM$ of derived stacks. Assume that the morphism $0_\bfH \colon \Spec(k) \to \bfH$ lifts to a morphism $0_{\bfH'} \colon \Spec(k) \to \bfH'$. We get a relative simplicial derived stack
\begin{align}
	u^\ell_\bullet \colon \calS^\ell_\bullet \bfFlagPerf_{\bfM', \bfH'}^{(m),\dagger}(\scrC;\sfV) \longrightarrow \calS_\bullet \bfPerf_{\bfH'}(\scrC) \ . 
\end{align}

We provide two lemmas that will be useful later on to check if the Assumption~\ref{item:algebra-(i)} of Corollary~\ref{cor:COHA} and the Assumption~\ref{item:COHA-representations-4} of Corollary~\ref{cor:COHA-representations} are satisfied.
\begin{lemma}\label{lem:induced-derived-lci}
	Assume that $\bfH$ and $\bfH'$ are geometric derived stacks locally of finite presentation over $k$ and that the following conditions are met:
	\begin{enumerate}\itemsep=0.2cm	
		\item \label{item:induced-derived-lci-1} the squares \eqref{eq:pullback-square-T} relative to both the stacks $\bfH$ and $\bfH'$ are pullback;
		\item \label{item:induced-derived-lci-2} the stack $\bfH$ satisfies Assumption~\ref{item:algebra-(i)} of Corollary~\ref{cor:COHA}.
	\end{enumerate}
	Then, the stack $\bfH'$ also satisfies Assumption~\ref{item:algebra-(i)} of Corollary~\ref{cor:COHA}.
\end{lemma}
\begin{proof}
	Consider the following ladder of commutative squares:
	\begin{align}
		\begin{tikzcd}[ampersand replacement = \&]
			\calS_2 \bfPerf_{\bfH'}(\scrC) \arrow{r} \arrow{d}{\partial_0 \times \partial_2} \& \calS_2 \bfPerf_\bfH(\scrC) \arrow{r} \arrow{d}{\partial_0 \times \partial_2} \& \calS_2 \bfPerfps(\scrC) \arrow{d}{\partial_0 \times \partial_2} \\
			\bfH' \times \bfH' \arrow{r} \& \bfH \times \bfH \arrow{r} \& \bfPerfps(\scrC) \times \bfPerfps(\scrC)
		\end{tikzcd} \ .
	\end{align}
	By assumption, the right square is a pullback, hence the middle vertical map is quasi-compact, universally finitely connected and derived lci. Also, by assumption, the outer square is a pullback. Thus, the left square is a pullback as well, and hence the leftmost vertical map inherits the same properties.
\end{proof}

\begin{corollary}\label{cor:induced_COHA-properness}
	Assume that $\bfH$, $\bfH'$, $\bfM$, and $\bfM'$ are quasi-separated geometric derived stacks locally almost of finite presentation over $k$ and that the following conditions are met:
	\begin{enumerate}\itemsep=0.2cm
		\item the pair $(\bfH, \bfH')$ satisfies the conditions of Lemma~\ref{lem:induced-derived-lci};
		
		\item the map
		\begin{align}
			\varpi_0 \colon \calS_1^\ell \bfFlagPerf^{(m),\dagger}_{\bfH,\bfM}(\scrC;\sfV) \longrightarrow \bfM
		\end{align}
		is locally rpas;
		
		\item \label{item:induced_COHA-representation-left-2} the map $\alpha \colon \bfH' \to \bfH$ is locally rpas; 
		
		\item \label{item:induced_COHA-representation-left-3} the square
		\begin{align}
			\begin{tikzcd}[ampersand replacement=\&]
				\calS_1^\ell \bfFlagPerf^{(m),\dagger}_{\bfH',\bfM'}(\scrC;\sfV) \arrow{r} \arrow{d}{u_1^\ell\times \varpi_0} \& \calS_1^\ell \bfFlagPerf^{(m),\dagger}_{\bfH,\bfM}(\scrC;\sfV) \arrow{d}{u_1^\ell\times \varpi_0} \\
				\bfH' \times \bfM' \arrow{r} \& \bfH \times \bfM
			\end{tikzcd}
		\end{align}		
		is a pullback (this happens for instance if both pairs $(\bfH,\bfM)$ and $(\bfH',\bfM')$ make the square \eqref{eq:diagram-0-u} into a pullback).
	\end{enumerate}
	Then the map
	\begin{align}
		\varpi_0 \colon \calS_1^\ell \bfFlagPerf^{(m),\dagger}_{\bfH',\bfM'}(\scrC;\sfV) \longrightarrow \bfM' 
	\end{align}
	is locally rpas.
\end{corollary}

\begin{proof}
	Consider the following commutative diagram:
	\begin{align}\label{eq:the_big_diagram}
		\begin{tikzcd}[ampersand replacement = \&]
			\calS_1^\ell \bfFlagPerf^{(m),\dagger}_{\bfH',\bfM'}(\scrC;\sfV) \arrow{r}{q_0} \arrow{d}{p_0} \& \scrX_{\bfM,\bfM',\bfH'} \arrow{r}{q_1} \arrow{d} \& \scrX_{\bfM,\bfM',\bfH} \arrow{r}{q_2} \arrow{d} \& \bfM' \arrow{d}{\mu}  \\
			\scrX_{\bfM', \bfM, \bfH'} \arrow{r} \arrow{d}{p_1} \& \scrX_{\bfM,\bfM,\bfH'} \arrow{r} \arrow{d} \& \calS_1^\ell \bfFlagPerf^{(m),\dagger}_{\bfH,\bfM}(\scrC;\sfV) \arrow{d}{u_1^\ell\times \varpi_1} \arrow{r}{\varpi_0} \& \bfM \\
			\bfH'\times \bfM'  \arrow{r}{\id_{\bfH'}\times \mu} \& \bfH' \times \bfM  \arrow{r}{\alpha\times \id_F} \& \bfH \times \bfM \&  
		\end{tikzcd}
	\end{align}
	where the objects $\scrX_{\ast,\ast,\ast}$ are defined by asking all the squares to be pullbacks. Then the map $\varpi_0\colon \calS_1^\ell \bfFlagPerf^{(m),\dagger}_{\bfH',\bfM'}(\scrC;\sfV) \longrightarrow \bfM'$ is canonically identified with the composite $q_2 \circ q_1 \circ q_0$. Although the upper left square may not be a pullback, condition~\eqref{item:induced_COHA-representation-left-3} implies that $q_0$ is an equivalence. Since $\varpi_0\colon \calS_1^\ell \bfFlagPerf^{(m),\dagger}_{\bfH,\bfM}(\scrC;\sfV) \longrightarrow \bfM$ is locally rpas by assumption, we see that $q_2$ is locally rpas. Since $\alpha$ is locally rpas by assumption, the same goes for $q_1$. The conclusion follows from the stability of $\mathsf{lrpas}$ under composition, see Lemma~\ref{lem:properties_locallly_rpas}.
\end{proof}

Similarly, by using again diagram~\eqref{eq:the_big_diagram}, one can prove the following.
\begin{corollary}\label{cor:induced_COHA-derived-lci}
	Assume that $\bfH$, $\bfH'$, $\bfM$, and $\bfM'$ are geometric derived stacks locally of finite presentation over $k$ and that the following conditions are met:
	\begin{enumerate}\itemsep=0.2cm
		\item the pair $(\bfH, \bfH')$ satisfies the conditions of Lemma~\ref{lem:induced-derived-lci};
		
		\item the map
		\begin{align}
			u_1^\ell\times \varpi_1 \colon \calS_1^\ell \bfFlagPerf^{(m),\dagger}_{\bfH,\bfM}(\scrC;\sfV) \longrightarrow \bfH\times \bfM
		\end{align}
		is quasi-compact, finitely connected and derived lci;
		
		\item the map $\mu \colon \bfM' \to \bfM$ is quasi-compact, finitely connected and derived lci;
		
		\item the square
		\begin{align}
			\begin{tikzcd}[ampersand replacement=\&]
				\calS_1^\ell \bfFlagPerf^{(m),\dagger}_{\bfH',\bfM'}(\scrC;\sfV) \arrow{r} \arrow{d}{u_1^\ell\times \varpi_0} \& \calS_1^\ell \bfFlagPerf^{(m),\dagger}_{\bfH,\bfM}(\scrC;\sfV) \arrow{d}{u_1^\ell\times \varpi_0} \\
				\bfH' \times \bfM' \arrow{r} \& \bfH \times \bfM
			\end{tikzcd}
		\end{align}		
		is a pullback (this happens for instance if both pairs $(\bfH,\bfM)$ and $(\bfH',\bfM')$ make the square \eqref{eq:diagram-0-u} into a pullback).
	\end{enumerate}
	Then the map
	\begin{align}
		u_1^\ell\times \varpi_1 \colon \calS_1^\ell \bfFlagPerf^{(m),\dagger}_{\bfH',\bfM'}(\scrC;\sfV) \longrightarrow \bfH'\times \bfM'
	\end{align}
	is quasi-compact, finitely connected and derived lci.
\end{corollary}

\begin{remark}
	Note that one can formulate ``right'' versions of Corollaries~\ref{cor:induced_COHA-properness} and \ref{cor:induced_COHA-derived-lci}. We leave to the interested reader to explicit formulate these statements.
\end{remark}

\subsection{The derived stack of $m$-flags as a linear stack} \label{subsec:linear-stack}

Let $\calF^{\mathsf{u}} \in \catQCoh(\bfPerfps(\scrC)) \otimes \scrC$ be the universal object on $\bfPerfps(\scrC)$. For $i = 1,2$ let
\begin{align}
	\mathsf{pr}_i \colon \bfPerfps(\scrC) \times \bfPerfps(\scrC) \longrightarrow \bfPerfps(\scrC) 
\end{align}
be the canonical projections and set
\begin{align}
	\calF_i^{\mathsf u} \coloneqq \mathsf{pr}_i^\ast(\calF^{\mathsf{u}}) \in \catQCoh( \bfPerfps(\scrC) \times \bfPerfps(\scrC) ) \otimes \scrC \ . 
\end{align}
Define 
\begin{align}
	\calE \coloneqq \Hom_{\catQCoh(\bfPerfps(\scrC) \times \bfPerfps(\scrC)) \otimes \scrC}( \calF_2^{\mathsf{u}}, \calF_1^{\mathsf{u}} ) \in \catQCoh(\bfPerfps(\scrC) \times \bfPerfps(\scrC)) \ .
\end{align}

\begin{remark}
	Let us briefly describe the relation with the commutative case, discussed in \cite[Construction~3.6 \& Proposition~3.7]{Porta_Sala_Hall}. In \textit{loc.\ cit.}, $\scrC = \catQCoh(Y)$. The tensor product
	\begin{align}
		\catQCoh(\bfPerfps(\scrC) \times \bfPerfps(\scrC)) \otimes \scrC 
	\end{align}
	plays the role here of $\catQCoh(\bfPerf(Y) \times \bfPerf(Y) \times Y)$ in \textit{loc.\ cit.}, and indeed, when $\catQCoh(Y)$ is compactly generated, the two formul{\ae} agree. The main difference is that, while $\catQCoh(\bfPerf(Y) \times \bfPerf(Y) \times Y)$ has a natural symmetric monoidal structure and the pullback functor along
	\begin{align}
		q \colon \bfPerf(Y) \times \bfPerf(Y) \times Y \longrightarrow \bfPerf(Y) \times \bfPerf(Y) 
	\end{align}
	is symmetric monoidal, in the noncommutative setting $\catQCoh(\bfPerfps(\scrC) \times \bfPerfps(\scrC)) \otimes \scrC$ only has the structure of a categorical module over $\catQCoh(\bfPerfps(\scrC) \times \bfPerfps(\scrC))$. In particular, we cannot define the internal hom $\calHom(\calF_2^{\mathsf{u}}, \calF_1^{\mathsf{u}})$, but we can define its pushforward, using the canonical enrichment of $\catQCoh(\bfPerfps(\scrC) \times \bfPerfps(\scrC)) \otimes \scrC$ over $\catQCoh(\bfPerfps(\scrC) \times \bfPerfps(\scrC))$.
\end{remark}

\begin{lemma}\label{lem:perfect-object}
	One has $\calE \in \catPerf(\bfPerfps(\scrC)\times \bfPerfps(\scrC))$.
\end{lemma}

\begin{proof}
	It is enough to prove that for every $A \in \CAlg$ and every pair $\calF_1, \calF_2 \in \bfPerfps(\scrC)(A)$ one has
	\begin{align}
		\Hom_{\scrC_A}(\calF_2, \calF_1) \in \catPerf(A) \ .
	\end{align}
	First, note that since $\scrC$ is of finite type, it is in particular smooth. Therefore, Proposition~\ref{prop:basics_smooth_and_proper}-\eqref{prop:basics_smooth_and_proper-6} implies that that $\calF_1 \in \scrC_A^\omega$. At this point the statement follows from the fact that $\calF_2$ is pseudo-perfect .
\end{proof}

Consider the linear stack
\begin{align}
	\V_{\bfPerfps(\scrC, \tau) \times \bfPerfps(\scrC)}(\calE^\vee[-1]) \coloneqq \Spec_{\bfPerfps(\scrC) \times \bfPerfps(\scrC)}( \Sym(\calE^\vee[-1]) )\  .
\end{align}
The following is a non-commutative analogue of \cite[Proposition~3.7]{Porta_Sala_Hall}:
\begin{proposition} \label{prop:stack_extensions_is_relatively_affine}
	There is a natural commutative diagram
	\begin{align}
		\begin{tikzcd}[column sep = small, ampersand replacement = \&]
			\calS_2\bfPerfps(\scrC) \arrow{dr}[swap]{\partial_0 \times \partial_2} \arrow{rr}{\phi} \& \& \V_{\bfPerfps(\scrC) \times \bfPerfps(\scrC)}( \calE^\vee[-1] ) \arrow{dl}{\pi} \\
			{} \& \bfPerfps(\scrC) \times \bfPerfps(\scrC)
		\end{tikzcd} \ ,
	\end{align}
	where $\phi$ is furthermore an equivalence. In particular, $\partial_0 \times \partial_2$ is quasi-compact and universally finitely connected.
\end{proposition}

\begin{proof}
	Fix $S = \Spec(A) \in \dAff_k$ and a map $x \colon S \to \bfPerfps(\scrC) \times \bfPerfps(\scrC)$ classifying a pair $(\calF_{12}, \calF_{01})$ of pseudo-perfect functors $\scrC\op \to \Modd_A$.
	By definition,
	\begin{align}
		\calF_{12} \simeq x^\ast(\calF_1^{\mathsf u}) \qquad \text{and} \qquad \calF_{01} \simeq x^\ast(\calF_2^{\mathsf{u}}) \ .  
	\end{align}
	Moreover,
	\begin{align}
		\Map(S, \calS_2 \bfPerfps(\scrC)) \times_{\Map(S, \bfPerfps(\scrC) \times \bfPerfps(\scrC))} \{x\} 
	\end{align}
	is the maximal $\infty$-groupoid contained consisting of fiber sequences of the form $\calF_{01} \to \calF \to \calF_{12}$. Unraveling the definition, we can review this as a full subgroupoid of the maximal $\infty$-groupoid contained inside
	\begin{align}
		\calS_2\FunR_k(\scrC\op, \Modd_A) \times_{\Modd_A \times_k \Modd_A} \{(\calF_{12}, \calF_{01})\} \ . 
	\end{align}
	Using the stability of $\scrC_A \simeq \FunR_k(\scrC\op,\Modd_A)$ and \cite[Corollary~4.3.2.16]{HTT}, we see that the above $\infty$-category is equivalent to
	\begin{align}
		\Map_{\scrC_A}( \calF_{12}, \calF_{01}[1] ) \simeq \Fun(\Delta^1, \FunR_k(\scrC\op, \Modd_A)) \times_{\Modd_A \times \Modd_A} \{(\calF_{12}, \calF_{01}[1])\} \ . 
	\end{align}
	In particular, it is already an $\infty$-groupoid.
	Moreover,
	\begin{align}
		\Map_{\scrC_A}( \calF_{12}, \calF_{01}[1] ) & \simeq \tau_{\geqslant 0}\Hom_A( A, \Hom_{\scrC_A}( \calF_{12}, \calF_{01}[1] ) ) \\
		& \simeq \tau_{\geqslant 0} \Hom_A( \Hom_{\scrC_A}(\calF_{01},\calF_{12})^\vee[-1], A ) \\
		& \simeq \tau_{\geqslant 0} \Hom_A( x^\ast \calE, A ) \\
		& \simeq \Map_{\CAlg}( \Sym_A(x^\ast(\calE^\vee[-1])), A ) \ ,
	\end{align}
	which gives the identification with the functor of points of $\V_{\bfPerfps(\scrC) \times \bfPerfps(\scrC)}( \calE^\vee[-1] )$. The last statement follows from Lemma~\ref{prop:properties_fconn}--\eqref{item:properties_fconn-3}.
\end{proof}

Now, we shall see that also $\calS_1^\ell \bfFlagPerfps^{(m),\dagger}(\scrC;\sfV)$ and $\calS_1^r \bfFlagPerfps^{(m),\dagger}(\scrC;\sfV)$ are linear stacks. First, consider the morphism in $\mathbf\Delta\op$ 
\begin{align}
	\ev_{m-1} \colon [m-1] \longrightarrow [0]
\end{align}
that selects $m-1$. It induces a natural transformation $\ev_{m-1} \star \id_{[n]} \colon [m-1] \star [n] \to [0] \star [n]$ which yields, via Construction~\ref{construction:Vflags_generalized}, a natural morphism
\begin{align}
	\lambda_\bullet^{\ell} \colon \calS^\ell_\bullet \bfFlagPerfps^{(m),\dagger}(\scrC;\sfV) \longrightarrow \calS^\ell_\bullet \bfFlagPerfps^{(1)}(\scrC) \ , 
\end{align}
compatible with the structural maps to $\calS_\bullet \bfPerfps(\scrC)$. With respect to the notation of Remarks~\ref{rem:flag_action_coordinates} and \ref{rem:flag_action_coordinates_II}, the morphism
\begin{align}
	\lambda_0^{\ell}\colon \bfFlagPerfps^{(m),\dagger}(\scrC;\sfV)  \longrightarrow \bfPerfps(\scrC)
\end{align}
which sends a diagram of the form \eqref{eq:flags} to $F_{0,m}$. On the other hand, thanks to the description of $ \calS^\ell_1 \bfFlagPerfps^{(m),\dagger}(\scrC;\sfV)$ as a pullback (cf.\ Equation~\eqref{eq:identifications-2}), the morphism
\begin{align}
	\lambda_1^{\ell} \colon  \calS^\ell_1 \bfFlagPerfps^{(m),\dagger}(\scrC;\sfV) \longrightarrow  \calS^\ell_1 \bfFlagPerfps^{(1)}(\scrC) \simeq \bfFlagPerfps^{(2)}(\scrC) \simeq \calS_2\bfPerfps(\scrC)\ , 	
\end{align}
can be described as sending a diagram of the form \eqref{eq:flag_action} to the sub-diagram
\begin{align}
	\begin{tikzcd}[ampersand replacement=\&]
		0 \arrow{r} \& F_{0,m} \arrow{r} \arrow{d} \& F_{0,m+1} \arrow{d} \\
		\& 0 \arrow{r} \& F_{m,m+1} \arrow{d} \\
		\& \& 0
	\end{tikzcd} \ .
\end{align}

\begin{proposition}\label{prop:pullback-extension}
	The square
	\begin{align}
		\begin{tikzcd}[column sep = large, ampersand replacement=\&]
			\calS_1^\ell \bfFlagPerfps^{(m),\dagger}(\scrC;\sfV) \arrow{r}{\lambda_1^{\ell}} \arrow{d}{u_1^\ell\times \varpi_1} 	\& 	\calS_1^\ell \bfFlagPerf^{(1)}_{\mathsf{ps}}(\scrC) \arrow{d}{u_1^\ell\times \varpi_1} \\
			\bfPerfps(\scrC)  \times \bfFlagPerf^{(m),\dagger}_{\mathsf{ps}}(\scrC;\sfV) \arrow{r}{\id\times \lambda_0^{\ell}} \& 		\bfPerfps(\scrC)\times \bfFlagPerf^{(1)}_{\mathsf{ps}}(\scrC) 
		\end{tikzcd}
	\end{align}
	is a pullback.
\end{proposition}		

\begin{proof}
	Unraveling the definitions, the claim follows from the fact that the Waldhausen construction $\calS_\bullet \bfPerfps(\scrC)$ is a $2$-Segal space (applying \cite[Proposition~2.3.2-(3)]{DK} with $n = m+1$, $i = 0$ and $j = n-1 = m$).
\end{proof}

Since the construction of linear stacks commutes with square pullbacks, we get:
\begin{corollary}\label{cor:linear-stack-l}
	There is a natural commutative diagram
	\begin{align}
		\begin{tikzcd}[column sep = tiny, ampersand replacement = \&]
			\calS_1^\ell \bfFlagPerfps^{(m),\dagger}(\scrC;\sfV)  \arrow{dr}[swap]{u_1^\ell\times \varpi_1} \arrow{rr}{\phi} \& [-35pt] \& \V_{\bfPerfps(\scrC) \times \bfFlagPerf^{(m),\dagger}_{\mathsf{ps}}(\scrC;\sfV)}\left( \Big(\id\times \lambda_0^{\ell}\Big)^\ast \calE \right) \arrow{dl}{\pi} \\
			{} \& \bfPerfps(\scrC) \times \bfFlagPerf^{(m),\dagger}_{\mathsf{ps}}(\scrC;\sfV)  
		\end{tikzcd} \ ,
	\end{align}
	where $\phi$ is furthermore an equivalence.
	In particular, $u_1^\ell \times \varpi_1$ is quasi-compact and universally finitely connected.
\end{corollary}

By using similar arguments, we obtain the following characterization of $\calS_1^r \bfFlagPerfps^{(m),\dagger}(\scrC;\sfV)$.
\begin{proposition}\label{prop:linear-stack-r}
	The square
	\begin{align}
		\begin{tikzcd}[column sep = large, ampersand replacement=\&]
			\calS_1^r \bfFlagPerfps^{(m),\dagger}(\scrC;\sfV) \arrow{r}{\lambda_0^{r}} \arrow{d}{\varpi_0 \times u_1^r} 	\& 	\calS_1^r \bfFlagPerf^{(1)}_{\mathsf{ps}}(\scrC) \arrow{d}{\varpi_0\times u_1^r} \\
			\bfFlagPerf^{(m),\dagger}_{\mathsf{ps}}(\scrC;\sfV) \times \bfPerfps(\scrC) \arrow{r}{\lambda_1^{r} \times \id} \& \bfFlagPerf^{(1)}_{\mathsf{ps}}(\scrC) \times \bfPerfps(\scrC)
		\end{tikzcd}
	\end{align}
	is a pullback, where $\lambda_n^r$ is induced by $\id_{[n]} \star \ev_0 \colon [n] \star [m-1] \to [n] \star [0]$. In particular, there is a natural commutative diagram
	\begin{align}
		\begin{tikzcd}[column sep = tiny, ampersand replacement = \&]
			\calS_1^r \bfFlagPerfps^{(m),\dagger}(\scrC;\sfV)  \arrow{dr}[swap]{\varpi_0\times u_1^r} \arrow{rr}{\phi} \& [-35pt] \& \V_{\bfFlagPerf^{(m),\dagger}_{\mathsf{ps}}(\scrC;\sfV) \times \bfPerfps(\scrC)}\left( \Big(\lambda_1^r \times \id\Big)^\ast \calE \right) \arrow{dl}{\pi} \\
			{} \& \bfFlagPerf^{(m),\dagger}_{\mathsf{ps}}(\scrC;\sfV) \times  \bfPerfps(\scrC)
		\end{tikzcd} \ ,
	\end{align}
	where $\phi$ is furthermore an equivalence. In particular, $\varpi_0 \times u_1^r$ is quasi-compact and universally finitely connected.
\end{proposition}

\subsubsection{The relative tangent complexes of $u_1^\ell\times \varpi_1$ and $\varpi_0\times u_1^r$} \label{subsec:m_flags_cotangent_complex}

In Corollary~\ref{cor:COHA-representations} (resp.\ Corollary~\ref{cor:COHA-representations-right}), we assumed that the map
\begin{align}
	u_1^\ell\times \varpi_1  \colon \calS_1^\ell \bfFlagPerf_{\bfH,\bfM}^{(m),\dagger}(\scrC;\sfV) \longrightarrow \bfH \times \bfM \quad \big(\text{resp.}\  \varpi_0\times u_1^r \colon \calS_1^r \bfFlagPerf_{\bfH,\bfM}^{(m),\dagger}(\scrC;\sfV) \longrightarrow \bfM\times \bfH \big)
\end{align}
is derived lci. We provide a method to compute the cotangent complex of this morphism in the simplest case, i.e., $\bfM = \bfFlagPerfps^{(m),\dagger}(\scrC;\sfV)$ and $\bfH = \bfPerfps(\scrC)$.

\begin{proposition}\label{prop:flag_tangent_complex}
	The relative tangent complex $\T_x$ of the map
	\begin{align}
		u_1^\ell\times \varpi_1 \colon \calS_1^\ell \bfFlagPerfps^{(m),\dagger}(\scrC;\sfV) \longrightarrow \bfPerfps(\scrC) \times \bfFlagPerf^{(m),\dagger}_{\mathsf{ps}}(\scrC;\sfV)  
	\end{align}
	at a point $x \colon \Spec(A) \to \calS_1^\ell \bfFlagPerfps^{(m),\dagger}(\scrC;\sfV)$ classifying a diagram of the form \eqref{eq:flag_action} fits into the following natural fiber sequence
	\begin{align}
		\T_x \longrightarrow \Hom_{\scrC_A}(F_{0,m+1},F_{0,m})[1] \oplus  \Hom_{\scrC_A}(F_{m,m+1}, F_{0,m+1})[1] \longrightarrow \Hom_{\scrC_A}(F_{0,m+1}, F_{0,m+1})[1] \ . 
	\end{align}
\end{proposition}

\begin{proof}
	Thanks to Proposition~\ref{prop:pullback-extension}, we apply \cite[Proposition~3.2]{Porta_Sala_Hall}\footnote{Note that in \textit{loc.\ cit.}, we computed the relative cotangent complex, while here we are computing the relative tangent complex. The two are related by duality.} to the map $u_1^\ell \times \varpi_1$, that implies that $\T_x$ is identified with the limit of the following diagram:
	\begin{align}
		\begin{tikzcd}[ampersand replacement = \&]
			\& \& 0 \arrow{d} \\
			\& \Hom_{\scrC_A}(F_{0,m+1},F_{0,m+1})[1] \arrow{r} \arrow{d} \& \Hom_{\scrC_A}(F_{0,m+1},F_{m,m+1})[1] \\
			0 \arrow{r} \& \Hom_{\scrC_A}(F_{0,m}, F_{0,m+1})[1]
		\end{tikzcd}\ .
	\end{align}
	Thus, we see that $\T_x$ fits in the following pullback square
	\begin{align}
		\begin{tikzcd}[ampersand replacement = \&]
			\T_x \arrow{r} \arrow{d} \& \Hom_{\scrC_A}(F_{0,m+1},F_{0,m})[1] \arrow{d} \\
			\Hom_{\scrC_A}(F_{m,m+1},F_{0,m+1})[1] \arrow{r} \& \Hom_{\scrC_A}(F_{0,m+1},F_{0,m+1})[1] 
		\end{tikzcd} \ ,
	\end{align}
	whence the conclusion.
\end{proof}

Similarly, one can prove the following.
\begin{proposition}
	The relative tangent complex $\T_x$ of the map
	\begin{align}
		\varpi_0\times u_1^r \colon \calS_1^r \bfFlagPerfps^{(m),\dagger}(\scrC;\sfV) \longrightarrow \bfFlagPerf^{(m),\dagger}_{\mathsf{ps}}(\scrC;\sfV) \times \bfPerfps(\scrC) 
	\end{align}
	at a point $x \colon \Spec(A) \to \calS_1^r \bfFlagPerfps^{(m),\dagger}(\scrC;\sfV)$ classifying a diagram of the form \eqref{eq:flag_action} fits into the following natural fiber sequence
	\begin{align}
		\T_x \longrightarrow \Hom_{\scrC_A}(F_{0,m+1},F_{1,m+1})[1] \oplus  \Hom_{\scrC_A}(F_{0,1}, F_{0,m+1})[1] \longrightarrow \Hom_{\scrC_A}(F_{0,m+1}, F_{0,m+1})[1] \ . 
	\end{align}
	\end{proposition}

\section{COHA and CatHA of a finite type category with a fixed $t$-structure}\label{sec:COHA-t-structure}

In this section, we construct a cohomological (in the motivic sense) and categorical Hall algebra associated to a stable $\infty$-category equipped with a $t$-structure, satisfying some natural conditions. We provide new examples of COHAs and CatHAs associated to \textit{2-Calabi-Yau completions}. In the quiver case, we recover Schiffmann-Vasserot's construction of COHAs via \textit{Hamiltonian reduction} and categorify it.

Fix a base ring $k$ and a $k$-linear compactly generated stable $\infty$-category $\scrC\in \PrLomega_k$ equipped with a $t$-structure $\tau = (\scrC_{\geqslant 0}, \scrC_{\leqslant 0})$ satisfying Assumption~\ref{assumption:t_structure_filtered_colimits}.

\subsection{Construction of the COHA and CatHA}

Let us make the following assumption on $(\scrC,\tau)$:
\begin{assumption}\label{assumption:existence-COHA}
	\hfill
	\begin{enumerate}[label=(C.\arabic*)]\itemsep0.2cm
		\item \label{assumption:existence-COHA-1} $\scrC$ is of finite type;
		\item  \label{assumption:existence-COHA-2} $(\scrC,\tau)$ satisfies Assumption~\ref{assumption:t_structure_filtered_colimits};
		\item \label{assumption:existence-COHA-3} the $t$-structure $\tau$ universally satisfies openness of flatness in the sense of Definition~\ref{def:openness-flatness};
		\item \label{assumption:existence-COHA-4} $\sfS_\scrC^![2]$ is $t$-exact with respect to $\tau$, where $\sfS_\scrC^!$ is the left Serre functor\footnote{The notion of Serre functors is introduced in Definition~\ref{def:Serre_functors}.} of $\scrC$.
	\end{enumerate}
\end{assumption}
The above assumption implies that the derived stack $\bfCohps(\scrC, \tau)$ of $\tau$-flat pseudo-perfect objects of $\scrC$ is a geometric derived stack locally of finite presentation over $k$, by Proposition~\ref{prop:openness}.

Setting $\bfH\coloneqq \bfCohps(\scrC, \tau)$ in the construction discussed in \S\ref{subsec:COHA-representations}, we obtain a simplicial stack
\begin{align}
	\calS_\bullet \bfCohps(\scrC, \tau)\coloneqq  \calS_\bullet \bfPerf_\bfH(\scrC)\ .
\end{align}

The following is a generalization of  \cite[Lemma~4.1]{Porta_Sala_Hall}, obtained by applying Proposition~\ref{prop:2-Segal-algebra}.
\begin{lemma}\label{lem:2-Segal-Coh}
	The square \eqref{eq:pullback-square-T} for $\bfCohps(\scrC, \tau)$ is a pullback. In particular, $\calS_\bullet \bfCohps(\scrC, \tau)$ is a 2-Segal derived stack.
\end{lemma}

\begin{proof}
	By the descriptions \eqref{eq:extensions} and \eqref{eq:extensions-T}, it is easy to see that if the extreme terms in a fiber sequence are $\tau$-flat, also the middle one is $\tau$-flat. This implies that the assumption of Proposition~\ref{prop:2-Segal-algebra} is satisfied. Hence, the claim holds.
\end{proof}

In addition, using Proposition~\ref{prop:stack_extensions_is_relatively_affine} and the same arguments of the proof of \cite[Proposition~3.10]{Porta_Sala_Hall}, we get:
\begin{lemma}\label{lem:derived-lci}
	Assume that there exists a nonzero positive integer $n_\scrC\in \N$ such that $\sfS_\scrC^![n_\scrC]$ is $t$-exact with respect to $\tau$, where $\sfS_\scrC^!$ is the left Serre functor of $\scrC$. Then, the relative cotangent complex of the map 
	\begin{align}
		\partial_0 \times \partial_2 \colon \calS_2\bfCohps(\scrC, \tau) \longrightarrow \bfCohps(\scrC, \tau)\times \bfCohps(\scrC, \tau)
	\end{align}
	is perfect and has tor-amplitude $[-1, n_\scrC - 1]$.
\end{lemma}

\begin{proof}
	This is essentially the content of \cite[Proposition~3.11]{Porta_Sala_Hall}, and it is a direct consequence of Sere duality (Proposition~\ref{prop:basics_Serre_functors}--\eqref{prop:basics_Serre_functors-2}). We briefly reproduce the proof for the convenience of the reader. Thanks to the analysis carried out in Proposition~\ref{prop:stack_extensions_is_relatively_affine}, it suffices to show that the tor-amplitude of
	\begin{align}
		\calE^\vee[-1] \simeq \big(\Hom_{\scrC \otimes \catQCoh(\bfCohps(\scrC,\tau) \times \bfCohps(\scrC,\tau)}(\calF_2^{\mathsf u}, \calF_1^{\mathsf u})\big)^\vee[-1] 
	\end{align}
	is concentrated in $[-1,n_\scrC-1]$. This can be checked on closed points, so we can fix a map
	\begin{align}
		x \colon \Spec(\C) \longrightarrow \bfCohps(\scrC,\tau) \times \bfCohps(\scrC,\tau) 
	\end{align}
	classifying two $\tau$-flat pseudo-perfect objects $\calF_1, \calF_2 \in \scrC$. Since $\scrC$ is smooth, Proposition~\ref{prop:basics_smooth_and_proper}--\eqref{prop:basics_smooth_and_proper-6} implies that $\calF_1$ and $\calF_2$ are compact. Using Proposition~\ref{prop:basics_Serre_functors}--\eqref{prop:basics_Serre_functors-2} we see that the left Serre functor realizes an equivalence
	\begin{align}\label{eq:lciness_Serre}
		\Hom_\scrC(\calF_2, \calF_1)^\vee[-1] \simeq \Hom_\scrC(\mathsf S^!_\scrC(\calF_1), \calF_2)[-1] \ .
	\end{align}
	Since $\calF_1$ and $\calF_2$ are pseudo-perfect and $\tau$-flat, we have
	\begin{align}
		\pi_i \Hom_\scrC(\calF_2, \calF_1) \simeq \Ext^{-i}(\calF_2, \calF_1) \simeq 0 
	\end{align}
	for every $i > 0$. Since $\Hom_\scrC(\calF_2,\calF_1) \in \catPerf(\C)$, it automatically follows that
	\begin{align}
		\pi_i\big( \Hom_\scrC(\calF_2, \calF_1)^\vee[-1] \big) \simeq 0 
	\end{align}
	for $i < -1$. On the other hand, since $\mathsf S^!_\scrC[n_\scrC]$ is $\tau$-exact, it follows that $\mathsf S^!_\scrC(\calF_2)[n_\scrC]$ is again pseudo-perfect and $\tau$-flat. Applying the same argument given above, we deduce that
	\begin{align}
		\pi_{i}\big( \Hom_\scrC(\mathsf S_\scrC^!(\calF_1), \calF_2)[-1] \big) \simeq \pi_{i-1+n_\scrC} \big( \Hom_\scrC(\mathsf S^!_\scrC(F_1)[n_\scrC], \calF_2) \big) 
	\end{align}
	vanishes for $i + 1 - n_\scrC > 0$. Thus, the conclusion follows from the equivalence \eqref{eq:lciness_Serre}.
\end{proof}

Lemmas~\ref{lem:2-Segal-Coh} and \ref{lem:derived-lci}, and Corollary~\ref{cor:COHA} yield the main theorem of this section:
\begin{theorem} \label{thm:COHA-Coh}
	If Assumption~\ref{assumption:existence-COHA} holds and the map
	\begin{align}\label{eq:partial-1-proper}
		\partial_1 \colon \calS_2 \bfCohps(\scrC, \tau) \longrightarrow\bfCohps(\scrC, \tau)
	\end{align}
	is locally rpas, then $\catCohb_{\mathsf{pro}}( \bfCohps(\scrC, \tau) )$ has the structure of an $\E_1$-monoidal stable pro-$\infty$-category, whose underlying tensor product is given by the composition
	\begin{align}\label{eq:multiplication}
		\begin{tikzcd}[ampersand replacement=\&]
			\catCohb_{\mathsf{pro}}( \bfCohps(\scrC, \tau) ) \times \catCohb_{\mathsf{pro}}( \bfCohps(\scrC, \tau) ) \ar{r}{\boxtimes} \& \catCohb_{\mathsf{pro}}( \bfCohps(\scrC, \tau) \times \bfCohps(\scrC, \tau) ) \\
			{} \ar{r}{(\partial_1)_\ast \circ (\partial_0 \times \partial_2)^\ast} \& \catCohb_{\mathsf{pro}}( \bfCohps(\scrC, \tau) )
		\end{tikzcd}\ .
	\end{align}
	
	Similarly, let $\bfD^\ast$ be a motivic formalism, and fix $\calA \in \CAlg(\bfD^\ast(\Spec(k)))$ and $\Gamma \subseteq \Pic(\bfD^\ast(\Spec(k)))$ such that Assumption~\ref{assumption:motivic_formalism} is satisfied. Then, the topological vector space $\HBMDGamma_0(\bfCohps(\scrC, \tau);\calA)$ becomes a unital associative algebra. In particular,
	\begin{align}
		G_0( \bfCohps(\scrC, \tau) )\quad \text{and} \quad \HBM_\ast( \bfCohps(\scrC, \tau) )
	\end{align}
	become unital associative algebras.
\end{theorem}		

\begin{remark}
	In \cite[\S5]{DHSM_BPS}, the authors constructed a cohomological Hall algebra associated to a moduli stack of objects of an admissible finite length abelian 2-Calabi-Yau subcategory of a fixed $\C$-linear finite type stable $\infty$-category. The construction is performed under three assumptions. Their assumption (1) corresponds to our assumption on the properness of the map \eqref{eq:partial-1-proper}.  Their assumption (2) ensures the existence of the pullback $(\partial_0\times \partial_2)^\ast$: in our case, the 2-Calabi-Yau property on the stable $\infty$-category yields straightforwardly the existence of $(\partial_0 \times \partial_2)^\ast$ thanks to Lemma~\ref{lem:perfect-object}, Proposition~\ref{prop:stack_extensions_is_relatively_affine}, and Lemma~\ref{lem:derived-lci}. Finally, their assumption (3) is needed in their setting to prove the associativity of the product: for us the associativity follows directly from the 2-Segal space structure of $\calS_\bullet \bfCohps(\scrC, \tau)$. Thus, Theorem~\ref{thm:COHA-Coh} provides a refinement and a categorification of the construction of the COHA in \cite{DHSM_BPS}.
\end{remark}

\subsection{COHA and CatHA associated to a Serre subcategory}

In this section we define the cohomological and categorical Hall algebra associated to an open substack of $\bfCohps(\scrC, \tau)$ which roughly speaking parametrizes objects belonging to a Serre subcategory of $\scrC^\heartsuit$. 

We need first the following preliminary result.
\begin{proposition} \label{prop:open_implies_fully_faithful}
	Let $f \colon F \to G$ be a morphism of derived stacks. If $f$ is representable by open immersions, then for every commutative algebra $A \in \CAlg$, the induced map $f(A) \colon F(A) \to G(A)$ is $(-1)$-truncated in the $\infty$-category $\calS$ of spaces, and hence it is fully faithful. Furthermore, $f$ is an equivalence if and only if for every field $k$, every morphism $\Spec(k) \to G$ factors through $f$.
\end{proposition}

\begin{proof}
	For the first half, it is enough to observe that the diagonal morphism $F \to F \times_G F$ is an equivalence, which is obvious from the fact that $f$ is representable by open immersions. For the second half, it is enough to prove that for every $A \in \CAlg$ the map $f_A\colon F_A \coloneqq \Spec(A) \times_G F \to \Spec(A)$ is an equivalence. By assumption $f_A$ is an open subscheme of $\Spec(A)$, and hence it is an isomorphism if and only if for every field $k$, every morphism $\Spec(k) \to \Spec(A)$ has image contained in $F_A$. The conclusion follows.
\end{proof}

\begin{corollary}
	Let $\bfU \to \bfCohps(\scrC,\tau)$ be a morphism representable by open immersions and $A \in \CAlg$. Then, $\bfU(A)$ is a full subgroupoid of $\bfCohps(\scrC,\tau)(A) = \catCohps(\scrC_A,\tau_A)^\simeq$.
\end{corollary}

\begin{definition}
	Let $A \in \CAlg$ and let $\scrC^\heartsuit_A$ for the heart of the induced $t$-structure on $\scrC_A$.
	We say that an object $F \in \scrC_A^\heartsuit$ is \textit{pseudo-compact} if for every $G \in \scrC_A^\omega$, one has $\Hom_{\scrC_A}(G,F) \in \catPerf(A)$.
	We denote by $\PsCpt(\scrC_A,\tau_A)$ the full subcategory of $\scrC_A^\heartsuit$ spanned by pseudo-compact objects.
\end{definition}

\begin{warning}
	Since short exact sequences in $\scrC^\heartsuit_A$ are fiber sequences in $\scrC_A$, we see that the category $\PsCpt(\scrC_A,\tau_A)$ has the following properties:
	\begin{enumerate}\itemsep=0.2cm
		\item it is closed under extensions in $\scrC_A^\heartsuit$;
		
		\item a morphism in $\PsCpt(\scrC_A,\tau_A)$ that is a monomorphism (resp.\ an epimorphism) in $\scrC_A$ admits a cokernel (resp.\ a kernel) in $\PsCpt(\scrC_A,\tau_A)$.
	\end{enumerate}
	Still, $\catCohps^\heartsuit(\scrC_A, \tau_A)$ might fail to be abelian, as it is not guaranteed that every morphism admits a (co)kernel. If the pseudo-perfect objects of $\scrC_A^\heartsuit$ are closed under subobjects (or, equivalently, under quotients), then $\catCohps^\heartsuit(\scrC_A,\tau_A)$ is automatically abelian.\footnote{We are very grateful to Enrico Lampetti for explaining to us the relevance of this assumption in the study of moduli of objects and in establishing the existence of good moduli spaces \cite{Lampetti_Good_moduli}.}
\end{warning}

In practice, we will only consider $\PsCpt(\scrC_A,\tau_A)$ when $A$ is either a field or a DVR. Notice that, in both cases, $A$ is regular.

\begin{example}
	Assume that $A$ is underived, regular and noetherian. We collect here a two notable cases in which $\PsCpt(\scrC_A,\tau_A)$ is closed under subobjects, and hence in which it is abelian.
	\begin{enumerate}\itemsep=0.2cm
		\item Let $X$ be a quasi-compact smooth scheme over $\Spec(k)$ and let $\scrC \coloneqq \catQCoh(X)$, equipped with the standard $t$-structure. Then an object $\calF \in \scrC_A^\heartsuit \simeq \catQCoh^\heartsuit(X \times \Spec(A))$ is pseudo-compact if and only if it is coherent and has proper support. Since the proper support condition is closed under monomorphism, we see that $\PsCpt(\scrC_A,\tau_A)$ is closed under subobjects and it is therefore abelian.
		
		\item Let $\scrD_0$ be a small category and set $\scrC \coloneqq \Fun(\scrD_0, \Modd_k)$, equipped with the standard (objectwise) $t$-structure. Then an object in
		\begin{align}
			\scrC_A^\heartsuit \simeq \Fun(\scrD_0, \Modd_A^\heartsuit) 
		\end{align}
		is pseudo-perfect if and only if it takes values in finitely generated $A$-modules. Since $A$ is noetherian and since a morphism in $\scrC_A^\heartsuit$ is a monomorphism if and only if it is a monomorphism objectwise, we conclude that $\PsCpt(\scrC_A,\tau)$ is closed under subobjects in this situation. Notice that this example covers quiver representations, local systems and constructible sheaves \cite{Porta_Teyssier_Exodromy,HPT_Beyond}. In \cite{Lampetti_Good_moduli}, it was proved that the perverse $t$-structure equally satisfies the above assumption.
	\end{enumerate}
\end{example}

In addition, note that when $A$ is underived, $\catCohps(\scrC_A,\tau_A) \subseteq \PsCpt(\scrC_A,\tau_A)$, with equality holding when $A$ is a field.

\begin{definition}\label{def:category_associated_to_open_substack}
	Let $\bfU \to \bfCohps(\scrC,\tau)$ be a morphism representable by open immersions and $A \in \CAlg$. We let $\PsCpt_\bfU(\scrC_A,\tau_A)$ be the full subcategory of $\PsCpt(\scrC_A,\tau_A)$ spanned by the objects that belong to the image of $\bfU(A)$.
\end{definition}

Now, we fix an open substack $\alpha\colon \bfT\to \bfCohps(\scrC, \tau)$. By setting $\bfH\coloneqq \bfT$ in the construction introduced in \S\ref{subsec:COHA-representations}, we obtain a simplicial stack
\begin{align}
	\calS_\bullet \bfCoh_\bfT(\scrC, \tau)\coloneqq  \calS_\bullet \bfPerf_\bfT(\scrC)\ .
\end{align}
The following lemma holds even when $\PsCpt(\scrC_\kappa,\tau_\kappa)$ is not abelian:
\begin{lemma}\label{lem:Serre_subcategory}
	\hfill
	\begin{enumerate}\itemsep0.2cm
		\item \label{item:Serre_subcategory-1} If for every field $\kappa$ the category $\PsCpt_\bfT(\scrC_\kappa,\tau_\kappa)$ is closed under extensions in $\PsCpt(\scrC_\kappa,\tau_\kappa)$, then the square
		\begin{align}
			\begin{tikzcd}[ampersand replacement=\&]
				\calS_2 \bfCoh_\bfT(\scrC,\tau) \arrow{r} \arrow{d}{\partial_0 \times \partial_2} \& \calS_2\bfCohps(\scrC,\tau) \arrow{d}{\partial_0 \times \partial_2} \\
				\bfT \times \bfT \arrow{r} \& \bfCohps(\scrC,\tau) \times \bfCohps(\scrC,\tau)
			\end{tikzcd} 
		\end{align}
		is a pullback.
		
		\item \label{item:Serre_subcategory-2} If for every field $\kappa$ the category $\PsCpt_\bfT(\scrC_\kappa,\tau_\kappa)$ is closed under subobjects and quotients inside $\PsCpt(\scrC_\kappa,\tau_\kappa)$, then the square
		\begin{align}
			\begin{tikzcd}[ampersand replacement=\&]
				\calS_2 \bfCoh_\bfT(\scrC,\tau) \arrow{d}{\partial_1} \arrow{r} \& \calS_2 \bfCohps(\scrC,\tau) \arrow{d}{\partial_1} \\
				\bfT \arrow{r} \& \bfCohps(\scrC,\tau)
			\end{tikzcd}
		\end{align}
		is a pullback. This holds in particular when $\PsCpt(\scrC_\kappa,\tau_\kappa)$ is abelian and $\PsCpt_\bfT(\scrC_\kappa,\tau_\kappa)$ is a Serre abelian subcategory of  $\PsCpt(\scrC_\kappa,\tau_\kappa)$.
	\end{enumerate}
\end{lemma}

\begin{proof}
	We give the argument for Statement~\eqref{item:Serre_subcategory-1}, as Statement~\eqref{item:Serre_subcategory-2} follows in an analogous way. Since $\bfT \to \bfCohps(\scrC,\tau)$ is representable by open immersions, Proposition~\ref{prop:open_implies_fully_faithful} shows that it is enough to prove that for every field $\kappa$ the square obtained by taking $\kappa$-points is a pullback. Proposition~\ref{prop:open_implies_fully_faithful} also guarantees that the bottom horizontal map is fully faithful, so it follows that the top horizontal one is fully faithful as well. To check essential surjectivity, recall that a $\kappa$-point of $\calS_2\bfCohps(\scrC,\tau)$ can be represented as a fiber sequence
	\begin{align}
		\mathsf E \coloneqq E_{01} \longrightarrow E_{02} \longrightarrow E_{12} 
	\end{align}
	in $\scrC_\kappa$, where in addition the $E_{ij}$ are pseudo-perfect and $\tau$-flat. Since $\kappa$ is a field, we have that $E_{ij} \in \scrC_\kappa^\heartsuit$, and therefore the $E_{ij}$ belong to the abelian category $\catCohps^\heartsuit(\scrC_\kappa,\tau_\kappa)$. This immediately implies that the map on the left is a monomorphism, and that the one on the right is an epimorphism in $\catCohps^\heartsuit(\scrC_\kappa,\tau_\kappa)$. By the definition of being closed under extensions, if $\partial_0(\mathsf E) = E_{01}$ and $\partial_2(\mathsf E) = E_{12}$ are in $\catCoh_\bfT(\scrC_k,\tau_k)$, the same goes for $\partial_1(\mathsf E) = E_{02}$. The conclusion follows.
\end{proof}

By using Corollary~\ref{cor:COHA} for $\bfH\coloneqq \bfT$ and Lemmas~\ref{lem:2-Segal-Coh}, \ref{lem:derived-lci},  and \ref{lem:Serre_subcategory}--\eqref{item:Serre_subcategory-1}, we obtain the following.
\begin{corollary}\label{cor:induced_COHA-Serre}
	Let $\alpha\colon \bfT\to \bfCohps(\scrC, \tau)$ be an open substack of $\bfCohps(\scrC, \tau)$ such that for every field $\kappa$ the category $\PsCpt_\bfT(\scrC_\kappa,\tau_\kappa)$ is closed under extensions in $\PsCpt(\scrC_\kappa,\tau_\kappa)$. If Assumption~\ref{assumption:existence-COHA} holds and the map
	\begin{align}
		\partial_1 \colon \calS_2 \bfCoh_{\bfT}(\scrC, \tau) \longrightarrow\bfT
	\end{align}
	is locally rpas, then $\catCohb_{\mathsf{pro}}( \bfT )$ has the structure of an $\E_1$-monoidal stable pro-$\infty$-category, whose underlying tensor product is induced by \eqref{eq:multiplication} with $\bfT$ instead of $\bfCohps(\scrC, \tau)$. 
	
	Similarly, let $\bfD^\ast$ be a motivic formalism, and fix $\calA \in \CAlg(\bfD^\ast(\Spec(k)))$ and $\Gamma \subseteq \Pic(\bfD^\ast(\Spec(k)))$ such that Assumption~\ref{assumption:motivic_formalism} is satisfied. Then, the topological vector space $\HBMDGamma_0(\bfT;\calA)$ becomes a unital associative algebra. In particular,
	\begin{align}
		G_0( \bfT )\quad \text{and} \quad \HBM_\ast( \bfT )
	\end{align}
	become unital associative algebras.
\end{corollary}

\begin{remark}
	If $\partial_1 \colon \calS_2 \bfCohps(\scrC, \tau) \longrightarrow\bfCohps(\scrC, \tau)$ is locally rpas and for every field $\kappa$ the category $\PsCpt_\bfT(\scrC_\kappa,\tau_\kappa)$ is closed under extensions in $\PsCpt(\scrC_\kappa,\tau_\kappa)$ (e.g.\ if the latter is abelian and the former is a Serre abelian subcategory), then the map $\partial_1 \colon \calS_2 \bfCoh_{\bfT}(\scrC, \tau) \longrightarrow\bfT$ is locally rpas as well.
	This follows combining Lemma~\ref{lem:Serre_subcategory}--\eqref{item:Serre_subcategory-2} and Lemma~\ref{lem:properties_locallly_rpas}.
\end{remark}

\subsection{COHAs and CatHAs, and Bridgeland stability conditions}

Let $X$ be a smooth projective complex variety and let $\scrD$ be a $\C$-linear, strong (in the sense of \cite[Definition~3.5]{BLMNPS_Stability}) semiorthogonal component of $\catPerf(X)$ of finite cohomological amplitude (in the sense of \cite[Definition~3.7]{BLMNPS_Stability}) equipped with a Serre functor $\sfS_\scrD\simeq [2]$. Let $\sigma$ be a stability condition on $\scrD$ with respect to a finite rank free abelian monoid $\Lambda$ (in the sense of \cite[Definition~21.15]{BLMNPS_Stability}). 

Denote by $\bfT$ either the moduli stack $\bfCohps(\scrC, \tau)$, where $\tau$ is the induced $t$-structure on $\scrC\coloneqq \Ind(\scrD)$ whose heart is $\Ind(\calP_\sigma((0,1]))$, where $\calP_\sigma$ is the slicing associated to $\sigma$, or the moduli stack $\bfCohps^{\sigma\textrm{-}\mathsf{ss}, \mu}(\scrC, \tau)$ of $\sigma$-semistable objects of fixed slope $\mu$ on $\scrC$, or the moduli stack
\begin{align}
	\bigsqcup_{\mathbf{v}\in \Z\calS}\, \bfCohps^{\sigma\textrm{-}\mathsf{ss}, \mu}(\scrC, \tau; \mathbf{v})\ ,
\end{align}
where $\bfCohps^{\sigma\textrm{-}\mathsf{ss},\mu}(\scrD, \tau; \mathbf{v})$ denotes the moduli stack of $\sigma$-semistable objects on $\scrD$ with slope $\mu$ and with Mukai vector $\mathbf{v}\in \Lambda$. Here, $\calS$ is a set of Mukai vectors spanning a sublattice of $\Lambda$. 
\begin{corollary}\label{cor:COHA-stability-condition}
	$\catCohb_{\mathsf{pro}}( \bfT )$ has the structure of an $\E_1$-monoidal stable pro-$\infty$-category, whose underlying tensor product is induced by \eqref{eq:multiplication}. 
	
	Similarly, let $\bfD^\ast$ be a motivic formalism, and fix $\calA \in \CAlg(\bfD^\ast(\Spec(k)))$ and $\Gamma \subseteq \Pic(\bfD^\ast(\Spec(k)))$ such that Assumption~\ref{assumption:motivic_formalism} is satisfied. Then, the topological vector space $\HBMDGamma_0(\bfT;\calA)$ becomes a unital associative algebra. In particular,
	\begin{align}
		G_0( \bfT )\quad \text{and} \quad \HBM_\ast( \bfT )
	\end{align}
	become unital associative algebras.
\end{corollary}

\begin{proof}
	First, note that the assumptions of Theorem~\ref{thm:COHA-Coh} holds. Indeed, $\scrC\coloneqq \Ind(\scrD)\subset \catQCoh(X)$ is a smooth and proper category. Moreover, (the induced $t$-structure) $\tau$ is compatible with filtered colimits and it is right complete. Notice now that, although the induced $t$-structure on $\scrC$ might fail to be left complete, Assumption~\ref{item-2-assumption:t_structure_filtered_colimits} is always satisfied: indeed, if $f \colon A \to B$ is a faithfully flat morphism in $\CAlg$, then we have a canonically commutative diagram
	\begin{align}
		\begin{tikzcd}[ampersand replacement=\&]
			\scrC_A \arrow{r}{f^\ast} \arrow[hook]{d} \& \scrC_B \arrow[hook]{d} \\
			\catQCoh(X \times \Spec(A)) \arrow{r}{f^\ast} \& \catQCoh(X \times \Spec(B)) 
		\end{tikzcd} \ .
	\end{align}
	Faithfully flat descent for $\catQCoh(X)$ implies that the bottom horizontal functor is conservative, and the two vertical functors are fully faithful by construction. Thus, the top horizontal functor is an equivalence as well. Moreover, the induced $t$-structure universally satisfies openness of flatness by \cite[Proposition~20.8]{BLMNPS_Stability}, while condition (3) of Definition~20.5, Corollary~20.10, and Proposition~11.11 in \textit{loc.\ cit.} yields the properness of the Quot spaces in this case, hence it ensures that the map $\partial_1$ is locally rpas. Finally, the subcategory of semistable objects of fixed slope is a Serre subcategory, and the corresponding moduli stack is open (cf.\ \cite[Lemma~21.12]{BLMNPS_Stability}). Thus, the claim follows.
\end{proof}

\begin{remark}
	The previous proposition can be applied when $\scrD$ is either $\catPerf(S)$, with $S$ a K3 surface, or the Kuznetsov component $\mathcal{K}u(X)$, with $X$ either a Fano 3fold of Picard rank one (different from the complete intersection of a quadric and a cubic in $\PP^5$), a smooth cubic 4fold in $\PP^5$, or a Gushel-Mukai variety. The existence of a stability condition on $\mathcal{K}u(X)$ is proved in the first two cases in \cite[Theorem~1.1]{BLMS_Kuznetsov}, while in the third case is proved in \cite[Theorem~1.2]{PPZ_Kuznetsov}.
\end{remark}

\section{COHAs, CatHAs, and their representations associated to a torsion pair}\label{sec:COHA_representations_torsion_pairs}

In this section we fix a base ring $k$ and a $k$-linear compactly generated presentable stable $\infty$-category $\scrC\in \PrLomega_k$ equipped with a $t$-structure $\tau$. Our main goal is to prove that a torsion pair on $\scrC^\heartsuit$ gives rise to several Hall algebras and their representations.

\subsection{Preliminaries and Notation}

Set $\bfH \coloneqq \bfCohps(\scrC,\tau)$.

Fix furthermore an $(m-1)$-flag $\sfV$ and let $\bfM\coloneqq \bfFlagCohps^{(m),\dagger}(\scrC,\tau;\sfV) $ be the open substack of $\bfFlagPerfps^{(m),\dagger}(\scrC;\sfV)$ parametrizing $\sfV$-flags of the form \eqref{eq:flags} where we ask $F_{i,m}$ to be $\tau$-flat for $i = 0, 1, \ldots, m-1$.
Then the assumptions of Proposition~\ref{prop:2-Segal-representation} and of its ``right'' version are satisfied, and we write
\begin{align}
	\calS_\bullet^\ell\bfFlagCohps^{(m),\dagger}(\scrC,\tau;\sfV) \longrightarrow \calS_\bullet \bfCohps(\scrC,\tau) \quad\text{and}\quad
	\calS_\bullet^r\bfFlagCohps^{(m),\dagger}(\scrC,\tau;\sfV) \longrightarrow \calS_\bullet \bfCohps(\scrC,\tau) 
\end{align}
for the resulting relative $2$-Segal spaces, respectively.

\subsection{Families of torsion pairs, their COHAs, and their representations}

We introduce the notion of pair of stacks which parameterize families of torsion pairs. Let $\bfT$ and $\bfF$ be two \textit{open} substacks of $\bfCohps(\scrC,\tau)$. In what follows, the notation introduced in Definition~\ref{def:category_associated_to_open_substack} is in use.

\medskip

We introduce the following assumption. The notation introduced in Definition~\ref{def:category_associated_to_open_substack} is in use.
\begin{assumption}\label{ass:torsion-pair}
	For every field $\kappa$, the category $\PsCpt(\scrC_\kappa,\tau_k) = \catCohps(\scrC_k,\tau_k)$ is abelian and the subcategories $(\PsCpt_\bfT(\scrC_\kappa,\tau_\kappa), \PsCpt_\bfF(\scrC_\kappa,\tau_\kappa))$ form a torsion pair on the $\catCohps(\scrC_k, \tau_k)$.
\end{assumption}

We think about $\bfT$ and $\bfF$ as parametrizing together \textit{families of torsion pairs} in $\scrC$.

\begin{notation}
	We will rather denote the relative $2$-Segal spaces
	\begin{align}
		\calS_\bullet^\ell \bfFlagPerf^{(1),\dagger}_{\bfT, \bfF}(\scrC) \longrightarrow \calS_\bullet \bfPerf_{\bfT}(\scrC) \quad\text{and}\quad
		\calS_\bullet^r \bfFlagPerf^{(1),\dagger}_{\bfF, \bfT}(\scrC) \longrightarrow \calS_\bullet \bfPerf_{\bfT}(\scrC) 
	\end{align}
	respectively by the notation
	\begin{align}
		\calS_\bullet^\ell\bfFlagCoh^{(1),\dagger}_{\bfT,\bfF}(\scrC,\tau) \longrightarrow \calS_\bullet \bfCoh_{\bfT}(\scrC,\tau) \quad\text{and}\quad
		\calS_\bullet^r\bfFlagCoh^{(1),\dagger}_{\bfF,\bfT}(\scrC,\tau) \longrightarrow \calS_\bullet \bfCoh_{\bfT}(\scrC,\tau) \ .
	\end{align}
\end{notation}

Lemma~\ref{lem:torsion_pairs_extensions_mono_epi} and the arguments of the proof of Lemma~\ref{lem:Serre_subcategory} immediately imply the following result.
\begin{lemma} \label{lem:torsion_pairs_key_squares}
	The squares
	\begin{align}
		\begin{tikzcd}[column sep = small, ampersand replacement=\&]
			\calS_2 \bfCoh_\bfT(\scrC,\tau) \arrow{r} \arrow{d} \& \calS_2\bfCohps(\scrC,\tau) \arrow{d}{\partial_0 \times \partial_2} \\
			\bfT \times \bfT \arrow{r} \& \bfCohps(\scrC,\tau) \times \bfCohps(\scrC,\tau)
		\end{tikzcd} \ \text{and} \ \begin{tikzcd}[column sep = small, ampersand replacement=\&]
			\calS_2 \bfCoh_\bfF(\scrC,\tau) \arrow{r} \arrow{d} \& \calS_2 \bfCohps(\scrC,\tau) \arrow{d}{\partial_0 \times \partial_2} \\
			\bfF \times \bfF \arrow{r} \& \bfCohps(\scrC,\tau) \times \bfCohps(\scrC,\tau)
		\end{tikzcd} 
	\end{align}	
	are pullback.
	
	Similarly, the squares
	\begin{align}
		\begin{tikzcd}[column sep = small, ampersand replacement=\&]
			\calS_1^r\bfFlagCoh^{(1)}_{\bfF,\bfT}(\scrC,\tau) \arrow{r} \arrow{d} \& \calS_1^r\bfFlagCohps^{(1)}(\scrC,\tau) \arrow{d}{\varpi_1\times u^r_1} \\
			\bfF \times \bfT \arrow{r} \& \bfFlagCohps^{(1)}(\scrC,\tau) \times \bfCohps(\scrC,\tau)
		\end{tikzcd}
	\end{align}
	and
	\begin{align}\label{eq:square-tor-tf}
		\begin{tikzcd}[column sep = small, ampersand replacement=\&]
			\calS_1^\ell\bfFlagCoh^{(1)}_{\bfT, \bfF}(\scrC,\tau) \arrow{r} \arrow{d} \& \calS_1^\ell\bfFlagCohps^{(1)}(\scrC,\tau) \arrow{d}{u_1^\ell\times \varpi_0} \\
			\bfT \times \bfF \arrow{r} \& \bfCohps(\scrC,\tau) \times \bfFlagCohps^{(1)}(\scrC,\tau) 
		\end{tikzcd}
	\end{align}
	are pullback.
\end{lemma}

\begin{rem}\label{rem:identification}
	Unraveling the definitions, we see that all the three derived stacks $\calS_2 \bfCohps(\scrC,\tau)$, $\calS_1^r\bfFlagCohps^{(1)}(\scrC,\tau)$ and $\calS_1^\ell\bfFlagCohps^{(1)}(\scrC,\tau)$ can be identified with the derived stack $\bfCohpsext(\scrC,\tau)$ parametrizing \textit{extensions} of $\tau$-flat pseudo-perfect objects in $\scrC$. However, we use the notation $\calS_2 \bfCohps(\scrC,\tau)$ to indicate the algebra structure in correspondences of $\bfCohps(\scrC,\tau)$, while we write $\calS_1^r\bfFlagCohps^{(1)}(\scrC,\tau)$ to denote the right action in correspondences of $\bfCohps(\scrC,\tau)$ on itself. Similarly, $\calS_1^\ell\bfFlagCohps^{(1)}(\scrC,\tau)$ corresponds to the left action of $\bfCohps(\scrC,\tau)$ on itself. \hfill $\triangle$
\end{rem}

Combining Proposition~\ref{prop:2-Segal-representation} and Lemma~\ref{lem:torsion_pairs_key_squares} we deduce that the morphisms
\begin{align}
	\calS^r_\bullet \bfFlagCoh^{(1)}_{\bfF,\bfT}(\scrC,\tau) \longrightarrow \calS_\bullet \bfCoh_\bfF(\scrC,\tau) \quad \text{and} \quad \calS^\ell_\bullet \bfFlagCoh^{(1)}_{\bfT,\bfF}(\scrC,\tau) \longrightarrow \calS_\bullet \bfCoh_\bfT(\scrC,\tau) 
\end{align}
are relative $2$-Segal spaces. 
\begin{proposition} \label{prop:torsion_pairs_left_right_actions}
	Fix a base ring $k$ and a $k$-linear compactly generated presentable stable $\infty$-category $\scrC\in \PrLomega_k$ equipped with a $t$-structure $\tau$ such that Assumption~\ref{assumption:existence-COHA} holds.
	We also fix a motivic formalism $\bfD^\ast$, coefficients $\calA \in \bfD^\ast(\Spec(k))$ and $\Gamma \subseteq \Pic(\bfD^\ast(\Spec(k)))$ such that Assumption~\ref{assumption:motivic_formalism} holds.
	
	Let $\bfT$ and $\bfF$ be two open substacks of $\bfCohps(\scrC,\tau)$ such that Assumption~\ref{ass:torsion-pair} holds. Then:
	\begin{enumerate}[font=\upshape, label=(Alg-T)]
		\item \label{item:algebra-torsion} If the morphism
		\begin{align}
			\partial_1 \colon \calS_2 \bfCoh_\bfT(\scrC,\tau) \longrightarrow \bfT
		\end{align}
		is locally rpas, then $\catCohb_{\mathsf{pro}}(\bfT)$ has an induced $\E_1$-monoidal structure. Similarly, the topological vector space $\HBMDGamma_0(\bfF;\calA)$ has the structure of a unital associative algebra. In particular,
		\begin{align}
			G_0( \bfT )\quad \text{and} \quad \HBM_\ast( \bfT )
		\end{align}
		have the structures of unital associative algebras.
	\end{enumerate}
	
	\begin{enumerate}[font=\upshape, label=(Alg-F)]
		\item \label{item:algebra-torsion-free} If the morphism
		\begin{align}
			\partial_1 \colon \calS_2 \bfCoh_\bfF(\scrC,\tau) \longrightarrow \bfF
		\end{align}
		is locally rpas, then $\catCohb_{\mathsf{pro}}(\bfF)$ has an induced $\E_1$-monoidal structure. Similarly, the topological vector space $\HBMDGamma_0(\bfF;\calA)$ has the structure of a unital associative algebra. In particular,
		\begin{align}
			G_0( \bfF )\quad \text{and} \quad \HBM_\ast( \bfF )
		\end{align}
		have the structures of unital associative algebras.
	\end{enumerate}
	
	\begin{enumerate}[font=\upshape, label=(Rep-F)]
		\item \label{item:rep-torsion-free} In addition to the assumption in \ref{item:algebra-torsion}, assume that the morphism
		\begin{align}
			\varpi_0 \colon \calS_1^\ell\bfFlagCoh^{(1),\dagger}_{\bfT,\bfF}(\scrC,\tau) \longrightarrow \bfF
		\end{align}
		is locally rpas. Then $\catCohb_{\mathsf{pro}}(\bfF)$ has an induced structure of a left categorical module over $\catCohb_{\mathsf{pro}}(\bfT)$. Similarly, the topological vector space $\HBMDGamma_0(\bfF;\calA)$ has the structure of a unital associative algebra. In particular,
		\begin{align}
			G_0( \bfF )\quad \text{and} \quad \HBM_\ast( \bfF )			
		\end{align}
		are left modules of $G_0( \bfT )$ and $\HBM_\ast( \bfT )$, respectively.
	\end{enumerate}
	
	\begin{enumerate}[font=\upshape, label=(Rep-T)]	
		\item \label{item:rep-torsion} In addition to the assumption in \ref{item:algebra-torsion-free}, assume that the morphism
		\begin{align}
			\varpi_1 \colon \calS_1^r\bfFlagCoh^{(1),\dagger}_{\bfF,\bfT}(\scrC,\tau) \longrightarrow \bfT
		\end{align}
		is locally rpas. Then $\catCohb_{\mathsf{pro}}(\bfT)$ has an induced structure of a right categorical module over $\catCohb_{\mathsf{pro}}(\bfF)$. Similarly, the topological vector space $\HBMDGamma_0(\bfF;\calA)$ has the structure of a unital associative algebra. In particular,
		\begin{align}
			G_0( \bfT )\quad \text{and} \quad \HBM_\ast( \bfT )			
		\end{align}
		are right modules of $G_0( \bfF )$ and $\HBM_\ast( \bfF )$, respectively.
	\end{enumerate}
\end{proposition}

\begin{proof}
	We first show that the theses of \ref{item:algebra-torsion} and \ref{item:rep-torsion-free} hold. For this, it is enough to show that, setting
	\begin{align}
		\bfH\coloneqq \bfT \quad \text{and} \quad \bfM\coloneqq \bfF\ ,
	\end{align}
	the assumptions of Corollaries~\ref{cor:COHA} and \ref{cor:COHA-representations} are satisfied.
	
	First, note that the square \eqref{eq:pullback-square-T} for $\bfCohps(\scrC, \tau)$ is a pullback by Lemma~\ref{lem:2-Segal-Coh}. Lemma~\ref{lem:torsion_pairs_key_squares} yields that the square \eqref{eq:pullback-square-T} for $\bfT$ is a pullback as well. Similarly, by Lemma~\ref{lem:derived-lci}, Condition~\ref{item:algebra-(i)} is satisfied by $\bfCohps(\scrC, \tau)$. By Lemma~\ref{lem:torsion_pairs_key_squares}, this condition is also satisfied by $\bfT$.
	
	Finally, by assumption, the map $\partial_1 \colon \calS_2 \bfCoh_\bfT(\scrC,\tau) \longrightarrow \bfT$ is locally rpas. Thus, we can apply Corollary~\ref{cor:COHA}. Therefore, the thesis of \ref{item:algebra-torsion} holds. 
	
	First, note that Condition~\eqref{item:2-Segal-representation-(2)} of Proposition~\ref{prop:2-Segal-representation} is satisfied by $\bfT, \bfF$ thanks to Lemma~\ref{lem:torsion_pairs_key_squares}. Now, note that the map
	\begin{align}
		u_1^\ell \times \varpi_1 \colon \calS_1^\ell\bfFlagCoh^{(1),\dagger}_{\bfT,\bfF}(\scrC,\tau) \longrightarrow \bfT\times \bfF
	\end{align}
	is quasi-compact, finitely connected and derived lci. Indeed, under the identification
	\begin{align}
		\calS_1^\ell\bfFlagCohps^{(1)}(\scrC,\tau) \simeq \calS_2 \bfCohps(\scrC,\tau) \ , 
	\end{align}
	the map
	\begin{align}
		u_1^\ell\times \varpi_1  \colon \calS_1^\ell\bfFlagCohps^{(1)}(\scrC,\tau) \longrightarrow \bfCohps(\scrC,\tau) \times \bfCohps(\scrC,\tau) 
	\end{align}
	corresponds to the map
	\begin{align}
		\partial_0\times \partial_2 \colon \calS_2\bfCohps(\scrC, \tau) \longrightarrow \bfCohps(\scrC, \tau)\times \bfCohps(\scrC, \tau)\ ,
	\end{align}
	where $\partial_0$ and $\partial_2$ are defined in Formula~\eqref{eq:partial}, which is guaranteed to be quasi-compact, finitely connected and derived lci by Proposition~\ref{prop:stack_extensions_is_relatively_affine}, Lemma~\ref{lem:derived-lci}, and Assumption~\ref{assumption:existence-COHA}. We therefore obtain a commutative square
	\begin{align} \label{eq:torsion_pairs_left_right_actions:lci}
		\begin{tikzcd}[ampersand replacement=\&]
			\calS_1^\ell\bfFlagCoh^{(1)}_{\bfT,\bfF}(\scrC,\tau) \arrow{r} \arrow{d} \& \calS_1\bfFlagCohps^{(1)}(\scrC,\tau) \arrow{d} \\
			\bfT\times \bfF \arrow{r} \& \bfCohps(\scrC,\tau) \times \bfCohps(\scrC,\tau) 
		\end{tikzcd} \ ,
	\end{align}
	where both horizontal arrows are open immersions and the right vertical map is quasi-compact, finitely connected and derived lci. It follows that the left vertical map is quasi-compact, finitely connected and derived lci as well. 
	
	Since, by assumption, the map
	\begin{align}
		\varpi_0 \colon \calS_1^\ell\bfFlagCoh^{(1),\dagger}_{\bfT,\bfF}(\scrC,\tau) \longrightarrow \bfM
	\end{align}
	is locally rpas, all the conditions in Corollary~\ref{cor:COHA-representations} are satisfied. Therefore, the thesis of \ref{item:rep-torsion-free} hold.
	
	Now, the theses of \ref{item:algebra-torsion-free} and \ref{item:rep-torsion} hold similarly, by using a Corollaries~\ref{cor:COHA} and \ref{cor:COHA-representations-right} for $\bfM\coloneqq \bfT$ and $\bfH\coloneqq \bfF$.
\end{proof}

\begin{remark}\label{rem:Serre}
	Note that if $\bfT$ (resp.\ $\bfF$) is such that it defines Serre subcategories of $\scrC_\kappa$ for every field $\kappa$ and the map $\partial_1 \colon \calS_2 \bfCohps(\scrC,\tau) \longrightarrow \bfCohps(\scrC,\tau)$ is locally rpas, then Corollary~\ref{cor:induced_COHA-Serre} yields the thesis of \ref{item:algebra-torsion} (resp.\ \ref{item:algebra-torsion-free}).
\end{remark}

\begin{remark}[$\Lambda$-gradings]
	Fix an abelian monoid $(\Lambda,+)$ and assume that $\calS_\bullet \bfCohps(\scrC,\tau)$ underlies a $\Lambda$-graded $2$-Segal derived stack.
	Then, using Corollary~\ref{cor:induced_Lambda_grading} and Theorem~\ref{thm:BM_Lambda_graded_functoriality}, one can prove that the algebra and module structures of Proposition~\ref{prop:torsion_pairs_left_right_actions} have compatible induced $\Lambda$-gradings.
\end{remark}

\subsection{Torsion COHAs and their torsion-free representations}\label{subsec:torsion-torsion-free}

In this section, we fix a torsion pair $\upsilon = (\scrT, \scrF)$ in $\scrC^\heartsuit$ which we assume to be open in the sense of Definition~\ref{def:open_torsion_pair}.
Recall from Construction~\ref{construction:torsion_pair_substacks} the two substacks $\bfCoh_\scrT(\scrC,\tau)$ and $\bfCoh_{\scrF}(\scrC,\tau)$ of $\bfCohps(\scrC,\tau)$.

\begin{notation}
	The maps $\bfCoh_\scrT(\scrC,\tau) \to \bfCohps(\scrC,\tau)$ and $\bfCoh_{\scrF}(\scrC,\tau) \to \bfCohps(\scrC,\tau)$ define $2$-Segal substacks of $\calS_\bullet \bfCohps(\scrC,\tau)$, $\calS_\bullet^\ell\bfFlagCoh^{(1)}_{\mathsf{ps}}(\scrC,\tau)$ and $\calS_\bullet^r\bfFlagCoh^{(1)}_{\mathsf{ps}}(\scrC,\tau)$ via Propositions~\ref{prop:2-Segal-algebra} and \ref{prop:2-Segal-representation}, and Lemma~\ref{lem:torsion_pairs_key_squares}. To keep the notation under control, we simply denote them by 
	\begin{align}
		\calS_\bullet \bfCoh_{\scrT} (\scrC,\tau) \ , \qquad \calS_\bullet^\ell\bfFlagCoh_{\scrT,\scrF}^{(1)}(\scrC,\tau) \ , \qquad \calS_\bullet^r\bfFlagCoh_{\scrF,\scrT}^{(1)}(\scrC,\tau) \ . 
	\end{align}
	Similarly, we use the notation
	\begin{align}
		\calS_\bullet \bfCoh_{\scrT} (\scrC,\tau_\upsilon) \ , \qquad \calS_\bullet^\ell\bfFlagCoh_{\scrT,\scrF[1]}^{(1)}(\scrC,\tau_\upsilon) \ , \qquad \calS_\bullet^r\bfFlagCoh_{\scrF[1],\scrT}^{(1)}(\scrC,\tau_\upsilon)
	\end{align}
	for the stacks obtained from the tilted $t$-structure $\tau_\upsilon$, using the maps
	$\bfCoh_\scrT(\scrC,\tau) \simeq \bfCoh_\scrT(\scrC,\tau_\upsilon) \to \bfCohps(\scrC,\tau_\upsilon)$ and $\bfCoh_{\scrF}(\scrC,\tau) \simeq \bfCoh_{\scrF[1]}(\scrC,\tau_\upsilon)\to \bfCohps(\scrC,\tau_\upsilon)$. 
\end{notation}

\begin{lemma}\label{lem:tilted-Assumption-C}
	Fix a base ring $k$ and a $k$-linear compactly generated stable $\infty$-category $\scrC\in \PrLomega_k$ equipped with a $t$-structure $\tau$ for which Assumptions~\ref{assumption:existence-COHA-1}, \ref{assumption:existence-COHA-2}, and \ref{assumption:existence-COHA-3} hold.
	
	Let $\upsilon = (\scrT, \scrF)$ be a torsion pair on $\scrC^\heartsuit$ such that
	\begin{enumerate}[label=(\roman*)]\itemsep0.2cm
		\item \label{item:tilted-Assumption-C-1} both $\scrT$ and $\scrC^\heartsuit$ are compactly generated,
		\item \label{item:tilted-Assumption-C-2} the inclusion $\scrT \hookrightarrow \scrC^\heartsuit$ preserves compact objects,
		\item \label{item:tilted-Assumption-C-3} $\upsilon = (\scrT, \scrF)$ is open in the sense of Definition~\ref{def:open_torsion_pair}.
	\end{enumerate}
	Then, the pair $(\scrC, \tau_\upsilon)$ satisfies Assumptions~\ref{assumption:existence-COHA-1}, \ref{assumption:existence-COHA-2}, and \ref{assumption:existence-COHA-3} as well.
\end{lemma}

\begin{proof}
	First, Assumption~\ref{assumption:existence-COHA-1} only depends on $\scrC$ (and not on the $t$-structure), so it is satisfied by $(\scrC,\tau_\upsilon)$. 
	
	Verifying Assumption~\ref{assumption:existence-COHA-2} means to check that $(\scrC,\tau_\upsilon)$ satisfies Assumption~\ref{assumption:t_structure_filtered_colimits}. To see this, first observe that our hypotheses~\ref{item:tilted-Assumption-C-1} and \ref{item:tilted-Assumption-C-2} together with the axioms of torsion pair imply that $\scrF \subset \scrC^\heartsuit$ is closed under filtered colimits. In turn, this implies that $(\scrC,\tau_\upsilon)$ is compatible with filtered colimits (see Remark~\ref{rem:assumption_A_for_tiltings}). Now observe that since $(\scrC,\tau_\upsilon)$ is by construction bounded with respect to $(\scrC,\tau)$, left completeness of the latter implies left completeness of the former. In other words, if Assumption~\ref{item-1-assumption:t_structure_filtered_colimits} holds for $(\scrC,\tau)$, then it holds for $(\scrC,\tau_\upsilon)$. On the other hand, Assumption~\ref{item-2-assumption:t_structure_filtered_colimits} only depends on $\scrC$ and not on the underlying $t$-structure. Thus, in conclusion, we see that $(\scrC,\tau_\upsilon)$ satisfies Assumption~\ref{assumption:t_structure_filtered_colimits} (i.e.\ Assumption~\ref{assumption:existence-COHA-2}).
	
	Finally our hypothesis~\ref{item:tilted-Assumption-C-3} allows to invoke Proposition~\ref{prop:openness-tilted}, which shows that Assumption~\ref{assumption:existence-COHA-3} is satisfied.
\end{proof}

\begin{theorem}\label{thm:left-right-action}
	Fix a base ring $k$ and a $k$-linear compactly generated stable $\infty$-category $\scrC\in \PrLomega_k$ equipped with a $t$-structure $\tau$ for which Assumption~\ref{assumption:existence-COHA} holds.
	
	Let $\upsilon = (\scrT, \scrF)$ be a torsion pair on $\scrC^\heartsuit$ such that
	\begin{enumerate}\itemsep0.2cm
		\item \label{item:torsion-pair-left-right-action-1} both $\scrT$ and $\scrC^\heartsuit$ are compactly generated,
		\item \label{item:torsion-pair-left-right-action-2} the inclusion $\scrT \hookrightarrow \scrC^\heartsuit$ preserves compact objects,
		\item \label{item:torsion-pair-left-right-action-3} $\sfS_\scrC^![2]$ is $t$-exact with respect to the tilted $t$-structure $\tau_\upsilon$,
		\item \label{item:torsion-pair-left-right-action-4} $\upsilon = (\scrT, \scrF)$ is open in the sense of Definition~\ref{def:open_torsion_pair}.
	\end{enumerate}
	In addition, assume that
	\begin{enumerate}[start=5]\itemsep=0.2cm
		
		\item  \label{item:left-right-action-1} the map
		\begin{align}
			\partial_1 \colon \calS_2\bfCoh_{\scrT}(\scrC,\tau) \longrightarrow \bfCoh_{\scrT}(\scrC,\tau)
		\end{align}
		is locally rpas. Equivalently, the map
		\begin{align}
			\partial_1 \colon \calS_2\bfCoh_{\scrT}(\scrC,\tau_\upsilon) \longrightarrow \bfCoh_{\scrT}(\scrC,\tau_\upsilon)
		\end{align}
		is locally rpas.
		
		\item \label{item:left-right-action-2} both maps
		\begin{align}
			\varpi_0 \colon \calS_1^\ell\bfFlagCoh^{(1),\dagger}_{\scrT,\scrF}(\scrC,\tau) \longrightarrow \bfCoh_{\scrF}(\scrC,\tau)
		\end{align}
		and
		\begin{align}
			\varpi_1 \colon \calS_1^r\bfFlagCoh^{(1),\dagger}_{\scrF[1],\scrT}(\scrC,\tau_\upsilon) \longrightarrow \bfCoh_{\scrF[1]}(\scrC,\tau)
		\end{align}
		are locally rpas.
	\end{enumerate}
	Then, 
	\begin{itemize}\itemsep0.2cm
		\item $\catCohb_{\mathsf{pro}}(\bfCoh_{\scrT}(\scrC,\tau))$ inherits the structure of a $\E_1$-monoidal pro-$\infty$-category, and
		
		\item $\catCohb_{\mathsf{pro}}(\bfCoh_{\scrF}(\scrC,\tau))$ has both the structure of a categorical left and of a categorical right module over $\catCohb_{\mathsf{pro}}(\bfCoh_{\scrT}(\scrC,\tau))$.
	\end{itemize}
	
	Similarly, let $\bfD^\ast$ be a motivic formalism, and fix $\calA \in \CAlg(\bfD^\ast(\Spec(k)))$ and $\Gamma \subseteq \Pic(\bfD^\ast(\Spec(k)))$ such that Assumption~\ref{assumption:motivic_formalism} is satisfied.
	Then,
	\begin{itemize}\itemsep0.2cm
		\item the topological vector space $\HBMDGamma_0(\bfCoh_{\scrT}(\scrC,\tau);\calA)$ has the structure of a unital associative algebra. In particular,
		\begin{align}
			G_0( \bfCoh_{\scrT}(\scrC,\tau) )\quad \text{and} \quad \HBM_\ast( \bfCoh_{\scrT}(\scrC,\tau) )
		\end{align}
		have the structures of unital associative algebras, and
		
		\item the topological vector space $\HBMDGamma_0(\bfCoh_{\scrF}(\scrC,\tau);\calA)$ has both the structure of a left and a right $\HBMDGamma_0(\bfCoh_{\scrT}(\scrC,\tau);\calA)$-module. In particular,
		\begin{align}
			G_0( \bfCoh_{\scrF}(\scrC,\tau) )\quad \text{and} \quad \HBM_\ast( \bfCoh_{\scrF}(\scrC,\tau) )			
		\end{align}
		have both the structures of a left and a right $G_0( \bfCoh_{\scrT}(\scrC,\tau) )$-module and $\HBM_\ast( \bfCoh_{\scrT}(\scrC,\tau) )$-module, respectively.
	\end{itemize}
\end{theorem}

\begin{proof}
	Let us begin by showing that both $(\scrC, \tau)$ and $(\scrC, \tau_\upsilon)$ satisfy Assumption~\ref{assumption:existence-COHA}. For $(\scrC,\tau)$ this is our hypothesis. For $(\scrC,\tau_\upsilon)$, notice that Assumptions~\ref{assumption:existence-COHA-1}, \ref{assumption:existence-COHA-2}, and \ref{assumption:existence-COHA-3} hold by Lemma~\ref{lem:tilted-Assumption-C}, while Assumption~\ref{assumption:existence-COHA-4} for $(\scrC,\tau_\upsilon)$ coincides with our hypothesis~\eqref{item:torsion-pair-left-right-action-3}. Thus, $(\scrC,\tau_\upsilon)$ satisfies Assumption~\ref{assumption:existence-COHA}. In particular, the derived stacks $\bfCoh_{\scrT}(\scrC,\tau)$ and $\bfCoh_{\scrT}(\scrC,\tau_\upsilon)$ are geometric and locally of finite presentation over $k$ by Proposition~\ref{prop:openness}.
	
	\medskip
	
	Observe now that the morphisms $\bfCoh_\scrT(\scrC,\tau) \to \bfCohps(\scrC,\tau)$ and $\bfCoh_\scrT(\scrC,\tau) \to\bfCohps(\scrC,\tau_\upsilon)$ extend to morphisms of $2$-Segal spaces thanks to Lemmas~\ref{lem:2-Segal-Coh} and \ref{lem:torsion_pairs_key_squares}, fitting in the following pullback square:
	\begin{align}
		\begin{tikzcd}[ampersand replacement=\&]
			\calS_\bullet \bfCoh_\scrT(\scrC,\tau) \arrow{r} \arrow{d} \& \calS_\bullet \bfCohps(\scrC,\tau) \arrow{d} \\
			\calS_\bullet \bfCohps(\scrC,\tau_\upsilon) \arrow{r} \& \calS_\bullet \bfPerfps(\scrC) 
		\end{tikzcd}\ .
	\end{align}
	In other words, we obtain a canonical equivalence of $2$-Segal spaces
	\begin{align}
		\calS_\bullet \bfCoh_\scrT(\scrC,\tau) \simeq \calS_\bullet \bfCoh_\scrT(\scrC,\tau_\upsilon) \ . 
	\end{align}
	
	Now, Assumption~\ref{ass:torsion-pair} is satisfied for both the stacks $\bfCoh_{\scrT}(\scrC,\tau)$ and $\bfCoh_{\scrT}(\scrF,\tau)$ because of hypothesis~\eqref{item:torsion-pair-left-right-action-3}, while it is satisfied by $\bfCoh_{\scrF}(\scrF,\tau)\simeq \bfCoh_{\scrF[1]}(\scrF,\tau_\upsilon)$ and $\bfCoh_{\scrT}(\scrC,\tau)\simeq \bfCoh_{\scrT}(\scrC,\tau_\upsilon)$ since the torsion pair $(\scrF[1], \scrT)$ is open by Proposition~\ref{prop:openness-tilted} and Remark~\ref{rem:torsion-pair-openness}. 
	
	Furthermore, Assumption~\eqref{item:left-right-action-1} implies that conditions~\ref{item:algebra-torsion} and \ref{item:algebra-torsion-free} of Proposition~\ref{prop:torsion_pairs_left_right_actions} are satisfied. Therefore, the existence of the algebra structure follows. Finally, the existence of the left and right actions follows from Proposition~\ref{prop:torsion_pairs_left_right_actions}, since Assumption~\eqref{item:left-right-action-2} implies that the conditions~\ref{item:rep-torsion-free} and \ref{item:rep-torsion} of \textit{loc.\ cit.} hold.
\end{proof}

Now we introduce generalizations in the present context of the notions of left and right Hecke patterns introduced in \cite[\S5.1]{KV_Hall} (see also \cite[Definition~6.1 and Remark~6.2]{MMSV}).
\begin{definition}\label{def:Hecke-patterns}
	Let $\bfT$ and $\bfF$ be quasi-separated geometric derived stacks locally of finite presentation over $k$, together with maps  $\bfT\to \bfCoh_{\scrT}(\scrC, \tau)$ and $\bfF\to \bfCoh_\scrF(\scrC, \tau)$. We say that
	\begin{enumerate}\itemsep0.2cm
		\item $\bfF$ is a \textit{left Hecke pattern for $\bfT$ with respect to the $t$-structure $\tau$} if the square
		\begin{align}\label{eq:left-Hecke}
			\begin{tikzcd}[column sep = small, ampersand replacement=\&]
				\calS_1^\ell\bfFlagCoh^{(1)}_{\bfT,\bfF}(\scrC,\tau) \arrow{r} \arrow{d} \& \calS_1^\ell\bfFlagCohps^{(1)}(\scrC,\tau) \arrow{d}{u_1^\ell\times \varpi_0} \\
				\bfT \times \bfF \arrow{r} \&\bfCohps(\scrC,\tau) \times \bfFlagCohps^{(1)}(\scrC,\tau)
			\end{tikzcd}
		\end{align}
		is a pullback; 
		
		\item $\bfF$ is a \textit{right Hecke pattern for $\bfT$ with respect to the $t$-structure $\tau$} if the square
		\begin{align}
			\begin{tikzcd}[column sep = small, ampersand replacement=\&]
				\calS_1^\ell\bfFlagCoh^{(1)}_{\bfT,\bfF}(\scrC,\tau) \arrow{r} \arrow{d} \& \calS_1^\ell\bfFlagCohps^{(1)}(\scrC,\tau) \arrow{d}{u_1^\ell\times \varpi_1} \\
				\bfT \times \bfF \arrow{r} \&\bfCohps(\scrC,\tau)\times \bfFlagCohps^{(1)}(\scrC,\tau)
			\end{tikzcd}
		\end{align}
		is pullback, and 
		
		\item $\bfF$ is a \textit{right Hecke pattern for $\bfT$ with respect to the $t$-structure $\tau_\upsilon$} if the square
		\begin{align}\label{eq:right-Hecke}
			\begin{tikzcd}[column sep = small, ampersand replacement=\&]
				\calS_1^r\bfFlagCoh^{(1)}_{\bfF, \bfT}(\scrC,\tau_\upsilon) \arrow{r} \arrow{d} \& \calS_1^r\bfFlagCohps^{(1)}(\scrC,\tau_\upsilon) \arrow{d}{\varpi_1\times u^r_1} \\
				\bfF \times \bfT \arrow{r} \& \bfFlagCohps^{(1)}(\scrC,\tau_\upsilon)\times \bfCohps(\scrC,\tau_\upsilon) 
			\end{tikzcd}
		\end{align}
		is a pullback.
	\end{enumerate}
	
	Furthermore, we say that $\bfF$ is a \textit{two-sided Hecke pattern for $\bfT$ with respect to the $t$-structure $\tau$} if it is a left and a right Hecke pattern for $\bfT$ with respect to the $t$-structure $\tau$.
\end{definition}

\begin{remark}
	It is easy to show that if $\bfF$ is a right Hecke pattern for $\bfT$ with respect to the $t$-structure $\tau$, then it is a right Hecke pattern for $\bfT$ with respect to the $t$-structure $\tau_\upsilon$. On the other hand, the vice versa is not true in general, for example, when $\bfT=\bfCoh_{\scrT}(\scrC, \tau)$ and $\bfF=\bfCoh_{\scrF}(\scrC, \tau)$.
\end{remark}

\begin{remark}
	Let $\bfF$ be a two-sided Hecke pattern for $\bfT$ with respect to $\tau$. Then, Corollary~\ref{cor:induced_COHA-derived-lci} holds also without the assumption that $\mu\colon \bfF=\bfM'\to \bfM=\bfCoh_{\scrF}(\scrC, \tau)$ is derived lci.
\end{remark}

\begin{corollary}\label{cor:induced_COHA_open}
	Keep the assumptions of Theorem~\ref{thm:left-right-action}. Let $\bfT$ and $\bfF$ quasi-separated geometric derived stacks locally of finite presentation over $k$, together with maps $\alpha\colon\bfT\to \bfCoh_{\scrT}(\scrC, \tau)$ and $\mu\colon \bfF\to \bfCoh_\scrF(\scrC, \tau)$ of derived stacks. Assume that
	\begin{enumerate}\itemsep0.2cm
		\item the square
		\begin{align}\label{eq:multiplication-T}
			\begin{tikzcd}[column sep = small, ampersand replacement=\&]
				\calS_2 \bfCoh_\bfT(\scrC,\tau) \arrow{r} \arrow{d} \& \calS_2\bfCohps(\scrC,\tau) \arrow{d}{\partial_0 \times \partial_2} \\
				\bfT \times \bfT \arrow{r} \& \bfCohps(\scrC,\tau) \times \bfCohps(\scrC,\tau)
			\end{tikzcd} 
		\end{align}
		is a pullback; 
		
		\item the map
		\begin{align}
			\partial_0 \times \partial_2\colon \calS_2 \bfCoh_\bfT(\scrC,\tau)\longrightarrow \bfT\times \bfT
		\end{align}
		is quasi-compact, finitely connected and derived lci;
		
		\item the map
		\begin{align}
			\partial_1\colon \calS_2 \bfCoh_\bfT(\scrC,\tau)\longrightarrow \bfT
		\end{align}
		is locally rpas;
		
		\item the map $\alpha \colon \bfT \to \bfCoh_\scrT(\scrC,\tau)$ is locally rpas;
		
		\item both maps
		\begin{align}
			\varpi_0 & \colon \calS_1^\ell\bfFlagCoh^{(1),\dagger}_{\scrT,\scrF}(\scrC,\tau) \longrightarrow \bfCoh_{\scrF}(\scrC,\tau) \ , \\
			\varpi_1 & \colon \calS_1^r\bfFlagCoh^{(1),\dagger}_{\scrF[1],\scrT}(\scrC,\tau_\upsilon) \longrightarrow \bfCoh_{\scrF[1]}(\scrC,\tau)
		\end{align}
		are locally rpas, and
		
		\item \label{item:Hecke} $\bfF$ is a left Hecke pattern for $\bfT$ with respect to $\tau$ and either 
		\begin{enumerate}\itemsep0.2cm
			\item $\mu$ is quasi-compact, finitely connected and derived lci, and $\bfF$ a right Hecke pattern for $\bfT$ with respect to $\tau_\upsilon$, or
			
			\item \label{item:two-sided} $\bfF$ a right Hecke pattern for $\bfT$ with respect to $\tau$ (i.e., a two-sided Hecke pattern for $\bfT$ with respect to $\tau$),		\end{enumerate}
	\end{enumerate}
	Then 
	\begin{itemize}\itemsep0.2cm
		\item $\catCohb_{\mathsf{pro}}(\bfT)$ has the structure of a $\E_1$-monoidal pro-$\infty$-category, and
		
		\item$\catCohb_{\mathsf{pro}}( \bfF )$ has both the structure of a categorical left and of a categorical right module over $\catCohb_{\mathsf{pro}}( \bfT )$.
	\end{itemize}
	Similarly, let $\bfD^\ast$ be a motivic formalism, and fix $\calA \in \CAlg(\bfD^\ast(\Spec(k)))$ and $\Gamma \subseteq \Pic(\bfD^\ast(\Spec(k)))$ such that Assumption~\ref{assumption:motivic_formalism} is satisfied. Then,
	\begin{itemize}\itemsep0.2cm
		\item the topological vector space $\HBMDGamma_0(\bfT;\calA)$ has the structure of a unital associative algebra. In particular,
		\begin{align}
			G_0( \bfT )\qquad \text{and} \qquad \HBM_\ast( \bfT )
		\end{align}
		have the structures of unital associative algebras, and
		
		\item the topological vector space $\HBMDGamma_0(\bfF;\calA)$ has both the structure of a left and a right $\HBMDGamma_0(\bfT;\calA)$-module. In particular,
		\begin{align}
			G_0( \bfF )\qquad \text{and} \qquad \HBM_\ast( \bfF )			
		\end{align}
		have both the structures of a left and a right $G_0( \bfT )$-module and $\HBM_\ast( \bfT )$-module, respectively.
	\end{itemize}
\end{corollary}

\begin{proof}
	It suffices to apply Corollaries~\ref{cor:induced_COHA-properness} and \ref{cor:induced_COHA-derived-lci} for $m = 1$ and
	\begin{align}
		\bfH \coloneqq \bfCoh_{\scrT}(\scrC, \tau) \ , \qquad \bfH' \coloneqq \bfT \ , \qquad \bfM \coloneqq \bfCoh_{\scrF}(\scrC,\tau) \ , \qquad  \bfM' \coloneqq \bfF \ ,
	\end{align}	
	to obtain the desired result.
\end{proof}

\begin{remark}
	Let $(\Lambda,+)$ be a monoid. Assume that $\calS_\bullet \bfCohps(\scrC,\tau)$ admits a $\Lambda^\star$-graded $2$-Segal structure. Then applying repeatedly Corollary~\ref{cor:induced_Lambda_grading} and Remark~\ref{rem:Lambda_graded_multiplication}, we deduce that the algebras $G_0(\bfT)$ and $\HBM_\ast(\bfT)$ (and their categorical counterpart $\catCohb_{\mathsf{pro}}(\bfT)$) inherits a canonical $\Lambda$-grading compatible with the Hall multiplication. Similarly, following Remark~\ref{rem:Lambda_graded_action}, we deduce that the left and right modules $G_0(\bfF)$ and $\HBM_\ast(\bfF)$ (and their categorical counterpart $\catCohb_{\mathsf{pro}}(\bfF)$) acquire canonical $\Lambda$-gradings.
\end{remark}

Let $(\bfT, \bfF)$ be two derived stacks satisfying the assumptions of Corollary~\ref{cor:induced_COHA_open}. We give now the main definition of this section:
\begin{defin}\label{def:Yangian-torsion}
	The \textit{categorical quantum loop algebra} $\scrH_{(\bfT, \bfF)}$ of the pair $(\bfT, \bfF)$ is the monoidal subcategory of the monoidal $\infty$-category of endofunctors $\End( \catCohb_{\mathsf{pro}}( \bfF ) )$ generated by the images of the two monoidal functors
	\begin{align}
		a_\ell&\colon \catCohb_{\mathsf{pro}}( \bfT )\longrightarrow\End( \catCohb_{\mathsf{pro}}( \bfF ) )\ ,\\[2pt]
		a_r&\colon \catCohb_{\mathsf{pro}}( \bfT )\longrightarrow\End( \catCohb_{\mathsf{pro}}( \bfF ) )\ ,
	\end{align}
	corresponding to the two module structures of $\catCohb_{\mathsf{pro}}( \bfF )$.
	
	The \textit{quantum loop algebra} $\scrU_{(\bfT, \bfF)}$ of the pair $(\bfT, \bfF)$ is the subalgebra of $\End( G_0( \bfF ) )$ generated by the images of the two maps of associative algebras
	\begin{align}
		a_\ell&\colon G_0( \bfT )\longrightarrow\End( G_0( \bfF ) )\ ,\\[2pt]
		a_r&\colon G_0( \bfT )\longrightarrow\End( G_0( \bfF ) )\ ,
	\end{align}
	corresponding to the two module structures of $G_0( \bfF )$. 
	
	The \textit{Yangian} $\scrY_{(\bfT, \bfF)}$ of the pair $(\bfT, \bfF)$ is the subalgebra of $\End( \HBM_\ast( \bfF ) )$ generated by the images of the two maps of associative algebras
	\begin{align}
		a_\ell&\colon \HBM_\ast( \bfT )\longrightarrow\End( \HBM_\ast( \bfF ) )\ ,\\[2pt]
		a_r&\colon \HBM_\ast( \bfT )\longrightarrow\End( \HBM_\ast( \bfF ) )\ ,
	\end{align}
	corresponding to the two module structures of $\HBM_\ast( \bfF )$. 
	
	Furthermore, if there is an action on $\bfT$ and $\bfF$ by a torus $T$, we introduce the \textit{equivariant} version of the notions above by replacing $\bfT$ and $\bfF$ by the quotient stacks $\bfT/T$ and $\bfF/T$, respectively.\footnote{Compare these definitions in the equivariant setting with \cite[\S4.3]{Porta_Sala_Hall}.}\hfill $\oslash$
\end{defin}

\section{Categorified commutators}\label{sec:categorical_commutators}

In this section, we study the relation between the compositions between the left and the right actions in the two possible orders. In particular, we provide a criterion to establish if the categoried commutator vanishes.

\subsection{Preliminaries}

Let us place ourselves again in the context of Theorem~\ref{thm:left-right-action}. Our goal is to provide a geometric interpretation of the commutator of the left and the right actions. For this, we first look at $\calS_3 \bfPerfps(\scrC)$, and observe that the simplicial identities guarantee that the diagrams
\begin{align}\label{eq:diagram-commutators-left}
	\begin{tikzcd}[ampersand replacement=\&]
		\bfPerfps(\scrC)^{\times 3} \arrow[equal]{d} \&[2em] \calS_2 \bfPerfps(\scrC) \times \bfPerfps(\scrC) \arrow{l}[swap]{(\partial_0 \times \partial_2) \times \id} \&[1em] \calS_3 \bfPerfps(\scrC) \arrow{l}[swap]{\partial_3\times \partial_{0,1}} \arrow[equal]{d} \\
		\bfPerfps(\scrC)^{\times 3} \& \bfPerfps(\scrC) \times \calS_2 \bfPerfps(\scrC) \arrow{l}[swap]{\id \times (\partial_0 \times \partial_2)} \& \calS_3 \bfPerfps(\scrC) \arrow{l}[swap]{\partial_{2,3}\times \partial_0}
	\end{tikzcd}
\end{align}
and	
\begin{align}\label{eq:diagram-commutators-right}
	\begin{tikzcd}[ampersand replacement=\&]
		\calS_3 \bfPerfps(\scrC) \arrow{r}{\partial_1} \arrow[equal]{d} \&[-1em] \calS_2 \bfPerfps(\scrC) \arrow{r}{\partial_1} \&[-1em] \bfPerfps(\scrC) \arrow[equal]{d} \\
		\calS_3 \bfPerfps(\scrC) \arrow{r}{\partial_2} \& \calS_2 \bfPerfps(\scrC) \arrow{r}{\partial_1} \& \bfPerfps(\scrC)
	\end{tikzcd}
\end{align}
commute, where the first horizontal lines in \eqref{eq:diagram-commutators-left} and \eqref{eq:diagram-commutators-right} correspond to the composition ``left action $\circ$ right action'', while the second horizontal lines in \eqref{eq:diagram-commutators-left} and \eqref{eq:diagram-commutators-right} corresponds to the composition ``right action $\circ$ left action''. 

Let us spell out the maps in the diagrams above. First, the stack $\calS_3 \bfPerf(\scrC)$ parametrizes flags of the form
\begin{align}\label{eq:commutator_3_flag}
	\begin{tikzcd}[ampersand replacement = \&]
		0 \arrow{r} \& F_{0,1} \arrow{r} \arrow{d} \& F_{0,2} \arrow{r} \arrow{d} \& F_{0,3} \arrow{d} \\
		\& 0 \arrow{r} \& F_{1,2} \arrow{r} \arrow{d} \& F_{1,3} \arrow{d} \\
		\& \& 0 \arrow{r} \& F_{2,3} \arrow{d} \\
		\& \& \& 0
	\end{tikzcd}\ .
\end{align}
In the first horizontal line of  \eqref{eq:diagram-commutators-left}, the map $\partial_3\colon \calS_3\bfPerfps(\scrC) \to \calS_2\bfPerfps(\scrC)$ sends the flag \eqref{eq:commutator_3_flag} to
\begin{align}
	\begin{tikzcd}[ampersand replacement = \&]
		0 \arrow{r} \& F_{0,1} \arrow{r} \arrow{d} \& F_{0,2} \arrow{d} \\
		\& 0 \arrow{r} \& F_{1,2} \arrow{d} \\
		\& \& 0
	\end{tikzcd}\ ,
\end{align}
while $\partial_{0,1}$ sends the flag \eqref{eq:commutator_3_flag} to $F_{2,3}$. The maps $\partial_0, \partial_2\colon \calS_2 \bfPerfps(\scrC)\to \bfPerfps(\scrC)$ have been introduced in \eqref{eq:partial}. On the other hand, in the first horizontal line of \eqref{eq:diagram-commutators-right}, the map $\partial_1\colon \calS_3 \bfPerfps(\scrC) \to \calS_2 \bfPerfps(\scrC)$ sends the flag \eqref{eq:commutator_3_flag} to 
\begin{align}
	\begin{tikzcd}[ampersand replacement = \&]
		0 \arrow{r} \& F_{0,2} \arrow{d} \arrow{r} \& F_{0,3} \arrow{d} \\
		\& 0 \arrow{r} \& F_{2,3} \arrow{d} \\
		\& \& 0
	\end{tikzcd}
\end{align}
while $\partial_1 \colon \calS_2\bfPerfps(\scrC)\to \bfPerfps(\scrC)$ has been introduced in \eqref{eq:partial} and sends the above flag to $F_{0,3}$.

In the second horizontal line of  \eqref{eq:diagram-commutators-left}, the map $\partial_0\colon \calS_3\bfPerfps(\scrC)\to \calS_2\bfPerfps(\scrC)$ sends the flag \eqref{eq:commutator_3_flag} to 
\begin{align}
	\begin{tikzcd}[ampersand replacement = \&]
		0 \arrow{r} \& F_{1,2} \arrow{r} \arrow{d} \& F_{1,3} \arrow{d} \\
		\& 0 \arrow{r} \& F_{2,3} \arrow{d} \\
		\& \& 0
	\end{tikzcd} 
\end{align}
while $\partial_{2,3}$ sends the flag \eqref{eq:commutator_3_flag} to $F_{0,1}$. In the second horizontal line of  \eqref{eq:diagram-commutators-right}, the map $\partial_2\colon \calS_3 \bfPerfps(\scrC) \to \calS_2\bfPerfps(\scrC)$ sends the flag \eqref{eq:commutator_3_flag} to
\begin{align}
	\begin{tikzcd}[ampersand replacement = \&]
		0 \arrow{r} \& F_{0,1} \arrow{d} \arrow{r} \& F_{0,3} \arrow{d} \\
		\& 0 \arrow{r} \& F_{1,3} \arrow{d} \\
		\& \& 0
	\end{tikzcd} \ ,
\end{align}
while $\partial_1 \colon \calS_2\bfPerfps(\scrC)\to \bfPerfps(\scrC)$ has been introduced in \eqref{eq:partial} and sends the above flag to $F_{0,3}$.

Define now an open substack $\calS_3 \bfFlagCoh^{(1)}_{(\scrT, (\scrF, \scrT))}(\scrC,\tau)$ of $\calS_3 \bfPerf(\scrC)$ parametrizing flags $\F$ of the form \eqref{eq:commutator_3_flag} satisfying
\begin{itemize}\itemsep=0.2cm
	\item $\partial_3(\F) \in \calS_1^r \bfFlagCoh^{(1)}_{\scrF, \scrT[-1]}(\scrC,\tau_{-\upsilon})$ and $\partial_{0,1}(\F) = F_{2,3} \in \bfCoh_{\scrT}(\scrC,\tau)$;
	
	\item $\partial_1(\F) \in \calS_1^\ell \bfFlagCoh^{(1)}_{\scrT, \scrF}(\scrC,\tau)$,
\end{itemize}
where $\tau_{-\upsilon} \coloneqq \tau_{\upsilon}[-1]$ denotes the anti-tilted $t$-structure.\footnote{In this subsection, we prefer to work with the torsion pair $(\scrF, \scrT[-1])$ rather than $(\scrF[1], \scrT)$, since we believe it makes our computations easier to follow.}
\begin{remark}
	Note that for the flags of the form \eqref{eq:commutator_3_flag} parametrized by $\calS_3 \bfFlagCoh^{(1)}_{((\scrT, \scrF), \scrT)}(\scrC,\tau)$ one has that $F_{1,3}$ is a pseudo-perfect object, which is $\tau$-flat, since it fits into the fiber sequence
	\begin{align}
		F_{0,2} \longrightarrow F_{1,3} \longrightarrow F_{2,3} \ .
	\end{align}
\end{remark}
The derived stack $\calS_3 \bfFlagCoh^{(1)}_{(\scrT, (\scrF, \scrT))}(\scrC,\tau)$ fits into the convolution diagram encoding the composition ``left action $\circ$ right action'':
\begin{align}
	\begin{tikzpicture}
		\node (A) at (0,1) {$\calS_3 \bfFlagCoh^{(1)}_{(\scrT, (\scrF, \scrT))}(\scrC,\tau)$};
		\node (B) at (4.5, 1) {$\calS_1^\ell \bfFlagCoh^{(1)}_{\scrT, \scrF}(\scrC,\tau)$};
		\node (C) at (8, 1) {$\bfCoh_{\scrF}(\scrC, \tau)$};
		\node (D) at (0, -0.5) {$\bfCoh_{\scrT}(\scrC, \tau)\times \calS_1^r \bfFlagCoh^{(1)}_{\scrF,\scrT[-1]}(\scrC,\tau_{-\upsilon})$};
		\node (E) at (0, -2) {$\bfCoh_{\scrT}(\scrC, \tau)\times \left(\bfCoh_{\scrF}(\scrC, \tau)\times \bfCoh_{\scrT[-1]}(\scrC, \tau_{-\upsilon}) \right)$};
		\draw (A) edge[->] node[above]{\small $\partial_1$} (B);
		\draw (B) edge[->] node[above]{\small $\partial_1$} (C);
		\draw (A) edge[->] node[left]{\small $\partial_{0,1}\times \partial_3$} (D);
		\draw (D) edge[->] node[left]{\small $\id\times (\partial_0 \times u_1^r)$} (E);
	\end{tikzpicture}\ .
\end{align}
Above, the map $\partial_0\times u_1^r$ corresponds to $\partial_2\times \partial_0$ of Formula~\eqref{eq:partial}.

We set
\begin{align}
	p_{\ell\circ r} \coloneqq (\id \times \partial_0 \times u_1^r) \circ (\partial_{0,1}\times \partial_3) \qquad \text{and} \qquad q^{\ell\circ r} \coloneqq \partial_1 \circ \partial_1 \ . 
\end{align}
Observe that $p_{\ell\circ r}$ is derived lci, while $q^{\ell\circ r}$ is locally rpas.

Define similarly an open substack $\calS_3 \bfFlagCoh^{(1)}_{((\scrT,\scrF), \scrT)}(\scrC,\tau)$ of $\calS_3 \bfPerf(\scrC)$ parametrizing flags $\F$ of the form \eqref{eq:commutator_3_flag} satisfying
\begin{itemize}\itemsep=0.2cm
	\item $\partial_0(\F) \in \calS_1^\ell \bfFlagCoh^{(1)}_{\scrT, \scrF}(\scrC,\tau)$ and $\partial_{2,3}(\F) = F_{0,1} \in \bfCoh_{\scrT[-1]}(\scrC,\tau_{-\upsilon})$;
	
	\item $\partial_2(\F) \in \calS_1^r \bfFlagCoh^{(1)}_{\scrF,\scrT[-1]}(\scrC,\tau_{-\upsilon})$.
\end{itemize}
\begin{remark}
	Note that for the flags of the form \eqref{eq:commutator_3_flag} parametrized by $\calS_3 \bfFlagCoh^{(1)}_{((\scrT,\scrF),\scrT)}(\scrC,\tau)$ one has that $F_{0,2}$ is a pseudo-perfect object, which is $\tau_\upsilon$-flat, since it fits into the fiber sequence
	\begin{align}
		F_{0,1}\longrightarrow F_{0,2}\longrightarrow F_{1,2} \ .
	\end{align}
\end{remark}
$\calS_3 \bfFlagCoh^{(1)}_{((\scrT,\scrF), \scrT)}(\scrC,\tau)$ fits into the convolution diagram encoding the composition ``right action $\circ$ left action'':
\begin{align}
	\begin{tikzpicture}
		\node (A) at (0,1) {$\calS_3 \bfFlagCoh^{(1)}_{((\scrT,\scrF), \scrT)}(\scrC,\tau)$};
		\node (B) at (5, 1) {$\calS_1^r \bfFlagCoh^{(1)}_{\scrF, \scrT[-1]}(\scrC,\tau_{-\upsilon})$};
		\node (C) at (9, 1) {$\bfCoh_{\scrF}(\scrC, \tau)$};
		\node (D) at (0, -0.5) {$\calS_1^\ell \bfFlagCoh^{(1)}_{\scrT, \scrF}(\scrC,\tau)\times\bfCoh_{\scrT[-1]}(\scrC, \tau_{-\upsilon})$};
		\node (E) at (0, -2) {$\left(\bfCoh_{\scrT}(\scrC, \tau)\times \bfCoh_{\scrF}(\scrC, \tau) \right)\times \bfCoh_{\scrT[-1]}(\scrC, \tau_{-\upsilon})$};
		\draw (A) edge[->] node[above]{\small $\partial_2$} (B);
		\draw (B) edge[->] node[above]{\small $\partial_1$} (C);
		\draw (A) edge[->] node[left]{\small $\partial_0\times \partial_{2,3}$} (D);
		\draw (D) edge[->] node[left]{\small $(u_1^\ell \times \partial_1)\times \id$} (E);
	\end{tikzpicture}\ .
\end{align}
Above, the map $u_1^\ell\times \partial_1$ corresponds to $\partial_2\times \partial_0$ of Formula~\eqref{eq:partial}.

We set
\begin{align}
	p_{r\circ \ell} \coloneqq (u_1^\ell \times \partial_1\times \id) \circ (\partial_0\times \partial_{2,3}) \qquad \text{and} \qquad q^{r\circ \ell} \coloneqq \partial_1 \circ \partial_2 \ . 
\end{align}
Observe that $p_{r\circ \ell}$ is derived lci, while $q^{r\circ \ell}$ is locally rpas.

Summarizing the discussion so far, we have diagrams
\begin{align}
	\begin{tikzpicture}
		\node (A) at (0,3.5) {$\calS_3 \bfFlagCoh^{(1)}_{((\scrT, \scrF), \scrT)}(\scrC,\tau) $};
		\node (B) at (0,2) {$\calS_3 \bfFlagCoh^{(1)}_{(\scrT, (\scrF, \scrT))}(\scrC,\tau) $};
		\node (C) at (-4,0) {$\bfCoh_\scrT(\scrC,\tau)\times \bfCoh_\scrF(\scrC,\tau) \times \bfCoh_{\scrT[-1]}(\scrC,\tau_{-\upsilon})$};
		\node(D) at (4,0) {$\bfCoh_\scrF(\scrC,\tau)$};
		\draw (A.west) edge[->, bend right = 15pt] node[xshift = -5pt, yshift = 5pt]{\tiny $p_{r\circ \ell}$} (C.100) ;
		\draw (B.west) edge[->,bend right = 15pt] node [xshift = 5pt,yshift = -5pt]{\tiny $p_{\ell\circ r}$} (C) ;
		\draw (A.east) edge[->, bend left = 25pt] node (E) [xshift = 5pt, yshift = 5pt]{\tiny $q^{r\circ \ell}$} (D.75) ;
		\draw (B.east) edge[->, bend left = 15pt] node[xshift=-3pt,yshift = -7pt]{\tiny $q^{\ell\circ r}$} (D) ;
	\end{tikzpicture}
\end{align}
The maps $p_{r\circ \ell}$ and $q^{r\circ \ell}$ induce after passing to $\catCohb_{\mathsf{pro}}$ the functor
\begin{align}
	\catCohb_{\mathsf{pro}}(\bfCoh_\scrT(\scrC,\tau)) \otimes \catCohb_{\mathsf{pro}}(\bfCoh_\scrF(\scrC,\tau)) \otimes \catCohb_{\mathsf{pro}}(\bfCoh_{\scrT[-1]}(\scrC,\tau_{-\upsilon}))  \longrightarrow \catCohb_{\mathsf{pro}}(\bfCoh_\scrF(\scrC,\tau)) 
\end{align}
that can be written
\begin{align}
	\calG' \otimes \calE \otimes \calG \mapsto q^{r\circ \ell}_\ast( p_{r\circ \ell}^\ast(\calG' \boxtimes \calE \boxtimes \calG) ) \simeq (\calG' \oast \calE) \oast \calG \ . 
\end{align}

On the other hand, the morphism $p_{\ell\circ r}$ and $q^{\ell\circ r}$ induce in the same way the functor
\begin{align}
	\catCohb_{\mathsf{pro}}(\bfCoh_\scrT(\scrC,\tau)) \otimes \catCohb_{\mathsf{pro}}(\bfCoh_\scrF(\scrC,\tau)) \otimes \catCohb_{\mathsf{pro}}(\bfCoh_{\scrT[-1]}(\scrC,\tau_{-\upsilon}))  \longrightarrow \catCohb_{\mathsf{pro}}(\bfCoh_\scrF(\scrC,\tau)) 
\end{align}
that can be written
\begin{align}
	\calG' \otimes \calE \otimes \calG \mapsto q^{\ell\circ r}_\ast( p_{\ell\circ r}^\ast(\calG' \boxtimes \calE \boxtimes \calG) ) \simeq \calG' \oast (\calE \oast \calG) \ . 
\end{align}

\subsection{Vanishing of the categorical commutator}

Now, we discuss when the two compositions coincide. First, we define $\calS_3 \bfFlagCoh^{\cup}_{\scrT,\scrF,\scrT}(\scrC,\tau)$ as the open substack of $\calS_3 \bfPerf(\scrC)$ parametrizing flags $\F$ of the form \eqref{eq:commutator_3_flag} for which $F_{0,1}[-1], F_{2,3}\in \bfCoh_{\scrT}(\scrC, \tau)$ and $F_{1,2}, F_{0,3}\in \bfCoh_{\scrF}(\scrC, \tau)$. Finally, we define $\calS_3 \bfFlagCoh^{\cap}_{\scrT,\scrF,\scrT}(\scrC,\tau)$ as the fiber pullback
\begin{align}
	\begin{tikzcd}[ampersand replacement=\&]
		\calS_3 \bfFlagCoh^{\cap}_{\scrT,\scrF,\scrT}(\scrC,\tau) \arrow{r}\arrow{d} \& \calS_3 \bfFlagCoh^{(1)}_{(\scrT, (\scrF, \scrT))}(\scrC,\tau) \arrow{d}\\
		\calS_3 \bfFlagCoh^{(1)}_{((\scrT, \scrF), \scrT)}(\scrC,\tau) \arrow{r} \& \calS_3 \bfFlagCoh^{\cup}_{\scrT,\scrF,\scrT}(\scrC,\tau)
	\end{tikzcd}\ .
\end{align}
Here, all the maps are open immersions. 
\begin{proposition}\label{prop:geometric_commutation_criterion}
	Assume that $\scrT$ is a Serre subcategory of $\scrC^\heartsuit$.
	\begin{enumerate}\itemsep=0.2cm		
		\item For every geometric point $y \colon \Spec(k) \to \calS_3 \bfFlagCoh^{(1)}_{(\scrT, (\scrF,\scrT))}(\scrC,\tau)$ corresponding to a flag $\F$ of the form \eqref{eq:commutator_3_flag}, the natural morphisms
		\begin{align}
			(F_{1,3})_{\mathsf{tor}} \longrightarrow F_{2,3} \qquad \text{and} \qquad (F_{1,3})_{\mathsf{tor}} \longrightarrow \calH^1_\tau(F_{0,1}) 
		\end{align}
		are monomorphisms in the heart $\scrC^\heartsuit$ of $\tau$.
		
		\item For every geometric point $y \colon \Spec(k) \to \calS_3 \bfFlagCoh^{(1)}_{((\scrT,\scrF), \scrT)}(\scrC,\tau)$ corresponding to a flag $\F$ of the form \eqref{eq:commutator_3_flag}, the natural morphisms
		\begin{align}
			\calH^1_\tau(F_{0,1}) \longrightarrow \calH^1_\tau(F_{0,2}) \qquad \text{and} \qquad \calH^0_\tau(F_{2,3}) \to \calH^1_\tau(F_{0,2}) 
		\end{align}
		are surjective in the heart $\scrC^\heartsuit$ of $\tau$.
	\end{enumerate}
\end{proposition}

\begin{proof}
	Let us prove the first assertion. Consider the following diagram in $\scrC^\heartsuit$
	\begin{align}
		\begin{tikzcd}[ampersand replacement=\&]
			0 \arrow{r} \& (F_{1,3})_{\mathsf{tor}} \arrow{r} \arrow{d}{\alpha} \& F_{1,3} \arrow{r} \arrow{d} \& (F_{1,3})_{\mathsf{t.f.}} \arrow{r} \arrow{d} \& 0 \\
			0 \arrow{r} \& (F_{2,3})_{\mathsf{tor}} \arrow{r}{\sim} \& F_{2,3} \arrow{r} \& 0 
		\end{tikzcd} \ ,
	\end{align}
	whose rows are exact. Applying the snake lemma we obtain that the exact sequence
	\begin{align}
		0 \longrightarrow \ker(\alpha) \longrightarrow F_{1,2} \longrightarrow (F_{1,3})_{\mathsf{t.f.}} \ . 
	\end{align}
	Since $F_{1,2}$ is torsion free by assumption, it follows that $\ker(\alpha)$ is torsion free as well. On the other hand, $\ker(\alpha)$ is a subobject of $(F_{1,3})_{\mathsf{tor}}$. Since $\scrT$ is a Serre subcategory of $\scrC^\heartsuit$, it follows that $\ker(\alpha)$ is torsion as well, and hence $\ker(\alpha) = 0$. This proves the injectivity of the first map. For the second map, we consider instead the following diagram:
	\begin{align}
		\begin{tikzcd}[ampersand replacement=\&]
			0 \arrow{r} \& (F_{1,3})_{\mathsf{tor}} \arrow{r} \arrow{d}{\alpha'} \& F_{1,3} \arrow{r} \arrow{d} \& (F_{1,3})_{\mathsf{t.f.}} \arrow{r} \arrow{d} \& 0 \\
			0 \arrow{r} \& \calH^1_\tau(F_{0,1})_{\mathsf{tor}} \arrow{r}{\sim} \& \calH^1_\tau(F_{0,1}) \arrow{r} \& 0 
		\end{tikzcd}\ .
	\end{align}
	Reasoning as above, we see that $\ker(\alpha')$ is a subobject of $F_{0,3}$, which is torsion-free by assumption. At the same time, $\ker(\alpha')$ is also a subobject of $(F_{1,3})_{\mathsf{tor}}$, and hence (since $\scrT$ is a Serre subcategory), we deduce that $\ker(\alpha') = 0$.
	
	We now prove the second assertion. Starting from the fiber sequence $F_{0,1} \to F_{0,2} \to F_{1,2}$, we obtain the long exact sequence
	\begin{align}
		\calH^1_\tau(F_{0,1}) \longrightarrow \calH^1_\tau(F_{0,2}) \longrightarrow \calH^1_\tau(F_{1,2}) = 0 \ , 
	\end{align}
	which implies the first surjectivity. On the other hand, starting from $F_{0,2} \to F_{0,3} \to F_{2,3}$, we obtain
	\begin{align}
		\calH^0_\tau(F_{2,3}) \longrightarrow \calH^1_\tau(F_{0,2}) \longrightarrow \calH^1_\tau(F_{0,3}) = 0 \ , 
	\end{align}
	whence the conclusion.
\end{proof}

\begin{remark}
	The above proposition provides a generalization of \cite[Claim~3.8]{Negut_Shuffle} and \cite[Lemma~6.5]{Toda_Hall_Categorical_DT}, where the authors dealt with the left and right actions of the stack of zero-dimensional sheaves on a smooth complex surface.
\end{remark}

\begin{corollary}\label{cor:support_commutation_criterion}
	Let $\bfT, \bfT'$ be two geometric $\Lambda$-graded derived stacks locally of finite presentation over $k$ and let $i \colon \bfT \to \bfCoh_{\scrT}(\scrC,\tau)$ and $i' \colon \bfT' \to \bfCoh_{\scrT[-1]}(\scrC,\tau_{-\upsilon})$ be two locally rpas morphisms for which the squares \eqref{eq:multiplication-T} and \eqref{eq:left-Hecke} are pullback for $\bfT$ and $\bfCoh(\scrC, \tau)$ and the squares \eqref{eq:multiplication-T} and \eqref{eq:left-Hecke} are pullback for $\bfT'$ and $\bfCoh(\scrC, \tau_{-\upsilon})$. Assume that the conditions of Theorem~\ref{thm:left-right-action} and the following condition holds:
	\begin{enumerate}\itemsep=0.2cm
		\item \label{item:condition-1} $\scrT$ is a Serre subcategory of $\scrC^\heartsuit$.
		
		\item \label{item:condition-2} for every geometric point $x \colon \Spec(k) \to \bfCohps(\scrC,\tau)$ classifying an object $M \in \scrC_k$, if there exist two surjections
		\begin{align}
			\calH^1_\tau(F') \longrightarrow M \qquad \text{and} \qquad \calH^0_\tau(F) \longrightarrow M  
		\end{align}
		in $\scrC_k^\heartsuit$, such that $F$ is classified by a point $y \in \mathbf T(k)$ via the embedding $i$, while $F'$ is classified by a point $y' \in \mathbf T'(k)$ via the embedding $i'$, then $M = 0$.
		
		\item \label{item:condition-3} for every geometric point $x \colon \Spec(k) \to \bfCohps(\scrC,\tau)$ classifying an object $M \in \scrC_k$, if there exist two monomorphisms
		\begin{align}
			M \longrightarrow \calH^1_\tau(F')  \qquad \text{and} \qquad M \longrightarrow \calH^0_\tau(F)   
		\end{align}
		in $\scrC_k^\heartsuit$, such that $F$ is classified by a point $y \in \mathbf T(k)$ via the embedding $i$, while $F'$ is classified by a point $y' \in \mathbf T'(k)$ via the embedding $i'$, then $M = 0$.
	\end{enumerate}
	In this case, the fiber product
	\begin{align}
		\left(\bfT \times \bfCoh_{\scrF}(\scrC,\tau) \times \bfT'\right) \times_{\bfCoh_{\scrT}(\scrC,\tau) \times \bfCoh_{\scrF}(\scrC,\tau) \times \bfCoh_{\scrT[-1]}(\scrC,\tau_{-\upsilon})} \calS_3 \bfFlagCoh^{(1)}_{(\scrT, (\scrF, \scrT))}(\scrC,\tau) 
	\end{align}
	is equivalent to
	\begin{align}
		\left(\bfT \times \bfCoh_{\scrF}(\scrC,\tau) \times \bfT' \right) \times_{\bfCoh_{\scrT}(\scrC,\tau) \times \bfCoh_{\scrF}(\scrC,\tau) \times \bfCoh_{\scrT[-1]}(\scrC,\tau_{-\upsilon})}  \calS_3 \bfFlagCoh^{(1)}_{((\scrT, \scrF), \scrT)}(\scrC,\tau) \ .
	\end{align}
	In particular, for any $\calG\in \catCohb_{\mathsf{pro}}(\bfT)$, $\calG'\in \catCohb_{\mathsf{pro}}(\bfT')$, and $\calE\in \catCohb_{\mathsf{pro}}(\bfCoh_{\scrF}(\scrC, \tau))$ we get
	\begin{align}
		(i_\ast'\calG' \oast \calE) \oast i_\ast\calG \simeq i_\ast'\calG' \oast (\calE \oast i_\ast\calG)\ .
	\end{align}
\end{corollary}

\begin{proof}
	The assumptions on $\bfT$ and $\bfT'$ ensure that $\catCohb_{\mathsf{pro}}(\bfT)$ and $\catCohb_{\mathsf{pro}}( \bfT' )$ have the structure of $\E_1$-monoidal pro-$\infty$-categories, and $\catCohb_{\mathsf{pro}}( \bfF )$ has the structure of a categorical left module over $\catCohb_{\mathsf{pro}}(\bfT)$ and of a categorical right module over $\catCohb_{\mathsf{pro}}( \bfT' )$. Indeed, it is enough to apply Corollaries~\ref{cor:COHA-representations}, \ref{cor:induced_COHA-properness}, and \ref{cor:induced_COHA-derived-lci} to $\bfM=\bfM'=\bfCoh_{\scrF}(\scrC, \tau)$, $\bfH=\bfCoh_{\scrT}(\scrC, \tau)$, and $\bfH'=\bfT$ and Corollaries~\ref{cor:COHA-representations-right} and the ``right'' version of Corollaries~\ref{cor:induced_COHA-properness} and \ref{cor:induced_COHA-derived-lci} to $\bfM=\bfM'=\bfCoh_{\scrF}(\scrC, \tau_\upsilon)$, $\bfH=\bfCoh_{\scrT[-1]}(\scrC, \tau_\upsilon)$, and $\bfH'=\bfT'$.
	
	We now prove the second statement.
	Consider the natural maps
	\begin{align}
		\begin{tikzcd}[column sep=small, ampersand replacement=\&]
			\& \calS_3 \bfFlagCoh^\cap_{\scrT,\scrF,\scrT}(\scrC,\tau) \arrow{dr} \arrow{dl} \\
			\calS_3 \bfFlagCoh^{(1)}_{((\scrT,\scrF),\scrT)}(\scrC,\tau) \& \& \calS_3 \bfFlagCoh^{(1)}_{(\scrT,(\scrF,\scrT))}(\scrC,\tau) 
		\end{tikzcd}\ .
	\end{align}
	We claim that our hypothesis \eqref{item:condition-1}, \eqref{item:condition-2}, and \eqref{item:condition-3} imply that they become equivalent after base-changing to $\bfT \times {\bfCoh_\scrF(\scrC,\tau)} \times \bfT'$.
	We argue for the left diagonal map, as the argument is symmetrical for the right one. Since the map in question is already known to be an open immersion, it suffices to show that it is surjective on closed points. Thus, fix a field $k$ and assume that we are given a flag $\F$ of the form \eqref{eq:commutator_3_flag}, which we reproduce here for the convenience of the reader:
	\begin{align}
		\begin{tikzcd}[ampersand replacement = \&]
			0 \arrow{r} \& F_{0,1} \arrow{r} \arrow{d} \& F_{0,2} \arrow{r} \arrow{d} \& F_{0,3} \arrow{d} \\
			\& 0 \arrow{r} \& F_{1,2} \arrow{r} \arrow{d} \& F_{1,3} \arrow{d} \\
			\& \& 0 \arrow{r} \& F_{2,3} \arrow{d} \\
			\& \& \& 0
		\end{tikzcd}\ .
	\end{align}
	If $\F$ classifies a $k$-point of
	\begin{align}
		(\bfT \times \bfCoh_\scrF(\scrC,\tau) \times \bfT') \times_{\bfCoh_{\scrT}(\scrC,\tau) \times \bfCoh_{\scrF}(\scrC,\tau) \times \bfCoh_{\scrT[-1]}} \calS_3 \bfFlagCoh^{(1)}_{((\scrT,\scrF),\scrT)}(\scrC,\tau) \ , 
	\end{align}
	we have that $F_{0,1}[-1]$ and $F_{2,3}$ belong to $\scrT_k$, and $F_{1,2}$, $F_{0,2}$ and $F_{0,3}$ belong to $\scrF_k$. Furthermore, $F_{0,1}$ classifies a $k$-point of $\mathbf T'$ and $F_{2,3}$ classifies a $k$-point of $\bfT$. To conclude that $\F$ classifies a $k$-point of $\calS_3 \bfFlagCoh^\cap_{\scrT,\scrF,\scrT}(\scrC,\tau)$, we need to argue that $F_{1,3}$ belongs to $\scrF_k$ as well.
	To begin with, the extension $\partial_0(\F)$ implies that $F_{1,3}$ belongs to $\scrC_k^\heartsuit$.
	Applying Proposition~\ref{prop:geometric_commutation_criterion}, we obtain monomorphisms
	\[ (F_{1,3})_{\mathsf{tor}} \longrightarrow F_{2,3} \qquad \text{and} \qquad (F_{1,3})_{\mathsf{tor}} \longrightarrow \calH^1_\tau(F_{0,1}) \ . \]
	This implies, via our hypothesis \eqref{item:condition-3} that $(F_{1,3})_{\mathsf{tor}} = 0$, and therefore that $F_{1,3} \in \scrF_k$.
\end{proof}

\subsection{Categorified commutator fitting in a fiber sequence}

Let $\bfT\to \bfCoh_{\scrT}(\scrC, \tau)$ be a morphism representable by open and closed immersions and  let $\bfF\to \bfCoh_\scrF(\scrC, \tau)$ be a morphism for which the squares \eqref{eq:multiplication-T}, \eqref{eq:right-Hecke}, \eqref{eq:left-Hecke} are pullbacks. Assume that $\bfF$ is a two-sided Hecke pattern for $\bfT$ with respect to $\tau$.

Define the stacks
\begin{align}
	\begin{tikzcd}[ampersand replacement=\&]
		\calS_3 \bfFlagCoh^{(1)}_{((\bfT, \bfF), \bfT)}(\scrC,\tau) \ar{r}\ar{d}\&\calS_3 \bfFlagCoh^{(1)}_{((\scrT, \scrF), \scrT)}(\scrC,\tau)\ar{d}\\
		\bfT\times \bfF\times \bfT[-1] \times \bfF\ar{r}\&\bfCoh_\scrT(\scrC,\tau)\times \bfCoh_\scrF(\scrC,\tau) \times \bfCoh_{\scrT[-1]}(\scrC,\tau_{-\upsilon})\times \bfCoh_\scrF(\scrC,\tau)
	\end{tikzcd}\ , \\
	\begin{tikzcd}[ampersand replacement=\&]
		\calS_3 \bfFlagCoh^{(1)}_{(\bfT, (\bfF, \bfT))}(\scrC,\tau) \ar{r}\ar{d}\&\calS_3 \bfFlagCoh^{(1)}_{(\scrT, (\scrF, \scrT))}(\scrC,\tau)\ar{d}\\
		\bfT\times \bfF\times \bfT[-1] \times \bfF\ar{r}\&\bfCoh_\scrT(\scrC,\tau)\times \bfCoh_\scrF(\scrC,\tau) \times \bfCoh_{\scrT[-1]}(\scrC,\tau_{-\upsilon})\times \bfCoh_\scrF(\scrC,\tau)
	\end{tikzcd}\ , \\
	\begin{tikzcd}[ampersand replacement=\&]
		\calS_3 \bfFlagCoh^{\cap}_{(\bfT, \bfF, \bfT)}(\scrC,\tau) \ar{r}\ar{d}\&\calS_3 \bfFlagCoh^{\cap}_{(\scrT, \scrF, \scrT)}(\scrC,\tau)\ar{d}\\
		\bfT\times \bfF\times \bfT[-1] \times \bfF\ar{r}\&\bfCoh_\scrT(\scrC,\tau)\times \bfCoh_\scrF(\scrC,\tau) \times \bfCoh_{\scrT[-1]}(\scrC,\tau_{-\upsilon})\times \bfCoh_\scrF(\scrC,\tau)
	\end{tikzcd}\ .
\end{align}
Under the current assumptions, $\calS_3 \bfFlagCoh^{(1)}_{((\bfT, \bfF), \bfT)}(\scrC,\tau)$ coincides with $\calS_3 \bfFlagCoh^{\cap}_{(\bfT, \bfF, \bfT)}(\scrC,\tau)$, since in this case the object $F_{0,2}$ appearing in the diagram of the form \eqref{eq:commutator_3_flag} is also $\tau$-flat. 

Let $\alpha\colon \calS_3 \bfFlagCoh^{\cap}_{\bfT,\bfF,\bfT}(\scrC,\tau)\to \calS_3 \bfFlagCoh^{(1)}_{(\bfT, (\bfF, \bfT))}(\scrC,\tau)$ be the canonical open embedding. Since $q^{\ell\circ r}$ and $q^{r\circ \ell}$ are locally rpas and $q^{r\circ\ell}=q^{\ell\circ r}\circ \alpha$, also $\alpha$ is locally rpas. Thus, there exists a morphism
\begin{align}
	\Psi\colon q^{\ell\circ r}_\ast\circ p_{\ell\circ r}^\ast \longrightarrow q^{r\circ \ell}_\ast\circ p_{r\circ \ell}^\ast\ . 
\end{align}

\begin{definition}
	The \textit{categorical commutator} between the two actions is the cofiber $\mathsf{cofib}(\Psi)$.
\end{definition}

\begin{remark}
	In the framework developed by Neguţ in \cite{Negut_Shuffle, Negut_Categorification}, Yu Zhao has computed explicitly the categorical commutator for the stacks $\bfT$ and $\bfF$ mentioned in the previous remark in \cite{Zhao_Commutators_Hilbert, Zhao_Commutators}.
	
	In our general framework, the explicit characterization of the categorical commutator seems to be an hard question to be addressed, also under the assumption of two-sided Hecke pattern property which allows us to encode the categorical commutator as a cofiber of $\Psi$.
\end{remark}

\section{Appendix: Categorical completions}\label{appendix:ind_objects}

We collect in this appendix some standard facts concerning categorical completions and ind-objects that are needed to define Borel-Moore homology in the non-quasi-compact case.

\medskip

Let $\scrK$ be a collection of $\infty$-categories and let $\scrC$ be a fixed $\infty$-category. Following \cite[\S5.3.6]{HTT}, we denote by $\scrP^{\scrK}(\scrC)$ the full subcategory of $\PSh(\scrC)$ closed under $\scrK$-colimits and containing all representable functors. We refer to $\scrP^{\scrK}(\scrC)$ as the \textit{free $\scrK$-cocompletion of $\scrC$}, a terminology justified by \cite[Corollary~5.3.5.4]{HTT}. We are mostly interested in the following two examples:

\begin{example}
	\hfill
	\begin{enumerate}\itemsep=0.2cm
		\item When $\scrK$ consists of the collection of small filtered diagrams, we write $\Ind(\scrC)$ instead of $\scrP^\scrK(\scrC)$, and we refer to it as the \textit{ind-completion} of $\scrC$.
		
		\item When $\scrK = \{\mathbb N\}$ (where $\mathbb N$ is seen as a discrete category), we write $\PSh^\sqcup(\scrC)$ instead of $\scrP^\scrK(\scrC)$, and we refer to it as the \textit{countable coproduct completion} of $\scrC$.
	\end{enumerate}
\end{example}

We represent an object $X \in \scrP^\scrK(\scrC)$ via the notation
\begin{align}
	X \simeq \fcolim_{i \in I} X_i \ , 
\end{align}
by which we mean that we are given a functor $I \to \scrC$ from an $\infty$-category $I \in \scrK$, and that we consider the colimit of this diagram as an object in $\scrP^\scrK(\scrC)$. In this case, we say that this is a \textit{presentation for $X$}. Mapping spaces in $\scrP^\scrK(\scrC)$ are characterized by the relation
\begin{align}
	\Map_{\scrP^\scrK(\scrC)}\big(\fcolim_{i \in I} X_i, \fcolim_{j \in J} Y_j\big) \simeq \lim_{i\in I} \colim_{j \in J} \Map_{\scrC}( X_i, Y_j ) \ . 
\end{align}
Notice that by definition $\scrP^\scrK(\scrC)$ is a full subcategory of $\PSh(\scrC)$. In particular, given a presheaf $F \colon \scrC\op \to \Spc$, the $\scrK$-structure on $F$, if it exists, is unique.

\medskip

In the main body of this paper, we are interested in taking $\scrK$ as in the above example, and $\scrC$ to be mild variations of the category of qcqs geometric stacks. However, we will not be interested in the whole $\scrP^{\scrK}(\scrC)$, but rather to specific subcategories cut out via the following procedure:
\begin{definition}\label{def:restricted_presentations}
	Let $\scrC_0 \subseteq \scrC$ be a (not necessarily full) subcategory of $\scrC$.
	\begin{enumerate}\itemsep=0.2cm
		\item We say that a presentation $\{X_i\}_{i \in I}$ for an ind-object $X \in \scrP^{\scrK}(\scrC)$ is a $\scrC_0$-presentation if the diagram $I \to \scrC$ factors through $\scrC_0$.
		
		\item We say that an object $X \in \scrP^{\scrK}(\scrC)$ is of \textit{$\scrC_0$-type} if it admits a $\scrC_0$-presentation.
	\end{enumerate}
	We let $\scrP^{\scrK}_{\scrC_0}(\scrC)$ denote the full subcategory of $\scrP^{\scrK}(\scrC)$ spanned by objects of $\scrC_0$-type.
\end{definition}

\begin{notation}
	Most frequently, the subcategory $\scrC_0$ will be determined by specifying a property $P$ for the objects or for the morphisms of $\scrC$ (in the latter case, $P$-morphisms should be closed under compositions, and all identities should satisfy $P$). This determines a subcategory $\scrC_P$ of $\scrC$, and we will say that an ind-object $\scrX$ is ind-$P$ rather than saying that it is of $\scrC_P$-type. We will also write $\Ind_P(\scrC)$ instead of $\Ind_{\scrC_P}(\scrC)$.
\end{notation}

When dealing with full subcategory, the following technical lemma is extremely useful in guaranteeing the existence of a restricted presentation.

\begin{lemma}\label{lem:Q_presentations}
	Let $\scrE$ be an $\infty$-category and let $\scrE_0$ be a \textit{full} subcategory of $\scrE$. Given $\scrX \in \Ind(\scrE)$, write $(\scrE_0)_{/\scrX} \coloneqq \scrE_0 \times_{\Ind(\scrE)} \Ind(\scrE)_{/\scrX}$. Then following statements are equivalent:
	\begin{enumerate}\itemsep=0.2cm
		\item \label{item:Q_presentations-1} the object $\scrX \in \Ind(\scrE)$ is of $\scrE_0$-type;
		
		\item \label{item:Q_presentations-2} the $\infty$-category $(\scrE_0)_{/\scrX}$ is filtered and the canonical functor $(\scrE_0)_{/\scrX} \to \scrE_{/\scrX}$ is colimit-final;
		
		\item \label{item:Q_presentations-3} for every morphism $u \colon \overline{X} \to \scrX$ where $\overline{X} \in \Ind(\scrE)^\omega$, there exists a factorization of $u$ as
		\begin{align}
			\overline{X} \longrightarrow X \longrightarrow \scrX \ , 
		\end{align}
		where $X$ belongs $\scrE_0$.
	\end{enumerate}
\end{lemma}

\begin{proof}
	The equivalence \eqref{item:Q_presentations-1} $\Leftrightarrow$  \eqref{item:Q_presentations-2} is standard and follows from Quillen's Theorem~A (see \cite[Theorem~4.1.3.1]{HTT}).
	
	The implication  \eqref{item:Q_presentations-1} $\Rightarrow$  \eqref{item:Q_presentations-3} simply follows unraveling the definitions.
	
	Let us prove that  \eqref{item:Q_presentations-3} $\Rightarrow$  \eqref{item:Q_presentations-2}. We will prove the finality statement; essentially the same argument proves that $(\scrE_0)_{/\scrX}$ is filtered as well. In virtue of Quillen's Theorem~A, we have to prove that for every $X \in \scrE_{/\scrX}$, the $\infty$-category of factorizations
	\begin{align}
		(\scrE_0)_{X/\!\!/\scrX} \coloneqq (\scrE_0)_{/\scrX} \times_{\scrE_{/\scrX}} \scrE_{X/\!\!/\scrX} 
	\end{align}
	is weakly contractible. For this, it is enough to prove that it is filtered. Let therefore $J$ be a finite category, and let $F \colon J \to (\scrE_0)_{X/\!\!/\scrX}$ be a diagram. To construct an extension to $J^{\triangleright}$, let $\overline{X}$ be the colimit of $F$ computed in $\scrE_{X/\!\!/\scrX}$. Since $J$ is finite, we see that $\overline{X}$ is a compact object in $\Ind(\scrE)$. Thus, assumption \eqref{item:Q_presentations-3} guarantees that the structural morphism $\overline{X} \to \scrX$ factors through an object $X$ satisfying $Q$. This provides the extension we were looking for.
\end{proof}

Fix now a presentable $\infty$-category $\scrT$ together with a property $P$ of objects in $\scrT$, and a property $Q$ of morphisms in $\scrT$. We will assume that $Q$-morphisms are stable under composition and that every identity of $\scrT$ satisfies $Q$. We write $\scrT_Q$ for the (non full) subcategory of $\scrT$ determined by the property $Q$, and $\scrT_{P \cap Q}$ for the full subcategory of $\scrT_Q$ spanned by objects satisfying $P$. The restricted Yoneda embedding gives rise to a functor
\begin{align}
	\Phi_P \colon \scrT \longrightarrow \PSh(\scrT_P) 
\end{align}
sending an object $X \in \scrT$ to the presheaf $\Phi_P(X)$ defined as
\begin{align}
	\Phi_P(X)(U) \coloneqq \Map_{\scrT}(U,X) \ . 
\end{align}
As observed previously, it makes sense to ask for which objects $X$, the presheaf $\Phi_P(X)$ defines an ind-object. To answer this, we introduce the following notion of admissibility:

\begin{defin}\label{def:PQ_admissibility}
	Let $X \in \scrT$.
	We say that $X$ is \textit{$(P,Q)$-admissible} if it admits a $\scrT_{P\cap Q}$-presentation $\{X_i\}_{i \in I}$ such that for every $S \in \scrT_P$ the canonical morphism
	\begin{align}
		\colim_{i \in I} \Map_{\scrT}(S, X_i) \longrightarrow \Map_{\scrT}(S, X) 
	\end{align}
	is an equivalence in $\Spc$. In this case, we say that $\{X_i\}$ is a \textit{$(P,Q)$-admissible presentation} for $X$. \hfill $\oslash$ 
\end{defin}

\begin{proposition}\label{prop:indization}
	Let $X \in \scrT$ be a $(P,Q)$-admissible object.
	Then:
	\begin{enumerate}\itemsep=0.2cm
		\item \label{prop-item:indization-1} the presheaf $\Phi_P(X) \in \PSh(\scrT_P)$ is an ind-object;
		
		\item \label{prop-item:indization-2} if $\{X_i\}_{i \in I}$ is a $(P,Q)$-presentation for $X$, then
		\begin{align}
			\Phi_P(X) \simeq \fcolim_{i \in I} X_i 
		\end{align}
		as objects in $\Ind(\scrT_P)$.
		
		\item \label{prop-item:indization-3} The functor
		\begin{align}
			\Phi_P \colon \scrT \longrightarrow \PSh(\scrT_P) 
		\end{align}
		commutes with limits and it is fully faithful once restricted to $(P,Q)$-admissible objects.
	\end{enumerate}
\end{proposition}

\begin{proof}
	Fix an admissible $(P,Q)$-presentation $\{X_i\}_{i \in I}$ for $X$. Applying \cite[Corollary~5.3.5.4(1)]{HTT}, we see that it is enough to check that for every object $S \in \scrT_P$, the canonical map
	\begin{align}
		\colim_{i \in I} \Map_{\scrT}(S, X_i) \longrightarrow \Map_{\scrT}\big(S, \colim_{i \in I} X_i \big) 
	\end{align}
	is an equivalence. Since $S$ satisfies $P$, and the diagram factors through $\scrT_{P \cap Q}$, this is true by assumption. This simultaneously proves points \eqref{prop-item:indization-1} and \eqref{prop-item:indization-2}. The Yoneda lemma guarantees that $\Phi_P$ commutes with limits. We are therefore left to check the second half of point \eqref{prop-item:indization-3}. Let therefore $X$ and $Y$ be two $(P,Q)$-admissible objects and fix presentations $\{X_i\}_{i \in I}$ and $\{Y_j\}_{j \in J}$. Then
	\begin{align}
		\Map_{\Ind(\scrT_P)}(\Phi_P(X), \Phi_P(Y)) & \simeq \Map_{\Ind(\scrT_P)}\big( \fcolim_{i \in I} X_i, \fcolim_{j \in J} Y_j \big) & \text{By Point \eqref{prop-item:indization-2}} \\
		& \simeq \lim_{i \in I} \colim_{j \in J} \Map_{\scrT_P}( X_i, Y_j ) \\
		& \simeq \lim_{i \in I} \colim_{j \in J} \Map_{\scrT}(X_i, Y_j) \\
		& \simeq \lim_{i \in I} \Map_{\scrT}\big( X_i, \colim_{j \in J} Y_j\big) & \text{By } (P,Q)\text{-admissibility of } Y \\
		& \simeq \Map_{\scrT}\big(\colim_{i \in I} X_i, \colim_{j \in J} Y_j\big) & \text{by Yoneda} \\
		& \simeq \Map_{\scrT}(X, Y) \ .
	\end{align}
	The conclusion follows.
\end{proof}

\newpage
\part{COHAs, CatHAs, and their representations: Applications}\label{part:application}

In this part, we apply the foundational results obtained in Parts~\ref{part:Segal} and \ref{part:foundations} to concrete examples. Specifically, given a smooth projective complex surface $S$, we introduce COHAs and CatHAs associated to \textit{torsion sheaves} on $S$ and (categorical) representations of them associated to \textit{torsion-free sheaves} on $S$. In particular, we categorify known constructions of usual cohomological and K-theoretical Hall algebras. Furthermore, we introduce representations of these COHAs and CatHAs associated to \textit{stable pairs} on $S$. 

Furthermore, we provide a categorification of the preprojective COHA of a quiver and the representation of it in terms of Nakajima quiver varieties.

\section{Perverse $t$-structures for surfaces}

In this section, we recall Deligne's perverse $t$-structures following \cite[\S3.1, \S7.1, and \S7.2]{Bayer_Large} (see also \cite{Bez_Perverse} and \cite{Kashiwara_Perverse}, where perverse $t$-structures were originally introduced). We also characterize those on the derived category $\catPerf(S)$ of perfect complexes on a smooth projective complex surface $S$.

\subsection{Recollection of Deligne's perverse $t$-structures}

Let $X$ be a smooth projective complex variety of dimension $n$. Fix a function $\sfp\colon \{0, 1, \ldots, n\}\to \Z$ and introduce the ``dual'' function $\overline{\sfp}\colon \{0, 1, \ldots, n\}\to \Z$ by $\overline{p}(d)=-d-p(d)$.
\begin{definition}
	We say that $\sfp$ is a \textit{perversity function} if the functions $\sfp$ and $\overline{\sfp}$ are monotone decreasing. In this is a case, we call $\overline{\sfp}$ the \textit{dual} perversity function.
\end{definition}

\begin{remark}
	Note that if $\sfp(d)\geq \sfp(d+1)\geq \sfp(d)-1$, then $\sfp$ is a perversity function.
\end{remark}
Fix a perversity function $\sfp$ and define the following increasing filtration of $\catCoh(X)$ by abelian subcategories:
\begin{align}
	\scrA^{\sfp, \leqslant k}\coloneqq \left\{ \calF\in \catCoh(X)\, \vert\, \sfp(\dim\mathsf{supp}(\calF))\geq -k\right\}\ .
\end{align}
Then the pair $(\catPerf(X)^{\sfp, \leqslant 0}, \catPerf(X)^{\sfp, \geqslant 0})$ given by
\begin{align}
	\catPerf(X)^{\sfp, \leqslant 0} &\coloneqq \Big\{ E\in \catPerf(X)\, \vert\, \calH^{-k}(E)\in \scrA^{\sfp, \leqslant k} \text{ for all } k\in \Z\Big\}\ , \\[4pt]
	\catPerf(X)^{\sfp, \geqslant 0} &\coloneqq \Big\{ E \in \catPerf(X)\, \vert\, \Hom_{\catPerf(X)} (A, E) = 0 \text{ for all } A\in \scrA^{\sfp, \leqslant k}[k+1] \text{ and for all } k\in \Z\Big\} 
\end{align}
defines a bounded $t$-structure $\tau^\sfp$ on $\catPerf(X)$. Objects in the heart $\scrA^\sfp \coloneqq \catPerf(X)^{\sfp, \leqslant 0}  \cap \catPerf(X)^{\sfp, \geqslant 0}$ are called \textit{perverse coherent sheaves}.

\subsubsection{Perverse $t$-structures and torsion pairs}

Fix $k=-\sfp(\ell)$ for some $0\leq \ell\leq n$. Define the function $\sfp^k\colon \{0, 1, \ldots, n\} \to \Z$ given by
\begin{align}\label{eq:perversity-k}
	\sfp^k(d)\coloneqq \begin{cases}
		\sfp(d) & \text{ if } \sfp(d)\geq -k\ , \\[2pt]
		\sfp(d)+1 & \text{ if } \sfp(d)< -k\ .
	\end{cases}
\end{align}
Then, $\sfp^k$ is a perversity function. 
\begin{lemma}\label{lem:perversity-torsion-pair}
	The $t$-structure $\tau^\sfp$ is the tilting of $\tau^{\sfp^k}$ with respect to the torsion pair $(\scrA^{\sfp^k}\cap \scrA^{\sfp}, \scrA^{\sfp^k}\cap \scrA^{\sfp}[1])$.
\end{lemma}

\begin{proof}
	Since $\scrA^{\sfp, \leqslant i}\subseteq \scrA^{\sfp^k, \leqslant i} \subseteq \scrA^{\sfp, \leqslant i+1}$, we get
	\begin{align}
		\catPerf(X)^{\sfp^k, \leqslant -1}\subseteq \catPerf(X)^{\sfp, \leqslant 0} \subseteq \catPerf(X)^{\sfp^k, \leqslant 0}\ .
	\end{align}
	By Remark~\ref{rem:relative_boundedness}, the claim follows.
\end{proof}

\subsubsection{Perverse $t$-structures and duality}

Let $\D\colon \catPerf(X)\op \to \catPerf(X)$ be the dualizing functor defined by
\begin{align}
	E\longmapsto \R\calH om(E, \omega_X[\dim(X)])\ ,
\end{align}
where $\omega_X$ is the canonical bundle of $X$. We recall the following result relating the $t$-structures associated to a perversity function and its dual perversity function.
\begin{proposition}\label{prop:duality}
	Let $\sfp\colon \{0, \ldots, n\} \to \Z$ be a perversity function and let $\overline{\sfp}$ be its dual perversity function.
	Then, the autoequivalence $\D\colon \catPerf(X)\op \to \catPerf(X)$ is $t$-exact with respect to the $t$-structure $(\tau^\sfp)\op$ on $\catPerf(X)\op$ and the $t$-structure $\tau^{\overline{\sfp}}$ on $\catPerf(X)$. 
\end{proposition}

\subsection{Perverse $t$-structures on surfaces}\label{subsec:perverse-t-structure-surfaces}

Let $S$ be a smooth projective irreducible complex surface. For any fixed $n\in \Z$, there are four possible perversity functions $\tensor*[_n]{\sfp}{_i}\colon \{0, 1, 2\} \to \Z$ for $i=1, 2, 3, 4$ given by
\begin{align}
	\tensor*[_n]{\sfp}{_1}(0) & = n = \tensor*[_n]{\sfp}{_1}(1) =  \tensor*[_n]{\sfp}{_1}(2) \ ,\\[4pt] 
	\tensor*[_n]{\sfp}{_2}(0) & = n = \tensor*[_n]{\sfp}{_2}(1) \ , \quad \tensor*[_n]{\sfp}{_2}(2)=n-1 \ , \\[4pt]
	\tensor*[_n]{\sfp}{_3}(0) & = n \ , \quad  \tensor*[_n]{\sfp}{_2}(1)=n-1=\tensor*[_n]{\sfp}{_2}(2) \ , \\[4pt]
	\tensor*[_n]{\sfp}{_4}(0) & = n \ , \quad \tensor*[_n]{\sfp}{_4}(1) = n-1 \ , \quad \tensor*[_n]{\sfp}{_4}(2) = n-2 \ .
\end{align}
Since the $t$-structure associated to $\tensor*[_n]{\sfp}{_i}$ is the shift by $-n$ of the $t$-structure associated to $\tensor*[_0]{\sfp}{_i}$, for $i=1, \ldots, 4$, we shall characterize explicitly only the latter. In particular, we set
\begin{align}
	\sfp_{\mathsf{std}}\coloneqq \tensor*[_0]{\sfp}{_1}\  , \quad \sfp_{1,2}\coloneqq \tensor*[_0]{\sfp}{_2}\ , \quad \sfp_{0, 1} \coloneqq \tensor*[_n]{\sfp}{_3}\ , \quad \sfp_{0,1,2} \coloneqq \tensor*[_n]{\sfp}{_4} \ .
\end{align}
Note that the perverse $t$-structure associated to $\sfp_{\mathsf{std}}$ is the standard $t$-structure $\tau_{\mathsf{std}}$ on $\catPerf(S)$. 

We shall denote by $\tau_\scrA$ the $t$-structure on $\catPerf(S)$ associated to $\sfp_{0,1}$ and by $\scrA$ the corresponding heart, while we shall denote by $\tau_\scrB$ the $t$-structure of $\catPerf(S)$ associated to $\sfp_{1,2}$ and by $\scrB$ the corresponding heart\footnote{The use of these notations is inspired by \cite[\S2]{Toda_Stable}.}. Finally, we denote by $\tau_{0,1,2}$ the $t$-structure on $\catPerf(S)$ associated to $\sfp_{0,1,2}$ and by $\catPerf(S)^{\heartsuit_{0,1,2}}$ the corresponding heart.

Note that $\sfp_{0,1,2}$ is the dual perversity function of $\sfp_{\mathsf{std}}$, while $\sfp_{0, 1}$ is the dual perversity function of $\sfp_{1, 2}$.

\begin{remark}\label{rem:large_volume}
	The perversity function $\sfp_{1,2}$ corresponds to the function 
	\begin{align}
		d\longmapsto -\left\lfloor \frac{d}{2}\right\rfloor\ .
	\end{align}
	For this reason, this is called the \textit{large volume perversity function} (cf.\ \cite[\S4]{Bayer_Large}). 
\end{remark}

Proposition~\ref{prop:duality} yields the following.
\begin{corollary}\label{cor:autoequivalences-perversities}
	Then autoequivalence
	\begin{align}
		\D\colon \catPerf(S)\op \longrightarrow \catPerf(S) 
	\end{align}
	restricts to equivalences
	\begin{align}
		\D \colon \scrA\op \xrightarrow{\sim} \scrB \quad \text{and}\quad
		\D \colon \catCoh(S)\op \xrightarrow{\sim} \catPerf(S)^{\heartsuit_{0,1,2}} \ .
	\end{align}
\end{corollary}

Now, let us introduce the following full subcategories of $\catCoh(S)$:
\begin{multline}
	\catCoh_{0}(S)\coloneqq \{ \calE\in \catCoh(S)\, \vert\, \dim\mathsf{Supp}(\calE)=0\}\\
	\shoveright{\subset\catCoh_{\mathsf{tor}}(S)\coloneqq \{\calE\in\catCoh(S)\, \vert\, \dim\mathsf{Supp}(\calE)\leq 1\}\ , }\\[2pt]
	\shoveleft{\catCoh_{\mathsf{t.f.}}(S)\coloneqq \{\calE\in\catCoh(S)\, \vert\, \calE\text{ is torsion-free}\}}\\ 
	\subset\catCoh_{\geqslant 1}(S) \coloneqq \{ \calE\in \catCoh(S) \, \vert \, \Hom_S(\calT, \calE)=0 \, \text{ for any }\, \calT\in \catCoh_{0}(S) \} \ .
\end{multline}
Since $\sfp_{1,2}^0=\sfp_{\mathsf{std}}$ and $\sfp_{0,1}^0=\sfp_{\mathsf{std}}$, by Lemma~\ref{lem:perversity-torsion-pair}, the $t$-structure $\tau_\scrA$ is the tilting of the standard $t$-structure of $\catPerf(S)$ by the torsion pair $(\catCoh_{0}(S), \catCoh_{\geqslant 1}(S))$\footnote{One can show that it is a torsion pair using the same arguments as in the proof of \cite[Lemma~2.10]{Toda_Stable}.}, while the $t$-structure $\tau_\scrB$ is the tilting of the standard $t$-structure by the torsion pair $(\catCoh_{\mathsf{tor}}(S), \catCoh_{\mathsf{t.f.}}(S))$. 

\begin{remark}\label{rem:tor_serre}
	Note that $\catCoh_{\mathsf{tor}}(S)$ is uniquely characterized as the subcategory of $\catCoh(S)$ consisting by rank zero sheaves. Thus, it is straightforward to see that $\catCoh_{\mathsf{tor}}(S)$ is Serre subcategory of $\catCoh(S)$, hence it is also a Serre subcategory of $\scrB$ by Lemma~\ref{lem:tor_Serre}.
\end{remark}

\subsection{A torsion pair on the heart $\scrA$}\label{subsec:torsion-pair-perverse}

In this section, we shall define a torsion pair on $\scrA$ which is ``dual'' to the torsion pair $(\catCoh_{\mathsf{tor}}(S), \catCoh_{\mathsf{t.f.}}(S))$ of $\catCoh(S)$.

Following \cite[\S2.3]{Toda_Stable}, let $\scrA_{\mathsf{tor}}$ be the smallest full subcategory of $\scrA$ closed under extensions and containing both $\calF[1]$ for pure one-dimensional sheaves $\calF$, and $\scrO_x$ for $x\in S$; we also set:
\begin{align} \label{eq:A12}
	\scrA_{\mathsf{t.f.}} \coloneqq \{ E\in \scrA\, \vert\, \Hom_{\scrA}(F, E) = 0 \text{ for any } F\in \scrA_{\mathsf{tor}} \}\subset \scrA \ .
\end{align}

The following lemma can be proven by an argument parallel to the one given in \cite[Lemma~2.16]{Toda_Stable}\footnote{One can also apply Lemma~\ref{lem:perversity-torsion-pair} since $\sfp_{0, 1, 2}^0=\sfp_{0, 1}$.}.
\begin{lemma}\label{lem:perverse_torsion_pair}
	The pair $(\scrA_{\mathsf{tor}},\scrA_{\mathsf{t.f.}} )$ is a torsion pair of $\scrA$.
\end{lemma}	

\begin{definition}
	We call the objects of $\scrA_{\mathsf{tor}}$ \textit{torsion perverse sheaves} on $S$, while the objects of $\scrA_{\mathsf{t.f.}}$ \textit{torsion-free perverse sheaves} on $S$.
\end{definition}

Now, we provide an explicit characterization of torsion perverse sheaves.
\begin{lemma} \label{lem:A_1-1}
	An object $E \in \scrA$ belongs to $\scrA_{\mathsf{tor}}$ if and only if $\mathsf{rk}(E) = 0$.
	In particular, 
	\begin{enumerate}\itemsep0.2cm
		\item $\calH^0(E)$ is zero-dimensional and $\calH^{-1}(E)$ is pure one-dimensional for any $E\in \scrA_{\mathsf{tor}}$,
		\item \label{item:A_1-1} $\scrA_{\mathsf{tor}}$ is a Serre subcategory of $\scrA$.
	\end{enumerate}
\end{lemma}

\begin{proof}
	Since the rank is additive, it immediately follows from the definition that every object $E$ in $\scrA_{\mathsf{tor}}$ satisfies $\mathsf{rk}(E) = 0$.
	Vice-versa, if $E \in \scrA$, then $\calH^0(E)$ is zero-dimensional by definition.
	If $\mathsf{rk}(E) = 0$, then we also have $\mathsf{rk}(\calH^{-1}(E)) = 0$, and therefore $\calH^{-1}(E) \in \catCoh_{\mathsf{tor}}(S)$.
	On the other hand, since $E \in \scrA$, we also have $\calH^{-1}(E) \in \catCoh_{\geqslant 1}(S)$, and therefore $\calH^{-1}(E)$ is pure one-dimensional.
	
	For statement \eqref{item:A_1-1}, since $\scrA_{\mathsf{tor}}$ is the torsion part of a torsion pair, it is automatically closed under quotients and extensions.
	If
	\begin{align}
		0 \longrightarrow T_1 \longrightarrow T \longrightarrow T_2 \longrightarrow 0 
	\end{align}
	is a short exact sequence in $\scrA$ and $T \in \scrA_{\mathsf{tor}}$, then $T_2 \in \scrA_{\mathsf{tor}}$, whence $\mathsf{rk}(T_1) = \mathsf{rk}(T) - \mathsf{rk}(T_2) = 0$.
	Thus, $T_1 \in \scrA_{\mathsf{tor}}$.
\end{proof}

Now, we provide a more insightful description of $\scrA_{\mathsf{t.f.}}$. 
\begin{theorem}\label{thm:duality-tor-torsion-free} 
	The autoequivalence
	\begin{align}
		\D \colon \scrA\op \xrightarrow{\sim} \scrB 
	\end{align}
	identifies the torsion pair $(\scrA_{\mathsf{t.f.}}\op, \scrA_{\mathsf{tor}}\op)$ on $\scrA\op$ with the torsion pair $(\catCoh_{\mathsf{t.f.}}(S)[1], \catCoh_{\mathsf{tor}}(S))$ on $\scrB$.
	In other words, $\D$ also induces equivalences
	\begin{align}
		\D \colon \scrA_{\mathsf{tor}}\op \xrightarrow{\sim} \catCoh_{\mathsf{tor}}(S) \quad \text{and}\quad 
		\D \colon \scrA_{\mathsf{t.f.}}\op \xrightarrow{\sim} \catCoh_{\mathsf{t.f.}}(S) [1] \ .
	\end{align}
\end{theorem}

\begin{proof}
	Since $\D$ induces an equivalence between $\scrA\op$ and $\scrB$ by Corollary~\ref{cor:autoequivalences-perversities}, it is enough to prove the following two statements:
	\begin{enumerate}\itemsep=0.2cm
		\item \label{item:duality-tor-torsion-free} an object $E \in \catPerf(S)$ belongs to $\scrA_{\mathsf{tor}}$ if and only if $\D(E)$ belongs to $\catCoh_{\mathsf{tor}}(S)$;
		
		\item an object $E \in \catPerf(S)$ belongs to $\scrA_{\mathsf{t.f.}}$ if and only if $\D(E)$ belongs to $\catCoh_{\mathsf{t.f.}}(S)[1]$.
	\end{enumerate}
	Moreover, since $\scrA_{\mathsf{t.f.}}\op$ is the left orthogonal to $\scrA_{\mathsf{tor}}\op$ inside $\scrA\op$, and similarly $\catCoh_{\mathsf{t.f.}}(S)[1]$ is the left orthogonal of $\catCoh_{\mathsf{tor}}(S)$ inside $\scrB$, it is actually enough to prove the first statement.
	
	We start by proving the ``only if'' part of statement \eqref{item:duality-tor-torsion-free}. Attached to any object $E$ of $\scrA_{\mathsf{tor}}$ there is a canonical fiber sequence
	\begin{align}
		\calH^{-1}(E)[1] \longrightarrow E \longrightarrow \calH^0(E) \ , 
	\end{align}
	Applying $\D(-)=(-)^\vee[2]\otimes_{\scrO_S} \omega_S$, we obtain
	\begin{align}
		\D(\calH^0(E)) \longrightarrow \D(E) \longrightarrow \D(\calH^{-1}(E)[1]) \ . 
	\end{align}
	Lemma~\ref{lem:A_1-1} guarantees that $\calH^{-1}(E)$ is a pure one-dimensional sheaf and $\calH^0(E)$ is zero-dimen\-sio\-nal. Hence \cite[Proposition~1.1.6--(i)]{HL_Moduli} provides us with canonical identifications
	\begin{align}
		\D(\calH^0(E)) \simeq \calExt^2_S(\calH^0(E),\omega_S) \qquad \text{and} \qquad \D(\calH^{-1}(E)[1]) \simeq \calExt^1_S(\calH^{-1}(E),\omega_S) \ . 
	\end{align}
	This shows at the same time that $\D(\calH^0(E))$ is zero-dimensional and that $\D(\calH^{-1}(E)[1])$ is purely $1$-dimensional.
	Thus, $\calH^i(\D(E)) = 0$ unless $i = 0$, and $\calH^0(\D(E))$ is a torsion sheaf.
		
	For the ``if'' of statement \eqref{item:duality-tor-torsion-free}, it is enough to prove that if $\calE \in \catCoh_{\mathsf{tor}}(S)$, then $\D(\calE)$ belongs to $\scrA_{\mathsf{tor}}$. Since $(\catCoh_0(S), \catCoh_{\geqslant 1}(S)$ is a torsion pair on $\catCoh(S)$, there is a canonical fiber sequence attached to $\calE$:
	\begin{align}
		\calE_0 \longrightarrow \calE \longrightarrow \calE_{\geqslant 1} \ ,  
	\end{align}
	where $\calE_0$ is the maximal zero-dimensional subsheaf of $\calE$, and $\calE_{\geqslant 1} \in \catCoh_{\geqslant 1}(S)$. Since $\calE \in \catCoh_{\mathsf{tor}}(S)$, $\calE_{\geqslant 1}$ is purely $1$-dimensional. Applying $\D(-)$ we obtain
	\begin{align}
		\D(\calE_{\geqslant 1}) \longrightarrow \D(\calE) \longrightarrow \D(\calE_0) \ . 
	\end{align}
	Then $\D(\calE_0) \simeq \calExt^2_S(\calE_0,\omega_S)$ is zero-dimensional, and $\D(\calE_{\geqslant 1}) \simeq \calExt^1_S(\calE_{\geqslant 1},\omega_S)[1]$ is purely $1$-dimensional. Thus, $\D(\calE) \in \scrA_{\mathsf{tor}}$.
	\end{proof}

We provide the following two structural results on $\scrA_{\mathsf{t.f.}}$ and $\scrA_{\mathsf{tor}}$.
\begin{corollary}\label{cor:duality-torsion-free}
	Let $E \in \scrA_{\mathsf{t.f.}}$.
	Then:
	\begin{enumerate}\itemsep=0.2cm
		
		\item  \label{item:duality-torsion-free-1} one has
		\begin{align}
			\mathsf{rk}(E) = - \mathsf{rk}(\D(E[1])) \ ;
		\end{align}
		
		\item \label{item:duality-torsion-free-2} the sheaf $\calH^{-1}(E)$ is locally free of rank $-\mathsf{rk}(E)$;
		
		\item  \label{item:duality-torsion-free-3} the sequence
		\begin{align}
			0\longrightarrow \D(E[1])\longrightarrow \D(\calH^{-1}(E))[-2]  \longrightarrow \D(\calH^0(E)) \longrightarrow 0
		\end{align}		
		is short exact and corresponds to the one associated to the inclusion of $\D(E[1])$ into its underived double dual.
	\end{enumerate}
\end{corollary}

\begin{proof}
	To begin with, Theorem~\ref{thm:duality-tor-torsion-free} allows to write $E \simeq \D(\calE[1])$, where $\calE \in \catCoh_{\mathsf{t.f.}}(S)$. Let $\calW$ be the underived double dual of $\calE$. Since $\calE$ is torsion free, $\calW$ is locally free, the canonical map $\calE \to \calW$ is injective and its cokernel $\calT$ is zero-dimensional. They fit into the short exact sequence
	\begin{align}\label{eq:torsion-free-double-dual}
		0\longrightarrow\calE\longrightarrow \calW\longrightarrow \calT\longrightarrow 0 \ .
	\end{align}
	Rotating, we obtain a fiber sequence
	\begin{align}
		\calT \longrightarrow \calE[1] \longrightarrow \calW[1] \ . 
	\end{align}
	Applying $\D(-)$, we obtain
	\begin{align}
		\D(\calW[1]) \longrightarrow E \longrightarrow \D(\calT) \ . 
	\end{align}
	Passing to the long exact sequence and using the fact that $\D(\calT)$ is concentrated in degree zero, we obtain
	\begin{align}
		\calH^{-1}(E) \simeq \calH^{-1}(\D(\calW[1]))\simeq \calW^\vee\otimes_{\scrO_S} \omega_S \ ,
	\end{align}
	and thus we deduce that $\calH^{-1}(E)$ is locally free. 
	
	Consider now the canonical fiber sequence
	\begin{align}
		\calH^{-1}(E)[1] \longrightarrow E \longrightarrow \calH^0(E) \ . 
	\end{align}
	Applying $\D(-)$ and rotating, we obtain the fiber sequence
	\begin{align}
		\calE \longrightarrow \D(\calH^{-1}(E))[-2] \longrightarrow \D(\calH^0(E))\ .
	\end{align}	
	Since $\D(\calH^{-1}(E))[-2]\simeq \calW$, and $\D(\calH^0(E))$ has amplitude $[0,\, 0]$, the above fiber sequence corresponds to the short exact sequence \eqref{eq:torsion-free-double-dual}. Hence
	\begin{align}
		\calT\simeq \D(\calH^{0}(E)) \simeq \calExt^2_S(\calH^{0}(E),\omega_S)\ .
	\end{align}
	Finally, since $\calT$ is zero-dimensional, it follows that
	\begin{align}
		\mathsf{rk}(\calE) = \mathsf{rk}(\D(\calH^{-1}(E))[-2]) = \mathsf{rk}(\calW)= \mathsf{rk}(\calH^{-1}(E)) = - \mathsf{rk}(E) \ . 
	\end{align}	
	The proof is therefore achieved.
\end{proof}

\begin{corollary}\label{cor:duality-torsion}
	\hfill
	\begin{enumerate}\itemsep=0.2cm
		\item Under the equivalence $\D \colon \catCoh_{\mathsf{tor}}(S) \simeq \scrA_{\mathsf{tor}}\op$, pure $1$-dimensional sheaves correspond to objects $F \in \scrA_{\mathsf{tor}}$ satisfying $\calH^0(F) = 0$.
		
		\item Under the equivalence $\D\colon \catCoh_{\mathsf{t.f.}}(S)[1] \simeq \scrA_{\mathsf{t.f.}}\op$, $\catVect(S)[1]$ corresponds to the full subcategory of $\scrA_{\mathsf{t.f.}}\op$ spanned by objects $E$ satisfying $\calH^0(E) = 0$. Here, $\catVect(S)\subset \catCoh(S)$ denotes the subcategory of locally free sheaves on $S$.
	\end{enumerate}
\end{corollary}

\begin{proof}
	Using \cite[Proposition~1.1.6]{HL_Moduli}, the first statement follows from Lemma~\ref{lem:A_1-1}, while the second one follows from Corollary~\ref{cor:duality-torsion-free}--\eqref{item:duality-torsion-free-2}.
\end{proof}

\section{COHA and CatHA of a surface, and their representation via torsion-free sheaves}\label{sec:action-torsion}

In this section, we introduce COHAs and CatHAs associated to \textit{torsion sheaves} on a smooth projective surface $S$ and a (categorical) representations associated to \textit{torsion-free sheaves} on $S$. In particular, we provide a categorification of known constructions of cohomological and K-theoretical Hall algebras of $S$.

\subsection{Algebra and representations}\label{subsec:torsion-pair}

Let $S$ be a smooth projective irreducible complex surface. Let $\tau_{\mathsf{std}}$ denotes the standard $t$-structure on $\catPerf(S)$. Furthermore, consider the torsion pair $(\scrT, \scrF)$ in $\catPerf(S)^\heartsuit=\catCoh(S)$ given by
\begin{align}
	\scrT\coloneqq\catCoh_{\mathsf{tor}}(S)\quad \text{and} \quad 
	\scrF\coloneqq \catCoh_{\mathsf{t.f.}}(S) \ .
\end{align}
As seen in \S\ref{subsec:perverse-t-structure-surfaces}, the $t$-structure obtained by tilting with respect to this torsion pair is $\tau_\scrB$.

Set $\scrC_S\coloneqq \catQCoh(S)=\mathsf{Ind}(\catPerf(S))$. The $t$-structure $\tau_{\mathsf{std}}$ and the torsion pair $(\scrT, \scrF)$ canonically lift to $\scrC_S$ by Proposition~\ref{prop:t-structure-ind-completion} and Corollary~\ref{cor:completion-torsion-pair}, respectively: we will keep the same notation also for these lifts.

Denote by $\bfCoh(S)=\bfCoh(\scrC_S, \tau_{\mathsf{std}})$ the derived moduli stack of coherent sheaves on $S$. Construction~\ref{construction:torsion_pair_substacks} yields two substacks
\begin{align}
	\bfCoh_{\mathsf{tor}}(S) \coloneqq \bfCoh_{\scrT}(\scrC_S, \tau_{\mathsf{std}}) \qquad \text{and} \qquad \bfCoh_{\mathsf{t.f.}}(S) \coloneqq \bfCoh_\scrF(\scrC_S, \tau_{\mathsf{std}}) 
\end{align}
of $\bfCoh(S)$ respectively parametrizing torsion and torsion-free sheaves on $S$. By \cite[Example~A.4]{AB_Moduli}, the torsion pair $(\scrT, \scrF)$ is open. By Proposition~\ref{prop:openness-tilted}, the tilted $t$-structure $\tau_\scrB$ universally satisfies openness of flatness and, by Proposition~\ref{prop:openness}, $\bfCoh(S,\tau_\scrB) \coloneqq \bfCoh(\scrC_S,\tau_\scrB)$ is a geometric derived stack locally of finite presentation over $\C$, which is open inside $\bfPerf(S)$.

Furthermore, Construction~\ref{construction:torsion_pair_substacks} immediately implies that the natural morphisms
\begin{align}
	\bfCoh_{\mathsf{tor}}(S) \longrightarrow \bfCoh(S,\tau_\scrB) \quad \text{and} \quad [1] \colon \bfCoh_{\mathsf{t.f.}}(S) \longrightarrow \bfCoh(S,\tau_\scrB) 
\end{align}
are representable by open immersions.

Recall that one has a decomposition of $\bfCoh(S)$ into open and closed substacks
\begin{align}
	\bfCoh(S) = \bigsqcup_{(r, \beta, n)\,\in\, \Z \oplus \mathsf{NS}(S) \oplus \Z} \bfCoh(S; r, \beta, n) \ , 
\end{align}
where each component corresponds to coherent sheaves $\calE$ on $S$ with rank $r$, first Chern class $\beta$, and Euler characteristic $n$. The stack $\bfCoh(S, \tau_\scrB)$ admits a similar decomposition. Thus, $\bfCoh_{\mathsf{tor}}(S)$ is equivalent to the connected component corresponding to zero rank objects of both $\bfCoh(S)$ and $\bfCoh(S, \tau_\scrB)$, and hence it is open and closed inside both $\bfCoh(S)$ and $\bfCoh(S, \tau_\scrB)$.

The main result of this section is the following:
\begin{theorem}\label{thm:action-torsion}
	The stable pro-$\infty$-category $\catCohb_{\mathsf{pro}}( \bfCoh_{\mathsf{tor}}(S) )$ has a $\E_1$-monoidal structure, and the stable pro-$\infty$-category $\catCohb_{\mathsf{pro}}( \bfCoh_{\mathsf{t.f.}}(S) )$ has the structure of a left and a right categorical module over $\catCohb_{\mathsf{pro}}( \bfCoh_{\mathsf{tor}}(S) )$. 
	
	Similarly, let $\bfD^\ast$ be a motivic formalism, and fix $\calA \in \CAlg(\bfD^\ast(\Spec(k)))$ and $\Gamma \subseteq \Pic(\bfD^\ast(\Spec(k)))$ such that Assumption~\ref{assumption:motivic_formalism} is satisfied. Then, the topological vector space $\HBMDGamma_0(\bfCoh_{\mathsf{tor}}(S);\calA)$ has the structure of a unital associative algebra. In particular,
	\begin{align}
		G_0( \bfCoh_{\mathsf{tor}}(S) )\quad \text{and} \quad \HBM_\ast( \bfCoh_{\mathsf{tor}}(S) )
	\end{align}
	have the structures of unital associative algebras, and the topological vector space $\HBMDGamma_0(\bfCoh_{\mathsf{t.f.}}(S);\calA)$ has both the structure of a left and a right $\HBMDGamma_0(\bfCoh_{\mathsf{tor}}(S);\calA)$-module. In particular,
	\begin{align}
		G_0( \bfCoh_{\mathsf{t.f.}}(S) )\quad \text{and} \quad \HBM_\ast( \bfCoh_{\mathsf{t.f.}}(S) )			
	\end{align}
	have both the structures of a left and a right $G_0( \bfCoh_{\mathsf{tor}}(S) )$-module and $\HBM_\ast( \bfCoh_{\mathsf{tor}}(S) )$-module, respectively.
\end{theorem}

\begin{proof}
	It is enough to apply Theorem~\ref{thm:left-right-action} to the pair $(\scrC_S, \tau_{\mathsf{std}})$ and the torsion pair $(\scrT, \scrF)$. First, $(\scrC_S, \tau_{\mathsf{std}})$ satisfies Assumption~\ref{assumption:existence-COHA}. Indeed, the category $\scrC_S$ is smooth and proper and it admits the Serre functor $\sfS(-)=(-)\otimes_{\scrO_S} \omega_S[2]$. We already argued that the torsion pair $(\scrT, \scrF)$ is open, moreover $\scrT$ is a Serre subcategory of $\catCoh(S)$ (cf.\ Remark~\ref{rem:tor_serre}). Therefore, conditions~\eqref{item:torsion-pair-left-right-action-3} and \eqref{item:torsion-pair-left-right-action-4} of Theorem~\ref{thm:left-right-action} are verified. Moreover, conditions~\eqref{item:torsion-pair-left-right-action-1} and \eqref{item:torsion-pair-left-right-action-2} of Theorem~\ref{thm:left-right-action} are trivially satisfied.
	
	The map $\partial_1\colon \calS_2\bfCoh(S)\to \bfCoh(S)$ is locally rpas since any connected component of the Quot scheme for the standard $t$-structure on $\scrC_S$ is proper (cf.\ \cite[Theorem~2.2.4]{HL_Moduli}). By Lemma~\ref{lem:torsion_pairs_key_squares}, condition~\eqref{item:left-right-action-1} of Theorem~\ref{thm:left-right-action} is satisfied.
	
	Finally, condition~\eqref{item:left-right-action-2} of Theorem~\ref{thm:left-right-action} is also satisfied. Indeed, the map
	\begin{align}\label{eq:properness-standard}
		\varpi_0 \colon \calS_1^\ell\bfFlagCoh^{(1),\dagger}_{\mathsf{tor},\mathsf{t.f.}}(S) \longrightarrow \bfCoh_{\mathsf{t.f.}}(S)
	\end{align}
	is locally rpas because again the connected components of the Quot scheme for the standard $t$-structure on $\scrC_S$ are proper and $\bfCoh_{\mathsf{tor}}(S)$ is closed in $\bfCoh(S)$, while the map
	\begin{align}\label{eq:properness-scrB}
		\varpi_1 \colon \calS_1^r\bfFlagCoh^{(1),\dagger}_{\mathsf{t.f.}[1],\mathsf{tor}}(S,\tau_\scrB) \longrightarrow \bfCoh_{\mathsf{t.f.}[1]}(S, \tau_\scrB)\simeq \bfCoh_{\mathsf{t.f.}}(S)
	\end{align}
	is locally rpas, since the connected components of the Quot scheme parametrizing quotients of an object of $\scrF[1]$ in the tilted heart are proper as explained in Proposition~\ref{prop:Kollar_husks} and $\bfCoh_{\mathsf{tor}}(S)$ is closed in $\bfCoh(S, \tau_\scrB)$. 
\end{proof}

\begin{remark}\label{rem:action-torsion}
	We have an decomposition of the stack $\bfCoh_{\mathsf{t.f.}}(S)$ of torsion-free sheaves on $S$:
	\begin{align}
		\bfCoh_{\mathsf{t.f.}}(S) = \bigsqcup_{(r, \beta,n)\in \Z\oplus N_{\leqslant 1}(S)} \bfCoh_{\mathsf{t.f.}}(S; r, \beta, n) \ . 
	\end{align}
	Set 
	\begin{align}
		\bfCoh_{\mathsf{t.f.}}(S; r) \coloneqq \bigsqcup_{(\beta,n)\in N_{\leqslant 1}(S)} \bfCoh_{\mathsf{t.f.}}(S; r, \beta, n)\ .
	\end{align}
	Then, we have a result analogous to Theorem~\ref{thm:action-torsion} with $\bfCoh_{\mathsf{t.f.}}(S; r)$ instead of $\bfCoh_{\mathsf{t.f.}}(S)$.  
\end{remark}

Let $\scrH_{\mathsf{tor}, \mathsf{t.f}}$, $\scrU_{\mathsf{tor}, \mathsf{t.f}}$, and $\scrY_{\mathsf{tor}, \mathsf{t.f}}$ be the algebras of Definition~\ref{def:Yangian-torsion}  associated to the pair $(\bfCoh_{\mathsf{tor}}(S),$ $\bfCoh_{\mathsf{t.f.}}(S))$. Similarly, let $\scrH_{\mathsf{tor}, \mathsf{t.f}}(r)$, $\scrU_{\mathsf{tor}, \mathsf{t.f}}(r)$, and $\scrY_{\mathsf{tor}, \mathsf{t.f}}(r)$ be the algebras associated to $(\bfCoh_{\mathsf{tor}}(S),$ $\bfCoh_{\mathsf{t.f.}}(S; r))$. The decomposition of $\bfCoh_{\mathsf{t.f.}}(S)$ with respect to the rank induces the decompositions:
\begin{align}
	\scrH_{\mathsf{tor}, \mathsf{t.f}} \simeq \bigoplus_r \scrH_{\mathsf{tor}, \mathsf{t.f}}(r)\; , \ \scrU_{\mathsf{tor}, \mathsf{t.f}} \simeq \bigoplus_r \scrU_{\mathsf{tor}, \mathsf{t.f}}(r)\quad \text{and}\quad \scrY_{\mathsf{tor}, \mathsf{t.f}} \simeq \bigoplus_r \scrY_{\mathsf{tor}, \mathsf{t.f}}(r)\ .
\end{align}

\begin{remark}
	We claim that the algebras above should be isomorphic for different nonzero $r, r'\in \N$. This is what happens in the geometric case described in the next section.
\end{remark}

\subsubsection{CatHA and COHA of zero-dimensional sheaves and their representations}

Let $\bfCoh_{0\textrm{-}\!\dim}(S)$ be the derived moduli stack of zero-dimensional torsion sheaves on $S$. Note that the square \eqref{eq:multiplication-T} is a pullback with $\bfT\coloneqq \bfCoh_{0\textrm{-}\!\dim}(S)$, since the subcategory $\catCoh_{0\textrm{-}\!\dim}(S)$ is Serre and $\bfCoh_{0\textrm{-}\!\dim}(S)$ is open and closed in $\bfCoh_{\mathsf{tor}}(S)$. Thus, we can obtain analogs of Theorem~\ref{thm:action-torsion} and Remark~\ref{rem:action-torsion} for the pairs $(\bfT\coloneqq \bfCoh_{0\textrm{-}\!\dim}(S), \bfF\coloneqq \bfCoh_{\mathsf{t.f.}}(S))$ and $(\bfT\coloneqq \bfCoh_{0\textrm{-}\!\dim}(S), \bfF\coloneqq \bfCoh_{\mathsf{t.f.}}(S; r))$, by applying Corollary~\ref{cor:induced_COHA_open}. Correspondingly, we can introduce the algebras of Definition~\ref{def:Yangian-torsion} associated to these pairs and have the following decompositions:
\begin{align}
	\scrH_{0\textrm{-}\!\dim, \mathsf{t.f}} \simeq \bigoplus_{r, c_1} \scrH_{0\textrm{-}\!\dim, \mathsf{t.f}}(r, c_1)\; , \qquad \scrU_{0\textrm{-}\!\dim, \mathsf{t.f}} \simeq \bigoplus_{r, c_1} \scrU_{0\textrm{-}\!\dim, \mathsf{t.f}}(r, c_1)
\end{align}
and
\begin{align}
	\scrY_{0\textrm{-}\!\dim, \mathsf{t.f}} \simeq \bigoplus_{r, c_1} \scrY_{0\textrm{-}\!\dim, \mathsf{t.f}}(r, c_1)\ ,
\end{align}
where the sums are over all possible ranks $r\in \Z_{\geq 1}$ and first Chern classes $c_1\in \mathsf{NS}(S)$ of torsion-free coherent sheaves on $S$.

\subsection{Vanishing of categorical commutators}

We now apply the study of the categorical commutators of \S\ref{sec:categorical_commutators} to provide a geometric criterion guaranteeing that two operator commute.
Fix two closed subschemes $Z_1$ and $Z_2$ of $S$ and let
\begin{align}
	\iota_1 \colon Z_1 \longrightarrow S \qquad \text{and} \qquad \iota_2 \colon Z_2 \longrightarrow S 
\end{align}
be the natural inclusions. There are induced morphisms
\begin{align}
	j_{1} \colon \bfCoh(Z_1) \longrightarrow \bfCoh_{\mathsf{tor}}(S) \qquad \text{and} \qquad j_{2} \colon \bfCoh(Z_2) \longrightarrow \bfCoh_{\mathsf{tor}}(S) \ , 
\end{align}
given by the pushforward along $\iota_1$ and $\iota_2$. It follows from that the morphisms $j_1$ and $j_2$ are closed immersions. In particular, for $i = 1,2$ we obtain functors
\begin{align}
	j_{i,\,\ast} \colon \catCohb_{\mathsf{pro}}(\bfCoh(Z_i)) \longrightarrow \catCohb_{\mathsf{pro}}(\bfCoh_{\mathsf{tor}}(S)) \ . 
\end{align}
Assume now that $Z_1 \cap Z_2 = \emptyset$. This guarantees that, if $\calF, \calF'$ are coherent sheaves on $Z_1$ and $Z_2$ respectively, and we are given surjections (resp.\ monomorphisms)
\begin{align}
	\iota_{1,\,\ast}(\calF) \twoheadrightarrow \calM \qquad &\text{and} \qquad \iota_{2,\,\ast}(\calF') \twoheadrightarrow \calM  \\
	\text{(resp.\ } \calM \hookrightarrow \iota_{1,\,\ast}(\calF) \qquad &\text{and}  \qquad \calM\hookrightarrow \iota_{2,\,\ast}(\calF')  \text{ )}
\end{align}
for some coherent sheaf $\calM$ on $S$, then $\mathsf{supp}(\calM) \subseteq Z_1 \cap Z_2$, and hence necessarily $\calM = 0$. Thus, the assumptions of Corollary~\ref{cor:support_commutation_criterion} are satisfied. Writing $\tau_{-\scrB} \coloneqq \tau_{\scrB}[-1]$ and setting $\bfCoh(S,\tau_{-\scrB}) \coloneqq \bfCoh(\scrC_S, \tau_{-\scrB})$, we obtain:

\begin{corollary}\label{cor:commutators_geometric_setting}
	Let $Z_1$ and $Z_2$ of $S$ be two disjoint closed subschemes. Denote by $\iota_i$ the corresponding inclusions of subschemes and by $j_i$ the inclusions of the corresponding moduli stacks of sheaves, for $i=1,2$. For $\calG_1 \in \catCohb_{\mathsf{pro}}(\bfCoh(Z_1))$ and $\calG_2 \in \catCohb_{\mathsf{pro}}(\bfCoh(Z_2))$, we have
	\begin{align}
		(j_{2,\,\ast}\calG_2 \oast \calE) \oast \Phi(j_{1,\,\ast}\calG_1) \simeq j_{2,\,\ast}\calG_2 \oast (\calE \oast \Phi(j_{1,\,\ast}\calG_1))
	\end{align}
	for any $\calE\in \catCohb_{\mathsf{pro}}( \bfCoh_{\mathsf{t.f.}}(S) )$. Therefore the operators induced by $\calG_1$ and $\calG_2$ on $\catCohb_{\mathsf{pro}}( \bfCoh_{\mathsf{t.f.}}(S) )$ commute. Here, we denote by $\Phi\colon  \catCohb_{\mathsf{pro}}( \bfCoh_{\mathsf{tor}}(S) ) \longrightarrow \catCohb_{\mathsf{pro}}( \bfCoh_{\scrT[-1]}(S, \tau_{-\scrB}) )$ induced by the equivalence $[-1] \colon \bfCoh_{\mathsf{tor}}(S) \to \bfCoh_{\scrT[-1]}(S, \tau_{-\scrB})$.
	
	Similar statements hold in motivic Borel-Moore homology.
\end{corollary}

\subsection{Decategorification and comparison with existing results}

In this section, we compare our construction with Negut's construction  \cite{Negut_Shuffle, Negut_Categorification}  of representations of the K-theoretical Hall algebra of zero-dimensional sheaves on $S$ and its categorification, and DeHority's construction \cite{DeHority_KM} of representations of a certain subalgebra of the cohomological Hall algebra of torsion sheaves on $S$.

\subsubsection{CatHA and COHA of zero-dimensional sheaves and moduli spaces of Gieseker-stable sheaves}\label{subsubsec:negut}

In this section we compare our approach to that of Negu\c{t} \cite{Negut_Shuffle, Negut_Categorification}. In \textit{loc. cit.}, Negu\c{t} introduced operators acting on the K-theory of moduli spaces of Gieseker-stable sheaves on $S$ arising from certain Hecke correspondences associated to zero-dimensional sheaves and proved that they generate an algebra satisfying relations resembling those of the elliptic Hall algebra. For such a reason, we call it the \textit{elliptic Hall algebra} of $S$.

Let $H$ be an ample divisor and let $\bfCoh_{\mathsf{t.f.}}^{H\textrm{-}(\mathsf{s})\mathsf{s}}(S)$ be the derived moduli stack of $H$-Gieseker-(semi-)stable torsion-free sheaves on $S$, which is an open substack of $\bfCoh_{\mathsf{t.f.}}(S)$. We have a decomposition into open and closed substacks
\begin{align}
	\bfCoh_{\mathsf{t.f.}}^{H\textrm{-}(\mathsf{s})\mathsf{s}}(S) \coloneqq \bigsqcup_{r, c_1, c_2} \bfCoh_{\mathsf{t.f.}}^{H\textrm{-}(\mathsf{s})\mathsf{s}}(S; r, c_1, c_2)
\end{align}
with respect to the rank $r$ and the Chern classs $c_1, c_2$ of stable torsion free sheaves. Fix $r, c_1$ and define
\begin{align}
	\bfCoh_{\mathsf{t.f.}}^{H\textrm{-}(\mathsf{s})\mathsf{s}}(S; r, c_1) \coloneqq \bigsqcup_{c_2} \bfCoh_{\mathsf{t.f.}}^{H\textrm{-}(\mathsf{s})\mathsf{s}}(S; r, c_1, c_2)\ .
\end{align}

Let $\calM_S(r, c_1)$ be the moduli space of the classical truncation $\trunc{\bfCoh_{\mathsf{t.f.}}^{H\textrm{-}\mathsf{s}}(S; r, c_1)}$ of $\bfCoh_{\mathsf{t.f.}}^{H\textrm{-}\mathsf{s}}(S; r, c_1)$. Then, $\trunc{\bfCoh_{\mathsf{t.f.}}^{H\textrm{-}\mathsf{s}}(S; r, c_1)}$ is a $\C^\ast$-gerbe over $\calM_S(r, c_1)$. 

Recall the construction of the \textit{Hilbert stack $\mathbf{Hilb}(S)$ of points on $S$} following \cite[\S7.1]{MMSV}. Let $\bfCoh_{\mathsf{t.f.}}(S; 1)$ be the the derived stack of torsion-free coherent sheaves on $S$ of rank one and let $\bfPic(S)$ be the derived Picard stack of $S$. Note that $\bfCoh_{\mathsf{t.f.}}(S; 1)\simeq \bfCoh_{\mathsf{t.f.}}^{H\textrm{-}\mathsf{s}}(S; 1)$ since any rank one torsion-free sheaf is automatically Gieseker-stable, with respect to any ample divisor $H$. There exists a determinant map $\bfCoh_{\mathsf{t.f.}}(S; 1)\to \bfPic(S)$ (cf.\ \cite{STV}), hence $\mathbf{Hilb}(S)$ is by definition the pullback
\begin{align}
	\begin{tikzcd}[ampersand replacement=\&]
		\mathbf{Hilb}(S) \ar{r} \ar{d}\& \bfCoh_{\mathsf{t.f.}}^{H\textrm{-}\mathsf{s}}(S; 1) \ar{d}\\
		B\G_m \ar{r}\& \bfPic(S)
	\end{tikzcd}\ ,
\end{align}
where the map $B\G_m\to \bfPic(S)$ is the closed embedding of the substack parametrizing trivial invertible sheaves. Thus, there is an induced closed embedding $\mathbf{Hilb}(S)\to \bfCoh_{\mathsf{t.f.}}^{H\textrm{-}\mathsf{s}}(S; 1, 0)$, where the latter is the derived stack parametrizing $H$-stable torsion-free sheaves on $S$ of rank one and vanishing first Chern class.

Now, let us fix $r$ and $c_1\in \mathsf{NS}(S)$ such that $\mathsf{g.c.d.}(r, c_1\cdot H)=1$ (cf.\ \cite[Assumption~\textit{A}]{Negut_Shuffle}). This implies that all semistable sheaves are automatically stable. We can apply Corollary~\ref{cor:induced_COHA_open} to $\bfT\coloneqq \bfCoh_{0\textrm{-}\!\dim}(S)$ and $\bfF$ equals either $\bfCoh_{\mathsf{t.f.}}^{H\textrm{-}\mathsf{s}}(S; r, c_1)$ or $\mathbf{Hilb}(S)$. Therefore, we obtain the following result.
\begin{theorem}\label{thm:action-zero-dim-stable}
	The stable pro-$\infty$-category $\catCohb_{\mathsf{pro}}( \bfCoh_{\mathsf{t.f.}}^{H\textrm{-}\mathsf{s}}(S; r, c_1) )$ has the structure of a left and a right categorical module over $\catCohb_{\mathsf{pro}}( \bfCoh_{0\textrm{-}\!\dim}(S) )$. 
	
	Similarly, let $\bfD^\ast$ be a motivic formalism, and fix $\calA \in \CAlg(\bfD^\ast(\Spec(k)))$ and $\Gamma \subseteq \Pic(\bfD^\ast(\Spec(k)))$ such that Assumption~\ref{assumption:motivic_formalism} is satisfied. Then, the topological vector space $\HBMDGamma_0(\bfCoh_{\mathsf{t.f.}}^{H\textrm{-}\mathsf{s}}(S; r, c_1);\calA)$ has both the structure of a left and a right $\HBMDGamma_0(\bfCoh_{0\textrm{-}\!\dim}(S);\calA)$-module. In particular,
	\begin{align}
		G_0( \bfCoh_{\mathsf{t.f.}}^{H\textrm{-}\mathsf{s}}(S; r, c_1) )\quad \text{and} \quad \HBM_\ast( \bfCoh_{\mathsf{t.f.}}^{H\textrm{-}\mathsf{s}}(S; r, c_1) )			
	\end{align}
	have both the structures of a left and a right $G_0( \bfCoh_{0\textrm{-}\!\dim}(S) )$-module and $\HBM_\ast( \bfCoh_{0\textrm{-}\!\dim}(S) )$-module, respectively.
	
	Similar statements hold after replacing $\bfCoh_{\mathsf{t.f.}}^{H\textrm{-}\mathsf{s}}(S; r, c_1)$ with $\mathbf{Hilb}(S)$.
\end{theorem}

\begin{proof}
	We have already discussed that the square \eqref{eq:multiplication-T} is a pullback with $\bfT\coloneqq \bfCoh_{0\textrm{-}\!\dim}(S)$, since the subcategory $\catCoh_{0\textrm{-}\!\dim}(S)$ is Serre, and $\bfCoh_{0\textrm{-}\!\dim}(S)$ is open and closed in $\bfCoh_{\mathsf{tor}}(S)$. Moreover, thanks to \cite[Proposition~5.5]{Negut_Shuffle}, the stacks $\bfCoh_{\mathsf{t.f.}}^{H\textrm{-}\mathsf{s}}(S; r, c_1)$ and $\mathbf{Hilb}(S)$ are two-sided Hecke patterns for $\bfCoh_{0\textrm{-}\!\dim}(S)$ (with respect to the standard heart). Thus, Corollary~\ref{cor:induced_COHA_open} applies to $\alpha\colon\bfCoh_{0\textrm{-}\!\dim}(S)\to \bfCoh_{\mathsf{tor}}(S)(\scrC, \tau)$ and $\mu\colon \bfCoh_{\mathsf{t.f.}}^{H\textrm{-}\mathsf{s}}(S; r, c_1)\to \bfCoh_{\mathsf{t.f.}}(S)$, and the first assertion follows.
	
	Now, let $\mu\colon \mathbf{Hilb}(S)\to \bfCoh_{\mathsf{t.f.}}(S)$ be the composition of the closed embedding $\mathbf{Hilb}(S)\to \bfCoh_{\mathsf{t.f.}}^{H\textrm{-}\mathsf{s}}(S; 1, 0)$ with the open embedding $\bfCoh_{\mathsf{t.f.}}^{H\textrm{-}\mathsf{s}}(S; 1, 0)\to \bfCoh_{\mathsf{t.f.}}(S)$. Corollary~\ref{cor:induced_COHA_open} still applies to $\alpha\colon\bfCoh_{0\textrm{-}\!\dim}(S)\to \bfCoh_{\mathsf{tor}}(S)(\scrC, \tau)$ and $\mu\colon \mathbf{Hilb}(S)\to \bfCoh_{\mathsf{t.f.}}(S)$ since $\mathbf{Hilb}(S)$ is a two-sided Hecke pattern for $\bfCoh_{0\textrm{-}\!\dim}(S)$, so it satisfies condition~\eqref{item:two-sided} in \textit{loc.\ cit.} Thus, the second statement follows as well.
\end{proof}

\begin{remark}
	Note that $\mathbf{Hilb}(S)$ coincides with its truncation and $\mathbf{Hilb}(S)\simeq \mathsf{Hilb}(S)\times B\G_m$, where $\mathsf{Hilb}(S)$ is the Hilbert scheme of points on $S$ (cf.\ \cite[Lemma~7.1]{MMSV}.
	
	Similarly under \cite[Assumption~\textit{S}]{Negut_Shuffle}, $\bfCoh_{\mathsf{t.f.}}^{H\textrm{-}\mathsf{s}}(S; r, c_1)$ coincides with its classical truncation. Let $\calM_S(r, c_1)$ be the fine moduli space of $\bfCoh_{\mathsf{t.f.}}^{H\textrm{-}\mathsf{s}}(S; r, c_1)$. Then, there exists a universal sheaf on $\calM_S(r, c_1)\times S$. Thus, the map $\bfCoh_{\mathsf{t.f.}}^{H\textrm{-}\mathsf{s}}(S; r, c_1)\to \calM_S(r, c_1)$ has a section. By \cite[Lemma~3.10]{Hein_Moduli}, one has $\bfCoh_{\mathsf{t.f.}}^{H\textrm{-}\mathsf{s}}(S; r, c_1)\simeq \calM_S(r, c_1) \times B\G_m$. 
\end{remark}

The previous remark yields the following result.
\begin{corollary}
	The stable pro-$\infty$-category 
	\begin{align}
		\bigoplus_{w\in \Z} \catCohb_{\mathsf{pro}}( \mathsf{Hilb}(S) )_w
	\end{align}
	has the structure of a left and a right categorical module over $\catCohb_{\mathsf{pro}}( \bfCoh_{0\textrm{-}\!\dim}(S) )$. Here, $\catCohb_{\mathsf{pro}}( \mathsf{Hilb}(S) )_w$ is the weight $w$ part of $\catCohb_{\mathsf{pro}}( \mathbf{Hilb}(S) )$.
	
	An analogous statement holds for motivic Borel–Moore homology, as well as after replacing $\mathsf{Hilb}(S)$ with $\calM_S(r, c_1)$.
\end{corollary}

Thanks to the theorem above, we can define the corresponding algebras following Definition~\ref{def:Yangian-torsion}: we denote them by  
\begin{align}
	\scrH_{0\textrm{-}\!\dim, \mathsf{st}}(r, c_1)\; , \ \scrU_{0\textrm{-}\!\dim, \mathsf{st}}(r, c_1)\quad \text{and} \quad \scrY_{0\textrm{-}\!\dim, \mathsf{st}}(r, c_1)\ .
\end{align}

\begin{remark}
	In the rank one and $c_1=0$ case, we can a priori define two families of algebras: the first one obtained by considering $\bfCoh_{\mathsf{t.f.}}^{H\textrm{-}\mathsf{s}}(S; 1, 0)$ to construct the representations, while the second one obtained by considering $\mathbf{Hilb}(S)$. We shall use the same notation for both of them when this does not create any confusion.
\end{remark}

\begin{remark}
	Since $\bfCoh_{\mathsf{t.f.}}^{H\textrm{-}\mathsf{(s)s}}(S; r, c_1)$ is open in $\bfCoh_{\mathsf{t.f.}}(S; r, c_1)$, we have corresponding pullbacks in Borel-Moore homology, K-theory, and bounded derived category. They induce a functor of monoidal categories and homomorphisms of associative algebras, respectively:
	\begin{multline}
		\scrH_{0\textrm{-}\!\dim, \mathsf{t.f.}}(r, c_1) \longrightarrow \scrH_{0\textrm{-}\!\dim, \mathsf{st}}(r, c_1)\; , \ \scrU_{0\textrm{-}\!\dim, \mathsf{t.f.}}(r, c_1) \longrightarrow \scrU_{0\textrm{-}\!\dim, \mathsf{st}}(r, c_1)\; ,\\
		\text{and} \quad \scrY_{0\textrm{-}\!\dim, \mathsf{t.f.}}(r, c_1) \longrightarrow \scrY_{0\textrm{-}\!\dim, \mathsf{st}}(r, c_1)\ .
	\end{multline}
\end{remark}

Now, let us introduce the algebras acting on the K-theory and the bounded derived category of $\calM_S(r, c_1)$, which is related to the elliptic Hall algebra of $S$ and its categorification. Let $\scrU_{\mathsf{pt}, \mathsf{st}}(r, c_1)$ be the subalgebra of $\scrU_{0\textrm{-}\!\dim, \mathsf{st}}(r, c_1)$ generated by $G_0( \bfCoh_{\mathsf{pt}}(S; 1) )$, where $\bfCoh_{\mathsf{pt}}(S; 1)$ is the closed substack of $\bfCoh_{0\textrm{-}\!\dim}(S)$ parametrizing zero-dimensional sheaves on $S$ scheme-theoretically supported at a single point and of length one. Similarly, let $\scrH_{\mathsf{pt}, \mathsf{st}}(r, c_1)$ be the monoidal subcategory of $\scrH_{0\textrm{-}\!\dim, \mathsf{st}}(r, c_1)$ generated by $\catCohb_{\mathsf{pro}}( \bfCoh_{\mathsf{pt}}(S; 1) )$. 

Let us recall \cite[Assumption \textit{S}]{Negut_Shuffle}: the canonical bundle $\calK_S$ of $S$ is either trivial or satisfies $c_1(\calK_S)\cdot H<0$. Under this assumption, $\calM_S(r, c_1, c_2)$ is a smooth projective variety. \cite[Theorem~1.2]{Negut_Shuffle} yields the following result.
\begin{proposition}
	Under Assumptions \textit{A} and \textit{S}, the algebra $\scrU_{\mathsf{pt}, \mathsf{st}}(r, c_1)$ is independent of $r, c_1$ and it coincides with the elliptic Hall algebra of $S$ introduced in \cite{Negut_Shuffle, Negut_Categorification}.
	
	Furthermore, the monoidal category $\scrH_{\mathsf{pt}, \mathsf{st}}(r, c_1)$ contains the categorification of the elliptic Hall algebra of $S$ introduced in \cite{Negut_Categorification}.
\end{proposition}

\subsubsection{CatHA and COHA of one-dimensional sheaves and moduli spaces of rank one stable sheaves}

Let $\bfCoh_{\mathsf{t.f.}}(S;1)$ be the moduli stack of rank one torsion-free sheaves on $S$. By applying Corollary~\ref{cor:induced_COHA_open} to $\bfF\coloneqq \bfCoh_{\mathsf{t.f.}}(S;1)$ and $\bfT\coloneqq  \bfCoh_{\mathsf{tor}}(S)$, we obtain the following:
\begin{corollary}\label{cor:action-torsion-rank-one}
	The stable pro-$\infty$-category $\catCohb_{\mathsf{pro}}( \bfCoh_{\mathsf{t.f.}}(S;1) )$ has the structure of a left and a right categorical module over $\catCohb_{\mathsf{pro}}( \bfCoh_{\mathsf{tor}}(S) )$. 
	
	Similarly, let $\bfD^\ast$ be a motivic formalism, and fix $\calA \in \CAlg(\bfD^\ast(\Spec(k)))$ and $\Gamma \subseteq \Pic(\bfD^\ast(\Spec(k)))$ such that Assumption~\ref{assumption:motivic_formalism} is satisfied. Then, the topological vector space $\HBMDGamma_0(\bfCoh_{\mathsf{t.f.}}(S; 1);\calA)$ has both the structure of a left and a right $\HBMDGamma_0(\bfCoh_{\mathsf{tor}}(S);\calA)$-module. In particular,
	\begin{align}
		G_0( \bfCoh_{\mathsf{t.f.}}(S;1) )\quad \text{and} \quad \HBM_\ast( \bfCoh_{\mathsf{t.f.}}(S;1) )			
	\end{align}
	have both the structures of a left and a right $G_0( \bfCoh_{\mathsf{tor}}(S) )$-module and $\HBM_\ast( \bfCoh_{\mathsf{tor}}(S) )$-module, respectively.
\end{corollary}

Thanks to the above corollary, we can define the algebras of Definition~\ref{def:Yangian-torsion} associated to the pair $(\bfCoh_{\mathsf{tor}}(S), \bfCoh_{\mathsf{t.f.}}(S;1))$. In particular, we denote by $\scrY_{\mathsf{tor}, \mathsf{rk}=1}$ the corresponding Yangian. We shall now show that our construction recovers \cite{DeHority_KM}.

From now on, let $S$ be a K3 surface such that $\mathsf{NS}(S)$ is generated by irreducible $-2$ curves and any pair of irreducible $-2$ curves on $S$ are either disjoint or intersect transversally at a single point. Let $\scrY_{\mathsf{NS(S)}, \mathsf{rk}=1}$ be the subalgebra of $\scrY_{\mathsf{tor}, \mathsf{rk}=1}$ generated by the fundamental classes of the substacks $\bfCoh_{\mathsf{tor}}(S; 0, c_1, c_2)$ parametrizing torsion sheaves with first Chern class given by a $-2$ curve and with arbitrary second Chern class. The main result of \cite{DeHority_KM} yields the following.
\begin{proposition}
	The algebra $\scrY_{\mathsf{NS(S)}, \mathsf{rk}=1}$ is isomorphic to the modified universal enveloping algebra $\widetilde{\sfU}(\widehat{\frakg}(\mathsf{NS}(S)))$ of the affine Lie algebra $\widehat{\frakg}(\mathsf{NS}(S))$, introduced in \cite{DeHority_KM}.
\end{proposition}
We can define the monoidal subcategory $\scrH_{\mathsf{NS(S)}, \mathsf{rk}=1}$ of $\scrH_{\mathsf{tor}, \mathsf{rk}=1}$ generated by the categories $\catCohb_{\mathsf{pro}}( \bfCoh_{\mathsf{tor}}(S; 0, c_1, c_2) )$'s, where $c_1, c_2$ vary as above. Then, $\scrH_{\mathsf{NS(S)}, \mathsf{rk}=1}$ contains a categorification of $\widetilde{\sfU}(\widehat{\frakg}(\mathsf{NS}(S)))$.

\section{Perverse COHA of a surface and its perverse torsion-free representation}\label{sec:action-torsion-perverse}

Let $S$ be a smooth projective irreducible complex surface. Consider the torsion pair $(\scrT, \scrF)$ in $\catCoh(S)$ given by
\begin{align}
	\scrT\coloneqq \catCoh_{0}(S)\quad \text{and} \quad \scrF\coloneqq \catCoh_{\geqslant 1}(S)\ .
\end{align}
As seen in \S\ref{subsec:perverse-t-structure-surfaces}, the $t$-structure obtained by tilting with respect to this torsion pair is $\tau_\scrA$. In \cite[Example~A.4--(1)]{AB_Moduli} it is shown that the torsion pair $(\scrT, \scrF)$ is open. By Proposition~\ref{prop:openness-tilted} the $t$-structure $\tau_\scrA$ universally satisfies openness of flatness. Thus, by Proposition~\ref{prop:openness} the derived moduli stack $\bfCoh(S,\tau_\scrA) \coloneqq \bfCoh(\scrC_S,\tau_{\scrA})$ parametrizing $\tau_\scrA$-flat families of perfect objects is geometric and locally of finite presentation over $\C$.\footnote{Here and in what follows, we do not distinguish between $t$-structures and torsion pairs in $\catPerf(S)$ and the induced ones in $\catQCoh(S)$.}

Let $(\scrA_{\mathsf{tor}},\scrA_{\mathsf{t.f.}})$ be the torsion pair on $\scrA$ introduced in \S\ref{subsec:torsion-pair-perverse}. Applying Construction~\ref{construction:torsion_pair_substacks} to this torsion pair, we obtain two substacks
\begin{align}
	\bfCoh_{\mathsf{tor}}(S, \tau_\scrA) \coloneqq \bfCoh_{\scrA_{\mathsf{tor}}}(\scrC_S,\tau_\scrA) \qquad \text{and} \qquad \bfCoh_{\mathsf{t.f.}}(S, \tau_\scrA) \coloneqq \bfCoh_{\scrA_{\mathsf{t.f.}}}(\scrC_S,\tau_\scrA) 
\end{align}
of $\bfCoh(S, \tau_\scrA)$.

Our goal in this section is to show that Theorem~\ref{thm:left-right-action} can be applied to the current setup. In order to do this, the first thing to do is to check that $(\scrA_{\mathsf{tor}}, \scrA_{\mathsf{t.f.}})$ is an open torsion pair. For the torsion part, this is not difficult to see:
\begin{lemma} \label{lem:A_1}
	$\bfCoh_{\mathsf{tor}}(S, \tau_\scrA)$ is both open and closed inside $\bfCoh(S,\tau_\scrA)$.
\end{lemma}

\begin{proof}
	Consider the group homomorphism
	\begin{align}
		\mathsf{cl}\colon \scrA \longrightarrow \Z \oplus \mathsf{NS}(S) \oplus \Z
	\end{align}
	sending an object $E$ to the pair $(\mathsf{rk}(E), \ch_1(E), \chi(E))$; it induces a decomposition of $\bfCoh(S, \tau_\scrA)$ into open and closed substacks
	\begin{align}\label{eq:decomposition_connected_components_A}
		\bfCoh(S, \tau_\scrA) = \bigsqcup_{v\,\in\, \Z \oplus \mathsf{NS}(S) \oplus \Z} \bfCoh(S, \tau_\scrA; v) \ , 
	\end{align}
	where each component corresponds to objects $E\in \scrA$ such that $\mathsf{cl}(E)=v$. Then Lemma~\ref{lem:A_1-1} shows that
	\begin{align}
		\bfCoh_{\mathsf{tor}}(S, \tau_\scrA) = \bigsqcup_{\mathsf{rk}(v)=0} \bfCoh(S, \tau_\scrA; v) \ .
	\end{align}
	Therefore $\bfCoh_{\mathsf{tor}}(S, \tau_\scrA)$ is both open and closed inside $\bfCoh(S,\tau_\scrA)$.
\end{proof}

On the other hand, to prove that $\bfCoh_{\mathsf{t.f.}}(S, \tau_\scrA)$ is open inside $\bfCoh(S,\tau_\scrA)$ we use the duality result proved in  \S\ref{subsec:torsion-pair-perverse}.
\begin{lemma}\label{lem:A_tf_open}
	The morphism of derived stacks
	\begin{align}
		\bfCoh_{\mathsf{t.f.}}(S, \tau_\scrA) \longrightarrow \bfCoh(S,\tau_{\scrA}) 
	\end{align}
	is representable by open immersions.
\end{lemma}

\begin{proof}
	Since we already argued that the $t$-structure $\tau_{\scrA}$ universally satisfies openness of flatness, we see that $\bfCoh(S,\tau_{\scrA})$ is open in $\bfPerf(S)$. It is therefore enough to prove that the same goes for $\bfCoh_{\mathsf{t.f.}}(S, \tau_\scrA)$. Let $j \colon \bfCoh_{\mathsf{t.f.}}(S, \tau_\scrA) \to \bfPerf(S)$ be the canonical morphism. Since $\D(-)[-1] \colon \catPerf(S)\op \simeq \catPerf(S)$ induces a self equivalence of the derived stack $\bfPerf(S)$ (which we still denote in the same way), we see that it is enough to prove that the map
	\begin{align}
		\D(-)[-1] \circ j \colon \bfCoh_{\mathsf{t.f.}}(S, \tau_\scrA) \longrightarrow \bfPerf(S) 
	\end{align}
	is an open inclusion.
		However, Theorem~\ref{thm:duality-tor-torsion-free} canonically identifies this map with the inclusion of $\bfCoh_{\mathsf{t.f.}}(S)$, which is indeed open in $\bfCoh(S)$, and hence, a posteriori, in $\bfPerf(S)$.
\end{proof}

\begin{theorem}\label{thm:action-perverse-torsion}
	The stable pro-$\infty$-category $\catCohb_{\mathsf{pro}}( \bfCoh_{\mathsf{tor}}(S, \tau_\scrA) )$ has a $\E_1$-monoidal structure, and the stable pro-$\infty$-category $\catCohb_{\mathsf{pro}}( \bfCoh_{\mathsf{t.f.}}(S, \tau_\scrA) )$ has the structure of a left and a right categorical module over $\catCohb_{\mathsf{pro}}( \bfCoh_{\mathsf{tor}}(S, \tau_\scrA) )$. 
	
	Similarly, let $\bfD^\ast$ be a motivic formalism, and fix $\calA \in \CAlg(\bfD^\ast(\Spec(k)))$ and $\Gamma \subseteq \Pic(\bfD^\ast(\Spec(k)))$ such that Assumption~\ref{assumption:motivic_formalism} is satisfied. Then, the topological vector space $\HBMDGamma_0(\bfCoh_{\mathsf{tor}}(S, \tau_\scrA);\calA)$ has the structure of a unital associative algebra. In particular,
	\begin{align}
		G_0( \bfCoh_{\mathsf{tor}}(S, \tau_\scrA) )\quad \text{and} \quad \HBM_\ast( \bfCoh_{\mathsf{tor}}(S, \tau_\scrA) )
	\end{align}
	have the structures of unital associative algebras, and the topological vector space $\HBMDGamma_0(\bfCoh_{\mathsf{t.f.}}(S, \tau_\scrA);\calA)$ has both the structure of a left and a right $\HBMDGamma_0(\bfCoh_{\mathsf{tor}}(S, \tau_\scrA);\calA)$-module. In particular,
	\begin{align}
		G_0( \bfCoh_{\mathsf{t.f.}}(S, \tau_\scrA) )\quad \text{and} \quad \HBM_\ast( \bfCoh_{\mathsf{t.f.}}(S, \tau_\scrA) )			
	\end{align}
	have both the structures of a left and a right $G_0( \bfCoh_{\mathsf{tor}}(S, \tau_\scrA) )$-module and $\HBM_\ast( \bfCoh_{\mathsf{tor}}(S, \tau_\scrA) )$-module, respectively.
\end{theorem}

\begin{proof}
	We shall apply Theorem~\ref{thm:left-right-action} to $(\scrC_S, \tau_\scrA)$ and the torsion pair $(\scrA_{\mathsf{tor}},\scrA_{\mathsf{t.f.}})$. First, we shall show that $(\scrC_S, \tau_\scrA)$ satisfies Assumption~\ref{assumption:existence-COHA}. By \cite[Example~A.4]{AB_Moduli}, the torsion pair $(\scrT, \scrF)$ is open. In particular, the hypotheses of Lemma~\ref{lem:tilted-Assumption-C} are satisfied, hence $(\scrC_S, \tau_\scrA)$ satisfies Assumptions~\ref{assumption:existence-COHA-1}, \ref{assumption:existence-COHA-2}, and \ref{assumption:existence-COHA-3}. Since $\sfS(-)=(-)\otimes_{\scrO_S} \omega_S[2]$ is the Serre functor of $\scrC_S$, we see that Assumption~\ref{assumption:existence-COHA-3} is satisfied as well for $\tau_\scrA$.
	
	Now, conditions~\eqref{item:torsion-pair-left-right-action-1} and \eqref{item:torsion-pair-left-right-action-2} of Theorem~\ref{thm:left-right-action} are trivially satisfied for $(\scrA_{\mathsf{tor}},\scrA_{\mathsf{t.f.}})$. Moreover, by Lemmas~\ref{lem:A_1} and \ref{lem:A_tf_open} the torsion pair $(\scrA_{\mathsf{tor}},\scrA_{\mathsf{t.f.}})$ is open. In particular, also conditions~\eqref{item:torsion-pair-left-right-action-3} and \eqref{item:torsion-pair-left-right-action-4} of Theorem~\ref{thm:left-right-action} are satisfied.
	
	Now, let $(\tau_{\scrA})_\upsilon$ be the tilting of $\tau_{\scrA}$ with respect to the torsion pair $\upsilon=(\scrA_{\mathsf{tor}},\scrA_{\mathsf{t.f.}})$, and let us denote by $\scrA^\upsilon$ the tilted heart.
	Since the tilting of $\tau_\scrB$ with respect to the torsion pair $(\catCoh_{\mathsf{t.f.}}(S)[1], \catCoh_{\mathsf{tor}}(S))$ coincides with the shift by $1$ of the standard $t$-structure, which satisfies again openness of flatness we immediately deduce from Theorem~\ref{thm:duality-tor-torsion-free} that $(\tau_{\scrA})_\upsilon$ also satisfies openness of flatness. Moreover, it follows that $(\scrA_{\mathsf{t.f.}}[1], \scrA_{\mathsf{tor}})$ is an open torsion pair on $\scrA^\upsilon$ and that $\scrA_{\mathsf{tor}}$ is a Serre subcategory of $\scrA^\upsilon$. 
	
	By combining Theorem~\ref{thm:duality-tor-torsion-free} together with the arguments in the first half of the proof of Theorem~\ref{thm:action-torsion}, we obtain that condition~\eqref{item:left-right-action-1} of Theorem~\ref{thm:left-right-action} is satisfied.  On the other hand, by combining Theorem~\ref{thm:duality-tor-torsion-free} together with the arguments in the second half of the proof of Theorem~\ref{thm:action-torsion}, we obtain that also condition~\eqref{item:left-right-action-2} of Theorem~\ref{thm:left-right-action} is satisfied. By applying Theorem~\ref{thm:left-right-action}, the claims follow.
\end{proof}

As in \S\ref{subsec:torsion-pair}, we can refine this action keeping track of the decomposition \eqref{eq:decomposition_connected_components_A}.
More precisely, setting
\begin{align}
	\bfCoh(S, \tau_\scrA; r)\coloneqq \bigsqcup_{\mathsf{rk}(v)=-r} \bfCoh(S, \tau_\scrA; v) \quad\text{and}\quad 	\bfCoh_{\mathsf{t.f.}}(S, \tau_\scrA; r)\coloneqq \bigsqcup_{\mathsf{rk}(v)=-r} \bfCoh_{\mathsf{t.f.}}(\scrA_{\mathsf{t.f.}}S, \tau_\scrA; v) \ ,
\end{align}
we obtain:
\begin{corollary}\label{cor:action-perverse-torsion}
	The stable pro-$\infty$-category $\catCohb_{\mathsf{pro}}( \bfCoh_{\mathsf{t.f.}}(S, \tau_\scrA; r) )$ has the structure of a left and a right categorical module over $\catCohb_{\mathsf{pro}}( \bfCoh_{\mathsf{tor}}(S, \tau_\scrA) )$. 
	
	Similarly, let $\bfD^\ast$ be a motivic formalism, and fix $\calA \in \CAlg(\bfD^\ast(\Spec(k)))$ and $\Gamma \subseteq \Pic(\bfD^\ast(\Spec(k)))$ such that Assumption~\ref{assumption:motivic_formalism} is satisfied. Then, the topological vector space $\HBMDGamma_0(\bfCoh_{\mathsf{t.f.}}(S, \tau_\scrA; r);\calA)$ has both the structure of a left and a right $\HBMDGamma_0(\bfCoh_{\mathsf{tor}}(S, \tau_\scrA);\calA)$-module. In particular,
	\begin{align}
		G_0( \bfCoh_{\mathsf{t.f.}}(S, \tau_\scrA; r) )\quad \text{and} \quad \HBM_\ast( \bfCoh_{\mathsf{t.f.}}(S, \tau_\scrA; r) )			
	\end{align}
	have both the structures of a left and a right $G_0( \bfCoh_{\mathsf{tor}}(S, \tau_\scrA) )$-module and $\HBM_\ast( \bfCoh_{\mathsf{tor}}(S, \tau_\scrA) )$-module, respectively.
\end{corollary}

As it is natural to expect, Theorem~\ref{thm:duality-tor-torsion-free} implies that the algebras and their representations constructed in Theorems~\ref{thm:action-torsion} and \ref{thm:action-perverse-torsion} coincide.
Indeed, Corollary~\ref{cor:autoequivalences-perversities} provides an identification of $2$-Segal spaces
\begin{align}
	\calS_\bullet \bfCoh(S, \tau_\scrA\op) \simeq \calS_\bullet \bfCoh(S, \tau_\scrB) \ , 
\end{align}
and compatibly an identification of relative $2$-Segal spaces
\begin{align}
	\calS_\bullet^\ell \bfCoh(S, \tau_\scrA\op) \simeq \calS_\bullet^\ell \bfCoh(S, \tau_\scrB) \qquad \text{and} \qquad \calS_\bullet^r \bfCoh(S, \tau_\scrA\op) \simeq \calS_\bullet^r \bfCoh(S, \tau_\scrB) \ . 
\end{align}
Paired with Theorem~\ref{thm:duality-tor-torsion-free} and with the results of \S\ref{subsec:torsion-pair} and \S\ref{sec:action-torsion-perverse}, this yields:
\begin{theorem}\label{thm:swap}
	\hfill
	\begin{enumerate}
		\item The equivalence
		\begin{align}\label{eq:D-tor}	
			\begin{tikzcd}[ampersand replacement=\&]
				\D(-) \colon \bfCoh_{\mathsf{tor}}(S, \tau_\scrA\op) \ar{r}{\sim} \&  \bfCoh_{\mathsf{tor}}(S)\
			\end{tikzcd}
		\end{align}
		induces an equivalence of $\E_1$-monoidal stable pro-$\infty$-categories:
		\begin{align}
			\begin{tikzcd}[ampersand replacement=\&]
				\Gamma \colon \catCohb_{\mathsf{pro}}( \bfCoh_{\mathsf{tor}}(S, \tau_\scrA\op) ) \ar{r}{\sim} \& \catCohb_{\mathsf{pro}}( \bfCoh_{\mathsf{tor}}(S) )\ .
			\end{tikzcd}\ .
		\end{align}
		Similarly, let $\bfD^\ast$ be a motivic formalism, and fix $\calA \in \CAlg(\bfD^\ast(\Spec(k)))$ and $\Gamma \subseteq \Pic(\bfD^\ast(\Spec(k)))$ such that Assumption~\ref{assumption:motivic_formalism} is satisfied. Then, the equivalence \eqref{eq:D-tor} induces an isomorphism of unital associative algebras
		\begin{align}
			\begin{tikzcd}[ampersand replacement=\&]
				\Gamma \colon \HBMDGamma_0(\bfCoh_{\mathsf{tor}}(S, \tau_\scrA);\calA)\op \ar{r}{\sim} \& \HBMDGamma_0( \bfCoh_{\mathsf{tor}}(S);\calA)
			\end{tikzcd}\ .
		\end{align}
		In particular, we have isomorphisms of unital associative algebras:
		\begin{align}
			\begin{tikzcd}[ampersand replacement=\&, row sep=tiny]
				\Gamma \colon G_0( \bfCoh_{\mathsf{tor}}(S, \tau_\scrA) )\op \ar{r}{\sim} \& G_0( \bfCoh_{\mathsf{tor}}(S)  ) \\
				\Gamma \colon \HBM_\ast( \bfCoh_{\mathsf{tor}}(S, \tau_\scrA) )\op  \ar{r}{\sim} \&  \HBM_\ast( \bfCoh_{\mathsf{tor}}(S)  ) 
			\end{tikzcd}\ . 
		\end{align}
				
		\item The equivalences
		\begin{align}\label{eq:D-torsion-torsion-free}
			\begin{tikzcd}[ampersand replacement=\&, row sep=tiny]
				\D(-)[-1]\colon  \bfCoh_{\mathsf{t.f.}}(S, \tau_\scrA\op; r)  \ar{r}{\sim} \&  \bfCoh_{\mathsf{t.f.}}(S; r)\\ 
				\D(-) \colon \bfCoh_{\mathsf{tor}}(S, \tau_\scrA\op) \ar{r}{\sim} \&  \bfCoh_{\mathsf{tor}}(S)
			\end{tikzcd}
		\end{align}
		induce an equivalence of left and right categorical modules over $\catCohb_{\mathsf{pro}}( \bfCoh_{\mathsf{tor}}(S) )$:
		\begin{align}
			\begin{tikzcd}[ampersand replacement=\&]
				\Psi\colon \catCohb_{\mathsf{pro}}( \bfCoh_{\mathsf{t.f.}}(S, \tau_\scrA\op; r) )  \ar{r}{\sim} \& \catCohb_{\mathsf{pro}}( \bfCoh_{\mathsf{t.f.}}(S; r) )
			\end{tikzcd}\ .
		\end{align}
		Similarly, let $\bfD^\ast$ be a motivic formalism, and fix $\calA \in \CAlg(\bfD^\ast(\Spec(k)))$ and $\Gamma \subseteq \Pic(\bfD^\ast(\Spec(k)))$ such that Assumption~\ref{assumption:motivic_formalism} is satisfied. Then, the equivalences \eqref{eq:D-torsion-torsion-free} induce an isomorphism of left and right $\HBMDGamma_0(\bfCoh_{\mathsf{tor}}(S);\calA)$-module
		\begin{align}
			\begin{tikzcd}[ampersand replacement=\&]
				\Phi\colon \HBMDGamma_0(\bfCoh_{\mathsf{t.f.}}(S, \tau_\scrA\op; r);\calA)  \ar{r}{\sim} \& \HBMDGamma_0(\bfCoh_{\mathsf{t.f.}}(S; r);\calA)
			\end{tikzcd}
		\end{align}
		In particular, we have isomorphisms
		\begin{align}
			\begin{tikzcd}[ampersand replacement=\&, row sep=tiny]
				G_0( \bfCoh_{\mathsf{t.f.}}(S, \tau_\scrA\op; r) )  \ar{r}{\sim} \&  G_0( \bfCoh_{\mathsf{t.f.}}(S; r) ) \\ 
				\HBM_\ast( \bfCoh_{\mathsf{t.f.}}(S, \tau_\scrA\op; r) )  \ar{r}{\sim} \&  \HBM_\ast( \bfCoh_{\mathsf{t.f.}}(S; r) )
			\end{tikzcd}
		\end{align}
		of left and right modules over $G_0( \bfCoh_{\mathsf{tor}}(S) )$ and of $\HBM_\ast( \bfCoh_{\mathsf{tor}}(S) )$, respectively.
			\end{enumerate}
\end{theorem}

\begin{remark}
	For every $\infty$-groupoid $K \in \scrS$, there is a functorial self-equivalence $\mathsf{inv}_K \colon K\op \simeq K$.
	This induces an equivalence of derived stacks
	\begin{align}
		\bfCoh(S, \tau_\scrA) \simeq \bfCoh(S, \tau_\scrA\op) \ , 
	\end{align}
	which nevertheless does not propagate through the natural $2$-Segal structure on both sides.
	Concretely, this means that $\E_1$-monoidal pro-$\infty$-category $\catCohb_{\mathsf{pro}}(\bfCoh(S, \tau_\scrA\op))$ has $\catCohb_{\mathsf{pro}}(\bfCoh(S, \tau_\scrA))$ as underlying pro-$\infty$-category, and it has the \textit{opposite tensor structure}.
\end{remark}

\section{Representations of the COHA of a surface via stable pairs}\label{sec:action-stable-pairs}

In this section, we construct representations of COHAs and CatHAs of torsion sheaves on a smooth projective surface $S$ associated to \textit{stable pairs}. We introduce such a notion in a great generality and make use of the formalism developed in Part~\ref{part:foundations} to construct the representations.

\subsection{$\calV$-stable (co-)pairs in a stable $\infty$-category} \label{subsubsec:abstract_stable_pairs}

We start by the following definition:
\begin{defin}
	A \textit{framed tilting datum} is a $4$-tuple $\bfC = (\scrC, \tau, \upsilon, \calV)$, where $\scrC$ is a stable $\infty$-category, $\tau$ is a $t$-structure on $\scrC$, $\upsilon=(\scrC^\heartsuit_{\mathsf{tor}}, \scrC^\heartsuit_{\mathsf{t.f.}})$ is a torsion pair on $\scrC^\heartsuit$, and $\calV \in \scrC$ is an object.~$\oslash$
\end{defin}

Given a framed tilting datum $\bfC$ as above, we refer to $\calV$ as the \textit{framing} and to the triple $(\scrC,\tau,\upsilon)$ as the \textit{tilting datum} underlying $\bfC$. We leave as an exercise to verify that framed tilting data can be naturally organized into an $\infty$-category $\mathbb T^{\mathsf{fr}}$. Concretely, given two framed tilting data $\bfC = (\scrC, \tau_\scrC, \upsilon_\scrC, \calV)$ and $\mathbf D = (\scrD, \tau_\scrD, \upsilon_\scrD, \calW)$ a morphism from $\bfC$ to $\mathbf D$ is a pair $(F, \alpha)$, where
\begin{align}
	F \colon \scrC \longrightarrow \scrD 
\end{align}
is an $t$-exact stable functor with the property that $F(\scrC^\heartsuit_{\mathsf{tor}}) \subseteq \scrD^\heartsuit_{\mathsf{tor}}$ and $F(\scrC^\heartsuit_{\mathsf{t.f.}}) \subseteq \scrD^\heartsuit_{\mathsf{t.f.}}$, and
\begin{align}
	\alpha \colon F(\calV) \xrightarrow{\sim} \calW 
\end{align}
is an equivalence in $\scrD$.

Recall that $\calS_2 \scrC$ parametrizes diagrams $\F$ of the form
\begin{align}
	\begin{tikzcd}[ampersand replacement=\&]
		0 \arrow{r} \& F_1 \arrow{r} \arrow{d} \& F_2 \arrow{d} \\
		\& 0 \arrow{r} \& F_3 \arrow{d} \\
		\& \& 0
	\end{tikzcd} \ , 
\end{align}
where the central square is asked to be a pullback. Moreover, we define
\begin{align}
	\partial_0(\F) = F_1\ , \qquad \partial_1(\F) = F_2\ , \qquad \partial_2(\F) = F_3 \ . 
\end{align}

\begin{definition}
	Let $\scrC$ be a stable $\infty$-category. Given an object $\calV \in \scrC$, we let $\scrC^\dagger(\calV)$ be the fiber product
	\begin{align}
		\begin{tikzcd}[ampersand replacement=\&]
			\scrC^\dagger(\calV) \arrow{r} \arrow{d} \& \calS_2 \scrC \arrow{d}{\partial_0} \\
			\ast \arrow{r}{\calV} \& \scrC 
		\end{tikzcd} \ .
	\end{align}
	We refer to $\scrC^\dagger(\calV)$ as \textit{the $\infty$-category of $\calV$-pairs}.
	
	Replacing $\partial_0$ with $\partial_2$ in the above diagram, we denote by $\scrC^\ddag(\calV)$ the resulting category, and we refer to it as the \textit{$\infty$-category of $\calV$-co-pairs.}
\end{definition}	

\begin{remark}
	Informally speaking, objects in $\scrC^\dagger(\calV)$ are fiber sequences of the form $\calV \to \calF \to E$, and morphisms are morphisms of fiber sequences that restrict to the identity on $\calV$. On the other hand, the objects in $\scrC^\ddag(\calV)$ are fiber sequences of the form $\calF \to E \to \calV$ and morphisms are morphisms of fiber sequences that restrict to the identity on $\calV$.
\end{remark}

\begin{definition}\label{def:abstract_pairs}
	Let $\bfC = (\scrC, \tau, \upsilon, \calV) \in \mathbb T^{\mathsf{fr}}$.
	\begin{enumerate}\itemsep=0.2cm
		\item A \textit{$\bfC$-pair} is a fiber sequence in $\scrC$ of the form
		\begin{align}
			\calV \longrightarrow \calF \longrightarrow E 
		\end{align}
		where $\calF[1] \in \scrC^\heartsuit_{\mathsf{tor}}$ and $E \in \scrC^\heartsuit_{\mathsf{t.f.}}$.
		We denote by $\sfP(\bfC)$ the full subcategory of $\scrC^\dagger(\calV)$ spanned by $\bfC$-pairs.
		
		\item A \textit{$\bfC$-co-pair} is a fiber sequence in $\scrC$ of the form
		\begin{align}
			\calF \longrightarrow E \longrightarrow \calV[1] 
		\end{align}
		where $\calF \in \scrC^\heartsuit_{\mathsf{t.f.}}$ and $E \in \scrC^\heartsuit_{\mathsf{tor}}$.
		We denote by $\sfP^{\mathsf c}(\bfC)$ the full subcategory of $\scrC^\ddag(\calV[1])$ spanned by $\bfC$-co-pairs.
	\end{enumerate}
\end{definition}

It is straightforward to check that the constructions $\sfP(\bfC)$ and $\sfP^{\mathsf c}(\bfC)$ give rise to functors
\begin{align}
	\sfP, \ \sfP^{\mathsf c} \colon \mathbb T^{\mathsf{fr}} \longrightarrow \Cat_\infty \ . 
\end{align}
In order to combine this construction with Theorem~\ref{thm:duality-tor-torsion-free}, we need to understand certain natural symmetries of $\mathbb T^{\mathsf{fr}}$, and how they interact with the previous two functors.

Let $t \textrm{-} \Cat_\infty^{\mathsf{st}}$ be the $\infty$-category of stable $\infty$-categories equipped with $t$-structures, and $t$-exact functors between them.
\begin{construction}
	\hfill
	\begin{enumerate}\itemsep=0.2cm
		\item Let $\mu_2 \coloneqq \{1, \zeta_2\} \subseteq \C^\times$.\footnote{Here, we write $\zeta_2$ instead of $-1$ to avoid confusion.}
		There is a canonical action of $\mu_2$ on $\Cat_\infty$ that takes an $\infty$-category $\scrC$ to its opposite $\scrC\op$.
		This action preserves stable $\infty$-categories, and therefore it lifts to an action on $\Cat_\infty^{\mathsf{st}}$. 
		The $\mu_2$-action lifts to $t \textrm{-} \Cat_\infty^{\mathsf{st}}$, observing that if $\tau = (\scrC^{\leqslant 0}, \scrC^{\geqslant 0})$ is a $t$-structure on $\scrC$, then $\tau\op \coloneqq \big((\scrC^{\geqslant 0})\op, (\scrC^{\leqslant 0})\op\big)$ is a $t$-structure on $\scrC\op$.
		
		\item The $\infty$-category $t\textrm{-}\Cat_\infty^{\mathsf{st}}$ also carries a canonical action of $\Z$, which is uniquely determined by the fact that $1$ acts on $(\scrC,\tau)$ by shifting the $t$-structure:
		\begin{align}
			1 \cdot (\scrC,\tau) \coloneqq (\scrC,\tau[1]) \ . 
		\end{align}
	\end{enumerate}	
\end{construction}

\begin{rem}
	The two actions do not commute: indeed the endo-functor $(-)\op \colon \Cat_\infty^{\mathsf{st}} \to \Cat_\infty^{\mathsf{st}}$ takes the self-equivalence $[1] \colon \scrC \simeq \scrC$ to the self-equivalence $[-1] \colon \scrC\op \simeq \scrC\op$.
	Thus,
	\begin{align}
		1 \cdot (\zeta_2 \cdot \tau) & = \tau\op[1] = \big((\scrC^{\geqslant 0})\op[1], (\scrC^{\leqslant 0})\op[1]\big) \\
		& = \big((\scrC^{\geqslant 0}[-1])\op, (\scrC^{\leqslant 0}[-1])\op\big) \\
		& = \big((\scrC^{\geqslant 1})\op,(\scrC^{\leqslant 1})\op\big) \ ,
	\end{align}
	while
	\begin{align}
		\zeta_2 \cdot( 1 \cdot \tau ) = (\tau[1])\op = \big((\scrC^{\geqslant 0}[1])\op, (\scrC^{\leqslant 0}[1])\op\big) = \big((\scrC^{\geqslant -1})\op, (\scrC^{\leqslant -1})\op\big) \ . 
	\end{align}
	The same computation immediately shows that
	\begin{align}
		\tau\op[-1] = (\tau[1])\op \ , 
	\end{align}
	which in turn implies that $t \textrm{-} \Cat_\infty^{\mathsf{st}}$ carries a natural $\mu_2 \ltimes \Z$-action, where $\mu_2$ acts on $\Z$ by multiplication by $-1$.\hfill $\triangle$
\end{rem}

Now, let $\mathbb T$ be the $\infty$-category of stable $\infty$-categories $\scrC$ equipped with a $t$-structure $\tau = (\scrC^{\leqslant 0}, \scrC^{\geqslant 0})$ and a torsion pair $\upsilon = (\scrC^\heartsuit_{\mathsf{tor}}, \scrC^\heartsuit_{\mathsf{t.f.}})$ on the heart $\scrC^\heartsuit$ of $\tau$.
\begin{construction}
	\hfill
	\begin{enumerate}\itemsep=0.2cm
		\item The $\mu_2 \ltimes \Z$-action on $t\textrm{-}\Cat_\infty^{\mathsf{st}}$ canonically lifts to $\mathbb T$ by setting
		\begin{align}
			1 \cdot \upsilon \coloneqq \big(\scrC^\heartsuit_{\mathsf{tor}}[1], \scrC^\heartsuit_{\mathsf{t.f.}}[1]\big) \qquad \text{and} \qquad \zeta_2\cdot\upsilon \coloneqq \upsilon\op \coloneqq \big((\scrC^\heartsuit_{\mathsf{t.f.}})\op, (\scrC^\heartsuit_{\mathsf{tor}})\op\big) \ .
		\end{align}
		The tilting operation provides an extension of this action to the bigger group $\mu_2 \ltimes \Z[\text{\sfrac{1}{2}}]$, where $\text{\sfrac{1}{2}}$ acts by
		\begin{align}
			\text{\sfrac{1}{2}}\cdot (\tau,\upsilon) \coloneqq (\tau_\upsilon, \upsilon_{\text{\sfrac{1}{2}}}) \ . 
		\end{align}
		Here, $\tau_{\upsilon}$ is the $t$-structure obtained by tilting $\tau$ with respect to $\upsilon$, and
		\begin{align}
			\upsilon_{\text{\sfrac{1}{2}}} \coloneqq \big(\scrC^\heartsuit_{\mathsf{t.f.}}[1],\scrC^\heartsuit_{\mathsf{tor}}\big) 
		\end{align}
		is the induced tilted torsion pair.
		
		\item The action of $\mu_2 \ltimes \Z[\text{\sfrac{1}{2}}]$ on $\mathbb T$ extends to $\mathbb T^{\mathsf{fr}}$, by setting
		\begin{align}
			\zeta_2 \cdot \calV \coloneqq \calV\op[-2] \qquad \textrm{and} \qquad \text{\sfrac{1}{2}} \cdot \calV \coloneqq \calV[1] \ . 
		\end{align}
		Here $\calV\op$ denote the object $\calV \in \scrC$ seen as an object in $\scrC\op$, and the shift $[2]$ is understood computed in $\scrC\op$.
		In other words:
		\begin{align}
			\calV\op[-2] \simeq (\calV[2])\op \ . 
		\end{align}
		Finally, observe that since $1 = \text{\sfrac{1}{2}} + \text{\sfrac{1}{2}}$, the above definition implies $1 \cdot \calV \coloneqq \calV[2]$.
	\end{enumerate}
\end{construction}
Summarizing, for any $\bfC = (\scrC, \tau, \upsilon, \calV) \in \mathbb T^{\mathsf{fr}}$ we get
\begin{align}
	1\cdot \bfC = (\scrC, \tau[1], \upsilon[1], \calV[2]) , \
	\sfrac{1}{2} \cdot \bfC = (\scrC, \tau_\upsilon, \upsilon_{\sfrac{1}{2}}, \calV[1]) , \
	\zeta_2\cdot \bfC = (\scrC, \tau\op, \upsilon\op, \calV\op[-2]) \ .
\end{align}

In what follows, we will exclusively need the notation introduced in this construction. For this reason, we allowed ourselves to ignore the technical details that would be needed to properly construct an $\infty$-categorical action of $\mu_2\ltimes \Z[\text{\sfrac{1}{2}}]$ on $\mathbb T^{\mathsf{fr}}$.

\begin{lemma}\label{lem:pair_copairs_symmetries}
	Let $\bfC = (\scrC,\tau,\upsilon,\calV) \in \mathbb T^{\mathsf{fr}}$. Then:
	\begin{enumerate}\itemsep=0.2cm
		\item \label{item:pair_copairs_symmetries-1} The shift-rotation self-equivalence $\rho \colon \calS_2 \scrC \to \calS_2 \scrC$ sending a fiber sequence $F_1 \to F_2 \to F_3$ to $F_2[1] \to F_3[1] \to F_1[2]$ induces natural equivalences
		\begin{align}
			\rho \colon \sfP^{\mathsf c}(\mathrm{\sfrac{1}{2}} \cdot \bfC) \simeq \sfP(\bfC) \quad \text{and} \quad \rho \colon \sfP(\mathrm{\sfrac{1}{2}} \cdot \bfC) \simeq \sfP^{\mathsf c}(1 \cdot \bfC)\ . 
		\end{align}
		
		\item \label{item:pair_copairs_symmetries-2} The canonical equivalence $\omega \colon (\calS_2 \scrC)\op \to \calS_2(\scrC\op)$ sending a fiber sequence $F_1 \to F_2\to F_3$ to $F_3\op \to F_2\op \to F_1\op$ induces a natural equivalence
		\begin{align}
			\rho^{-1} \circ \omega \colon \sfP( \zeta_2 \cdot \bfC) \simeq \sfP(\bfC)\op \ , 
		\end{align}
		where $\rho^{-1}$ is the inverse to the shift-rotation functor $\rho$ considered in the previous point.
	\end{enumerate}
\end{lemma}

\begin{proof}
	First observe that the shift-rotation functor induces a canonical equivalence
	\begin{align}
		\begin{tikzcd}[column sep = small, ampersand replacement=\&]
			\rho \colon \scrC^\dagger(\calV) \arrow{r}{\sim} \& \scrC^\ddag(\calV[2]) \ .
		\end{tikzcd}
	\end{align}
	Now, a $\calV$-pair $\calV \to \calF \to E$ belongs to $\sfP(\bfC)$ if and only if $\calF[1] \in \scrC^\heartsuit_{\mathsf{tor}}$ and $E \in \scrC^\heartsuit_{\mathsf{t.f.}}$. After applying $\rho$, we obtain the fiber sequence
	\begin{align}
		\calF[1] \longrightarrow E[1] \longrightarrow \calV[2] \ , 
	\end{align}
	and $\calF[1] \in \scrC^\heartsuit_{\mathsf{tor}} = (\scrC^\heartsuit_{\upsilon})_{\mathsf{t.f.}}$, while $E[1] \in \scrC^\heartsuit_{\mathsf{t.f.}}[1] = (\scrC^\heartsuit_{\upsilon})_{\mathsf{tor}}$. Therefore, the above sequence is a co-pair for $\text{\sfrac{1}{2}} \cdot \bfC = (\scrC,\tau_\upsilon, \upsilon_{\text{\sfrac{1}{2}}}, \calV[1])$. Thus, $\rho$ induces a well defined functor
	\begin{align}
		\rho \colon \sfP(\bfC) \longrightarrow \sfP^{\mathsf c}(\text{\sfrac{1}{2}} \cdot \bfC) \ , 
	\end{align}
	and it is straightforward to check that it is an equivalence. Similarly, one can show that $\rho$ induces the equivalence $\sfP(\text{\sfrac{1}{2}}\cdot \bfC) \simeq \sfP^{\mathsf c}(1 \cdot \bfC)$. This proves statement \eqref{item:pair_copairs_symmetries-1}.
	
	We now prove statement \eqref{item:pair_copairs_symmetries-2}. To begin with, the natural equivalence $\omega$ tautologically induces an equivalence
	\begin{align}
		\begin{tikzcd}[column sep = small, ampersand replacement=\&]
			\omega \colon \big(\scrC^\dagger(\calV)\big)\op \arrow{r}{\sim} \& (\scrC\op)^\ddag(\calV\op) \ ,
		\end{tikzcd} 
	\end{align}
	which sends a $\calV$-pair $\calV \to \calF \to E$ to $E\op \to \calF\op \to \calV\op$.
	Further applying
	\begin{align}
		\rho^{-1} \colon (\scrC\op)^\ddag(\calV\op) \longrightarrow (\scrC\op)^\dagger(\calV\op[-2]) \ , 
	\end{align}
	we obtain
	\begin{align}
		\calV\op[-2] \longrightarrow E\op[-1] \longrightarrow \calF\op[-1] \ . 
	\end{align}
	We now observe that $E\op[-1][1] = \calE\op \in (\scrC^\heartsuit_{\mathsf{t.f.}})\op = (\scrC\op)^\heartsuit_{\mathsf{tor}}$, while $\calF\op[-1] \simeq (\calF[1])\op$ belongs to $(\scrC^\heartsuit_{\mathsf{tor}})\op = (\scrC\op)^\heartsuit_{\mathsf{t.f.}}$. Thus, $\omega$ induces a well defined functor
	\begin{align}
		\rho^{-1} \circ \omega \colon \bfP(\bfC)\op \longrightarrow \bfP(\zeta_2 \cdot \bfC) \ , 
	\end{align}
	which is readily checked to be an equivalence.
\end{proof}

Now, we are ready to state the main result of this section, which follows from Theorem~\ref{thm:duality-tor-torsion-free}.
\begin{theorem}\label{thm:duality_stable_pairs_categorical}
	Let $S$ be a smooth projective irreducible complex surface and let $\calV$ be a locally free sheaf on $S$. Set $\upsilon_{\scrB}\coloneqq (\catCoh_{\mathsf{tor}}(S), \catCoh_{\mathsf{t.f.}}(S))$ and $\upsilon_{\scrA}\coloneqq (\scrA_{\mathsf{tor}}, \scrA_{\mathsf{t.f.}})$. 
	
	Let $\bfC \coloneqq (\catPerf(S), \tau_{\mathsf{std}}, \upsilon_{\scrB}, \D(\calV[2])[-1])$ and $\mathbf A \coloneqq (\catPerf(S), \tau_{\scrA}, \upsilon_{\scrA}, \calV)$.
	Then there are canonical equivalences
	\begin{align}
		\D(-)\colon \sfP(\mathbf A)\op \simeq \sfP(\mathsf{\sfrac{1}{2}} \cdot \bfC) \qquad \text{and} \qquad \D(-)\colon\sfP^{\mathsf c}(\mathsf{\sfrac{1}{2}} \cdot \mathbf A)\op \simeq \sfP^{\mathsf c}(1 \cdot \bfC) \ . 
	\end{align}
		The former can informally be described as sending a fiber sequence $\calV \to \calF \to E$ to $\D(\calV)[-2] \to \D(E)[-1] \to \D(\calF)[-1]$, and the latter can be described as sending a fiber sequence $\calF \to E \to \calV[2]$ to $\D(E)[1] \to \D(\calF)[1] \to \D(\calV[2])[2]$.
\end{theorem}

\begin{proof}
	To begin with, we observe that the functor
	\begin{align}
		\D(-) \colon \catPerf(S)\op \longrightarrow \catPerf(S) 
	\end{align}
	takes $\calV\op[-2] = (\calV[2])\op$ to $\calV^\vee[-2][2] \simeq \calV^\vee$. Thus, Theorem~\ref{thm:duality-tor-torsion-free} can be restated saying that the anti-equivalence $\D(-)$ induces an equivalence
	\begin{align}
		\D(-) \colon \zeta_2 \cdot \mathbf A \simeq \text{\sfrac{1}{2}} \cdot \bfC
	\end{align}
	in $\mathbb T^{\mathsf{fr}}$.
	At this point, applying repeatedly Lemma~\ref{lem:pair_copairs_symmetries}, we find:
	\begin{align}
		\sfP(\mathbf A)\op \simeq \sfP(\zeta_2 \cdot \mathbf A) \simeq \sfP(\text{\sfrac{1}{2}} \cdot \bfC) \ . 
	\end{align}
	Similarly,
	\begin{align}
		\sfP^{\mathsf c}(\text{\sfrac{1}{2}} \cdot \mathbf A)\op \simeq \sfP(\mathbf A)\op \simeq \sfP(\text{\sfrac{1}{2}}\cdot \bfC) \simeq \sfP^{\mathsf c}(1 \cdot \bfC) \ . 
	\end{align}
	The explicit formulas are obtained by unwinding the definitions of these equivalences.
\end{proof}

\subsection{Stable (Co-)pairs on a smooth surface $S$}

Let $S$ be a smooth projective irreducible complex surface and let $\calV$ be a locally free sheaf on $S$. 
\begin{definition}\label{def:stable_pairs}
	A \textit{$\calV$-stable pair on $S$} is a $(\catPerf(S), \tau_{\scrA}, \upsilon_{\scrA}, \calV)$-pair.\footnote{Our viewpoint on stable pairs resembles that in \cite{Bridgeland_Hall} for Pandharipande-Thomas stable pairs on Calabi-Yau threefolds.} When $\calV = \scrO_S^{\oplus r}$, we refer to $\calV$-stable pairs as \textit{rank $r$ stable pairs}. When $r = 1$, we simply refer to them as \textit{stable pairs}.
\end{definition}
We set $\sfP(S;\calV)\coloneqq \sfP( (\catPerf(S), \tau_{\scrA}, \upsilon_{\scrA}, \calV) )$. When $\calV = \scrO_S^{\oplus r}$, we write $\sfP(S;r)$ instead of $\sfP(S;\scrO_S^{\oplus r})$.

When $r=1$, our notion coincides with the usual one of \textit{Pandharipande-Thomas stable pairs} and \textit{Bradlow pairs} on $S$, as implied by the following proposition.
\begin{proposition} \label{prop:stable_pairs_different_formulations}
	Let $\calV$ be a locally free sheaf on $S$ of finite rank.
	For a fiber sequence
	\begin{align}\label{eq:stable_pairs_as_triples}
		\begin{tikzcd}[column sep = 16pt, ampersand replacement = \&]
			\calV \arrow{r}{s} \& \calF \arrow{r} \& E
		\end{tikzcd}
	\end{align}
	in $\catPerf(S)$, the following statements are equivalent:
	\begin{enumerate}\itemsep=0.2cm
		\item \label{item:stable_pairs_different_formulations-1} \eqref{eq:stable_pairs_as_triples} corresponds to a Pandharipande-Thomas $\calV$-stable pair on $S$, i.e., $\calF$ is purely one-dimensional coherent sheaf and the morphism $s \colon \calV \to \calF$ has zero-dimensional cokernel;
		
		\item \label{item:stable_pairs_different_formulations-2} \eqref{eq:stable_pairs_as_triples} is a $(\catPerf(S), \tau_{\scrA}, \upsilon_{\scrA}, \calV)$-pair, i.e., $E$ belongs to $\scrA_{\mathsf{t.f.}}$ and $\calF[1]$ belongs to $\scrA_{\mathsf{tor}}$;
		
		\item \label{item:stable_pairs_different_formulations-3} $E$ belongs to $\scrA_{\mathsf{t.f.}}$, $\calF[1]$ belongs to $\scrA$ and $\mathsf{rk}(E) = - \mathsf{rk}(\calV)$.
	\end{enumerate}
	In this case, we further have $\ch_1(E) = \ch_1(\calF)$ if and only if $\ch_1(\calV)=0$.
\end{proposition}

\begin{proof}
	We first prove that \eqref{item:stable_pairs_different_formulations-1} $\Rightarrow$ \eqref{item:stable_pairs_different_formulations-2}.
	
	Since $\calF$ is purely one-dimensional, we have $\calF[1] \in \scrA_{\mathsf{tor}} \subseteq \scrA$ by definition. So we only have to show that $E \in \scrA_{\mathsf{t.f.}}$. Taking the long exact sequence of cohomology sheaves associated to \eqref{eq:stable_pairs_as_triples} we obtain
	\begin{align}\label{eq:stable_pairs_long_exact_sequence}
		0 = \calH^{-1}(\calF) \longrightarrow \calH^{-1}(E) \longrightarrow \calH^{0}(\calV) \longrightarrow \calH^{0}(\calF) \longrightarrow \calH^0(E) \longrightarrow 0 \ .
	\end{align}
	The central terms are canonically identified with $\calV$ and $\calF$, respectively. Thus $\calH^{-1}(E)$ is torsion-free, and our assumption guarantees that $\calH^0(E)$ is zero-dimensional. In other words, $E \in \scrA$. 
		Rotating the sequence \eqref{eq:stable_pairs_as_triples} we obtain the fiber sequence
	\begin{align}
		\begin{tikzcd}[column sep = 16pt, ampersand replacement = \&]
			E \arrow{r} \& \calV[1] \arrow{r}{s[1]} \& \calF[1] 
		\end{tikzcd}\ ,
	\end{align}
	where the three terms belong to $\scrA$. It follows that this is in fact a short exact sequence in $\scrA$ and therefore that the map $E \to \calV[1]$ is injective. Since $\calV$ is locally free, we have that $\calV^\vee$ is again locally free, and in particular $\D(\calV[1]) \simeq \calV^\vee[1]\otimes_{\scrO_S}\omega_S$ belongs to $\catCoh_{\mathsf{t.f.}}(S)[1]$. Hence, Theorem~\ref{thm:duality-tor-torsion-free} guarantees that $\calV[1]$ belongs to $\scrA_{\mathsf{t.f.}}$. Thus, Lemma~\ref{lem:torsion_pairs_extensions_mono_epi}--\eqref{item:torsion_pairs_extensions_mono_epi-3} guarantees that $E \in \scrA_{\mathsf{t.f.}}$. Finally, since $\calF$ is purely one-dimensional coherent sheaf, we get that $\mathsf{rk}(\calF[1]) = 0$, and hence that
	\begin{align}
		\mathsf{rk}(E) = \mathsf{rk}(\calV[1]) = - \mathsf{rk}(\calV) \ . 
	\end{align}
	
	We now prove that \eqref{item:stable_pairs_different_formulations-2} $\Rightarrow$ \eqref{item:stable_pairs_different_formulations-1}.
	
	Since $\calF[1] \in \scrA_{\mathsf{tor}}$, Lemma~\ref{lem:A_1-1} guarantees that $\mathsf{rk}(\calF[1]) = 0$, that $\calH^0(\calF) \simeq \calH^{-1}(\calF[1])$ is purely one-dimensional and that $\calH^1(\calF) \simeq \calH^0(\calF[1])$ is zero-dimensional. Passing to the long exact sequence of cohomology sheaves associated to \eqref{eq:stable_pairs_as_triples}, we obtain
	\begin{align}
		0 = \calH^1(\calV) \longrightarrow \calH^1(\calF) \longrightarrow \calH^1(E) = 0 \ . 
	\end{align}
	Thus, $\calF \simeq \calH^{-1}(\calF[1])$ is purely one-dimensional. Finally, the sequence \eqref{eq:stable_pairs_long_exact_sequence} canonically identifies the cokernel of $\calV \to \calF$ with $\calH^0(E)$, which is zero-dimensional because $E \in \scrA$.
	
	Finally, the equivalence \eqref{item:stable_pairs_different_formulations-3} $\Leftrightarrow$ \eqref{item:stable_pairs_different_formulations-2} follows directly from Lemma~\ref{lem:A_1-1}: indeed, if $\calF[1]$ belongs to $\scrA$, then $\calF[1] \in \scrA_{\mathsf{tor}}$ if and only if $\mathsf{rk}(\calF[1]) = 0$, and this is equivalent to ask $\mathsf{rk}(E) = - \mathsf{rk}(\calV)$.
\end{proof}

\begin{remark}\label{rem:duality_stable_pairs_categorical}
	In light of Proposition~\ref{prop:stable_pairs_different_formulations}, Theorem~\ref{thm:duality_stable_pairs_categorical} can be stated saying that for a fiber sequence $\calV \to \calF \to E$ the following statements are equivalent:
	\begin{enumerate}\itemsep=0.2cm
		\item the fiber sequence is a $\calV$-stable pair, i.e., $\calF$ is purely $1$-dimensional and $E \in \scrA_{\mathsf{t.f.}}$;
		
		\item all the terms of the associated fiber sequence $\D(\calV)[-2] \to \D(E)[-1] \to \D(\calF)[-1]$ belong to $\catCoh^\heartsuit(S)$, and $\D(E)[-1]$ is torsion free, while $\D(\calF)[-1]$ is purely $1$-dimensional.
	\end{enumerate}
\end{remark}

\begin{remark}
	Let $X$ be a Calabi-Yau threefold. In \cite[Lemma~4.5]{Toda_Stable} it is shown that the forgetful functor
	\begin{align}
		\partial_0 \colon \sfP(X;1) \longrightarrow \scrA_{\mathsf{t.f.}} 
	\end{align}
	sending a stable pair $\calV \xrightarrow{s} \calF \to E$ to $E = \mathsf{cofib}(s)$ is an equivalence. This is not the case for surfaces, and indeed the proof given in \textit{loc.\ cit.} uses as essential ingredient the fact that for any $E \in \scrA_{\mathsf{t.f.}}$ there is a canonical isomorphism $\Hom_{\scrA}(E, \scrO_X[1]) \simeq \C$. This is obtained as follows: start with the fiber sequence
	\begin{align}
		\calH^{-1}(E)[1] \longrightarrow E \longrightarrow \calH^0(E) \ ;
	\end{align}
	applying $(-)^\vee[1]$ and passing to the long exact sequence of cohomology groups we obtain
	\begin{align}
		\Ext^1_X(\calH^0(E), \scrO_X) \longrightarrow \Ext^1_X(E,\scrO_X) \longrightarrow \Ext^0_X(\calH^{-1}(E),\scrO_X) \longrightarrow \Ext^2_X(\calH^0(E), \scrO_X) \ . 
	\end{align}
	Since $X$ is a threefold, the extremes vanish, and therefore the central morphism is an isomorphism. Finally, since $\calH^{-1}(E)$ is the ideal sheaf of a curve, we have
	\begin{align}
		\Ext^0_X(\calH^{-1}(E),\scrO_X) \simeq \C \ . 
	\end{align}
	On the other hand, in the surface case, we have by Serre duality
	\begin{align}
		\Ext^2_S(\calH^0(E), \scrO_X) \simeq \sfH^0(\calH^0(E) \otimes \omega_S)^\vee \ , 
	\end{align}
	which is typically non-zero. Moreover, in the rank one case, $\calH^{-1}(E)$ is a line bundle, and thus we see that $\Hom_S(\calH^{-1}(E), \scrO_S)$ is in general not isomorphic to $\C$. Thus, one cannot prove the existence of a canonical morphism $E\to \scrO_S[1]$. 
\end{remark}

As before, let us denote by $\scrA^\upsilon$ the tilted heart with respect to the torsion pair $(\scrA_{\mathsf{tor}},\scrA_{\mathsf{t.f.}})$. Recall that $(\scrA_{\mathsf{t.f.}}[1], \scrA_{\mathsf{tor}})$ is again a torsion pair on $\scrA^\upsilon$. We write as usual
\begin{align}
	\scrA^\upsilon_{\mathsf{tor}} \coloneqq \scrA_{\mathsf{t.f.}}[1] \qquad \text{and} \qquad \scrA^\upsilon_{\mathsf{t.f.}} \coloneqq \scrA_{\mathsf{tor}} \ . 
\end{align}
We start with the following definition.
\begin{definition}\label{def:stable_copair}
	Let $\calV$ be a locally free sheaf of finite rank on $S$. A \textit{$\calV$-stable co-pair} is a $(\sfrac{1}{2}\cdot (\catPerf(S), \tau_{\scrA}, \upsilon_{\scrA}, \calV))$-pair, i.e., a fiber sequence
	\begin{align}
		\calF \longrightarrow E \longrightarrow \calV[2] \ , 
	\end{align}
	where $\calF \in \scrA^\upsilon_{\mathsf{t.f.}}$ and $E \in \scrA^\upsilon_{\mathsf{tor}}$.
\end{definition}

The following lemma is tautological.
\begin{lemma} \label{lem:costable_stable}
	Let $\calV$ be a locally free sheaf of finite rank on $S$.
	A fiber sequence
	\begin{align}
		\calF \longrightarrow E \longrightarrow \calV[2] 
	\end{align}
	is a $\calV$-stable co-pair if and only if the shifted-rotated sequence
	\begin{align}
		\calV \longrightarrow \calF[-1] \longrightarrow E[-1] 
	\end{align}
	is a $\calV$-stable pair.
\end{lemma}

\begin{corollary} \label{cor:costable_pair_numerical}
	Let $\calV$ be a locally free sheaf of finite rank on $S$. Let
	\begin{align}
		\sfE = (\calF \longrightarrow E \longrightarrow \calV[2])
	\end{align}
	be a fiber sequence in $\catPerf(S)$, where $\calF$ and $E$ belong to $\scrA^\upsilon$. Then, the following statements are equivalent:
	\begin{enumerate}\itemsep=0.2cm
		\item $\sfE$ is a $\calV$-stable co-pair;
		
		\item $E \in \scrA^\upsilon_{\mathsf{tor}}$ and $\mathsf{rk}(E) = - \mathsf{rk}(\calV)$;
		
		\item $E \in \scrA^\upsilon_{\mathsf{tor}}$ and $\mathsf{rk}(\calF) = 0$.
	\end{enumerate}
\end{corollary}

\begin{proof}
	This follows at once combining Proposition~\ref{prop:stable_pairs_different_formulations} and Lemma~\ref{lem:costable_stable}.
\end{proof}

\subsection{Extension of stable pairs by torsion perverse sheaves}

The goal of this section is to show that the various Hall algebras attached to $\bfCoh_{\mathsf{tor}}(S,\tau_\scrA)$ naturally act on stable pairs. As usual, the fundamental mechanism is given by extensions. In the case at hand, it takes the following form:
\begin{definition}\label{def:extension_stable_pairs}
	Let $\sfE \coloneqq (\calV \to \calF \to E)$ be a $\calV$-stable pair and let $G \in \scrA_{\mathsf{tor}}$.
	An \textit{extension of $\sfE$ by $G$} is a commutative diagram $\E$ in $\catPerf(S)$
	\begin{align} \label{eq:stable_pairs_extensions}
		\begin{tikzcd}[ampersand replacement = \&]
			0 \arrow{r} \& \calV \arrow{r} \arrow{d} \& \calF \arrow{r} \arrow{d} \& \calF' \arrow{d} \\
			\& 0 \arrow{r} \& E \arrow{r} \arrow{d} \& E' \arrow{d} \\
			\& \& 0 \arrow{r} \& G \arrow{d} \\
			\& \& \& 0
		\end{tikzcd} \ ,
	\end{align}
	all of whose squares are pullback, and where we further ask that $\calV \to \calF' \to E'$ is a $\calV$-stable pair.
\end{definition}

\begin{notation} \label{notation:extension_stable_pairs}
	Compatibly with the simplicial notation of Part~\ref{part:Segal} (cf.\ Remark~\ref{rem:flag_action_coordinates_II}), in the situation of the above definition we set
	\begin{align}
		\varpi_1(\E) \coloneqq ( \calV \to \calF \to E ) \ , \quad \varpi_0(\E) \coloneqq (\calV \to \calF' \to E') \ , \quad u_1^\ell(\E) \coloneqq G \ .
	\end{align}
\end{notation}

\begin{proposition}\label{prop:stable_pairs_extensions}
	Let $\E$ be a diagram in $\catPerf(S)$ of the form \eqref{eq:stable_pairs_extensions}, where $G \in \scrA_{\mathsf{tor}}$. Assume that $\varpi_0(\E) = (\calV \to \calF' \to E')$ is a $\calV$-stable pair. Then $\varpi_1(\E) = (\calV \to \calF \to E)$ is a $\calV$-stable pair if and only if $E \in \scrA$.
\end{proposition}

\begin{proof}
	If $\varpi_1(\E)$ is a $\calV$-stable pair, then $E \in \scrA_{\mathsf{t.f.}}$ by definition, and in particular it belongs to $\scrA$. On the other hand, if $E \in \scrA$, then $E \to E' \to G$ is a fiber sequence in $\catPerf(S)$ and all its terms belong to $\scrA$. Thus, it is a short exact sequence in $\scrA$. In particular, $E \to E'$ is a monomorphism and therefore Lemma~\ref{lem:torsion_pairs_extensions_mono_epi}--\eqref{item:torsion_pairs_extensions_mono_epi-3} guarantees that $E \in \scrA_{\mathsf{t.f.}}$. Moreover, since $\mathsf{rk}(G) = 0$, we also obtain $\mathsf{rk}(E) = \mathsf{rk}(E')=-\mathsf{rk}(\calV)$. Finally, consider the fiber sequence
	\begin{align}
		G \longrightarrow \calF[1] \longrightarrow \calF'[1] \ . 
	\end{align}
	Since its extremes belong to $\scrA$, we see that $\calF[1] \in \scrA$ as well. Thus, the conclusion follows from Proposition~\ref{prop:stable_pairs_different_formulations}.
\end{proof}

\subsection{Right representations of $0$-dimensional sheaves via stable pairs} \label{subsec:stable_pairs_left_rep}

Fix a locally free sheaf $\calV$ of finite rank on $S$. To keep the notation under control, we let
\begin{align}
	\bfPerfps^\dagger(S;\calV) \coloneqq \bfFlagPerfps^{(2),\dagger}(\scrC_S;\calV) \ , 
\end{align}
where the latter is the stack of $\calV$-flags of length $2$ (see Definition~\ref{def:V_flags}). Intuitively, it is the derived stack parametrizing fiber sequences of the form
\begin{align}
	\sfE \coloneqq ( \calV \longrightarrow \calF \longrightarrow E ) \ . 
\end{align}
Compatibly with the simplicial notation of Part~\ref{part:Segal}, we set
\begin{align}
	\partial_0(\sfE) \coloneqq E \qquad \textrm{and} \qquad \partial_1(\sfE) \coloneqq \calF \ . 
\end{align}
We set
\begin{align}
	\calS^\ell_\bullet \bfPerfps^\dagger(S;\calV) \coloneqq \calS^\ell_\bullet \bfFlagPerfps^{(2),\dagger}(\scrC_S;\calV) \ ,
\end{align}
where the latter is the simplicial derived stack of \eqref{eq:universal_m_flag_representation}, obtained via Construction~\ref{construction:Vflags_generalized}. This simplicial stack comes equipped with a natural morphism
\begin{align} \label{eq:stable_pairs_action_I}
	u_\bullet^\ell \colon \calS^\ell_\bullet \bfPerfps^\dagger(S;\calV) \longrightarrow \calS_\bullet \bfPerfps(S) \ ,
\end{align}
which is a relative $2$-Segal space.

We are going to show that this action induces an action at the level of $\calV$-stable pairs. To begin with, following Definition~\ref{def:stable_pairs}, we introduce the moduli stack of $\calV$-stable pairs as the fiber product
\begin{align}
	\begin{tikzcd}[ampersand replacement=\&]
		\bfP(S;\calV) \arrow{r} \arrow{d} \& \bfPerfps^\dagger(S;\calV) \arrow{d}{\partial_1[1] \times \partial_0} \\
		\bfCoh_{\mathsf{tor}}(S,\tau_\scrA) \times \bfCoh_{\mathsf{t.f.}}(S,\tau_\scrA) \arrow{r} \& \bfPerfps(S) \times \bfPerfps(S)
	\end{tikzcd} \ .
\end{align}
When $\calV = \scrO_S^{\oplus r}$, we write $\bfP(S;r)$ instead of $\bfP(S;\calV)$.

Specializing the construction in \S\ref{subsec:COHA-representations} with $\bfH \coloneqq \bfCoh_{\mathsf{tor}}(S,\tau_\scrA)$ and $\bfM \coloneqq \bfP(S;\calV)$, we obtain a new relative simplicial derived stack that we denote
\begin{align}\label{eq:stable_pairs_action_II}
	u^\ell_\bullet \colon \calS^\ell_\bullet \bfP(S;\calV) \longrightarrow \calS_\bullet \bfCoh_{\mathsf{tor}}(S,\tau_\scrA) \ .
\end{align}
In order to prove that this is a relative $2$-Segal space and prove that it gives rise to actual representations, it is convenient to provide a different description. 

Introduce the auxiliary derived stack $\bfPerf^\dagger_{0\textrm{-}\scrA}(S;\calV)$ defined via the fiber product
\begin{align}
	\begin{tikzcd}[ampersand replacement=\&]
		\bfPerf^\dagger_{0\textrm{-}\scrA}(S;\calV) \arrow{r} \arrow{d} \& \bfPerfps^\dagger(\scrC;\calV) \arrow{d}{\partial_0} \\
		\bfCoh_{\mathsf{t.f.}}(S,\tau_\scrA) \arrow{r} \& \bfPerfps(S)
	\end{tikzcd}  \ ,
\end{align}
which can informally be described as the derived stack parametrizing $\calV$-extensions of the form $\calV \to \calF \to E$, where $E \in \scrA_{\mathsf{t.f.}}$, but where no condition is put on $\calF$.

Specializing the construction in \S\ref{subsec:COHA-representations} with $\bfH \coloneqq \bfCoh_{\mathsf{tor}}(S,\tau_\scrA)$ and $\bfM \coloneqq \bfPerf^\dagger_{0\textrm{-}\scrA}(S;\calV)$ we obtain the relative simplicial derived stack
\begin{align}\label{eq:stable_pairs_action_auxiliary}
	u^\ell_\bullet \colon \calS^\ell_\bullet \bfPerf^\dagger_{0\textrm{-}\scrA}(S;\calV) \longrightarrow \calS_\bullet \bfCoh_{\mathsf{tor}}(S,\tau_\scrA) \ .
\end{align}
\begin{lemma} \label{lem:stable_pairs_action_auxiliary}
	The square
	\begin{align}
		\begin{tikzcd}[ampersand replacement=\&]
			\calS^\ell_1 \bfPerf^\dagger_{0\textrm{-}\scrA}(S;\calV) \arrow{r} \arrow{d}{u_1^\ell\times \partial_1} \& \calS^\ell_1 \bfPerfps^\dagger(S;\calV) \arrow{d}{u_1^\ell\times \partial_1} 	\\
			\bfCoh_{\mathsf{tor}}(S,\tau_\scrA)\times \bfPerf^\dagger_{0\textrm{-}\scrA}(S;\calV) \arrow{r} \& \bfPerfps(S)\times \bfPerfps^\dagger(S;\calV)
		\end{tikzcd} 
	\end{align}
	is a pullback. In particular, \eqref{eq:stable_pairs_action_auxiliary} is a relative $2$-Segal space.
\end{lemma}

\begin{proof}
	The second half follows directly from the first one and Proposition~\ref{prop:2-Segal-algebra}. It is thus enough to prove the first statement. Unwinding the definitions, we first see that both horizontal morphisms are representable by open immersions. Therefore, it is enough to prove that we indeed have a pullback square after evaluating at geometric points. In this case, we are called to prove that given a diagram of the form \eqref{eq:stable_pairs_extensions}, if $E \in \scrA$ and $G \in \scrA_{\mathsf{tor}}$, then $E' \in \scrA$ as well. However, since $E \to E' \to G$ is a fiber sequence, the conclusion follows directly from the long exact sequence of homotopy groups associated to the $t$-structure $\tau_\scrA$.
\end{proof}

We now observe that the derived moduli stack of stable pairs $\bfP(S;\calV)$ maps canonically to $\bfPerf^\dagger_{0\textrm{-}\scrA}(S;\calV)$. We therefore obtain a commutative diagram $\rho$
\begin{align}
	\rho\coloneqq \begin{tikzcd}[ampersand replacement=\&]
		\bfP(S;\calV) \arrow{r} \arrow{d} \& \Spec(\C) \arrow{r} \arrow[equal]{d} \& \bfCoh_{\mathsf{tor}}(S,\tau_\scrA) \arrow[equal]{d} \\
		\bfPerf^\dagger_{0\textrm{-}\scrA}(S;\calV) \arrow{r} \& \Spec(\C) \arrow{r} \& \bfCoh_{\mathsf{tor}}(S,\tau_{\scrA}) 
	\end{tikzcd}\ ,
\end{align}
which defines a boundary datum for the relative $2$-Segal space \eqref{eq:stable_pairs_action_auxiliary} in the sense of Definition~\ref{def:Hecke_datum}.

Applying Construction~\ref{construction:Hecke_pattern}, we obtain once again the stable pair action \eqref{eq:stable_pairs_action_II}.
The usefulness of this description of \eqref{eq:stable_pairs_action_II} is explained by the following:
\begin{proposition} \label{prop:stable_pair_action_2-Segal}
	The square
	\begin{align}
		\begin{tikzcd}[ampersand replacement=\&]
			\calS_1^\ell \bfP(S;\calV) \arrow{r} \arrow{d}{\partial_0} \& \calS_1^\ell \bfPerf^\dagger_{0\textrm{-}\scrA}(S;\calV) \arrow{d}{\partial_0} \\
			\bfP(S;\calV) \arrow{r} \& \bfPerf^\dagger_{0\textrm{-}\scrA}(S;\calV)
		\end{tikzcd}
	\end{align}
	is a pullback. In particular, \eqref{eq:stable_pairs_action_II} is a relative $2$-Segal space.
\end{proposition}

\begin{proof}
	Observe first that the horizontal maps are representable by open immersions. Therefore, it is enough to check that the statement is true after evaluating on geometric points. In this case, we are called to show that given a diagram $\sfE$ of the form \eqref{eq:stable_pairs_extensions}, if $\calV \to \calF' \to E'$ is a $\calV$-stable pair, $G \in \scrA_{\mathsf{tor}}$ and $E \in \scrA$, then $\calV \to \calF \to E$ is also a $\calV$-stable pair. This is guaranteed to be true by Proposition~\ref{prop:stable_pairs_extensions}. For the second part of the statement, observe that what we just showed implies immediately that the square
	\begin{align}
		\begin{tikzcd}[ampersand replacement=\&]
			\calS_1^\ell \bfP(S;\calV) \arrow{r} \arrow{d}{u_1^\ell\times \partial_0} \& \calS_1^\ell \bfPerf^\dagger_{0\textrm{-}\scrA}(S;\calV) \arrow{d}{u_1^\ell\times \partial_0} \\
			\bfCoh_{\mathsf{tor}}(S,\tau_\scrA)\times \bfP(S;\calV) \arrow{r} \& \bfCoh_{\mathsf{tor}}(S,\tau_\scrA)\times \bfPerf^\dagger_{0\textrm{-}\scrA}(S;\calV)
		\end{tikzcd} 
	\end{align}
	is a pullback, so that the conclusion follows from Lemma~\ref{lem:stable_pairs_action_auxiliary} and Corollary~\ref{cor:Hecke_pattern}.
\end{proof}

Having constructed a $2$-Segal action of $\bfCoh_{\mathsf{tor}}(S,\tau_\scrA)$ on $\bfP(S;\calV)$, we proceed to check the conditions of Corollary~\ref{cor:COHA-representations}. We begin with properness, and we need to fix the following notation.
\begin{notation} \label{notation:stable_pairs_to_perverse_tf}
	Let
	\begin{align}
		\pi^{\mathsf{t.f.}}_0 \coloneqq \partial_0 \colon \bfP(S;\calV) \longrightarrow \bfCoh_{\mathsf{t.f.}}(S,\tau_\scrA) 
	\end{align}
	be the morphism that sends a $\calV$-stable pair $\sfE = (\calV \to \calF \to E)$ to its $\scrA_{\mathsf{t.f.}}$-component $E = \partial_0(\sfE)$.
	It can be promoted to a morphism of simplicial derived stacks
	\begin{align}
		\pi^{\mathsf{t.f.}}_\bullet \colon \calS^\ell_\bullet \bfP(S;\calV) \longrightarrow \calS^\ell_\bullet \bfFlagCoh^{(1)}_{\mathsf{t.f.},\mathsf{tor}}(S,\tau_\scrA) \ , 
	\end{align}
	and the level $1$ component $\pi^{\mathsf{t.f.}}_1$ can be explicitly described as the map sending an extension of stable pairs of the form \eqref{eq:stable_pairs_extensions} to the sub-extension $E \to E' \to G$.
\end{notation}

\begin{lemma}\label{lem:stable_pairs_action_properness}
	The commutative square
	\begin{align}
		\begin{tikzcd}[ampersand replacement=\&]
			\calS^\ell_1 \bfP(S;\calV) \arrow{r}{\pi^{\mathsf{t.f.}}_1} \arrow{d}{\partial_0} \& \calS^\ell_1 \bfFlagCoh^{(1)}_{\mathsf{t.f.},\mathsf{tor}}(S,\tau_\scrA) \arrow{d}{\partial_0} \\
			\bfP(S;\calV) \arrow{r}{\pi^{\mathsf{t.f.}}_0} \& \bfCoh_{\mathsf{t.f.}}(S,\tau_\scrA)
		\end{tikzcd}
	\end{align}
	is a pullback square. In particular, the left vertical morphism is locally rpas.
\end{lemma}

\begin{proof}
	The right vertical morphism is locally rpas by combining Theorem~\ref{thm:duality-tor-torsion-free} and the proof that the map \eqref{eq:properness-scrB} in the proof of Theorem~\ref{thm:action-torsion} is locally rpas.
	
	Thus all we have to do is to prove that the square of the statement is a pullback. Unwinding the definitions we see that we have to show that given any solid diagram of the form
	\begin{align}
		\begin{tikzcd}[ampersand replacement=\&]
			0 \arrow{r} \& \calV \arrow[dashed]{r} \arrow[bend left = 25pt]{rr} \arrow{d} \& \calF \arrow[dashed]{r} \arrow[dashed]{d} \& \calF' \arrow{d} \\
			\& 0 \arrow{r} \& E \arrow{r} \arrow{d} \& E' \arrow{d} \\
			\& \& 0 \arrow{r} \& G \arrow{d} \\
			\& \& \& 0
		\end{tikzcd} 
	\end{align}
	where $G \in \scrA_{\mathsf{tor}}$, $E \in \scrA$ and $\calV \to \calF' \to E'$ is a $\calV$-stable pair, then there exists a unique way (up to a contractible space of choices) of finding an object $\calF$ together with the dashed arrows making the above diagram into an extension of stable pairs in the sense of Definition~\ref{def:extension_stable_pairs}. Since the Waldhausen construction $\calS_\bullet \bfPerfps(S)$ is a $2$-Segal space, we see that setting
	\begin{align}
		\calF \coloneqq \mathsf{fib}(\calF' \to G ) 
	\end{align}
	provides a unique way of filling the above diagram inside $\calS_3 \bfPerfps(S)$, and therefore inside $\calS_1^\ell \bfFlagPerfps^{(2),\dagger}(S;\calV)$. Since $\calS_1^\ell \bfP(S;\calV)$ is open inside $\calS_1^\ell \bfFlagPerfps^{(2),\dagger}(S;\calV)$, all we are left to do is to argue that in this situation $\calV \to \calF \to E$ is a $\calV$-stable pair. As this is guaranteed by Proposition~\ref{prop:stable_pairs_extensions}, the conclusion follows.
\end{proof}

We now look at the map
\begin{align}\label{eq:stable_pairs_action_lci}
	u_1^\ell\times \partial_1  \colon \calS^\ell_1 \bfP(S;\calV) \longrightarrow \bfCoh_{\mathsf{tor}}(S,\tau_\scrA)  \times \bfP(S;\calV) \ .
\end{align}
As consequence of Proposition~\ref{prop:flag_tangent_complex}, we obtain the following.
\begin{corollary}\label{cor:relativecotangentcomplex-stable-pairs}
	The tangent complex $\T$ of the map \eqref{eq:stable_pairs_action_lci} at a point $x \colon \Spec(A) \to \calS^\ell_1 \bfP(S;\calV)$ classifying an extension of stable pairs of the form \eqref{eq:stable_pairs_extensions} fits in the following natural fiber sequence:
	\begin{align}
		\T \longrightarrow  \R\Hom_{S_A}(\calF',\calF)[1]\oplus  \R\Hom_{S_A}(\calG,\calF')[1] \longrightarrow \R\Hom_{S_A}(\calF',\calF')[1] \ . 
	\end{align}
	In particular, letting $T \coloneqq \D(G)$ and $P' \coloneqq \D(\calF')[-1]$, we obtain an isomorphism
	\begin{align}
		\sfH^2(\T) \simeq \Ext^2_{S_A}(P',T) \ . 
	\end{align}
	Thus, if $T$ is a zero-dimensional coherent sheaf and $A$ is a field, we have $\sfH^2(\T) = 0$.
\end{corollary}

\begin{proof}
	The existence of the given fiber sequence is the exact consequence of Proposition~\ref{prop:flag_tangent_complex}. Passing to the long exact sequence of cohomology groups, we obtain
	\begin{align}
		\Ext^2_{S_A}(\calF', \calF)\oplus 	\Ext^2_{S_A}(G, \calF') \xlongrightarrow{\phi} & \Ext^2_{S_A}(\calF', \calF') \longrightarrow \sfH^2(\T) \longrightarrow \\
		& \longrightarrow \Ext^3_{S_A}(\calF',\calF)\oplus \Ext^3_{S_A}(G, \calF') \longrightarrow \Ext^3_{S_A}(\calF',\calF') \ .
	\end{align}
	Setting $P \coloneqq \D(\calF)[-1]$ and applying Theorem~\ref{thm:duality-tor-torsion-free}, we obtain
	\begin{align}
		\Ext^3_{S_A}(\calF', \calF') \simeq \Ext^3_{S_A}(P',P') \simeq 0 \qquad \text{and} \qquad \Ext^3_{S_A}(\calF', \calF) \simeq \Ext^3_{S_A}(P,P') \simeq 0 \ . 
	\end{align}
	Similarly,
	\begin{align}
		\Ext^3_{S_A}(G, \calF') \simeq \Ext^2_{S_A}(P',T) \ . 
	\end{align}
	The conclusion therefore follows if we can prove that the map $\phi$ is surjective.
	
	After applying Theorem~\ref{thm:duality-tor-torsion-free}, it takes the following form:
	\begin{align}
		\Ext^2_{S_A}(P,P')\oplus \Ext^1_{S_A}(P',T)  \longrightarrow \Ext^2_{S_A}(P',P') \ . 
	\end{align}
	We claim that the first component is surjective. To see this, start from the fiber sequence $\calF \to \calF' \to G$. Applying $\D(-)$ and rotating it, it induces the fiber sequence $P' \to P \to T$. Applying $\R\Hom_{S_A}(-,P')$ and passing to the long exact sequence of cohomology groups, we find
	\begin{align}
		\Ext^2_{S_A}(T,P') \longrightarrow \Ext^2_{S_A}(P,P') \longrightarrow \Ext^2_{S_A}(P',P') \longrightarrow \Ext^3_{S_A}(T,P') \simeq 0 \ , 
	\end{align}
	whence the conclusion.
	
	Finally, the last statement is an obvious consequence of Serre duality: since
	\begin{align}
		\Ext^2_{S_A}(P',T)^\vee \simeq \Hom_{S_A}(T, P' \otimes \omega_{S_A / A}) \ , 
	\end{align}
	we see that if $T$ is $0$-dimensional the purity of $P'$ implies the vanishing of the above group.
\end{proof}

\begin{theorem} \label{thm:stable_pairs_left_action}
	Let $\bfCoh_{0\textrm{-}\!\dim}(S)$ be the moduli stack of zero-dimensional coherent sheaves on $S$. Then, the pro-$\infty$-category $\catCohb_{\mathsf{pro}}(\bfP(S;\calV))$ carries a right categorical module structure over $\catCohb_{\mathsf{pro}}( \bfCoh_{0\textrm{-}\!\dim}(S) )$.
	
	Similarly, let $\bfD^\ast$ be a motivic formalism, and fix $\calA \in \CAlg(\bfD^\ast(\Spec(k)))$ and $\Gamma \subseteq \Pic(\bfD^\ast(\Spec(k)))$ such that Assumption~\ref{assumption:motivic_formalism} is satisfied. Then, the topological vector space $\HBMDGamma_0(\bfP(S;\calV);\calA)$ has the structure of a right $\HBMDGamma_0( \bfCoh_{0\textrm{-}\!\dim}(S);\calA)$-module. In particular,
	\begin{align}
		G_0( \bfP(S;\calV) )\quad \text{and} \quad \HBM_\ast( \bfP(S;\calV) )			
	\end{align}
	have the structures of a right $G_0( \bfCoh_{0\textrm{-}\!\dim}(S) )$-module and $\HBM_\ast( \bfCoh_{0\textrm{-}\!\dim}(S) )$-module, respectively.
\end{theorem}

\begin{proof}
	Consider 
	\begin{align}
		\bfH\coloneqq \bfCoh_{0\textrm{-}\!\dim}(S)=\bfCoh_{0\textrm{-}\!\dim}(S, \tau_{\mathsf{std}}\op) \qquad\textit{and}\qquad \bfM\coloneqq \bfP(S; \calV)\ .
	\end{align}
	First, note that the category of zero-dimensional sheaves on $S$ is a Serre subcategory of both $\catCoh(S)$ and $\scrA$ (the latter via duality - Theorem~\ref{thm:duality-tor-torsion-free}). Moreover, we will interpret $\bfCoh_{0\textrm{-}\!\dim}(S)$ as the open and closed substack of both $\bfCoh_{\mathsf{tor}}(S)$ and $\bfCoh_{\mathsf{tor}}(S, \tau_\scrA)$ defined by the condition that the torsion objects have zero first Chern class.
	
	The property of being a Serre subcategory implies that the assumption in Proposition~\ref{prop:2-Segal-algebra} is satisfied by $\bfH$. Furthermore, thanks to Proposition~\ref{prop:stable_pairs_extensions}, also condition \eqref{item:2-Segal-representation-(2)}  of Proposition~\ref{prop:2-Segal-representation} is satisfied by $\bfH, \bfM$.
	
	Now note that $\bfCoh_{0\textrm{-}\!\dim}(S)$ fits into the following pullback squares:
	\begin{align}
		\begin{tikzcd}[ampersand replacement=\&]
			\calS_2\bfPerf_\bfH(\scrC_S) \arrow{r} \arrow{d}{\partial_0 \times \partial_2} \& \calS_2\bfPerf_{\scrA_{\mathsf{tor}}}(\scrC_S) \arrow{d}{\partial_0 \times \partial_2} \\
			\bfCoh_{0\textrm{-}\!\dim}(S, \tau_{\mathsf{std}}\op)\times \bfCoh_{0\textrm{-}\!\dim}(S, \tau_{\mathsf{std}}\op) \arrow{r} \& \bfCoh_{\mathsf{tor}}(S,\tau_\scrA)\times \bfCoh_{\mathsf{tor}}(S,\tau_\scrA)
		\end{tikzcd}
	\end{align}
	and
	\begin{align} 
		\begin{tikzcd}[ampersand replacement=\&]
			\calS_2\bfPerf_\bfH(\scrC_S) \arrow{r} \arrow{d}{\partial_1} \& \calS_2\bfPerf_{\scrA_{\mathsf{tor}}}(\scrC_S) \arrow{d}{\partial_1} \\
			\bfCoh_{0\textrm{-}\!\dim}(S, \tau_{\mathsf{std}}\op)\arrow{r} \& \bfCoh_{\mathsf{tor}}(S,\tau_\scrA)
		\end{tikzcd}\ .
	\end{align}
	Here $\calS_\bullet \bfPerf_{\scrA_{\mathsf{tor}}}(\scrC_S)\coloneqq  \calS_\bullet \bfPerf_{\bfCoh_{\mathsf{tor}}(S,\tau_\scrA)}(\scrC_S)$ (cf.\ Remark~\ref{rem:relative-2-Segal}).  
	
	The right vertical maps are quasi-compact, finitely connected and derived lci, and locally rpas, respectively, as shown in the proof of Theorem~\ref{thm:action-perverse-torsion}. Then, the left vertical maps satisfy the same properties, hence the assumptions of Corollary~\ref{cor:COHA} are satisfied by $\bfCoh_{0\textrm{-}\!\dim}(S)$.
	
	Finally, the map
	\begin{align}
		u_1^\ell\times \varpi_1 \colon \calS^\ell_1 \bfFlagPerf_{\bfH, \bfM}^{(m),\dagger}(\scrC_S;\calV) \longrightarrow \bfCoh_{0\textrm{-}\!\dim}(S)  \times \bfP(S;\calV)  
	\end{align}
	fits into the pullback square
	\begin{align}
		\begin{tikzcd}[ampersand replacement=\&]
			\calS^\ell_1 \bfFlagPerf_{\bfH, \bfM}^{(m),\dagger}(\scrC_S;\calV) \arrow{d}{u_1^\ell\times \varpi_1} \arrow{r} \& \calS^\ell_1 \bfP(S;\calV)\arrow{d}{u_1^\ell \times \varpi_1}\\
			\bfCoh_{0\textrm{-}\!\dim}(S, \tau_{\mathsf{std}}\op) \times \bfP(S;\calV) \arrow{r}\&\bfCoh_{\mathsf{tor}}(S,\tau_\scrA)\times  \bfP(S;\calV)
		\end{tikzcd}\ .
	\end{align}
	Thanks to the computation of the relative cotangent complex of the right vertical map in Corollary~\ref{cor:relativecotangentcomplex-stable-pairs}, we see that the left vertical map is derived lci. Moreover, it is quasi-compact and finitely connected by Corollary~\ref{cor:linear-stack-l}. Using a similar argument and Lemma~\ref{lem:stable_pairs_action_properness}, one shows that 
	\begin{align}
		\varpi_0 \colon \calS^\ell_1 \bfFlagPerf_{\bfH, \bfM}^{(m),\dagger}(\scrC_S;\calV) \longrightarrow \bfP(S;\calV)
	\end{align}
	is locally rpas. 
	
	Thus, all the assumptions of Corollary~\ref{cor:COHA-representations} are satisfied. The assertion follows noticing that the duality (Theorem~\ref{thm:duality-tor-torsion-free}) exchanges left and right representations (a phenomenon seen already in Theorem~\ref{thm:swap}).
\end{proof}

\subsection{Extension of stable co-pairs by torsion perverse sheaves}

We now introduce the notion of extensions of stable co-pairs by torsion perverse sheaves on $S$.
\begin{definition}\label{def:extension_copair}
	Let $\calV$ be a locally free sheaf of finite rank on $S$. Let $\sfE = (\calF \to E \to \calV[2])$ be a $\calV$-stable co-pair and let $G \in \scrA^\upsilon_{\mathsf{t.f.}}$. An \textit{extension of $\sfE$ by $G$} is a commutative diagram $\E$ in $\catPerf(S)$
	\begin{align}\label{eq:stable_copair_extensions}
		\begin{tikzcd}[ampersand replacement = \&]
			0 \arrow{r} \& G \arrow{r} \arrow{d} \& \calF' \arrow{r} \arrow{d} \& E' \arrow{d} \\
			\& 0 \arrow{r} \& \calF \arrow{r} \arrow{d} \& E \arrow{d} \\
			\& \& 0 \arrow{r} \& \calV[2] \arrow{d} \\
			\& \& \& 0
		\end{tikzcd} \ ,
	\end{align}
	where every square is asked to be a pullback and where we further ask $\calF' \to E' \to \calV[2]$ to be a $\calV$-stable co-pair.
\end{definition}

\begin{notation}
	Let $\E$ be an extension of $\calV$-stable co-pairs of the form \eqref{eq:stable_copair_extensions}. Compatibly with the simplicial notation of the Part~\ref{part:Segal}, we set
	\begin{align}
		\varpi_0(\E) \coloneqq ( \calF \longrightarrow E \longrightarrow \calV[2] ) \qquad \text{and} \qquad \varpi_1(\E) \coloneqq (\calF' \longrightarrow E' \longrightarrow \calV[2] ) \ . 
	\end{align}
	We equally set
	\begin{align}
		u_1^r(\E) \coloneqq G \ . 
	\end{align}
\end{notation}

Using Corollary~\ref{cor:costable_pair_numerical} instead of Proposition~\ref{prop:stable_pairs_different_formulations}, the same proof of Proposition~\ref{prop:stable_pairs_extensions} yields:
\begin{proposition} \label{prop:stable_copairs_extensions}
	Let $\E$ be an extension of $\calV$-stable co-pairs of the form \eqref{eq:stable_copair_extensions}. Assume that $G \in \scrA_{\mathsf{tor}}$ and that $\varpi_1(\E) = (\calF' \to E' \to \calV[2])$ is a $\calV$-stable co-pair. Then $\varpi_0(\E) = (\calF \to E \to \calV[2])$ is a $\calV$-stable co-pair if and only if $E \in \scrA^\upsilon$.
\end{proposition}

\subsection{Left representations of torsion sheaves via stable co-pairs}\label{subsec:right_representation_stable_pair}

Recall that a $\calV$-stable pair is a fiber sequence $\calV \to \calF \to \calE$, where $\calF[1] \in \scrA_{\mathsf{tor}}$ and $\calE \in \scrA_{\mathsf{t.f.}}$. This gives rise to the projection map
\begin{align}
	\pi^{\mathsf{t.f.}}_0 \colon \bfP(S;\calV) \longrightarrow \bfCoh_{\mathsf{t.f.}}(S,\tau_\scrA) 
\end{align}
of Notation~\ref{notation:stable_pairs_to_perverse_tf}. Recall also that $\bfCoh_{\mathsf{tor}}(S,\tau_\scrA)$ acts both from the left and from the right on $\bfCoh_{\mathsf{t.f.}}(S,\tau_\scrA)$.

In the previous section we saw that the left action restricted to the substack $\bfCoh_{0\textrm{-}\!\dim}(S, \tau_{\mathsf{std}}\op) \subseteq \bfCoh_{\mathsf{tor}}(S,\tau_{\scrA})$ of zero-dimensional coherent sheaves, can be lifted to an action at the level of stable pairs. We are now going to analyze the right action; the main result of the section is that the same lifting is possible, and this time it will take place at the level of the whole $\bfCoh_{\mathsf{tor}}(S,\tau_\scrA)$.

We now introduce just as in the previous subsection the derived stack $\bfP^{\mathsf{c}}(S;\calV)$ of stable co-pairs. As the procedure is formally identical, we limit ourselves to explain the broad steps. To begin with, applying the \textit{left} version of Construction~\ref{construction:Vflags_generalized} to the $2$-Segal space $\calS_\bullet \bfPerfps(S)$ with $\sfV\coloneqq \calV[2]$, we obtain a relative $2$-Segal space
\begin{align}\label{eq:stable_copair_right_action_I}
	u_\bullet^r \colon \calS^r_\bullet \bfPerfps^\ddag(S;\calV)  \coloneqq \calS^r_\bullet \bfFlagPerfps^{(2),\dagger}(\scrC_S;\calV[2])  \longrightarrow \calS_\bullet \bfPerfps(S) \ .
\end{align}
Unraveling the definitions, we see that
\begin{align}
	\bfPerfps^\ddag(S;\calV) = \calS^r_0 \bfPerfps^\ddag(S;\calV) 
\end{align}
parametrizes fiber sequences of the form
\begin{align}
	\calF \longrightarrow E \longrightarrow \calV[2] \ .
\end{align}

On the other hand, $\calS^r_1 \bfPerfps^\ddag(S;\calV)$ parametrizes extensions of the form \eqref{eq:stable_copair_extensions}, where there is still no restriction whatsoever on the terms $\calG$, $\calF$, $\calF'$, $E$ and $E'$. We now define the derived moduli stack $\bfP^{\mathsf{c}}(S;\calV)$ of stable co-pairs as the fiber product
\begin{align}
	\begin{tikzcd}[ampersand replacement=\&]
		\bfP^{\mathsf{c}}(S;\calV) \arrow{r} \arrow{d} \& \bfPerfps^\ddag(S;\calV) \arrow{d}{\partial_2 \times \partial_1} \\
		\bfCoh_{\mathsf{t.f.}}(S,\tau_\scrA^\upsilon) \times \bfCoh_{\mathsf{tor}}(S,\tau_\scrA^\upsilon) \arrow{r} \& \bfPerfps(S) \times \bfPerfps(S) 
	\end{tikzcd} \ .
\end{align}

Then, Lemma~\ref{lem:costable_stable} immediately implies the following:
\begin{lemma}\label{lem:pair_copair_rotating}
	The self-equivalence
	\begin{align}
		\rho \colon \calS_2 \bfPerfps(S) \longrightarrow \calS_2 \bfPerfps(S) 
	\end{align}
	that sends a fiber sequence $E_1 \to E_2 \to E_3$ to the shifted-rotated sequence $E_2[1] \to E_3[1] \to E_1[2]$
	induces an equivalence
	\begin{align}
		\rho \colon \bfP(S;\calV) \simeq \bfP^{\mathsf c}(S;\calV) \ . 
	\end{align}
\end{lemma}

Applying Construction~\ref{construction:Hecke_pattern} to the relative $2$-Segal space \eqref{eq:stable_copair_right_action_I} with the choice of the \textit{boundary datum} given by the diagram
\begin{align}
	\begin{tikzcd}[ampersand replacement=\&]
		\bfP^{\mathsf{c}}(S;\calV) \arrow{r} \arrow{d} \& \Spec(\C) \arrow[equal]{d} \arrow{r}{0} \& \bfCoh_{\mathsf{t.f.}}(S,\tau_\scrA^\upsilon) \arrow{d} \\
		\bfPerfps^{\ddag}(S;\calV) \arrow{r} \& \Spec(\C) \arrow{r} \& \bfPerfps(S)
	\end{tikzcd} \ , 
\end{align}
we obtain a relative simplicial derived stack
\begin{align}\label{eq:stable_copairs_action_II}
	u^r_\bullet \colon \calS^r_\bullet \bfP^{\mathsf{c}}(S;\calV) \longrightarrow \calS_\bullet \bfCoh_{\mathsf{t.f.}}(S,\tau_{\scrA}^\upsilon) \ .
\end{align}
The same procedure that leads to Proposition~\ref{prop:stable_pair_action_2-Segal} yields:
\begin{proposition}
	The morphism \eqref{eq:stable_copairs_action_II} is a relative $2$-Segal space.
\end{proposition}
In order to convey the information contained in the $2$-Segal property into an actual right representation, we need to verify the assumptions of Corollary~\ref{cor:COHA-representations-right}. We start with the property of being locally rpas.

\begin{notation}
	Let
	\begin{align}
		\pi^{\upsilon\textrm{-}\mathsf{t.f.}}_0 \colon \bfP^{\mathsf{c}}(S;\calV) \longrightarrow \bfCoh_{\mathsf{t.f.}}(S,\tau_\scrA^\upsilon) 
	\end{align}
	be the morphism that sends a stable co-pair $\sfE = (\calF \to \calE \to \calV[2])$ to $\partial_2(\sfE) = \calF$. It can be promoted to a morphism of simplicial derived stacks
	\begin{align}
		\pi^{\upsilon\textrm{-}\mathsf{t.f.}}_\bullet \colon \calS^r_\bullet \bfP^{\mathsf{c}}(S;\calV) \longrightarrow \calS_\bullet \bfCoh_{\mathsf{t.f.}}(S,\tau_\scrA^\upsilon) \ , 
	\end{align}
	and the level $1$ component $\pi^{\upsilon\textrm{-}\mathsf{t.f.}}_1$ can be explicitly described as the map sending an extension of $\calV$-stable co-pairs of the form \eqref{eq:stable_copair_extensions} to the sub-extension $\calG \to \calF' \to \calF$.
\end{notation}

\begin{lemma}\label{lem:stable-co-pairs-properness}
	The square
	\begin{align}
		\begin{tikzcd}[ampersand replacement=\&]
			\calS^r_1 \bfP^{\mathsf{c}}(S;\calV) \arrow{r}{\pi_1^{\upsilon\textrm{-}\mathsf{t.f.}}} \arrow{d}{\varpi_1} \& \calS^r_1 	\bfFlagCoh^{(1)}_{\mathsf{t.f.},\mathsf{tor}}(S,\tau_\scrA^\upsilon) \arrow{d}{\varpi_1} \\
			\bfP^{\mathsf{c}}(S;\calV) \arrow{r}{\pi_1^{\upsilon\textrm{-}\mathsf{t.f.}}} \& \bfCoh_{\mathsf{t.f.}}(S,\tau_{\scrA}^\upsilon)
		\end{tikzcd} 
	\end{align}
	is a pullback. In particular, the left vertical map is locally rpas.
\end{lemma}

\begin{proof}
	The right vertical morphism is locally rpas by combining Theorem~\ref{thm:duality-tor-torsion-free} and the proof that the map \eqref{eq:properness-standard} in the proof of Theorem~\ref{thm:action-torsion} is locally rpas.
	
	It is then enough to prove the first half of the statement, i.e., that the given square is a pullback. Unwinding the definitions, we have to prove that given a solid diagram of the form
	\begin{align}
		\begin{tikzcd}[ampersand replacement=\&]
			0 \arrow{r} \& G \arrow{r} \arrow{d} \& \calF' \arrow{r} \arrow{d} \& E' \arrow[dashed]{d} \arrow[bend left = 30pt]{dd} \\
			\& 0 \arrow{r} \& \calF \arrow[dashed]{r} \arrow{d} \& E \arrow[dashed]{d} \\
			\& \& 0 \arrow{r} \& \calV[2] \arrow{d} \\
			\& \& \& 0
		\end{tikzcd} \ , 
	\end{align}
	where $G \in \scrA_{\mathsf{t.f.}}^\upsilon$, $\calF' \to E' \to \calV[2]$ is a $\calV$-stable co-pair and $\calF \in \scrA^\upsilon$, there is a unique way (up to a contractible space of choices) of defining $E$ and the dashed arrows making the above diagram into an extension of $\calV$-stable co-pairs. Since the Waldhausen construction $\calS_\bullet \bfPerfps(S)$ satisfies the $2$-Segal property, we see that setting
	\begin{align}
		E \coloneqq \mathsf{cofib}( G \longrightarrow E' ) 
	\end{align}
	provides the unique filling inside $\calS_3 \bfPerfps(S)$ (and hence in $\calS_1^r \bfFlagPerfps^\ddag(S;\calV)$) we were looking for. To complete the proof, we only have to check that in this situation $\calF \to E \to \calV[2]$ is automatically a $\calV$-stable co-pair. To see this, observe first that, since both $\calF$ and $\calV[2]$ belong to $\scrA^\upsilon$, we automatically have that $E \in \scrA^\upsilon$ as well. The fiber sequence $G \to E' \to E$ is therefore a short exact sequence in $\scrA^\upsilon$, and since $E' \in \scrA^\upsilon_{\mathsf{tor}}$, it follows that $E \in \scrA^{\upsilon}_{\mathsf{tor}}$ as well. Finally, since $\mathsf{rk}(G) = 0$, we deduce $\mathsf{rk}(E) = \mathsf{rk}(E') = - \mathsf{rk}(\calV)$, so the conclusion follows from Proposition~\ref{prop:stable_copairs_extensions}.
\end{proof}

We now turn to the map
\begin{align}\label{eq:stable-co-pairs-pullback}
	\varpi_0\times u_1^r \colon \calS^r_1 \bfP^{\mathsf{c}}(S;\calV) \longrightarrow \bfP^{\mathsf{c}}(S;\calV)\times \bfCoh_{\mathsf{t.f.}}(S,\tau_\scrA^\upsilon)  \ . 
\end{align}
Contrary to what happens in the previous subsection, we have the following.
\begin{proposition}\label{prop:stable-co-pairs-derived-lci}
	The square
	\begin{align}
		\begin{tikzcd}[ampersand replacement=\&]
			\calS^r_1 \bfP^{\mathsf c}(S;\calV) \arrow{r} \arrow{d}{\varpi_0 \times u_1^r} \arrow{r} \& \calS^r_1 \bfFlagCoh_{\mathsf{t.f.},\mathsf{tor}}^{(1)}(S,\tau_\scrA^\upsilon) 	\arrow{d}{\varpi_0\times u_1^r} \\
			\bfP^{\mathsf c}(S;\calV) \times \bfCoh_{\mathsf{t.f.}}(S,\tau_\scrA^\upsilon) \arrow{r} \& \bfCoh_{\mathsf{tor}}(S,\tau_\scrA^\upsilon)\times \bfCoh_{\mathsf{t.f.}}(S,\tau_\scrA^\upsilon)
		\end{tikzcd} 
	\end{align}
	is a pullback. In particular, the left vertical map is quasi-compact, finitely connected and derived lci.
\end{proposition}

\begin{proof}
	It has already been verified in Theorem~\ref{thm:left-right-action} that the right vertical map is quasi-compact, finitely connected and derived lci. We are therefore left to check that the above square is a pullback. Equivalently, we have to check that the square 
	\begin{align}
		\begin{tikzcd}[ampersand replacement=\&]
			\calS^r_1 \bfP^{\mathsf c}(S;\calV) \arrow{r} \arrow{d}{\varpi_0} \arrow{r} \& \calS^r_1 \bfFlagCoh_{\mathsf{t.f.},\mathsf{tor}}^{(1)}(S,\tau_\scrA^\upsilon) \arrow{d}{\varpi_0} \\
			\bfP^{\mathsf c}(S;\calV) \arrow{r} \& \bfCoh_{\mathsf{tor}}(S,\tau_\scrA^\upsilon)
		\end{tikzcd}
	\end{align}
	is a pullback. Unraveling the definitions, we are called to show that given a solid diagram of the form
	\begin{align}
		\begin{tikzcd}[ampersand replacement=\&]
			0 \arrow{r} \& G \arrow[dashed]{r} \arrow[bend left = 25pt]{rr} \arrow{d} \& \calF' \arrow[dashed]{r} \arrow[dashed]{d} \& E' \arrow{d} \\
			\& 0 \arrow{r} \& \calF \arrow{r} \arrow{d} \& E \arrow{d} \\
			\& \& 0 \arrow{r} \& \calV[2] \arrow{d} \\
			\& \& \& 0
		\end{tikzcd} \ , 
	\end{align}
	where $G \in \scrA_{\mathsf{t.f.}}^\upsilon$, $\calF \to E \to \calV[2]$ is a $\calV$-stable co-pair and $E' \in \scrA_{\mathsf{tor}}^\upsilon$, there is a unique way (up to a contractible space of choices) to define $\calF'$ together with the dashed morphisms making the above diagram into an extension of $\calV$-stable co-pairs in the sense of Definition~\ref{def:extension_copair}.
	As usual, the $2$-Segal property satisfied by the Waldhausen construction $\calS_\bullet \bfPerfps(S)$ shows that the only possible choice is to take
	\begin{align}
		\calF' \coloneqq \mathsf{fib}( E' \to \calV[2] ) \ , 
	\end{align}
	and we are left to check that in this situation the fiber sequence $\calF' \to E' \to \calV[2]$ is a $\calV$-stable co-pair. Since both $G$ and $\calF$ belong to $\scrA^\upsilon_{\mathsf{t.f.}} = \scrA_{\mathsf{tor}}$, the same goes for $\calF'$. Moreover, $E'$ belongs to $\scrA_{\mathsf{tor}}^\upsilon$ by assumption and since $\mathsf{rk}(G) = 0$, we have
	\begin{align}
		\mathsf{rk}(E') = \mathsf{rk}(E) = - \mathsf{rk}(\calV) \ . 
	\end{align}
	Thus, the conclusion follows from Proposition~\ref{prop:stable_copairs_extensions}.
\end{proof}

The following is the main result of this section.
\begin{theorem}\label{thm:stable_copairs_left_action}
	The pro-$\infty$-category $\catCohb_{\mathsf{pro}}( \bfP^{\mathsf{c}}(S;\calV) )$  has the structure of a right categorical module over the $\E_1$-monoidal $\infty$-category $\catCohb_{\mathsf{pro}}( \bfCoh_{\mathsf{t.f.}}(\scrA^\upsilon)  )$.
	
	Similarly, let $\bfD^\ast$ be a motivic formalism, and fix $\calA \in \CAlg(\bfD^\ast(\Spec(k)))$ and $\Gamma \subseteq \Pic(\bfD^\ast(\Spec(k)))$ such that Assumption~\ref{assumption:motivic_formalism} is satisfied. Then, the topological vector space $\HBMDGamma_0( \bfP^{\mathsf{c}}(S;\calV);\calA)$ has the structure of a right $\HBMDGamma_0( \bfCoh_{\mathsf{t.f.}}(\scrA^\upsilon);\calA)$-module. In particular,
	\begin{align}
		G_0( \bfP^{\mathsf{c}}(S;\calV) )\quad \text{and} \quad \HBM_\ast( \bfP^{\mathsf{c}}(S;\calV) )			
	\end{align}
	have the structures of a right $G_0( \bfCoh_{\mathsf{t.f.}}(\scrA^\upsilon) )$-module and $\HBM_\ast( \bfCoh_{\mathsf{t.f.}}(\scrA^\upsilon) )$-module, respectively.
\end{theorem}

\begin{proof}
	We need to check the assumptions of Corollary~\ref{cor:COHA-representations-right} for 
	\begin{align}
		\bfH\coloneqq \bfCoh_{\mathsf{t.f.}}(S, \tau_\scrA^\upsilon) \quad \text{and} \quad \bfM\coloneqq \bfP^{\mathsf{c}}(S;\calV) \ .
	\end{align}
	First, note that the assumptions of Corollary~\ref{cor:COHA} are satisfied by $\bfH$, since in Theorem~\ref{thm:left-right-action} we observed that there is a canonical equivalence
	\begin{align}
		\calS_\bullet \bfCoh_{\mathsf{tor}}(\scrA) \simeq \calS_\bullet \bfCoh_{\mathsf{t.f.}}(\scrA^\upsilon) \ ,
	\end{align}
	and the above conditions are satisfied by $\bfCoh_{\mathsf{tor}}(\scrA)$, as already shown in Theorem~\ref{thm:action-perverse-torsion}. Moreover, the map
	\begin{align}
		\varpi_0 \times u_1^r \colon \calS^r_1 \bfP^{\mathsf{c}}(S;\calV) \longrightarrow  \bfP^{\mathsf{c}}(S;\calV)\times \bfCoh_{\mathsf{t.f.}}(S,\tau_\scrA^\upsilon) \ . 
	\end{align}
	is quasi-compact, finitely connected and derived lci by Proposition~\ref{prop:stable-co-pairs-derived-lci}, while the map 
	\begin{align}
		\varpi_1 \colon \calS^r_1 \bfP^{\mathsf{c}}(S;\calV)  \longrightarrow \bfP^{\mathsf{c}}(S;\calV) 
	\end{align}
	is locally rpas, as shown in Lemma~\ref{lem:stable-co-pairs-properness}. Thus, all assumptions of Corollary~\ref{cor:COHA-representations-right} are verified and the assertion follows.
\end{proof}

\begin{corollary}\label{cor:action-torsion-stable-pairs} 
	The pro-$\infty$-category $\catCohb_{\mathsf{pro}}(\bfP(S;\calV))$ carries a right categorical module structure over $\catCohb(\bfCoh_{\mathsf{tor}}(\scrA))$.
	
	Similarly, let $\bfD^\ast$ be a motivic formalism, and fix $\calA \in \CAlg(\bfD^\ast(\Spec(k)))$ and $\Gamma \subseteq \Pic(\bfD^\ast(\Spec(k)))$ such that Assumption~\ref{assumption:motivic_formalism} is satisfied. Then, the topological vector space $\HBMDGamma_0( \bfP(S;\calV);\calA)$ has the structure of a right $\HBMDGamma_0( \bfCoh_{\mathsf{tor}}(\scrA);\calA)$-module. In particular,
	\begin{align}
		G_0( \bfP(S;\calV) )\quad \text{and} \quad \HBM_\ast( \bfP(S;\calV) )			
	\end{align}
	have the structures of a right $G_0( \bfCoh_{\mathsf{tor}}(\scrA) )$-module and $\HBM_\ast( \bfCoh_{\mathsf{tor}}(\scrA) )$-module, respectively.
\end{corollary}

\begin{proof}
	First, Theorem~\ref{thm:left-right-action} yields a canonical equivalence
	\begin{align}
		\calS_\bullet \bfCoh_{\mathsf{tor}}(\scrA) \simeq \calS_\bullet \bfCoh_{\mathsf{t.f.}}(\scrA^\upsilon) \ ,
	\end{align}
	while Lemma~\ref{lem:pair_copair_rotating} yields an equivalence of derived stacks
	\begin{align}
		\rho \colon \bfP(S;\calV) \longrightarrow \bfP^{\mathsf{c}}(S;\calV) \ . 
	\end{align}
	Via this equivalence, the right $2$-Segal action of $\bfCoh_{\mathsf{tor}}(\scrA)$ on $\bfP^{\mathsf{c}}(S;\calV)$ is transferred to a right $2$-Segal action of $\bfCoh_{\mathsf{tor}}(\scrA)$ on $\bfP(S;\calV)$. The statements for the categorical, $G$-theoretical and Borel-Moore homology actions follow automatically.
\end{proof}

Since $\bfCoh_{\mathsf{tor}}(S)\op\simeq \bfCoh_{\mathsf{tor}}(\scrA)$, we obtain:
\begin{corollary}\label{cor:stable-pairs-left-action}
	The pro-$\infty$-category $\catCohb_{\mathsf{pro}}(\bfP(S;\calV))$ carries a left categorical module structure over $\catCohb(\bfCoh_{\mathsf{tor}}(S))$.
	
	Similarly, let $\bfD^\ast$ be a motivic formalism, and fix $\calA \in \CAlg(\bfD^\ast(\Spec(k)))$ and $\Gamma \subseteq \Pic(\bfD^\ast(\Spec(k)))$ such that Assumption~\ref{assumption:motivic_formalism} is satisfied. Then, the topological vector space $\HBMDGamma_0( \bfP(S;\calV);\calA)$ has the structure of a left $\HBMDGamma_0( \bfCoh_{\mathsf{tor}}(S);\calA)$-module. In particular,
	\begin{align}
		G_0( \bfP(S;\calV) )\quad \text{and} \quad \HBM_\ast( \bfP(S;\calV) )			
	\end{align}
	have the structures of a left $G_0( \bfCoh_{\mathsf{tor}}(S) )$-module and $\HBM_\ast( \bfCoh_{\mathsf{tor}}(S) )$-module, respectively.
\end{corollary}

Theorem~\ref{thm:stable_pairs_left_action} and Corollary \ref{cor:stable-pairs-left-action} yield the following.
\begin{corollary}\label{cor:action-zero-stable-pairs}
	The pro-$\infty$-category $\catCohb_{\mathsf{pro}}( \bfP(S;\calV) )$  has the structure of a left and right categorical module over the $\E_1$-monoidal $\infty$-category $\catCohb_{\mathsf{pro}}( \bfCoh_{0\textrm{-}\!\dim}(S) )$.
	
	Similarly, let $\bfD^\ast$ be a motivic formalism, and fix $\calA \in \CAlg(\bfD^\ast(\Spec(k)))$ and $\Gamma \subseteq \Pic(\bfD^\ast(\Spec(k)))$ such that Assumption~\ref{assumption:motivic_formalism} is satisfied. Then, the topological vector space $\HBMDGamma_0( \bfP(S;\calV);\calA)$ has both the structure of a left and a right $\HBMDGamma_0( \bfCoh_{0\textrm{-}\!\dim}(S);\calA)$-module. In particular,
	\begin{align}
		G_0( \bfP(S;\calV) )\quad \text{and} \quad \HBM_\ast( \bfP(S;\calV) )			
	\end{align}
	have both the structures of a left and a right $G_0( \bfCoh_{0\textrm{-}\!\dim}(S) )$-module and $\HBM_\ast( \bfCoh_{0\textrm{-}\!\dim}(S) )$-module, respectively.
\end{corollary}

\begin{proof}
	Using Corollary~\ref{cor:action-torsion-stable-pairs} we obtain a left categorical module structure of the pro-$\infty$-category $\catCohb_{\mathsf{pro}}( \bfP(S;\calV) )$ over $\catCohb_{\mathsf{pro}}( \bfCoh_{0\textrm{-}\!\dim}(S) )$ via the inclusion $\bfCoh_{0\textrm{-}\!\dim}(S, \tau_{\mathsf{std}}\op)\subset \bfCoh_{\mathsf{tor}}(S, \tau_\scrA)$, while the right (categorical) structure is constructed in Theorem~\ref{thm:stable_pairs_left_action}. Similar arguments yield the claim for motivic Borel-Moore homology.
\end{proof}	

\begin{remark}
	In the local surface case, Toda constructed a right categorical module structure of $\catCohb_{\mathsf{pro}}$ of Pandharipande-Thomas moduli spaces of stable pairs over the categorical Hall algebra of zero-dimensional sheaves (cf.\ \cite[\S4]{Toda_Hall_Categorical_DT}). In this case, there is no left categorical module structure because of a wall-crossing phenomenon which does not appear in our two-dimensional case.
\end{remark}

\begin{rem}
	In \cite{Toda_Hall_Categorical_DT}, Toda constructed a categorical right and left actions of the CatHA generated by $\catCohb_{\mathsf{pro}}( \bfCoh_{\mathsf{pt}}(S;1))$ (introduced in \S\ref{subsubsec:negut}) and computed the commutation relations between the two actions in K-theory (cf.\ \cite[Proposition~6.7]{Toda_Hall_Categorical_DT}).\hfill $\triangle$
\end{rem}

Consider $\calV=\scrO_S$. Then $\trunc{\bfP(S; \scrO_S)}$ is equivalent to the moduli space $\calP(S)$ of Pandharipande-Thomas stable pairs on $S$, which is a projective scheme.\footnote{\cite[Proposition~5.2]{PT_BPS} shows that, for a K3 surface $S$, the moduli space of stable pairs of fixed Chern class $\beta$ and Euler characteristic $n$ is smooth for $\beta\in \mathsf{NS}(S)$ irreducible.} We have the following:
\begin{corollary}\label{cor:action-zero-PT-stable-pairs}
	The stable pro-$\infty$-category $\catCohb_{\mathsf{pro}}( \bfP(S) )$ has the structure of a left and a right categorical module over $\catCohb_{\mathsf{pro}}( \bfCoh_{0\textrm{-}\!\dim}(S) )$. 
	
	Similarly, let $\bfD^\ast$ be a motivic formalism, and fix $\calA \in \CAlg(\bfD^\ast(\Spec(k)))$ and $\Gamma \subseteq \Pic(\bfD^\ast(\Spec(k)))$ such that Assumption~\ref{assumption:motivic_formalism} is satisfied. Then, the topological vector space $\HBMDGamma_0( \calP(S);\calA)$ has both the structure of a left and a right $\HBMDGamma_0( \bfCoh_{0\textrm{-}\!\dim}(S);\calA)$-module. In particular,
	\begin{align}
		G_0( \calP(S) )\quad \text{and} \quad \HBM_\ast( \calP(S) )			
	\end{align}
	have both the structures of a left and a right $G_0( \bfCoh_{0\textrm{-}\!\dim}(S) )$-module and $\HBM_\ast( \bfCoh_{0\textrm{-}\!\dim}(S) )$-module, respectively.
\end{corollary}

\subsubsection{Duality and representations}\label{subsec:duality_and_representations}

Theorem~\ref{thm:duality_stable_pairs_categorical} has an immediate counterpart at the level of moduli stacks. Indeed, mimicking the construction of $\bfP(S;\calV)$ performed in \S\ref{subsec:stable_pairs_left_rep}, we introduce the derived stack $\bfP_\scrB(S;\calV)$ as the fiber product
\begin{align}
	\begin{tikzcd}[ampersand replacement=\&]
		\bfP_\scrB(S;\D(\calV)[-2]) \arrow{r} \arrow{d} \& \bfPerfps^\dagger(S;\D(\calV)[-2]) \arrow{d}{\partial_1 \times \partial_0} \\
		\bfCoh_{\mathsf{t.f.}}(S) \times \bfCoh_{\mathsf{tor}}(S) \arrow{r} \& \bfPerfps(S) \times \bfPerfps(S) 
	\end{tikzcd} \ .
\end{align}

\begin{remark}\label{rem:1_pure}
	Let $\bfCoh_{1\textrm{-}\mathsf{pure}}(S)$ be the open substack of $\bfCoh_{\mathsf{tor}}(S)$ parametrizing pure $1$-dimen\-sio\-nal coherent sheaves. Combining Proposition~\ref{prop:stable_pairs_different_formulations}, Theorem~\ref{thm:duality-tor-torsion-free} and Remark~\ref{rem:duality_stable_pairs_categorical} we deduce that the natural projection $\bfP_\scrB(S;\D(\calV)[-2]) \to \bfCoh_{\mathsf{tor}}(S)$ factors through $\bfCoh_{1\textrm{-}\mathsf{pure}}(S)$. This leads to the following alternative description of $\bfP_\scrB(S;\calV)$: let $\bfP_{1\textrm{-}\mathsf{pure}}(S)$ be the fiber product
	\begin{align}
		\begin{tikzcd}[ampersand replacement=\&]
			\bfP_{1\textrm{-}\mathsf{pure}}(S;\D(\calV)[-2]) \arrow{r} \arrow{d} \& \bfPerfps^\dagger(S;\D(\calV)[-2]) \arrow{d}{\partial_0} \\
			\bfCoh_{1\textrm{-}\mathsf{pure}}(S) \arrow{r} \& \bfPerfps(S) 
		\end{tikzcd} \ .
	\end{align}
	Concretely, $\bfP_{1\textrm{-}\mathsf{pure}}(S;\D(\calV)[-2])$ parametrizes extensions of the form $\D(\calV)[-2] \to \calE \to \calF$, where $\calF$ is purely $1$-dimensional. Observe that this automatically implies that $\calE \in \bfCoh(S)$. At this point, it follows that
	\begin{align}
		\begin{tikzcd}[ampersand replacement=\&]
			\bfP_\scrB(S;\D(\calV)[-2]) \arrow{r} \arrow{d} \& \bfP_{1\textrm{-}\mathsf{pure}}(S;\D(\calV)[-2]) \arrow{d}{\partial_1} \\
			\bfCoh_{\mathsf{t.f.}}(S) \arrow{r} \& \bfPerfps(S)
		\end{tikzcd} 
	\end{align}
	is a fiber product, thus realizing $\bfP_\scrB(S;\D(\calV)[-2])$ as an open substack inside $\bfP_{1\textrm{-}\mathsf{pure}}(S;\D(\calV)[-2])$.
	\end{remark}

\begin{remark}\label{rem:relative_Hilbert}
	When $\calV = \scrO_S$, we can provide a more explicit description of (the underived truncation of) $\bfP_\scrB(S;\D(\calV)[-2])$. Let $\sfH_{1\textrm{-}\mathsf{pure}}(S)$ be the Hilbert scheme parametrizing pure one-dimensional subschemes $C \subset S$. Let $\calC\subset S\times \sfH_{1\textrm{-}\mathsf{pure}}(S)$ be the universal curve and consider the (underived) relative \textit{Hilbert scheme} $\sfHilb(\calC/\sfH_{1\textrm{-}\mathsf{pure}}(S))$ of $\sfH_{1\textrm{-}\mathsf{pure}}(S)$-flat families of zero-dimen\-sio\-nal quotients of $\scrO_{\calC}$. Then \cite[Proposition~B.8]{PT_BPS} yields a canonical identification
	\begin{align}
		\trunc{\bfP_\scrB(S;\scrO_S)} \simeq \sfHilb(\calC/\sfH_{1\textrm{-}\mathsf{pure}}(S)) \ . 
	\end{align}
\end{remark}

At this point, Theorem~\ref{thm:duality_stable_pairs_categorical} immediately implies:
\begin{theorem}\label{thm:hilbert-action}
	There is a canonical equivalence
	\begin{align}
		\bfP(S;\calV) \simeq \bfP_\scrB(S;\D(\calV)[-2]) \ . 
	\end{align}
	Thus, $\catCohb_{\mathsf{pro}}( \bfP_\scrB(S;\D(\calV)[-2]) )$ is a right categorical module over $\catCohb_{\mathsf{pro}}( \bfCoh_{\mathsf{tor}}(S) )$ and a left categorical module over $\catCohb_{\mathsf{pro}}( \bfCoh_{0\textrm{-}\!\dim}(S) )$.
	Similar statements hold in motivic Borel-Moore homology.
\end{theorem}

\begin{proof}
	For every affine derived scheme $T = \Spec(A)$, let $p_S \colon T \times S \to S$ be the canonical projection. There is a canonical equivalence
	\begin{align}
		\catPerf^\dagger(T \times S; p_S^\ast(\calV))\op \simeq \catPerf^\dagger(T \times S; p_S^\ast(\D(\calV)[-2])) \ , 
	\end{align}
	that sends a $p_S^\ast(\calV)$-extension $p_S^\ast(\calV) \to \calF \to \calE$ to the $p_S^\ast(\D(\calV)[-2])$-extension $\D(\calV)[-2] \to \D(\calE)[-1] \to \D(\calF)[-1]$. Passing to the maximal $\infty$-groupoids and using the canonical equivalence
	\begin{align}
		\mathsf{inv} \colon \catPerf^\dagger(T \times S; p_S^\ast(\calV))^\simeq \simeq \big( \catPerf^\dagger(T \times S; p_S^\ast(\calV))\op \big)^\simeq \ , 
	\end{align}
	this yields a canonical equivalence of derived stacks
	\begin{align}
		\bfPerf^\dagger(S;\calV) \simeq \bfPerf^\dagger(S;\D(\calV)[-2]) \ . 
	\end{align}
	Theorem~\ref{thm:duality_stable_pairs_categorical} and Remark~\ref{rem:duality_stable_pairs_categorical} imply that this equivalence restricts to an equivalence $\bfP(S;\calV) \simeq \bfP_\scrB(S;\D(\calV)[-2])$. The existence of the actions at the categorical (resp.\ motivic Borel-Moore homology) level is then a direct consequence of Theorems~\ref{thm:stable_pairs_left_action} and \ref{thm:stable_copairs_left_action}.
\end{proof}

\begin{rem}
	Rather than transferring the left and the right action via the equivalence $\bfP(S;\calV) \simeq \bfP_\scrB(S;\D(\calV)[-2])$, it would be possible to independently construct, mimicking what was done in \S\ref{subsec:stable_pairs_left_rep} and \S\ref{subsec:right_representation_stable_pair} for the $t$-structure $\tau_\scrB$ and the torsion pair $(\catCoh_{\mathsf{t.f.}}(S)[1], \catCoh_{\mathsf{tor}}(S))$. In this case, the former theorem would be upgraded, with no additional cost, to state that the equivalence $\bfP(S;\calV) \simeq \bfP_\scrB(S;\D(\calV)[-2])$ is compatible with the natural actions on both sides.\hfill $\triangle$
\end{rem}

\begin{corollary}\label{cor:hilbert-action}
	With respect to the notations of Remark~\ref{rem:relative_Hilbert}, $G_0( \sfHilb(\calC/\sfH_{1\textrm{-}\mathsf{pure}}(S)) )$ is a right module over $G_0( \bfCoh_{\mathsf{tor}}(S) )$ and a left module over $G_0( \bfCoh_{0\textrm{-}\!\dim}(S) )$.
	Similar statements hold for Borel-Moore homology.
\end{corollary}

\begin{proof}
	Since $G_0$ and Borel-Moore homology are insensitive to the derived structure, this immediately follows combining Theorem~\ref{thm:hilbert-action} with Remark~\ref{rem:relative_Hilbert}.
\end{proof}

\section{Preprojective COHA of a quiver, its categorification, and its representation via quiver varieties}\label{sec:COHA-quiver-Yangian}

In this section, we categorify the preprojective COHA of a quiver within our formalism. Moreover, we define a (categorical) representation in terms of Nakajima quiver varieties.

\subsection{Preprojective algebras}

Fix an algebraically closed field $k$ of characteristic zero. Let $\calQ$ be a \textit{quiver}, i.e., an oriented graph, with a finite vertex set $\calQ_0$ and a finite arrow set $\calQ_1$. For any arrow $a\in \calQ_1$, we denote by $\sfs(a)$ the \textit{source} of $a$ and by $\sft(a)$ the \textit{target} of $a$. The \textit{path algebra} $k \calQ$ of the quiver $\calQ$ is the associative algebra with basis all possible paths of length $\ell\geq 0$ of $\calQ$, endowed with the multiplication given by concatenation of paths, whenever possible, otherwise zero.

The \textit{double} $\calQ^{\mathsf{db}}$ of $\calQ$  is the quiver that has the same vertex set as $\calQ$ and whose set of arrows $\calQ_1^{\mathsf{db}}$ is a disjoint union of the set $\calQ_1$ of arrows of $\calQ$ and of the set 
\begin{align}
	\calQ_1^{\mathsf{opp}}\coloneqq \{a^\ast\,\vert\, a\in \calQ_1\}
\end{align}
consisting of an arrow $a^\ast$ for any arrow $a\in \calQ_1$, with the reverse orientation (i.e., $\sfs(a^\ast)=\sft(a)$ and $\sft(a^\ast)=\sfs(a)$).

\begin{definition}[{cf.\ \cite[\S4.1.4]{BCS_Quiver}}]\label{def:deformed-preprojective-algebra}
	The \textit{derived preprojective algebra} $\Pi_\calQ$ is the derived push-out 
	\begin{align}
		\Pi_\calQ\coloneqq \scrG_2(k\calQ)= k\calQ^{\mathsf{db}}\overset{\LL}{\underset{k[x]}{\coprod}} k\ ,
	\end{align}
	where the morphism $k[x]\to k\calQ^{\mathsf{db}}$ sends $x$ to 
	\begin{align}\label{eq:preprojective}
		\sum_{a\in \calQ_1}\, (a a^\ast - a^\ast a)\in k\calQ^{\mathsf{db}}\ .
	\end{align}
	
	The \textit{preprojective algebra} $\trunc{\Pi}_\calQ$ of $\calQ$ is the $0$-th cohomology of $\Pi_\calQ$, i.e., the quotient of the path algebra $k\calQ^{\mathsf{db}}$ by the two-sided ideal generated by the element \eqref{eq:preprojective}.
\end{definition}

\begin{remark}
	The notion of derived preprojective algebra was originally introduced by Ginzburg \cite{Ginzburg_CY}. It is equivalent to the notion of \textit{2-Calabi-Yau completion of $k\calQ$} by Keller \cite{Keller_Deformed}. As pointed out for example in \textit{loc.cit.}, when $\calQ$ is not a finite Dynkin quiver, $\Pi_\calQ$ is quasi-isomorphic to its preprojective algebra $\trunc{\Pi}_\calQ$. 
\end{remark}

\subsection{COHA of the preprojective algebra of a quiver and its representation via quiver varieties}

\begin{definition}[{\cite{CB_Moment}}]\label{def:CB-quiver}
	Let $w\in \N^{\calQ_0}$. The \textit{Crawley-Boevey quiver} associated to $\calQ$ is the quiver $\calQ^w$ whose set of vertices is given by $\calQ_0\sqcup \{\infty\}$, where $\infty$ is a new additional vertex, and whose set of arrows is the disjoint union of $\calQ_1$ and the set of $w_i$ additional edges of the form $\beta_i\colon i\to \infty$ for each vertex $i\in \calQ_0$.
\end{definition}

Fix a quiver $\calQ$ and fix $w\in \N^{\calQ_0}$ such that $\calQ^w$ is not a finite Dynkin quiver. Consider the abelian category $\Modd(\Pi_{\calQ^w})$ of finitely generated right modules (i.e., \textit{representations}) of $\Pi_{\calQ^w}$. Following e.g. \cite[\S5]{Ginzburg_Quiver}, a representation $\overline{M}$ corresponds to a pair of vector spaces over $k$
\begin{align}
	M\coloneqq \bigoplus_{i\in \calQ_0} M_i \quad \text{and} \quad M_\infty \ ,
\end{align}
and a quadruple $(x, y, \alpha, \beta)$ of collections of $k$-linear maps
\begin{align}
	x_a \colon M_{\sfs(a)} \longrightarrow M_{\sft(a)} \; ,\ y_a \colon M_{\sft(a)} \longrightarrow M_{\sfs(a)}  \; , \ \alpha_i \colon M_\infty\longrightarrow M_i \ \text{ and } \ \beta_i\colon M_i \to M_\infty
\end{align}
satisfying the preprojective relations. Consider the full subcategory $\scrT$ of $\Modd(\Pi_{\calQ^w})$ consisting of those representations $\overline{M}$ for which $M_\infty = 0$. This is a Serre subcategory. Moreover, it is equivalent to the category $\Modd(\Pi_\calQ)$ of representations of $\Pi_\calQ$. 
\begin{lemma}
	Let $\scrF\coloneqq \scrT^\perp$. Then, $\upsilon=(\scrT, \scrF)$ is a torsion pair of $\Modd(\Pi_{\calQ^w})$.
\end{lemma}

\begin{proof}
	One has to show that the pair $\upsilon=(\scrT, \scrF)$ satisfies the two conditions in Definition~\ref{def:torsion_pair}. The first one holds by definition of $\scrF$. 
	
	Now we prove that the second condition in the definition of torsion pairs hold for $\upsilon$. Let $\overline{M}$ be a representation of $\Pi_{\calQ^w}$. We need to show that it fits into a short exact sequence of the form
	\begin{align}
		0\longrightarrow \overline{T}\longrightarrow \overline{M} \longrightarrow \overline{F} \longrightarrow 0 
	\end{align}
	with $\overline{T}\in \scrT$ and $\overline{F}\in \scrF$. 
	
	Let $\overline{T}_{\mathsf{tor}}$ be the maximal submodule of $\overline{M}$ such that $T_\infty=0$ and set $\overline{Q}\coloneqq \overline{M}/\overline{T}$. If $\overline{Q}\notin\scrF$, there exists $\overline{T}'\in \scrT$ such that $\Hom_{\Pi_{\calQ^w}}(\overline{T}', \overline{Q})\neq 0$. In particular, we can assume that $\overline{T}'\subset \overline{Q}$. We claim that there exists a maximal $\overline{T}'\subset \overline{Q}$ such that $\overline{T}'\in \scrT$. Indeed, we can consider a sequence $\overline{T}_1\subset \overline{T}_2\subset \cdots \overline{T}_n\subset \cdots \subset \overline{Q}$ with $\overline{T}_i \in \scrT$ for any $i$. In particular, it must stabilize.
	
	Now, we set $\overline{F}\coloneqq \overline{Q}/\overline{T}'$, where $\overline{T}'\subset \overline{Q}$ is the maximal submodule such that $\overline{T}'\in \scrT$. Then $\overline{F} \in \scrF$. Therefore, set $\overline{T}\coloneqq \ker (\overline{E}\to \overline{F})$. Since $\overline{T}$ fits into the short exact sequence
	\begin{align}
		0\longrightarrow \overline{T}_{\mathsf{tor}}\longrightarrow \overline{T} \longrightarrow \overline{T}' \longrightarrow 0
	\end{align}
	and $\scrT$ is closed under extensions, $\overline{T}\in \scrT$. Thus, the claim follows.
\end{proof}

Set $\scrC_{\calQ^w}\coloneqq \Pi_{\calQ^w}\Mod$. Recall that the compact objects are exactly the finitely generated ones, hence $\scrC_{\calQ^w}=\Ind( \Pi_{\calQ^w}\Mod^{\mathsf{fg}} )$, where $\Pi_{\calQ^w}\Mod^{\mathsf{fg}}$ is the stable $\infty$-category of finitely generated right modules over $\Pi_{\calQ^w}$. Since $\Pi_{\calQ^w}$ is concentrated in homologically nonnegative degrees, the heart of the standard $t$-structure $\tau_{\mathsf{std}}$ of $\scrC_{\calQ^w}$ is the abelian category of right modules over the corresponding preprojective algebra $\Pi_{\calQ^w}$. The category $\scrC_{\calQ^w}$ is a $k$-linear compactly generated presentable stable $\infty$-category. Since the $\infty$-category $k\calQ^w\Mod$ is of finite type, also $\scrC_{\calQ^w}$ is of finite type by \cite[Theorem~5.9]{BCS_Quiver}. Moreover, the standard $t$-structure $\tau_{\mathsf{std}}$ satisfies Assumption~\ref{assumption:t_structure_filtered_colimits}. We shall use the notation $\upsilon=(\scrT, \scrF)$ also for the canonical induced torsion pair on the abelian category of right modules of $\Pi_{\calQ^w}$, whose existence is assured by Corollary~\ref{cor:completion-torsion-pair}. 

Let $\bfRep(\Pi_{\calQ^w})\coloneqq\bfCohps(\scrC_{\calQ^w}, \tau_{\mathsf{std}})$ be the derived stack of $\tau_{\mathsf{std}}$-flat pseudo-perfect objects of $\scrC_{\calQ^w}$. The truncation $\trunc{\bfRep(\Pi_{\calQ^w})}$ of $\bfRep(\Pi_{\calQ^w})$ is a classical geometric stack of finite presentation over $k$ (in particular, $\trunc{\bfRep(\Pi_{\calQ^w})}$ is a quotient stack). The $t$-structure $\tau_{\mathsf{std}}$ universally satisfies openness of flatness, thus $\bfRep(\Pi_{\calQ^w})$ is a geometric derived stack of finite presentation over $k$ by Proposition~\ref{prop:openness}. 

The stack $\bfRep(\Pi_{\calQ^w})$ admits a decomposition into open and closed substacks depending on the dimension vectors $v$ and $w_\infty$ of $M$ and $M_\infty$, respectively:
\begin{align}
	\bfRep(\Pi_{\calQ^w}) = \bigsqcup_{\genfrac{}{}{0pt}{}{v\in \N^{\calQ_0}\ ,}{w_\infty\in \N}} \bfRep(\Pi_{\calQ^w}; v, w_\infty)\ . 
\end{align}

Now, we note that
\begin{align}
	\bfRep_{\scrT}(\Pi_{\calQ^w})\coloneqq \bfCoh_{\scrT}(\scrC_{\calQ^w}, \tau_{\mathsf{std}}) \simeq \bigsqcup_{v\in \N^{\calQ_0}} \bfCohps(\scrC_{\calQ^w}, \tau_{\mathsf{std}}; v, 0)\ .
\end{align}
Hence, $\bfRep_{\scrT}(\Pi_{\calQ^w})$ is an open and closed substack of $\bfRep(\Pi_{\calQ^w})$. Moreover, 
\begin{align}
	\bfRep_{\scrT}(\Pi_{\calQ^w}) \simeq \bfCohps(\scrC_{\calQ}, \tau_{\mathsf{std}})\eqqcolon \bfRep(\Pi_{\calQ})\ ,
\end{align}
where $\scrC_{\calQ}\coloneqq \Pi_\calQ\Mod$. 

Consider the derived stack $\bfRep_{\scrF}(\Pi_{\calQ^w})\coloneqq \bfCoh_{\scrF}(\scrC_{\calQ^w}, \tau_{\mathsf{std}})$ parametrizing $\tau_{\mathsf{std}}$-flat pseudo-perfect objects belonging to $\scrF$. We can characterize it thanks to the following proposition.

Let us first introduce some notation. Fix $\theta_\infty\in \Q$, set 
\begin{align}
	\theta\coloneqq (1, \ldots, 1)\in \Q^{\calQ_0}\quad\text{and} \quad \overline{\theta}\coloneqq(\theta, \theta_\infty)\in \Q^{\calQ_0}\oplus\Q\ .
\end{align}
Define the $\overline{\theta}$-slope of a finite-dimensional representation $\overline{M}=(M, M_\infty, x, y, \alpha, \beta)$ of $\Pi_{\calQ^w}$ as
\begin{align}
	\mu_{\overline{\theta}}(\overline{M})\coloneqq \frac{\sum_{i\in \calQ_0}\dim M_i+\theta_\infty \dim M_\infty}{\sum_{i\in \calQ_0}\dim M_i+\dim M_\infty} \ .
\end{align}
Fix $\theta_\infty <1$. Then, for any finite-dimensional representation $\overline{M}$ of $\Pi_{\calQ^w}$ we get $\mu_{\overline{\theta}}(\overline{M})\leq 1$. The following is evident.
\begin{proposition}\label{prop:conditions-quiver}
	Let $\overline{M}$ be a a finite-dimensional representation of $\Pi_{\calQ^w}$. Then
	\begin{enumerate}\itemsep0.2cm
		\item $\overline{M}\in \scrT$ if and only if $\overline{M}$ is $\overline{\theta}$-semistable of slope one.
		\item \label{item:condition-torsion-free} $\overline{M}\in \scrF$ if and only if $\mu_{\overline{\theta}\textrm{-}\!\max}(\overline{M})<1$.\footnote{Here, $\mu_{\overline{\theta}\textrm{-}\!\max}(\overline{M})$ is the $\overline{\theta}$-slope of the maximal destabilizing sub-representation of $M$.} 
	\end{enumerate}
\end{proposition}
Condition~\eqref{item:condition-torsion-free} above is open, hence $\bfRep_{\scrF}(\Pi_{\calQ^w})$ is an open substack of $\bfRep(\Pi_{\calQ^w})$. 

Also, the stack $\bfRep_{\scrF}(\Pi_{\calQ^w})$ admits a decomposition into open and closed substacks:
\begin{align}
	\bfRep_{\scrF}(\Pi_{\calQ^w}) = \bigsqcup_{\genfrac{}{}{0pt}{}{v\in \N^{\calQ_0}\ ,}{w_\infty\in \N, \ w_\infty\neq 0}} \bfRep_{\scrF}(\Pi_{\calQ^w}; v, w_\infty)\ . 
\end{align}
Set
\begin{align}
	\bfRep_{\scrF}(\Pi_{\calQ^w}; w_\infty) = \bigsqcup_{v\in \N^{\calQ_0}} \bfRep_{\scrF}(\Pi_{\calQ^w}; v, w_\infty)
\end{align}
for any $w_\infty\in \N$, $w_\infty\neq 0$.

Now, note that by construction $\scrC_{\calQ^w}$ is 2-Calabi-Yau (cf.\ \cite[Theorem~5.9]{BCS_Quiver}). In addition, the map \eqref{eq:partial-1-proper} for $\bfRep(\Pi_{\calQ^w})$ is locally rpas (see e.g. \cite[\S3.2]{YZ_COHA}). Thus, $\catCohb_{\mathsf{pro}}( \bfRep(\Pi_{\calQ^w}) )$ has the structure of an $\E_1$-monoidal stable pro-$\infty$-category by Theorem~\ref{thm:COHA-Coh}. A similar statement holds at the level of motivic Borel-Moore homology groups.

Moreover, the conditions \ref{item:algebra-torsion} and \ref{item:rep-torsion-free} of Proposition~\ref{prop:torsion_pairs_left_right_actions} hold for $\bfRep(\Pi_\calQ)$. Thus, we get the following result.
\begin{theorem}\label{thm:action-quiver}
	The stable pro-$\infty$-category $\catCohb_{\mathsf{pro}}( \bfRep(\Pi_\calQ) )$ has a $\E_1$-monoidal structure. Moreover, the stable $\infty$-pro-category $\catCohb_{\mathsf{pro}}( \bfRep_{\scrF}(\Pi_{\calQ^w})  )$ has the structure of a left categorical module over $\catCohb_{\mathsf{pro}}( \bfRep(\Pi_\calQ) )$. 
	
	Similarly, let $\bfD^\ast$ be a motivic formalism, and fix $\calA \in \CAlg(\bfD^\ast(\Spec(k)))$ and $\Gamma \subseteq \Pic(\bfD^\ast(\Spec(k)))$ such that Assumption~\ref{assumption:motivic_formalism} is satisfied. Then, the topological vector space $\HBMDGamma_0(\bfRep(\Pi_\calQ);\calA)$ has the structure of a unital associative algebra. In particular,
	\begin{align}
		G_0( \bfRep(\Pi_\calQ) )\quad \text{and} \quad \HBM_\ast( \bfRep(\Pi_\calQ) )
	\end{align}
	have the structures of unital associative algebras, and the topological vector space $\HBMDGamma_0(\bfRep_{\scrF}(\Pi_{\calQ^w});\calA)$ has the structure of a left $\HBMDGamma_0(\bfRep(\Pi_\calQ);\calA)$-module. In particular,
	\begin{align}
		G_0( \bfRep_{\scrF}(\Pi_{\calQ^w}) )\quad \text{and} \quad \HBM_\ast( \bfRep_{\scrF}(\Pi_{\calQ^w}) )			
	\end{align}
	have the structures of a left $G_0( \bfRep(\Pi_\calQ) )$-module and $\HBM_\ast( \bfRep(\Pi_\calQ) )$-module, respectively.
	
	For each dimension vector $w_\infty$, the same result holds for $\bfRep_{\scrF}(\Pi_{\calQ^w}; w_\infty)$ instead of $\bfRep_{\scrF}(\Pi_{\calQ^w})$. Moreover, the same results hold equivariantly with respect to the torus $T$ introduced in \cite[\S3.3]{SV_generators}.
\end{theorem}
Note that the first part of the above theorem, i.e., the constructions of the two-dimensional cohomological Hall algebra of a quiver and its categorification have been already given in \cite{SV_Cherednik, SV_Elliptic, SV_generators, YZ_COHA, VV_KHA}. The left (categorical) action constructed above has not been studied before in the literature.
\begin{remark}
	Note that we cannot apply directly Theorem~\ref{thm:left-right-action} since the map 
	\begin{align}
		\calS_2\bfCohps(\scrC_{\calQ^w},\tau_\upsilon) \longrightarrow \bfCohps(\scrC_{\calQ^w},\tau_\upsilon) 
	\end{align}
	is not locally rpas in general. Here, $\tau_\upsilon$ is the tilted $t$-structure with respect to $(\scrT, \scrF)$.
\end{remark}

\begin{remark}
	Let $\Lambda^0_\calQ$ and $\Lambda^1_\calQ$ be the Lagrangian substacks introduced in \cite{SV_generators}. Since they are closed substacks of $\trunc{\bfRep(\Pi_\calQ)}$, a version of Theorem~\ref{thm:action-quiver} holds for $\Lambda^0_\calQ$ and $\Lambda^1_\calQ$ instead of $\bfRep(\Pi_\calQ)$, providing a categorification of \cite[Theorem~B--(a) and (b)]{SV_generators}.
\end{remark}

Now, fix $w_\infty=1$. Then, the condition, appearing in Proposition~\ref{prop:conditions-quiver}--\eqref{item:condition-torsion-free}, that $\mu_{\overline{\theta}\textrm{-}\!\max}(\overline{M})<1$ for a finite-dimensional representation $\overline{M}$ is equivalent to the requirement that all sub-represen\-ta\-ti\-ons of $\overline{M}$ have a nonzero vector space associated to the vertex $\infty$. Following \cite[Page~261]{CB_Moment}, if $\overline{M}$ satisfies this condition, we say that is  \textit{$\infty$-co-generated}. 

First, $\bfRep_{\scrF}(\Pi_{\calQ^w}; v, 1)$ coincides with its classical truncation. Moreover, as shown in \textit{loc. cit.}, the stack $\bfRep_{\scrF}(\Pi_{\calQ^w}; v, 1)$ is a $\G_m$-gerbe over the \textit{Nakajima quiver variety} $\calM_{\calQ, \theta}(v, w)$ of $\theta$-stable representations of the preprojective algebra of the framed quiver of $\calQ$ with dimension vectors $v$ and $w_\infty$. Since, $\calM_{\calQ, \theta}(v, w)$ admits a universal sheaf, there is a section $\calM_{\calQ, \theta}(v, w) \to \bfRep_{\scrF}(\Pi_{\calQ^w}; v, 1)$, hence by \cite[Lemma~3.10]{Hein_Moduli} we get that
\begin{align}
	\bfRep_{\scrF}(\Pi_{\calQ^w}; v, 1) \simeq \calM_{\calQ, \theta}(v, w) \times B\G_m\ .
\end{align}
\begin{theorem}\label{thm:action-Nakajima-quiver-variety}
	The stable $\infty$-pro-category 
	\begin{align}
		\bigoplus_{d\in \Z} \catCohb_{\mathsf{pro}}( \calM_{\calQ, \theta}(v, w) )_d
	\end{align}
	has the structure of a left categorical module over $\catCohb_{\mathsf{pro}}( \bfRep(\Pi_\calQ) )$. Here, $\catCohb_{\mathsf{pro}}( \calM_{\calQ, \theta}(v, w) )_d$ is the weight $d$ part of $\catCohb_{\mathsf{pro}}( \bfRep_{\scrF}(\Pi_{\calQ^w}; v, 1) )$.
	
	A similar statement holds at the level of motivic Borel-Moore homology and after replacing $\bfRep(\Pi_\calQ)$ with either $\Lambda^0_\calQ$ or $\Lambda^1_\calQ$. Furthermore, the same results hold equivariantly with respect to the torus $T$ introduced in \cite[\S3.3]{SV_generators}.
\end{theorem}
The left action of the two-dimensional (nilpotent) cohomological Hall algebra of a quiver on the cohomology of Nakajima quiver varieties has been constructed in \cite{YZ_COHA, SV_generators}: the above theorem provides a categorification of such an action.

\newpage

\newcommand{\etalchar}[1]{$^{#1}$}
\providecommand{\bysame}{\leavevmode\hbox to3em{\hrulefill}\thinspace}
\providecommand{\MR}{\relax\ifhmode\unskip\space\fi MR }
\providecommand{\MRhref}[2]{%
	\href{http://www.ams.org/mathscinet-getitem?mr=#1}{#2}
}
\providecommand{\href}[2]{#2}

\end{document}